\documentclass[10pt]{amsart}

\usepackage{latexsym,exscale,enumerate,amsfonts,amssymb, ulem, xparse, mathtools}
\usepackage{amsmath,amsthm,amsfonts,amssymb,amscd, stmaryrd,textcomp,mathscinet}
\usepackage{hyperref}
\usepackage{ulem}
\usepackage{bbold}
\usepackage{tocvsec2}
\usepackage{wasysym}

\usepackage{libertine}
\usepackage[usenames,dvipsnames]{xcolor}
\colorlet{green}{black!30!green} 

\usepackage[all]{xy}
\SelectTips{cm}{}

\usepackage{graphicx}

\usepackage{tikz}
\usetikzlibrary{calc}
\usetikzlibrary{decorations.markings}
\usetikzlibrary{decorations.pathreplacing}
\usetikzlibrary{arrows,shapes,positioning}
\usetikzlibrary{patterns}
\tikzstyle directed=[postaction={decorate,decoration={markings,
    mark=at position #1 with {\arrow{>}}}}]
\tikzstyle rdirected=[postaction={decorate,decoration={markings,
    mark=at position #1 with {\arrow{<}}}}]
\tikzset{anchorbase/.style={baseline={([yshift=-0.5ex]current bounding box.center)}}}

\tikzset{
    partial ellipse/.style args={#1:#2:#3}{
        insert path={+ (#1:#3) arc (#1:#2:#3)}
    }
}

\usepackage{tkz-fct}

\newcommand{\rep}{\Rep(\glnn{N})}

\newcommand{\sessAWeb}{N\cat{A}\cat{Web}_\mathrm{s}^{+,\mathrm{ess}}}
\newcommand{\repp}{\Rep^{+}(\glnn{N})}
\newcommand{\reph}{\Rep(\h)}
\newcommand{\rephp}{\Rep^{+}(\h)}
\newcommand{\Sym}{\mathrm{Sym}}\newcommand{\upcev}[1]{\scriptsize\reflectbox{\ensuremath{\vec{\reflectbox{\ensuremath{#1}}}}}}

\def\X{{\mathbbm X}}
\def\h{{\mathfrak h}}

\newcommand{\ep}{\underline{\epsilon}}
\newcommand{\onen}{{\mathbf 1}_{n}}
\newcommand{\onenn}[1]{{\mathbf 1}_{#1}}
\newcommand{\onenp}{{\mathbf 1}_{n'}}
\newcommand{\onenpp}{{\mathbf 1}_{n''}}
\newcommand{\onem}{{\mathbf 1}_{m}}
\newcommand{\onel}{{\mathbf 1}_{\lambda}}
\newcommand{\onell}[1]{{\mathbf 1}_{#1}}
\newcommand{\onea}{{\mathbf 1}_{{\bf a}}}
\newcommand{\oneaa}[1]{{\mathbf 1}_{#1}}
\newcommand{\onelp}{{\mathbf 1}_{\lambda'}}
\newcommand\rE{{\sf{E}}}
\newcommand\rF{{\sf{F}}}
\newcommand\sE{{\cal{E}}}
\newcommand\sF{{\cal{F}}}
\def\cal#1{\mathcal{#1}}%

\newcommand{\bV}{\textstyle{\bigwedge_q}}
\newcommand{\Alt}[2]{\textstyle{\bigwedge_{#1}^{#2}}}
\newcommand{\Alts}[2]{ \textstyle{\scriptstyle\bigwedge_{#1}^{#2}}}

\def\1{\mathbbm{1}}%
\newcommand\E{{\sf{E}}}
\newcommand\F{{\sf{F}}}
\def\l{\lambda}

\newcommand{\lb}{\text{\textlbrackdbl}}
\newcommand{\rb}{\text{\textrbrackdbl}}

\newcommand{\END}{{\rm END}}
\newcommand{\Gr}{\cat{Flag}_{N}}
\newcommand{\Grn}[1]{\cat{Flag}_{#1}}
\newcommand{\catRep}{{\mathsf{Rep}}}

\theoremstyle{plain}
\newtheorem{theorem}{Theorem}
\newtheorem{theorem*}[theorem]{Theorem*}
\newtheorem{theorem**}[theorem]{Theorem*}
\numberwithin{theorem}{section}

\newtheorem{conjecture}[theorem]{Conjecture}

\newtheorem{corollary*}[theorem]{Corollary*}
\newtheorem{corollary**}[theorem]{Corollary*}

\newtheorem{lemma}[theorem]{Lemma}
\newtheorem{lemma*}[theorem]{Lemma*}
\newtheorem{lemma**}[theorem]{Lemma*}

\newtheorem{notation*}[theorem]{Notation*}

\newtheorem{proposition*}[theorem]{Proposition*}
\newtheorem{proposition**}[theorem]{Proposition*}

\theoremstyle{definition}
\newtheorem{definition}[theorem]{Definition}
\newtheorem{definition*}[theorem]{Definition*}
\newtheorem{definition**}[theorem]{Definition*}

\theoremstyle{remark}

\newcommand{\UU}{{\bf U}}
\newcommand{\Uq}{{\bf U}_q(\mathfrak{sl}_2)}
\newcommand{\U}{\dot{{\bf U}}}
\newcommand{\Ucat}{\cal{U}}
\newcommand{\Ucatc}{\check{\cal{U}}}
\newcommand{\Ucatq}{\cal{U}_q}
\newcommand{\UcatD}{\dot{\cal{U}}}
\newcommand{\UcatDq}{\dot{\cal{U}}_q}
\newcommand{\B}{\dot{\mathbb{B}}}
\newcommand{\Bnm}{{_m\dot{\cal{B}}_n}}
\newcommand{\UA}{{_{\cat{A}}\dot{{\bf U}}}}
\newcommand{\SL}{\mathrm{SL}}
\newcommand{\sln}{\mf{sl}_n}
\newcommand{\slm}{\mf{sl}_m}
\newcommand{\slnn}[1]{\mf{sl}_{#1}}
\newcommand{\gln}{\mf{gl}_n}
\newcommand{\glm}{\mf{gl}_m}
\newcommand{\glnn}[1]{\mf{gl}_{#1}}
\newcommand{\gl}{\mf{gl}}

\newcommand{\und}[1]{\underline{#1}}

\newcommand{\sh}{w}
\newcommand{\cat}[1]{\mathchoice
  {\ensuremath{\mbox{\bfseries {\upshape {#1}}}}}
  {\ensuremath{\mbox{\bfseries {\upshape {#1}}}}}
  {\scalebox{.7}{\ensuremath{\mbox{\bfseries {\upshape {#1}}}}}}
  {\scalebox{.5}{\ensuremath{\mbox{\bfseries {\upshape {#1}}}}}}%
  }

\newcommand{\diagrep}{\phi} 
\newcommand{\pTr}{\mathrm{pTr}}
\newcommand{\bVn}{\bigwedge}
\newcommand{\Kar}{\operatorname{Kar}}
\newcommand{\sym}{\mathrm{Sym}}
\newcommand{\Kh}{\cat{Kh}}
\newcommand{\Bl}{\cat{Bl}}
\newcommand{\APS}{\cat{APS}}
\newcommand{\KhR}{\cat{KhR}}

\newcommand{\TL}{\cat{TL}}
\newcommand{\Web}[1][]{#1\cat{Web}}
\newcommand{\Webp}[1][]{#1\cat{Web}^{+}}
\newcommand{\Webq}[1][]{#1\cat{Web}_{q}}

\newcommand{\wrap}{D}
\newcommand{\wrapi}{D^{-1}}
\newcommand{\A}{\cat{A}}
\newcommand{\ATL}{\cat{A}\cat{TL}}
\newcommand{\AWeb}[1][]{#1\cat{A}\cat{Web}}
\newcommand{\AWebq}[1][]{#1\cat{A}\cat{Web}_{q}}
\newcommand{\AWebp}[1][]{#1\AWeb^+}

\newcommand{\essAWeb}[1][]{#1\AWeb^{\mathrm{ess}}}
\newcommand{\essAWebp}[1][]{#1\AWeb^{\mathrm{ess},+}}

\newcommand{\essbAWeb}[1][]{\overline{#1\AWeb}^{\mathrm{ess}}}
\newcommand{\essbAWebp}[1][]{\overline{#1\AWeb}^{\mathrm{ess},+}}

\newcommand{\Afoam}[1][]{#1\cat{A}\cat{Foam}}
\newcommand{\AFoam}[1][]{#1\cat{A}\cat{Foam}}

\def\mC{\cat{C}}
\def\basis{\mathrm{B}_T}
\def\basisJW{\mathrm{B}_{S}}
\def\basisstd{\mathrm{B}}

\newcommand{\T}{\cat{T}}
\newcommand{\TWebq}[1][]{#1\T\cat{Web}_{q}}
\newcommand{\TWebA}[1][]{#1\T\cat{Web}_{A}}
\newcommand{\Tfoam}[1][]{#1\T\cat{Foam}}
\newcommand{\essTfoam}[1][]{#1\Tfoam^{ess}}
\newcommand{\Tlink}{\T\cat{Link}}
\newcommand{\Ttanweb}{\T\cat{TanWeb}}
\newcommand{\Tfoamred}[1][]{#1\T\cat{Foam}^{\mathrm{red}}}

\newcommand{\Stwo}{\cat{S}^\mathbf{2}}
\newcommand{\Su}{\cat{S}}

\newcommand{\Sk}{\rm Sk}
\newcommand{\SWebq}[1][]{#1\Su\cat{Web}_{q}}
\newcommand{\SWebA}[1][]{#1\Su\cat{Web}_{A}}
\newcommand{\Sfoam}[1][]{#1\Su\cat{Foam}}
\newcommand{\Sskel}[1][]{#1\Su\cat{Skel}}

\newcommand{\Sfoamred}[1][]{#1\Sfoam^{\mathrm{red}}}
\newcommand{\Sfoamor}[1][]{#1\Sfoam^{\mathrm or}}
\newcommand{\SCob}{\Su\cat{Cob}}
\newcommand{\Slink}{\Su\cat{Link}}
\newcommand{\Slinko}{\Su\cat{Link}^{\circ}}
\newcommand{\Stanweb}{\Su\cat{TanWeb}}
\newcommand{\Stanwebo}{\Su\cat{TanWeb}^{\circ}}

\newcommand{\repr}{\mathrm{K_0}(\Rep(\glnn{2}))}
\newcommand{\reprs}{\mathrm{K_0}(\Rep(\slnn{2}))}

\newcommand{\cev}[1]{\reflectbox{\ensuremath{\vec{\reflectbox{\ensuremath{#1}}}}}}
\newcommand{\stwo}{*_{\wedge}}
\newcommand{\w}[1]{\wedge^{(#1)}}

\newcommand{\foam}[3][N]{#2\cat{Foam}_{#3}(#1)}
\newcommand{\Bfoam}[3][N]{#2\cat{BFoam}_{#3}(#1)}
\newcommand{\nWeb}{\cat{nWeb}}
\newcommand{\Foam}[2]{#1\cat{Foam}_{#2}}
\newcommand{\twoFoam}{\cat{Foam}}
\newcommand{\Vect}{\cat{Vect}}

\newcommand{\colTlink}{\T\cat{Link}_{1,2}}
\newcommand{\Tfoamr}{\T\cat{Foam}_r}
\newcommand{\TfoamrZ}{\T\cat{Foam}_r^\Z}
\newcommand{\Tfoamtwo}{\T\cat{Foam}_2}
\newcommand{\HC}{\operatorname{K}}

\newcommand{\Kom}{\operatorname{Kom}}
\newcommand{\Tr}{\operatorname{Tr}}
\newcommand{\vTr}{\operatorname{vTr}}
\newcommand{\hTr}{\operatorname{hTr}}
\newcommand{\grvTr}{\widetilde{\vTr}}
\newcommand{\Rep}{\cat{Rep}}
\newcommand{\grRep}{\cat{grRep}}
\newcommand{\SH}{\operatorname{SH}}
\newcommand{\SHz}{\operatorname{SH}_{t=0}}
\newcommand{\hSHz}{\widehat{\operatorname{SH}}_{t=0}}
\newcommand{\BN}[2]{{_{#1}\mathcal{BN}_{#2}}}
\newcommand{\saKh}{\mathsf{saKh}}
\newcommand{\saKhR}{\mathsf{saKhR}}

\newcommand{\numroman}{\renewcommand{\labelenumi}{\roman{enumi})}}
\newcommand{\numarabic}{\renewcommand{\labelenumi}{\arabic{enumi})}}
\newcommand{\numAlph}{\renewcommand{\labelenumi}{\Alph{enumi}.}}

\newcommand{\To}{\Rightarrow}
\newcommand{\TO}{\Rrightarrow}
\newcommand{\Hom}{{\rm Hom}}
\newcommand{\HOM}{{\rm HOM}}
\renewcommand{\to}{\rightarrow}
\newcommand{\maps}{\colon}
\newcommand{\op}{{\rm op}}
\newcommand{\co}{{\rm co}}
\newcommand{\iso}{\cong}
\newcommand{\id}{{\rm id}}
\newcommand{\bigb}[1]{
\begin{pspicture}(0,0)
 \rput(0,0){\psframebox[framearc=.5,fillstyle=solid]{\small $#1$}}
\end{pspicture}}
\newcommand{\del}{\partial}
\newcommand{\Res}{{\rm Res}}
\newcommand{\End}{{\rm End}}
\newcommand{\Aut}{{\rm Aut}}
\newcommand{\im}{{\rm im\ }}
\newcommand{\coim}{{\rm coim\ }}
\newcommand{\chr}{{\rm char\ }}
\newcommand{\coker}{{\rm coker\ }}
\newcommand{\Span}{{\rm Span}}
\newcommand{\spann}{{\rm span}}
\newcommand{\rk}{{\rm rk\ }}
\def\bigboxtimes{\mathop{\boxtimes}\limits}

\newcommand{\scs}{\scriptstyle}

\def\Res{{\mathrm{Res}}}
\def\Ind{{\mathrm{Ind}}}
\def\lra{{\longrightarrow}}
\def\dmod{{\mathrm{-mod}}}   
\def\fmod{{\mathrm{-fmod}}}   
\def\pmod{{\mathrm{-pmod}}}  
\def\rk{{\mathrm{rk}}}
\def\NH{{\mathrm{NH}}}
\def\pseq{{\mathrm{Seqd}}}
\def\Id{\mathrm{Id}}
\def\mc{\mathcal}
\def\mf{\mathfrak}
\def\Af{{_{\mc{A}}\mathbf{f}}}    
\def\primef{{'\mathbf{f}}}    
\def\shuffle{\,\raise 1pt\hbox{$\scriptscriptstyle\cup{\mskip
               -4mu}\cup$}\,}
\newcommand{\define}{\stackrel{\mbox{\scriptsize{def}}}{=}}

\numberwithin{equation}{section}

\def\OC#1{\textcolor[rgb]{1.00,0.00,0.00}{[Orientability Comment: #1]}}%
\def\BC#1{\textcolor[rgb]{1.00,0.00,0.00}{[Boerner Relation Comment: #1]}}%
\def\comm#1{}%
\def\new#1{\b #1\e}%

\def\emph#1{{\sl #1\/}}
\def\ie{{\sl i.e. \/}}
\def\eg{{\sl e.g. \/}}
\def\Eg{{\sl E.g.\/}}
\def\etc{{\sl etc.\/}}
\def\cf{{\sl c.f.\/}}
\def\etal{\sl{et al.\/}}%
\def\vs{\sl{vs.\/}}%

\let\hat=\widehat
\let\tilde=\widetilde

\let\phi=\varphi
\let\theta=\vartheta
\let\epsilon=\varepsilon

\usepackage{bbm}
\def\C{{\mathbbm C}}
\def\N{{\mathbbm N}}
\def\R{{\mathbbm R}}
\def\Z{{\mathbbm Z}}
\def\Q{{\mathbbm Q}}
\def\H{{\mathbbm H}}
\def\P{{\mathbbm P}}
\newcommand{\Bb}{\mathbb{B}}
\newcommand{\I}{\mathbb{I}}
\newcommand{\Ss}{\mathbb{S}}
\newcommand{\D}{\mathbb{D}}

\def\cal#1{\mathcal{#1}}%
\def\1{\mathbbm{1}}%
\def\ev{\mathrm{ev}}%
\def\coev{\mathrm{coev}}%
\def\tr{\mathrm{tr}}%
\def\st{\mathrm{st}}%
\def\pullback#1#2#3{%
  \,\mbox{\raisebox{-.8ex}{$\scriptstyle #1$}}%
  \!\prod\!
  \mbox{\raisebox{-.8ex}{$\scriptstyle #3$}}\,}%
\def\nn{\notag}
\newcommand{\ontop}[2]{\genfrac{}{}{0pt}{2}{\scriptstyle #1}{\scriptstyle #2}}

\def\la{\langle}
\def\ra{\rangle}

\newcommand{\dotsheet}[2][1]{
\xy
(0,0)*{
\begin{tikzpicture} [fill opacity=0.2,decoration={markings, mark=at position 0.6 with {\arrow{>}}; }, scale=#1]
	\draw [very thick, fill=red] (1,1) -- (-1,2) -- (-1,-1) -- (1,-2) -- cycle;
	\node [opacity=1] at (0,0) {$\bullet^{#2}$};
\end{tikzpicture}};
\endxy}

\newcommand{\sphere}[1][1]{
\xy
(0,0)*{
\begin{tikzpicture} [fill opacity=0.2, scale=#1]
	\path [fill=red] (0,0) circle (1);
	\draw (-1,0) .. controls (-1,-.4) and (1,-.4) .. (1,0);
	\draw[dashed] (-1,0) .. controls (-1,.4) and (1,.4) .. (1,0);
	\draw[very thick] (0,0) circle (1);
\end{tikzpicture}};
\endxy}

\newcommand{\dottedsphere}[2][1]{
\xy
(0,0)*{
\begin{tikzpicture} [fill opacity=0.2, scale=#1]
	\path [fill=red] (0,0) circle (1);
	\draw (-1,0) .. controls (-1,-.4) and (1,-.4) .. (1,0);
	\draw[dashed] (-1,0) .. controls (-1,.4) and (1,.4) .. (1,0);
	\draw[very thick] (0,0) circle (1);
	\node [opacity=1]  at (0,0.6) {$\bullet$};
	\node [opacity=1] at (-0.3,0.6) {$#2$};
\end{tikzpicture}};
\endxy}

\newcommand{\cylinder}[1][1]{
\xy
(0,0)*{
\begin{tikzpicture} [fill opacity=0.2,  decoration={markings, mark=at position 0.6 with {\arrow{>}};  }, scale=#1]
	\path [fill=red] (0,4) ellipse (1 and 0.5);
	\path [fill=red] (0,0) ellipse (1 and 0.5);
	\path[fill=red, opacity=.3] (-1,4) .. controls (-1,3.34) and (1,3.34) .. (1,4) --
		(1,0) .. controls (1,.66) and (-1,.66) .. (-1,0) -- cycle;
	\draw [very thick] (0,4) ellipse (1 and 0.5);
	\draw [very thick] (0,0) ellipse (1 and 0.5);
	\draw[very thick] (1,4) -- (1,0);
	\draw[very thick] (-1,4) -- (-1,0);
\end{tikzpicture}};
\endxy
}

\newcommand{\slthncfour}[1][1]{
\xy
(0,0)*{
\begin{tikzpicture} [fill opacity=0.2,  decoration={markings, mark=at position 0.6 with {\arrow{>}};  }, scale=#1]
	\path [fill=red] (0,4) ellipse (1 and 0.5);
	\path [fill=red, opacity=.3] (1,4) .. controls (1,3.34) and (-1,3.34) .. (-1,4) --
		(-1,4) .. controls (-1,2) and (1,2) .. (1,4);
	\path [fill=red] (0,0) ellipse (1 and 0.5);
	\path [fill=red, opacity=.3] (1,0) .. controls (1,.66) and (-1,.66) .. (-1,0) --
		(-1,0) .. controls (-1,2) and (1,2) .. (1,0);
	\draw [very thick] (0,4) ellipse (1 and 0.5);
	\draw [very thick] (0,0) ellipse (1 and 0.5);
	\draw[very thick] (-1,0) .. controls (-1,2) and (1,2) .. (1,0);
	\draw [very thick] (-1,4) .. controls (-1,2) and (1,2) .. (1,4);
	\node[opacity=1] at (0,1) {$\bullet$};
\end{tikzpicture}};
\endxy
}

\newcommand{\slthncfive}[1][1]{
\xy
(0,0)*{
\begin{tikzpicture} [fill opacity=0.2,  decoration={markings, mark=at position 0.6 with {\arrow{>}};  }, scale=#1]
	\path [fill=red] (0,4) ellipse (1 and 0.5);
	\path [fill=red, opacity=.3] (1,4) .. controls (1,3.34) and (-1,3.34) .. (-1,4) --
		(-1,4) .. controls (-1,2) and (1,2) .. (1,4);
	\path [fill=red] (0,0) ellipse (1 and 0.5);
	\path [fill=red, opacity=.3] (1,0) .. controls (1,.66) and (-1,.66) .. (-1,0) --
		(-1,0) .. controls (-1,2) and (1,2) .. (1,0);
	\draw [very thick] (0,4) ellipse (1 and 0.5);
	\draw [very thick] (0,0) ellipse (1 and 0.5);
	\draw[very thick] (-1,0) .. controls (-1,2) and (1,2) .. (1,0);
	\draw [very thick] (-1,4) .. controls (-1,2) and (1,2) .. (1,4);
	\node[opacity=1] at (0,3) {$\bullet$};
\end{tikzpicture}};
\endxy
}

\newcommand{\capFEfoam}[3][.5]{
\xy
(0,0)*{
\begin{tikzpicture} [scale=#1,fill opacity=0.2]
\path[fill=red] (2.5,2) to (-2.5,2) to (-2.5,-1) to (2.5,-1);
\path[fill=blue] (1.5,-2) to [out=90,in=0] (0,.5) to [out=180,in=90] (-1.5,-2) to
	(-.5,-1) to [out=90,in=180] (0,-.25) to [out=0,in=90] (.5,-1);
\path[fill=red] (2,1) to (-2,1) to (-2,-2) to (2,-2);
\draw[very thick, directed=.5] (2.5,2) to (-2.5,2);
\draw[very thick] (2.5,2) to (2.5,-1);
\draw[very thick] (-2.5,2) to (-2.5,-1);
\draw[very thick, directed=.5] (2.5,-1) to (-2.5,-1);
\draw[very thick, red, directed=.75] (-.5,-1) to [out=90,in=180] (0,-.25)
	to [out=0,in=90] (.5,-1);
\node[red, opacity=1] at (1.75,1.5) {$_{#2}$};
\draw[very thick, directed=.5] (1.5,-2) to (.5,-1);
\draw[very thick, directed=.5] (-.5,-1) to (-1.5,-2);
\draw[very thick, directed=.5] (2,1) to (-2,1);
\draw[very thick] (2,1) to (2,-2);
\draw[very thick] (-2,1) to (-2,-2);
\draw[very thick, directed=.5] (2,-2) to (-2,-2);
\draw[very thick, red, directed=.65] (1.5,-2) to [out=90,in=0] (0,.5)
	to [out=180,in=90] (-1.5,-2);
\node[red, opacity=1] at (1.25,.5) {$_{#3}$};
\end{tikzpicture}};
\endxy
}

\newcommand{\torus}[2]{
\draw [green, thick, directed=.25]  (0,0) to (#1,0);
\draw [green, thick, directed=.25] (0,#2)to (#1,#2);
\draw [red, thick,rdirected=.3, rdirected=.25] (#1,#2) to (#1,0);
\draw [red, thick, rdirected=.3, rdirected=.25] (0,#2) to (0,0);}

\newcommand{\torusd}[2]{
\draw [green, thick, directed=.55]  (0,0)to (#1,0) ;
\draw [green, thick, directed=.55] (0-0.5*#2,#2) to  (#1-0.5*#2,#2) ;
\draw [red, thick, rdirected=.5, rdirected=.55] (#1-0.5*#2,#2) to (#1,0);
\draw [red, thick, rdirected=.5, rdirected=.55] (0-0.5*#2,#2) to (0,0);}

\newcommand{\torusfront}[3]{
\draw [green, thick, directed=.25] (0,0)to(#1,0) ;
\draw [red, thick, rdirected=.3, rdirected=.25] (#1+0.5*#2,#2) to (#1,0);
\draw [green, thick, directed=.25] (0,0+#3) to (#1,0+#3);
\draw [green, thick, directed=.25] (0+0.5*#2,#2+#3) to (#1+0.5*#2,#2+#3);
\draw [red, thick, rdirected=.3, rdirected=.25] (#1+0.5*#2,#2+#3) to (#1,0+#3);
\draw [red, thick, rdirected=.3, rdirected=.25] (0+0.5*#2,#2+#3) to (0,0+#3);
\draw (0,0) to (0,0+#3);
\draw(#1,0) to(#1,0+#3);
}

\newcommand{\torusback}[3]{
\draw [opacity=.5] (0.5*#2,#2)  to  (0.5*#2,#2+#3); 
\draw (#1+0.5*#2,#2) to (#1+0.5*#2,#2+#3);
\draw [green,opacity=.5, thick, directed=.25] (0+0.5*#2,#2)to (#1+0.5*#2,#2) ;
\draw [red,opacity=.5, thick, rdirected=.3, rdirected=.25] (0+0.5*#2,#2) to (0,0);
}

\newcommand{\annfrontd}[3]{
\draw (0,0)to(#1,0) ;
\draw [red, thick, rdirected=.3, rdirected=.25] (#1+0.5*#2,#2) to (#1,0);
\draw [dashed] (0,0+#3) to (#1,0+#3);
\draw [dashed] (0+0.5*#2,#2+#3) to (#1+0.5*#2,#2+#3);
\draw [dashed,red, thick, rdirected=.3, rdirected=.25] (#1+0.5*#2,#2+#3) to (#1,0+#3);
\draw [dashed, red, thick, rdirected=.3, rdirected=.25] (0+0.5*#2,#2+#3) to (0,0+#3);
\draw (0,0) to (0,0+#3);
\draw(#1,0) to(#1,0+#3);
}

\newcommand{\annfront}[3]{
\draw (0,0)to(#1,0) ;
\draw [red, thick, rdirected=.3, rdirected=.25] (#1+0.5*#2,#2) to (#1,0);
\draw (0,0+#3) to (#1,0+#3);
\draw (0+0.5*#2,#2+#3) to (#1+0.5*#2,#2+#3);
\draw [red, thick, rdirected=.3, rdirected=.25] (#1+0.5*#2,#2+#3) to (#1,0+#3);
\draw [red, thick, rdirected=.3, rdirected=.25] (0+0.5*#2,#2+#3) to (0,0+#3);
\draw (0,0) to (0,0+#3);
\draw(#1,0) to(#1,0+#3);
}

\newcommand{\annback}[3]{
\draw [opacity=.5] (0.5*#2,#2)  to  (0.5*#2,#2+#3); 
\draw (#1+0.5*#2,#2) to (#1+0.5*#2,#2+#3);
\draw [opacity=.5] (0+0.5*#2,#2)to (#1+0.5*#2,#2) ;
\draw [red,opacity=.5, thick, rdirected=.3, rdirected=.25] (0+0.5*#2,#2) to (0,0);
}

\begin{document}
\allowdisplaybreaks

\title[Functoriality and skein positivity]{$\gl_2$ foam functoriality and skein positivity}

\author{Hoel Queffelec}
 \address{IMAG\\ Univ. Montpellier\\ CNRS \\ Montpellier \\ France}
\email{hoel.queffelec@umontpellier.fr}

\begin{abstract}
We prove full functoriality of Khovanov homology for tangled framed $\gl_2$ webs. We use this functoriality result to prove a strong positivity result for (orientable) surface skein algebras. The argument goes categorical and consists in proving that so-called linear complexes are stable under superposition.
\end{abstract}

\maketitle
\maxtocdepth{subsection}
\tableofcontents

\section{Introduction}

The definition of skein modules goes back to Przytycki~\cite{Prz1} and Turaev~\cite{Tur}, in the spirit of Conway's interpretation of knot invariants. This definition provides a diagrammatic extension of the Jones polynomial~\cite{Jones} to 3-manifolds, alternative to the original extension of Witten, Reshetikhin and Turaev~\cite{Witten3,RT}. In the case where the 3-manifold is a thickened surface $\Su\times [0,1]$, the skein module $\Sk(\Su):=\Sk(\Su\otimes [0,1])$ is endowed with an algebra structure by stacking.

These algebra have been extensively studied, starting from Frohman-Gelca's work~\cite{FG}. There they found a basis constructed using Chebyshev polynomials, that enjoys positive and particularly easy structure constants. Relations to character varieties have been explored for example by Bullock~\cite{Bul}, and Bonahon-Wong beautifully described the relation to quantum Teichmüller theory~\cite{BW,BW1} (see also Le's work~\cite{Le_teichm}). Skein algebras for surfaces also enjoy deep connections with the theory of cluster algebras. These connections were explored by Muller~\cite{Muller} and further investigated by Le~\cite{Le_triangular} in relation to Bonahon-Wong's quantum trace. Such connections to the theory of cluster algebras also underlie Thurston's work on positivity~\cite{Thu}, reformulating a conjecture of Fock-Goncharov~\cite{FoG} that has since been also extensively studied by Le, for example in~\cite{LeThurstonYu}.

This conjecture, central in skein theory, states that the Chebyshev basis introduced by Frohman-Gelca can be extended to any surface and yields a positive basis in general. It is known to hold for any surface once the quantum variable $q$ is evaluated to $1$~\cite{Thu}, or at generic $q$ for the cases of the torus~\cite{FG}, the punctured torus and the 4-punctured sphere~\cite{Bousseau}. The latter two results of Bousseau surprisingly go through Gromov-Witten invariants.

The goal of this paper is to exploit the structures offered by categorification in skein theory in order to prove positivity of surface skein algebras. Indeed, just as the Jones polynomial was categorified by Khovanov~\cite{Kh1}, one can lift the skein setting one level up. The most direct way to go is to extend Bar-Natan's reformulation of Khovanov homology from~\cite{BN2}, and consider a category of cobordisms over the surface, up to the very same local relations that controlled Khovanov homology. This strategy, however, does not yield very handy invariants, and attempts have been made to incorporate surface-specific data to the invariants. In particular, the explicit invariants defined by Asaeda, Przytycki and Sikora~\cite{APS} (and Roberts in the annular case~\cite{Roberts}) can be reinterpreted as the degree-zero part of the extension of Bar-Natan's category for a surface-specific non-negative grading (see joint work with Paul Wedrich \cite{QW_SkeinCat}, inspired by the work of Grigsby-Licata-Wehrli in the annular case~\cite{GLW}).

One can then go on and classify simples in the degree-zero part of these cobordism categories, and prove that the situation somewhat differs depending whether $\Su$ is the torus or another surface. Indeed, it is proven in \cite{QW_SkeinCat} that the categories are semi-simple with simples:
\begin{itemize}
\item if $\Su\neq\Ss^1\times \Ss^1$, obtained from Jones-Wenzl projectors put upon multicurves, decategorifying to the skein basis obtained by using Chebyshev polynomials of the $2^{\rm nd}$ kind;
\item if $\Su=\Ss^1\times \Ss^1$ (and after suitable reduction of the category), obtained from highest weight projectors put upon multicurves , decategorifying to the Frohman-Gelca skein basis obtained by using Chebyshev polynomials of the $1^{\rm st}$ kind (see~\cite{QW}).
\end{itemize}

This observation led us to reformulate the Fock-Goncharov conjecture into:
\begin{conjecture}
  The skein basis obtained from the multicurve basis by application of Chebyshev
  polynomials of the second kind is positive provided
  $\Su\neq \Ss^1\times \Ss^1$.
\end{conjecture}

A proof for this conjecture is the main result of the present paper (Theorem~\ref{thm:positivity}). Note that this does not mean anything for the positivity of the other basis, which is at the other extreme of the range of possibilities (see~\cite{LTY}). One might hope though that since the Frohman-Gelca basis is finer than the one used in this paper, the latter being positive might imply that the former is as well. Categorically, this would go by adding more idempotents to break down the simples into smaller ones.

The main issue when trying to prove a categorical analog of it lies in defining a stacking operation at the categorical level. Indeed, to be able to state a precise definition, one needs to go to a setting where Khovanov homology is functorial. Recall that in the original setting, one assigns the homotopy class of a complex to a knot or link, but one can also translate a cobordism between two knots into a chain map between the complexes associated to the source and target links. However, this assignment varies by a sign when performing isotopies.

Fixes to the functoriality gap have been provided by Clark, Morrison and Walker~\cite{CMW}, Caprau~\cite{Cap} and Blanchet~\cite{Blan}. The latter fix makes use of so-called foams instead of cobordisms. Unfortunately, passing from cobordisms to foams also invites to extend the objects in the category, from links to knotted webs. In order to properly define the stacking operation, we need a functoriality result for the category of knotted webs, which so far was missing.

In a recent collaboration with K. Walker~\cite{QWalker}, a definition for framed foams was proposed, with a description of their categories by generators and relations (movie moves). This extends work of Carter and Saito~\cite{CS_book} in the case of links, and of Carter~\cite{Carter_foams} for foams (see also~\cite{CarterKamada}). We specialize this work to state Theorem~\ref{thm:MM}, providing a full list of movie moves for $\gl_2$ foams.

Then proving that Khovanov homology is fully functorial in the web and foam setting is an exhausting but almost straightforward computation. We build upon previous results in the link case by Blanchet~\cite{Blan} and Ehrig-Tubbenhauer-Wedrich~\cite{ETW}, and make use of the spectral sequence to fully deformed homology to compare constants when projective functoriality is ensured by Bar-Natan's trick~\cite{BN2}.

Finally, with full functoriality in hand, one can stack webs and foams and
define a bifunctor on the category of webs and foams over a surface. Although this bifunctor does not a priori extend to a monoidal structure on the homotopy category, it still allows us to lift the skein product categorically.

Positivity rephrases as follows: applying Khovanov's functor to stacked Jones-Wenzl idempotents over multicurves yields a (so-called linear) complex whose objects have balanced homological and quantum shifts (this describes the heart of a $t$-structure). Using this reformulation allows us to very much reduce the cases to check in Lemma~\ref{lem:reduc}, using the idea that at the categorical level:
\[
  A\oplus B\;\text{linear}\;\Rightarrow A\;\text{linear}\;\rm{and}\; B\;\text{linear}.
\]
Note that an analogous decategorified statement is completely wrong: it is not because $A+B$ is positive that $A$ and $B$ are positive!

Then the proof goes by contradiction: from a minimal pair yielding a non-linear complex, one can build a smaller pair with the same property. Thus no such pair can exist.

\settocdepth{section}
\subsection*{Acknowledgements}
This work originates in a collaboration with Paul Wedrich. He largely contributed to the present paper and I would like to warmly thank him for the many hours of discussions! I also would like to thank Thang L\^{e} for several helpful conversations over the years, as well as Kevin Walker for the interest he showed in this work and Dylan Thurston for helpful discussions.

\subsection*{Funding}
This work was partly funded by the CNRS-MSI partnership {\it FuMa}, the ANR grants QUANTACT and CATORE and a CNRS PEPS. I am currently funded by the European Union’s Horizon 2020 research and innovation programme under the Marie Sklodowska-Curie grant agreement No 101064705.

\settocdepth{subsection}

\section{Foams and movie moves}
\label{sec:foams}

In this section we briefly review the notions of webs and foams. We rely on joint work with K. Walker~\cite{QWalker} to give a presentation of the category of foams up to isotopy. This is akin to Carter and Saito's list of movie moves~\cite{CS_book,Carter_foams}.

We follow the definitions and settings of~\cite{QWalker} for the topological aspects (in particular everything is smooth). This is equivalent to adopting Blanchet's $\glnn{2}$ foam formalism~\cite{Blan}, with an extra framing structure as in~\cite{QWalker} and in the smooth setting.

The main tool we will use in the faithfulness part of this paper is the following theorem, which is the $\glnn{2}$ specialization of the last theorem in~\cite{QWalker}: we have omitted all moves that contain a 6-valent vertex. Indeed, such a situation does not arise in the context of $\glnn{2}$ foams. Before stating the theorem, let us simply recall that $\glnn{2}$ webs are oriented labeled trivalent graphs. Edges are labeled $1$ or $2$ (then we often draw them doubled) and trivalent vertices satisfy the obvious flow condition. Webs are not meant to be connected and can contain loops. Furthermore, we assume that trivalent vertices arrange as follows:
    
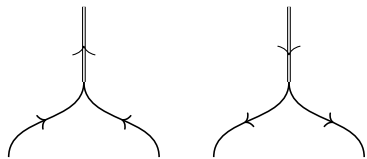
\[
  \begin{tikzpicture}[anchorbase,decoration={markings, mark=at position 0.5 with {\arrow{>}}; }]
    \draw [semithick, postaction={decorate}] (0,0) to [out=90,in=-90] (1,1);
    \draw [semithick, postaction={decorate}] (2,0) to [out=90,in=-90] (1,1);
    \draw [double, postaction={decorate}] (1,1) to [out=90,in=-90] (1,2);
  \end{tikzpicture}
  \qquad
  \begin{tikzpicture}[anchorbase,decoration={markings, mark=at position 0.5 with {\arrow{<}}; }]
    \draw [semithick, postaction={decorate}] (0,0) to [out=90,in=-90] (1,1);
    \draw [semithick, postaction={decorate}] (2,0) to [out=90,in=-90] (1,1);
    \draw [double, postaction={decorate}] (1,1) to [out=90,in=-90] (1,2);
  \end{tikzpicture}
\]

Then foams are natural cobordisms between webs. They are realized by singular surfaces. In the $\glnn{2}$ case, the singularities are simply singular lines around which three sheets meet. Just as webs, facets of a foam carry a label that matches the one on the boundary.

The following theorem thus gives a presentation by generators and relations of the $\glnn{2}$ foam category. It is presented in terms of movies and movie moves: generators are successions of frames that realize an elementary cobordism (see for example the cap morphism represented at the beginning of Equation~\eqref{eq:gen1}). Labels and orientations have been omitted, and any labeling and orientation compatible with the rules stated above should be considered. For examples, the cap can be $1$ or $2$ labeled and the circle in the first frame can have any orientation. The blister at the end of Equation~\eqref{eq:gen1} on the other only has choice of labels (and two orientations can appear but they are related by a $180^\circ$ rotation):
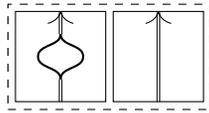
\[
  \begin{tikzpicture}[anchorbase]
    \draw [dashed] (-.7,-.7) rectangle (2,.7);
    \node at (0,0) {
      \begin{tikzpicture}[anchorbase, scale=.3,rotate=90]
        \draw [double] (0,2) -- (1,2);
        \draw [thick] (1,2) to [out=0,in=180] (2,3) to [out=0,in=180] (3,2);
        \draw [thick] (1,2) to [out=0,in=180] (2,1) to [out=0,in=180] (3,2);
        \draw [double,->] (3,2) -- (4,2);
        \draw (0,0) rectangle (4,4);
      \end{tikzpicture}};
    \node at (1.3,0) {\begin{tikzpicture}[anchorbase, scale=.3,rotate=90]
        \draw [double,->] (0,2) -- (4,2);
        \draw (0,0) rectangle (4,4);
      \end{tikzpicture}};
  \end{tikzpicture}
\]

We again refer to~\cite{QWalker} for detailed explanations regarding framing. Let us simply state that the symbol $\LEFTcircle$ that appears on strands represents a full twist in framing. Assuming that we like the framing to be in a standard position, a twist is a place where this is not realized and it goes all around the strand, in one or the other direction (hence the sign next to the twist). Often, one resolves this by inserting a curl. We insist that this is not what we want to do. Again, reasons for doing so are fully detailed in~\cite{QWalker}.

\begin{theorem} \label{thm:MM} 
  Framed foams between oriented framed web-tangles with preferred diagrams admit the following movie generators, in addition to identity movies over web-tangles (each picture is a representative of a family, obtained by crossing changes, changing half-twists types and signs, orientations, reversing frame ordering or taking planar symmetries):

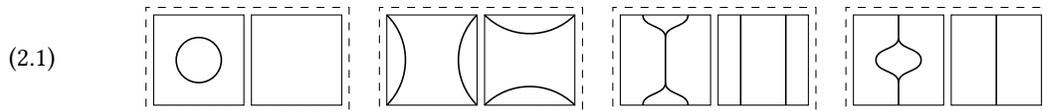
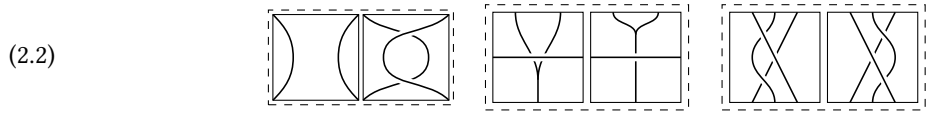
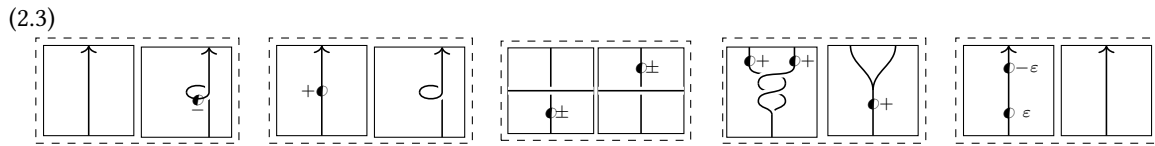
\begin{equation} \label{eq:gen1}

\end{equation}

Isotopies of framed foams induce sequences of movie moves from the following list:

\noindent {\bf framed version classical movie moves}

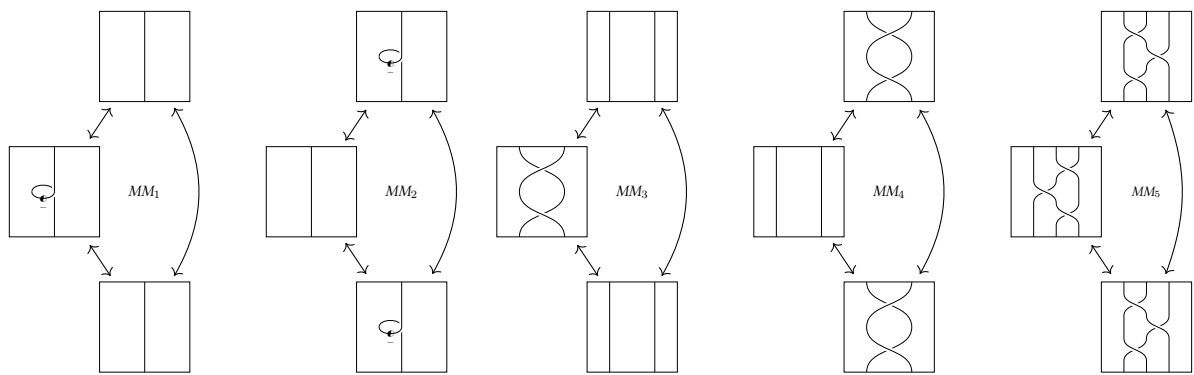
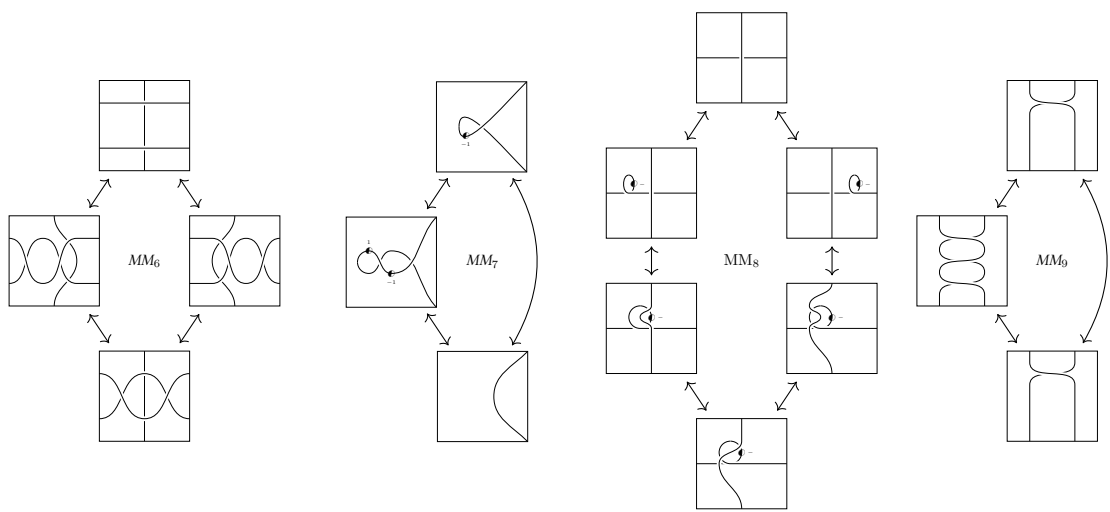
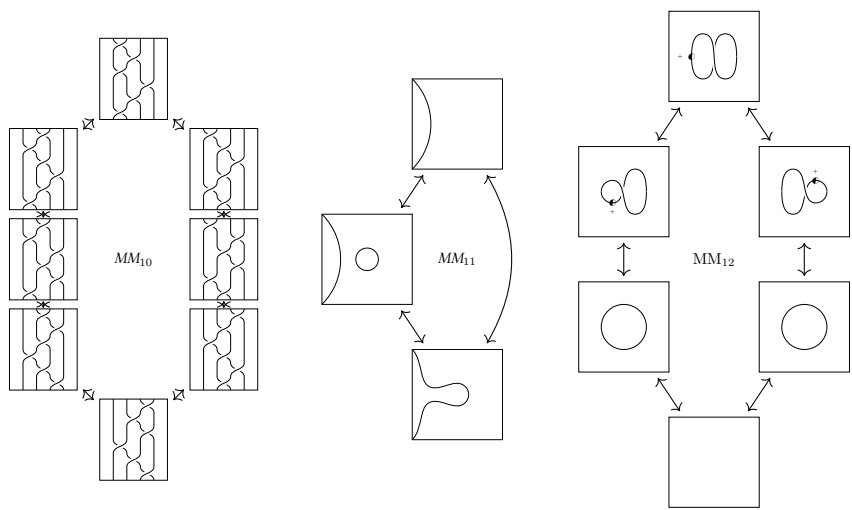
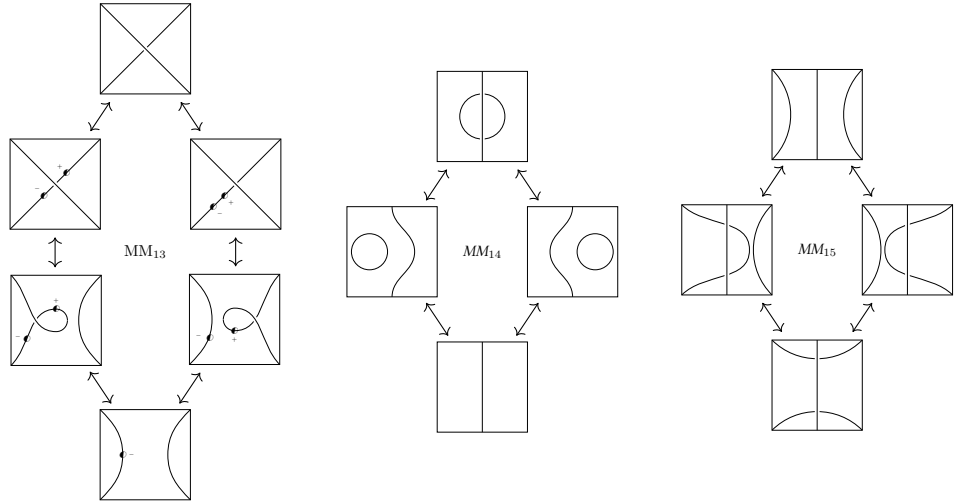
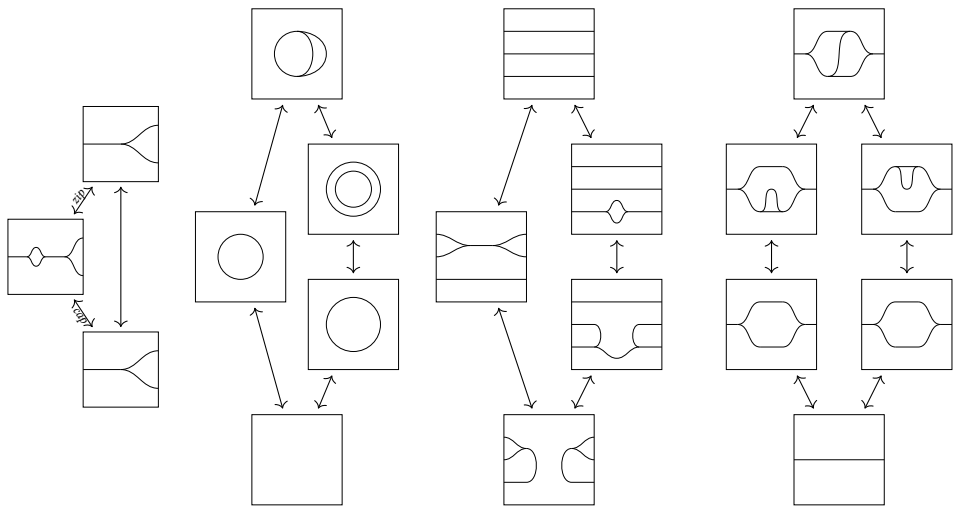
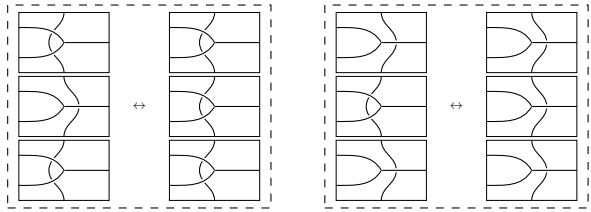
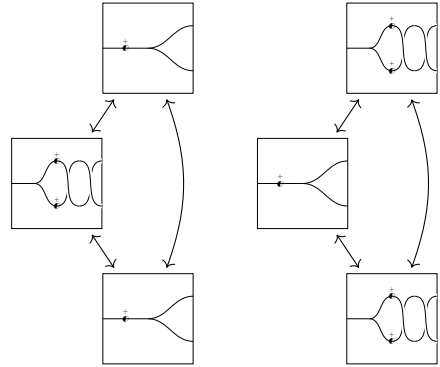
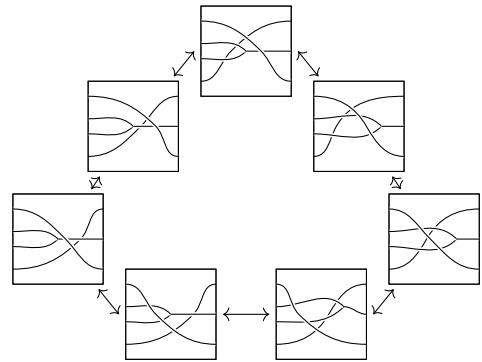
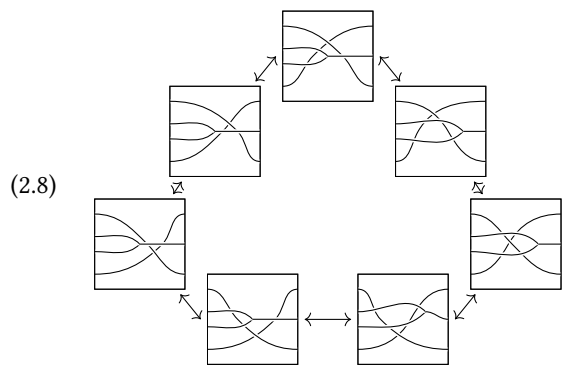
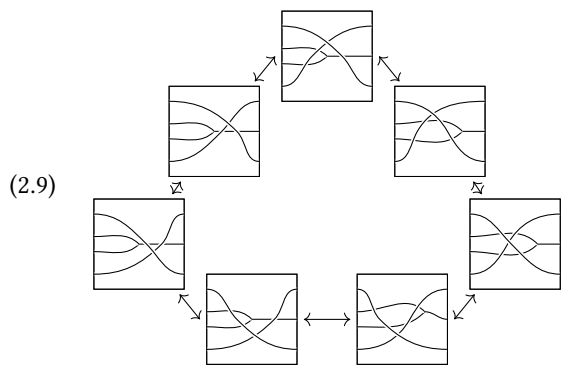
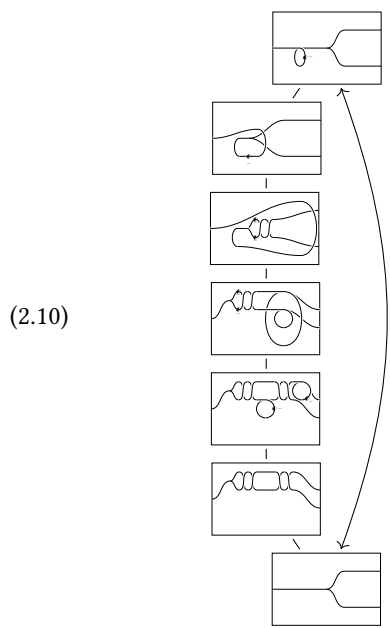
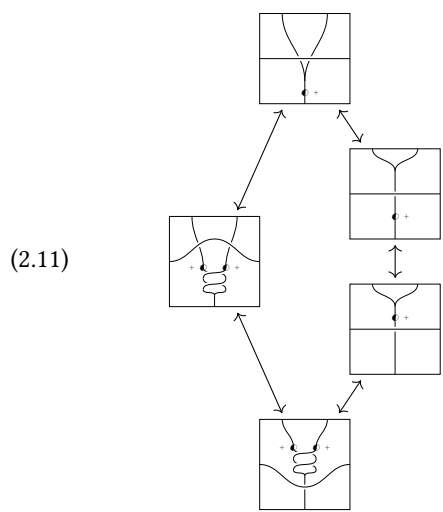
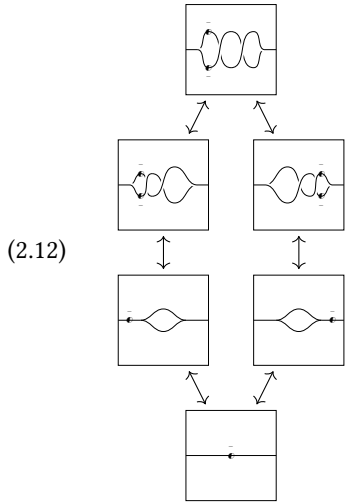
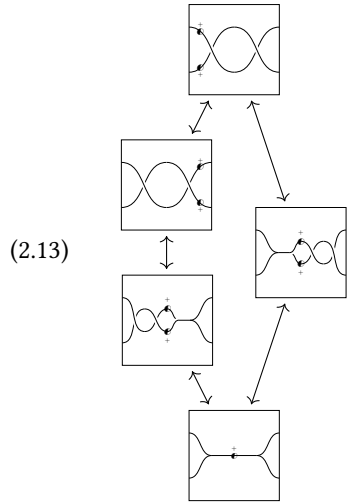
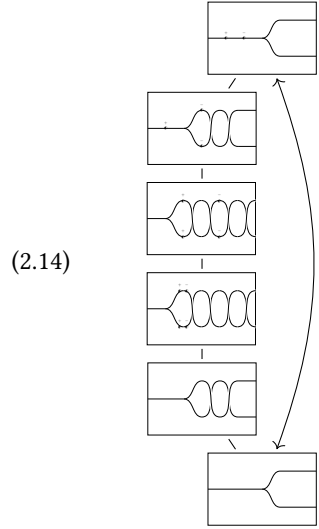
\[
  %
\]

\noindent {\bf Pure framing moves}

\[
  %
\]

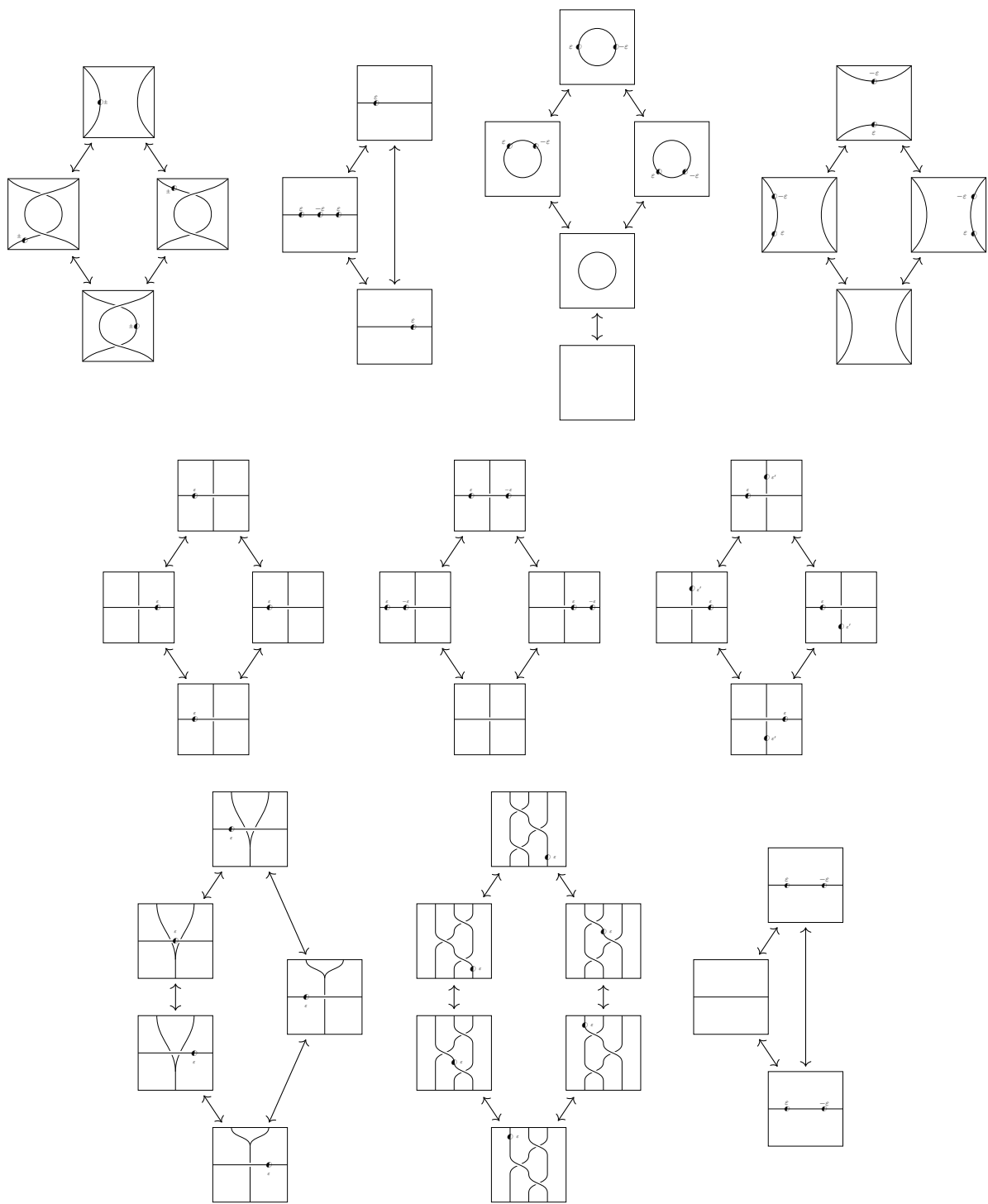
\noindent {\bf Invertibility of new generators}

\[
  %
  };
  \draw [-] (A) -- (B1);
  \draw [-] (B1) -- (C1);
  \draw [-] (C1) -- (D1);
  \draw [-] (D1) -- (E1);
  \draw [-] (E1) -- (F);
  \draw [<->] (A) to [out=-70,in=70] (F);
\end{tikzpicture}
\end{equation}
\end{minipage}
\]

\end{theorem}

\subsection{Over surfaces}

The analysis from~\cite{QWalker} can be ran without any change for webs embedded in a thickened surface $\Su\times [0,1]$ (with projection down to the surface). Indeed, the whole process is run locally. Thus, Theorem \ref{thm:MM} holds for framed tangled foams embedded in $\Su\times [0,1]\times [0,1]$.

\begin{definition}
  We denote by $\Stanwebo$ the category of tangled webs and foams over the surface $S$ presented by the previous theorem.
\end{definition}

\section{Full functoriality for  Khovanov--Rozansky homology}

\subsection{\texorpdfstring{$\slnn{2}$}{sl2}-foams}

The version of Khovanov homology we will use is based on Blanchet's $\glnn{2}$-foams~\cite{Blan}, that are obtained by considering foams as above modulo the following relations. See for example~\cite{LQR} for more historical details. Many of the pictures below are extracted from this paper or from~\cite{QW_SkeinCat}.

\begin{equation}\label{sl2closedfoam}

\quad .
\end{equation}
The neck-cutting relation in  \eqref{sl2closedfoam} implies that a $1$-handle on a $1$-labeled facet simplifies to twice a dot. 

\begin{equation} \label{sl2neckcutting_enh_2lab}
%
\end{equation}

The following relations are consequences of the previous ones and will prove useful.

\begin{equation} \label{sl2bubble12}
\xy
{0,0}*{
%
};
\endxy
\end{equation}

The relations below, valid for the Lee-Rasmussen deformation (see~\cite{Lee,Ras} for the original definition, and \cite{Blan} for the extension to foams), will also prove useful in several proofs. In a nutshell, the Lee-Rasmussen deformation is a variant of Khovanov homology that appears as the last page of a spectral sequence whose second page is the ordinary Khovanov homology. When trying to compare coefficients of some maps, it is often useful to go all the way to the last page: if both maps survive, it is often much easier to compare them then.

The new rule for dots is: $\bullet^2=1$, and facets of the foam can now be decorated by idempotents $p_+=\frac{1+\bullet}{2}$ and $p_{-}=\frac{1-\bullet}{2}$. All relations above remain valid, and one can deduce from them and from the rule $\bullet^2=1$ the following additional relations. A facet decorated with a $+$ stands for a facet carrying the idempotent $p_+$, and similarly with minus signs. To help distinguish between classical sheets and sheets decorated by idempotents, the latter are often colored blue.

\begin{gather}
  \label{deformedbb1}

  \end{equation}

\subsection{Khovanov homology}

In this section, we want to prove that Khovanov homology induces a functor from the category of framed tangled webs and foams to the homotopy category of framed webs and foams:
\[
  \Su\Kh\colon \Stanwebo \to \HC(\Sfoam)
\]
This will answer (positively) Conjecture 4.8 from \cite{QW_SkeinCat}.

We first define the map, by assigning locally to crossings pieces of complexes described below. This definition is very close to Blanchet's one \cite{Blan} and includes in the larger framework of foam-based exterior colored homologies (see \cite{KhR,MSV,QR} and references therein).

\begin{gather}
  \Kh\left( \;
    \begin{tikzpicture}[anchorbase, scale=.5]
      \draw [very thick, ->] (2,1) to [out=180,in=0] (0,0);
      \draw [white,line width=.15cm] (2,0) to [out=180,in=0] (0,1) ;
      \draw [very thick, ->] (2,0) to [out=180,in=0] (0,1);
    \end{tikzpicture}
    \;\right)
  \;\;= \;\;
  \begin{tikzpicture}[anchorbase, scale=.5]
    \draw [very thick, ->] (2,1) to (0,1);
    \draw [very thick, ->] (2,0) to (0,0);
  \end{tikzpicture}
  \;\;\to\;\;
  q^{-1} t \;
  \begin{tikzpicture}[anchorbase, scale=.5]
    \draw [very thick] (2,0) to[out=180,in=315] (1.3,.5);
    \draw [very thick] (2,1) to[out=180,in=45] (1.3,.5);
    \draw [double] (1.3,.5) -- (.7,.5);
    \draw [very thick, ->] (.7,.5) to[out=135,in=0]  (0,1);
    \draw [very thick, ->] (.7,.5) to[out=225,in=0] (0,0);
  \end{tikzpicture}
  \quad,\quad
  \Kh\left( \;
    \begin{tikzpicture}[anchorbase, scale=.5]
      \draw [very thick, ->] (2,0) to [out=180,in=0] (0,1);
      \draw [white,line width=.15cm] (2,1) to [out=180,in=0] (0,0) ;
      \draw [very thick, ->] (2,1) to [out=180,in=0] (0,0);
    \end{tikzpicture}
    \;\right)
  \;\;=\;\;
  q t^{-1}  \;
  \begin{tikzpicture}[anchorbase, scale=.5]
    \draw [very thick] (2,0) to[out=180,in=315] (1.3,.5);
    \draw [very thick] (2,1) to[out=180,in=45] (1.3,.5);
    \draw [double] (1.3,.5) -- (.7,.5);
    \draw [very thick, ->] (.7,.5) to[out=135,in=0]  (0,1);
    \draw [very thick, ->] (.7,.5) to[out=225,in=0] (0,0);
  \end{tikzpicture}
  \;\;\to\;\;
  \begin{tikzpicture}[anchorbase, scale=.5]
    \draw [very thick, ->] (2,1) to (0,1);
    \draw [very thick, ->] (2,0) to (0,0);
  \end{tikzpicture}\nonumber
  \\
  \label{eq:thickcrossing-2}
  \Kh\left(\begin{tikzpicture}[anchorbase, scale=.5]
      \draw [very thick, ->] (2,1) to [out=180,in=0] (0,0);
      \draw [white,line width=.15cm] (2,0) to [out=180,in=0] (0,1) ;
      \draw [double, ->] (2,0) to [out=180,in=0] (0,1);
    \end{tikzpicture}
  \right)
  = 
  q^{-1} t \;
  \begin{tikzpicture}[anchorbase, scale=.5]
    \draw [double] (2,0) -- (1.4,0);
    \draw [very thick, ->] (1.4,0) -- (0,0);
    \draw [very thick] (2,1) -- (.6,1);
    \draw [double, ->] (0.6,1) -- (0,1);
    \draw [very thick] (.6,1) -- (1.4,0);
  \end{tikzpicture}
  \;\;,\;\;
  \Kh\left(\begin{tikzpicture}[anchorbase, scale=.5]
      \draw [double, ->] (2,1) to [out=180,in=0] (0,0);
      \draw [white,line width=.15cm] (2,0) to [out=180,in=0] (0,1) ;
      \draw [very thick, ->] (2,0) to [out=180,in=0] (0,1);
    \end{tikzpicture}
  \right)
  = 
  q^{-1} t \;
  \begin{tikzpicture}[anchorbase, scale=.5]
    \draw [double] (2,1) -- (1.4,1);
    \draw [very thick, ->] (1.4,1) -- (0,1);
    \draw [very thick] (2,0) -- (.6,0);
    \draw [double, ->] (0.6,0) -- (0,0);
    \draw [very thick] (.6,0) -- (1.4,1);
  \end{tikzpicture}
  \;\;,\;\;
  \Kh\left(\begin{tikzpicture}[anchorbase, scale=.5]
      \draw [double, ->] (2,1) to [out=180,in=0] (0,0);
      \draw [white,line width=.15cm] (2,0) to [out=180,in=0] (0,1) ;
      \draw [double, ->] (2,0) to [out=180,in=0] (0,1);
    \end{tikzpicture}
  \right)
  = q^{-2} t^{2} \;
  \begin{tikzpicture}[anchorbase, scale=.5]
    \draw [double, ->] (2,0) -- (0,0);
    \draw [double, ->] (2,1) -- (0,1);
  \end{tikzpicture}
  \\ \nonumber
\Kh\left( 
    \begin{tikzpicture}[anchorbase, scale=.5]
      \draw [very thick, ->] (2,0) to [out=180,in=0] (0,1);
      \draw [white,line width=.15cm] (2,1) to [out=180,in=0] (0,0) ;
      \draw [double, ->] (2,1) to [out=180,in=0] (0,0);
    \end{tikzpicture}
  \right)
  =
  q t^{-1}  \;
  \begin{tikzpicture}[anchorbase, scale=.5]
    \draw [double] (2,1) -- (1.4,1);
    \draw [very thick, ->] (1.4,1) -- (0,1);
    \draw [very thick] (2,0) -- (.6,0);
    \draw [double, ->] (0.6,0) -- (0,0);
    \draw [very thick] (.6,0) -- (1.4,1);
  \end{tikzpicture}
  \;\;,\;\;
  \Kh\left(
    \begin{tikzpicture}[anchorbase, scale=.5]
      \draw [double, ->] (2,0) to [out=180,in=0] (0,1);
      \draw [white,line width=.15cm] (2,1) to [out=180,in=0] (0,0) ;
      \draw [very thick, ->] (2,1) to [out=180,in=0] (0,0);
    \end{tikzpicture}
  \right)
  = 
  q t^{-1} \;
  \begin{tikzpicture}[anchorbase, scale=.5]
    \draw [double] (2,0) -- (1.4,0);
    \draw [very thick, ->] (1.4,0) -- (0,0);
    \draw [very thick] (2,1) -- (.6,1);
    \draw [double, ->] (0.6,1) -- (0,1);
    \draw [very thick] (.6,1) -- (1.4,0);
  \end{tikzpicture}
  \;\;,\;\;
  \Kh\left( 
    \begin{tikzpicture}[anchorbase, scale=.5]
      \draw [double, ->] (2,0) to [out=180,in=0] (0,1);
      \draw [white,line width=.15cm] (2,1) to [out=180,in=0] (0,0) ;
      \draw [double, ->] (2,1) to [out=180,in=0] (0,0);
    \end{tikzpicture}
  \right)
  = q^{2} t^{-2} \;
  \begin{tikzpicture}[anchorbase, scale=.5]
    \draw [double, ->] (2,0) -- (0,0);
    \draw [double, ->] (2,1) -- (0,1);
  \end{tikzpicture}
\end{gather}

Above, $q$ stands for the quantum degree while $t$ is the homological one.

The twists are a part that is usually not considered. Here are the rules we adopt for them. They simply consist of shifts that allow for the $R_I$ moves to be grading-preserving. 
\begin{gather*}
\Kh\left(
    \begin{tikzpicture}[anchorbase,scale=.5]
      \draw [very thick, ->] (2,0) -- (0,0);
      \node at (1,0) {\small $\LEFTcircle$};
      \node at (1,.4) {\tiny $+$};
      \node at (1,-.4) { \vphantom{\tiny $+$}};
      \end{tikzpicture}
    \right)
    =
    q\;
        \begin{tikzpicture}[anchorbase,scale=.5]
      \draw [very thick, ->] (2,0) -- (0,0);
    \end{tikzpicture}
    \;\;,\;\;
      \Kh\left(
    \begin{tikzpicture}[anchorbase,scale=.5]
      \draw [very thick, ->] (2,0) -- (0,0);
      \node at (1,0) {\small $\LEFTcircle$};
      \node at (1,.4) {\tiny $-$};
      \node at (1,-.4) { \vphantom{\tiny $+$}};
      \end{tikzpicture}
    \right)
    =
    q^{-1}\;
        \begin{tikzpicture}[anchorbase,scale=.5]
      \draw [very thick, ->] (2,0) -- (0,0);
    \end{tikzpicture}
\\
  \Kh\left(
    \begin{tikzpicture}[anchorbase,scale=.5]
      \draw [double, ->] (2,0) -- (0,0);
      \node at (1,0) {\small $\LEFTcircle$};
      \node at (1,.4) {\tiny $+$};
      \node at (1,-.4) {\vphantom{\tiny $+$}};
      \end{tikzpicture}
    \right)
    =
    q^{-2}t^{2}\;
        \begin{tikzpicture}[anchorbase,scale=.5]
      \draw [double, ->] (2,0) -- (0,0);
    \end{tikzpicture}
    \;\;,\;\;
  \Kh\left(
    \begin{tikzpicture}[anchorbase,scale=.5]
      \draw [double, ->] (2,0) -- (0,0);
      \node at (1,0) {\small $\LEFTcircle$};
      \node at (1,.4) {\tiny $-$};
      \node at (1,-.4) {\vphantom{\tiny $+$}};
      \end{tikzpicture}    \right)
    =
    q^{2}t^{-2}\;
        \begin{tikzpicture}[anchorbase,scale=.5]
      \draw [double, ->] (2,0) -- (0,0);
    \end{tikzpicture}
 \end{gather*}

To extend these local rules to a knot diagram, we use Blanchet's process and consider that terms in the complex are tensor products of the web depicted above with $\C\simeq \bigwedge^{N_1}(\otimes_{c \in S_1} \C c_i)$. $S_1$ stands for the set of crossings that involve two 1-labeled strands and are smoothed as a split-merge web, or crossings involving a $1$ and a $2$-labeled strand (notice that these are exactly the crossings that yield an odd homological shift on some of the terms of the complex). The differentials in the complexes above are thus zip or unzip cobordisms for the web part, together with $\wedge c$ or $\langle -,c\rangle$ depending whether one increases or decreases the cardinality of $S_1$. We adopt the following convention for the evaluation map:
\[
\langle x\wedge c,c\rangle=x\quad \text{and more generally}\quad \langle x\wedge c,y\wedge c\rangle=\langle x,y\rangle.
\]
For example:
\[
  \langle c_1\wedge c_2,c_1\wedge c_2\rangle =1.
\]

In order to make Khovanov homology into a functor, we want to assign to elementary cobordisms between tangled webs preferred maps between the source and target complexes. We build upon~\cite{Kh1,Kh2,BN2,Blan,ETW}, though surprising differences occur: for example, it seems that the normalizations for $R_3$ moves that appear in~\cite{ETW} are not suitable for tangled webs.

Caps, cups, zips, unzips and saddle cobordisms are just naively sent onto themselves. For Reidemeister moves and forkslide moves, we refer to the list given in Appendix~\ref{app:maps}. We insist that these conventions agree with the ones from~\cite{Blan} and with those from~\cite{ETW} to the exception of the Reidemeister 3 moves involving two double strands. This will be discussed in the proof of Theorem~\ref{thm:functoriality}.

Continuing in Equation~\eqref{eq:thmMM3}, sliding a twist through a crossing is assigned an identity cobordism, and canceling a pair of opposite twists is also assigned an identity cobordism. Finally, it remains to explain our conventions for a twist passing through a trivalent vertex.

It relies on the following lemma. 

\begin{lemma} \label{lem:trivalenttwist}
  The following maps form a homotopy equivalence.
  
  \[

\]
\end{lemma}

Forgetting the $2$-labeled part above, one sees two $R_1$ moves being applied successively, giving inspiration for the shape of the above map.

\begin{proof}
  It suffices to exhibit homotopies:

\[
  h^1_{sw}=0,\quad h^1_{nw}=
      %
      \otimes \langle \cdot,c_1\rangle
    \]
  \end{proof}

  One can work out similar homotopies with orientation reversed, or with the other sign for the twists. Then for each of these generating moves, there is a choice of coefficient $\pm 1$. We adopt the conventions listed in the Section~\ref{sec:trivtwists}.

We now come to the central result of this section.

\begin{theorem}
  \label{thm:functoriality}
  Khovanov homology induces a well-defined functor:
  \[
      \Su\Kh\colon \Stanwebo \to \HC(\Sfoam)
  \]
\end{theorem}

\begin{proof}
  The proof goes by checking that moves from Theorem~\ref{thm:MM} are satisfied. The choice of generator for the framed version of $R_1$ is compatible with Blanchet's choice, and we preserved backtrack compatibility with the conventions of~\cite{ETW}, to the exception of Reidemeister 3 moves involving two double strands, where we introduce an extra sign. These moves always come in pairs however in classical movie moves, so the results are unchanged. This ensures that the moves $\rm{MM}_1$ to $\rm{MM}_{15}$ all hold.

  We thus turn our attention to the new moves.

Crossingless moves hold trivially by construction, as well as pure framing moves given our choice for the image of twists.

Forkslide invertibility is built in from the choice of representatives for these moves, while the invertibility of a twist moving through a trivalent vertex is ensured by Lemma~\ref{lem:trivalenttwist}.

\noindent {\bf Strand around vertex~\eqref{StrandAroundVertex}}:  We start with the check for \eqref{StrandAroundVertex}, as this move will serve later on to reduce the number of checks for other moves.

  We first consider the case of a 2-labeled single strand. There are a priori 8 situations to consider, and for each of them the single strand can pass over or under the trivalent vertex:
  
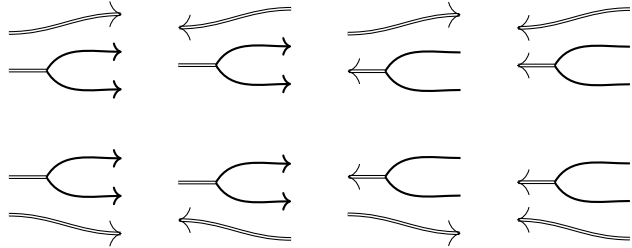
\begin{gather*}
\begin{tikzpicture}[anchorbase,scale=.5]
  \draw [double, ->] (0,1.5) to [out=0,in=180]  (3,2);
  \draw [double] (0,.5) to [out=0,in=180] (1,.5);
  \draw [thick,->] (1,.5) to [out=60,in=180] (3,1);
  \draw [thick,->] (1,.5) to [out=-60,in=180] (3,0);
\end{tikzpicture}
\qquad
\begin{tikzpicture}[anchorbase,scale=.5]
  \draw [double, <-] (0,1.5) to [out=0,in=180]  (3,2);
  \draw [double] (0,.5) to [out=0,in=180] (1,.5);
  \draw [thick,->] (1,.5) to [out=60,in=180] (3,1);
  \draw [thick,->] (1,.5) to [out=-60,in=180] (3,0);
\end{tikzpicture}
\qquad
\begin{tikzpicture}[anchorbase,scale=.5]
  \draw [double, ->] (0,1.5) to [out=0,in=180]  (3,2);
  \draw [double,<-] (0,.5) to [out=0,in=180] (1,.5);
  \draw [thick] (1,.5) to [out=60,in=180] (3,1);
  \draw [thick] (1,.5) to [out=-60,in=180] (3,0);
\end{tikzpicture}
\qquad
\begin{tikzpicture}[anchorbase,scale=.5]
  \draw [double, <-] (0,1.5) to [out=0,in=180]  (3,2);
  \draw [double,<-] (0,.5) to [out=0,in=180] (1,.5);
  \draw [thick] (1,.5) to [out=60,in=180] (3,1);
  \draw [thick] (1,.5) to [out=-60,in=180] (3,0);
\end{tikzpicture}
\\
\vspace{.5cm} \\
\begin{tikzpicture}[anchorbase,scale=.5]
  \draw [double, ->] (0,.5) to [out=0,in=180]  (3,0);
  \draw [double] (0,1.5) to [out=0,in=180] (1,1.5);
  \draw [thick,->] (1,1.5) to [out=60,in=180] (3,2);
  \draw [thick,->] (1,1.5) to [out=-60,in=180] (3,1);
\end{tikzpicture}
\qquad
\begin{tikzpicture}[anchorbase,scale=.5]
  \draw [double, <-] (0,.5) to [out=0,in=180]  (3,0);
  \draw [double] (0,1.5) to [out=0,in=180] (1,1.5);
  \draw [thick,->] (1,1.5) to [out=60,in=180] (3,2);
  \draw [thick,->] (1,1.5) to [out=-60,in=180] (3,1);
\end{tikzpicture}
\qquad
\begin{tikzpicture}[anchorbase,scale=.5]
  \draw [double, ->] (0,.5) to [out=0,in=180]  (3,0);
  \draw [double,<-] (0,1.5) to [out=0,in=180] (1,1.5);
  \draw [thick] (1,1.5) to [out=60,in=180] (3,2);
  \draw [thick] (1,1.5) to [out=-60,in=180] (3,1);
\end{tikzpicture}
\qquad
\begin{tikzpicture}[anchorbase,scale=.5]
  \draw [double, <-] (0,.5) to [out=0,in=180]  (3,0);
  \draw [double,<-] (0,1.5) to [out=0,in=180] (1,1.5);
  \draw [thick] (1,1.5) to [out=60,in=180] (3,2);
  \draw [thick] (1,1.5) to [out=-60,in=180] (3,1);
\end{tikzpicture}
\end{gather*}

Starting with the first case, one computes for the left-hand side of the move, using~\eqref{R214} and~\eqref{FS14} :
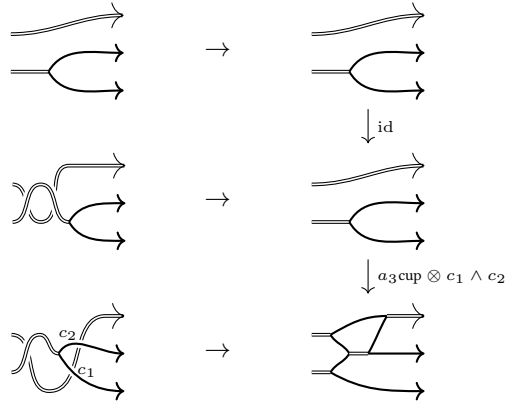
\begin{gather*}
  \begin{tikzpicture}[anchorbase]
    \node (topleft) at (0,0) {
      \begin{tikzpicture}[anchorbase,scale=.5]
        \draw [double, ->] (0,1.5) to [out=0,in=180]  (3,2);
        \draw [double] (0,.5) to [out=0,in=180] (1,.5);
        \draw [thick,->] (1,.5) to [out=60,in=180] (3,1);
        \draw [thick,->] (1,.5) to [out=-60,in=180] (3,0);
      \end{tikzpicture}
    };
    \node at (2,0) {$\rightarrow$};
    \node (topright) at (4,0) {
      \begin{tikzpicture}[anchorbase,scale=.5]
        \draw [double, ->] (0,1.5) to [out=0,in=180]  (3,2);
        \draw [double] (0,.5) to [out=0,in=180] (1,.5);
        \draw [thick,->] (1,.5) to [out=60,in=180] (3,1);
        \draw [thick,->] (1,.5) to [out=-60,in=180] (3,0);
      \end{tikzpicture}
    };
        \node (midleft) at (0,-2) {
      \begin{tikzpicture}[anchorbase,scale=.5]
        \draw [double, ->] (0,1.5) to [out=0,in=180] (.75,.5) to [out=0,in=180] (1.5,2) --  (3,2);
        \draw [white, line width=4] (0,.5) to [out=0,in=180] (.75,1.5) to [out=0,in=180] (1.5,.5);
        \draw [double] (0,.5) to [out=0,in=180] (.75,1.5) to [out=0,in=180] (1.5,.5);
        \draw [thick,->] (1.5,.5) to [out=60,in=180] (3,1);
        \draw [thick,->] (1.5,.5) to [out=-60,in=180] (3,0);
      \end{tikzpicture}
    };
    \node at (2,-2) {$\rightarrow$};
    \node (midright) at (4,-2) {
      \begin{tikzpicture}[anchorbase,scale=.5]
        \draw [double, ->] (0,1.5) to [out=0,in=180]  (3,2);
        \draw [double] (0,.5) to [out=0,in=180] (1,.5);
        \draw [thick,->] (1,.5) to [out=60,in=180] (3,1);
        \draw [thick,->] (1,.5) to [out=-60,in=180] (3,0);
      \end{tikzpicture}
    };
        \node (bottomleft) at (0,-4) {
      \begin{tikzpicture}[anchorbase,scale=.5]
        \draw [double, ->] (0,1.5) to [out=0,in=180] (1,0) to [out=0,in=180] (2.5,2) --  (3,2);
        \draw [white, line width=4] (0,.5) to [out=0,in=180] (.75,1.5) to [out=0,in=180] (1.25,1);
        \draw [white, line width=3] (1.25,1) to [out=60,in=180] (3,1);
        \draw [white, line width=3] (1.25,1) to [out=-60,in=180] (3,0);
        \draw [double] (0,.5) to [out=0,in=180] (.75,1.5) to [out=0,in=180] (1.25,1);
        \draw [thick,->] (1.25,1) to [out=60,in=180] (3,1);
        \draw [thick,->] (1.25,1) to [out=-60,in=180] (3,0);
        \node at (2,.5) {\tiny $c_1$};
        \node at (1.5,1.5) {\tiny $c_2$};
      \end{tikzpicture}
    };
    \node at (2,-4) {$\rightarrow$};
    \node (bottomright) at (4,-4) {
      \begin{tikzpicture}[anchorbase,scale=.5]
        \draw [double] (0,1.5) to (.5,1.5);
        \draw [double] (0,.5) to (.5,.5);
        \draw [thick,->] (.5,.5) to [out=-30,in=180] (3,0);
        \draw [thick] (.5,.5) to [out=60,in=-120] (1,1);
        \draw [thick] (.5,1.5) to [out=-60,in=120] (1,1);
        \draw [double] (1,1) -- (1.5,1);
        \draw [thick] (1.5,1) to [out=60,in=-120] (2,2);
        \draw [thick] (.5,1.5) to [out=30,in=180] (2,2);
        \draw [thick,->] (1.5,1) -- (3,1);
        \draw [double,->] (2,2) -- (3,2);
      \end{tikzpicture}
    };
    \draw [->] (topright) -- (midright) node[midway,right] {\tiny $\id$};
    \draw [->] (midright) -- (bottomright) node[midway,right] {\tiny $a_3\text{cup}\otimes c_1\wedge c_2$};
  \end{tikzpicture}
\end{gather*}

The right-hand side translates into:
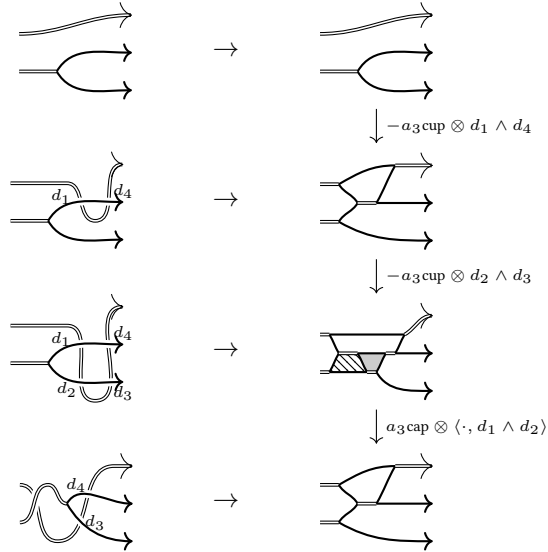
\begin{gather*}
  \begin{tikzpicture}[anchorbase]
    \node (topleft) at (0,0) {
      \begin{tikzpicture}[anchorbase,scale=.5]
        \draw [double, ->] (0,1.5) to [out=0,in=180]  (3,2);
        \draw [double] (0,.5) to [out=0,in=180] (1,.5);
        \draw [thick,->] (1,.5) to [out=60,in=180] (3,1);
        \draw [thick,->] (1,.5) to [out=-60,in=180] (3,0);
      \end{tikzpicture}
    };
    \node at (2,0) {$\rightarrow$};
    \node (topright) at (4,0) {
      \begin{tikzpicture}[anchorbase,scale=.5]
        \draw [double, ->] (0,1.5) to [out=0,in=180]  (3,2);
        \draw [double] (0,.5) to [out=0,in=180] (1,.5);
        \draw [thick,->] (1,.5) to [out=60,in=180] (3,1);
        \draw [thick,->] (1,.5) to [out=-60,in=180] (3,0);
      \end{tikzpicture}
    };
        \node (midleft+) at (0,-2) {
      \begin{tikzpicture}[anchorbase,scale=.5]
        \draw [double, ->] (0,1.5) -- (1.5,1.5) to [out=0,in=180] (2.25,.5) to [out=0,in=180]  (3,2);
        \draw [white, line width=3] (1,.5) to [out=60,in=180] (3,1);
        \draw [white, line width=3] (1,.5) to [out=-60,in=180] (3,0);
        \draw [double] (0,.5) -- (1,.5);
        \draw [thick,->] (1,.5) to [out=60,in=180] (3,1);
        \draw [thick,->] (1,.5) to [out=-60,in=180] (3,0);
        \node at (1.35,1.15) {\tiny $d_1$};
        \node at (3,1.3) {\tiny $d_4$};
      \end{tikzpicture}
    };
    \node at (2,-2) {$\rightarrow$};
    \node (midright+) at (4,-2) {
      \begin{tikzpicture}[anchorbase,scale=.5]
        \draw [double] (0,1.5) to (.5,1.5);
        \draw [double] (0,.5) to (.5,.5);
        \draw [thick,->] (.5,.5) to [out=-30,in=180] (3,0);
        \draw [thick] (.5,.5) to [out=60,in=-120] (1,1);
        \draw [thick] (.5,1.5) to [out=-60,in=120] (1,1);
        \draw [double] (1,1) -- (1.5,1);
        \draw [thick] (1.5,1) to [out=60,in=-120] (2,2);
        \draw [thick] (.5,1.5) to [out=30,in=180] (2,2);
        \draw [thick,->] (1.5,1) -- (3,1);
        \draw [double,->] (2,2) -- (3,2);
      \end{tikzpicture}
    };
        \node (midleft-) at (0,-4) {
      \begin{tikzpicture}[anchorbase,scale=.5]
        \draw [double, ->] (0,1.5) -- (1.5,1.5) to [out=0,in=180] (2.25,-.5) to [out=0,in=180]  (3,2);
        \draw [white, line width=3] (1,.5) to [out=60,in=180] (3,1);
        \draw [white, line width=3] (1,.5) to [out=-60,in=180] (3,0);
        \draw [double] (0,.5) -- (1,.5);
        \draw [thick,->] (1,.5) to [out=60,in=180] (3,1);
        \draw [thick,->] (1,.5) to [out=-60,in=180] (3,0);
        \node at (1.35,1.15) {\tiny $d_1$};
        \node at (1.5,-.2) {\tiny $d_2$};
        \node at (3,-.3) {\tiny $d_3$};
        \node at (3,1.3) {\tiny $d_4$};
      \end{tikzpicture}
    };
    \node at (2,-4) {$\rightarrow$};
    \node (midright-) at (4,-4) {
      \begin{tikzpicture}[anchorbase,scale=.5]
        \fill [fill=black, opacity=.2] (1.25,.5) -- (1.5,.5) to [out=60,in=-120] (1.75,1) -- (1,1) to [out=-60,in=120] (1.25,.5);
        \fill[pattern=north west lines] (.25,.5) -- (1.25,.5) to [out=120,in=-60] (1,1) -- (.5,1) to [out=-120,in=60] (.25,.5);
        \draw [double] (0,1.5) to (.25,1.5);
        \draw [double] (0,.5) to (.25,.5);
        \draw [thick] (.25,.5) to [out=60,in=-120] (.5,1);
        \draw [thick] (.25,1.5) to [out=-60,in=120] (.5,1);
        \draw [double] (.5,1) -- (1,1);
        \draw [thick] (1,1) to [out=-60,in=120] (1.25,.5);
        \draw [thick] (.25,.5) -- (1.25,.5);
        \draw [double] (1.25,.5) -- (1.5,.5);
        \draw [thick] (1.5,.5) to [out=60,in=-120] (1.75,1);
        \draw [thick] (1,1) -- (1.75,1);
        \draw [double] (1.75,1) -- (2,1);
        \draw [thick, ->] (1.5,.5) to [out=-60,in=180] (3,0);
        \draw [thick,->] (2,1) -- (3,1);
        \draw [thick] (2,1) to [out=60,in=-120] (2.25,1.5);
        \draw [thick] (.25,1.5) -- (2.25,1.5);
        \draw [double,->] (2.25,1.5) to [out=30,in=180] (3,2);
      \end{tikzpicture}
    };
        \node (bottomleft) at (0,-6) {
      \begin{tikzpicture}[anchorbase,scale=.5]
        \draw [double, ->] (0,1.5) to [out=0,in=180] (1,0) to [out=0,in=180] (2.5,2) --  (3,2);
        \draw [white, line width=4] (0,.5) to [out=0,in=180] (.75,1.5) to [out=0,in=180] (1.25,1);
        \draw [white, line width=3] (1.25,1) to [out=60,in=180] (3,1);
        \draw [white, line width=3] (1.25,1) to [out=-60,in=180] (3,0);
        \draw [double] (0,.5) to [out=0,in=180] (.75,1.5) to [out=0,in=180] (1.25,1);
        \draw [thick,->] (1.25,1) to [out=60,in=180] (3,1);
        \draw [thick,->] (1.25,1) to [out=-60,in=180] (3,0);
        \node at (2,.5) {\tiny $d_3$};
        \node at (1.5,1.5) {\tiny $d_4$};
      \end{tikzpicture}
    };
    \node at (2,-6) {$\rightarrow$};
    \node (bottomright) at (4,-6) {
      \begin{tikzpicture}[anchorbase,scale=.5]
        \draw [double] (0,1.5) to (.5,1.5);
        \draw [double] (0,.5) to (.5,.5);
        \draw [thick,->] (.5,.5) to [out=-30,in=180] (3,0);
        \draw [thick] (.5,.5) to [out=60,in=-120] (1,1);
        \draw [thick] (.5,1.5) to [out=-60,in=120] (1,1);
        \draw [double] (1,1) -- (1.5,1);
        \draw [thick] (1.5,1) to [out=60,in=-120] (2,2);
        \draw [thick] (.5,1.5) to [out=30,in=180] (2,2);
        \draw [thick,->] (1.5,1) -- (3,1);
        \draw [double,->] (2,2) -- (3,2);
      \end{tikzpicture}
    };
    \draw [->] (topright) -- (midright+) node[midway,right] {\tiny $-a_3\text{cup}\otimes d_1\wedge d_4$};
    \draw [->] (midright+) -- (midright-) node[midway,right] {\tiny $-a_3\text{cup}\otimes d_2\wedge d_3$};
    \draw [->] (midright-) -- (bottomright) node[midway,right] {\tiny $a_3\text{cap}\otimes \langle \cdot, d_1\wedge d_2\rangle$};
  \end{tikzpicture}
\end{gather*}

Above, the second cup is the one that creates the shaded rectangle, while the hatched rectangle is capped of at step 3. The foam induced by the composition of these two steps is simply an identity foam. The composition of all three steps thus induces the identity foam, with coefficient given by:
\[(-a_3)^2a_3\langle d1\wedge d_4\wedge d_2\wedge d_3,d_1\wedge d_2\rangle = -a_3d_4\wedge d_3=a_3d_3\wedge d_4. \]
Under the identification of $c_1$ with $d_3$ and $c_2$ with $d_4$, this matches the previous computation and proves that the move holds.

The same version of the move with the 2-labeled free strand passing over the trivalent vertex holds for a similar reason: on the level of foams, the computation goes unchanged. One can notice that the first step on the LHS always induces an identity, and that steps 1 and 2 on the RHS will always produce the same coefficients, that square to 1. Thus one can focus only on forkslide moves and track coefficients. Actually, the same goes true for all of the following moves:
  
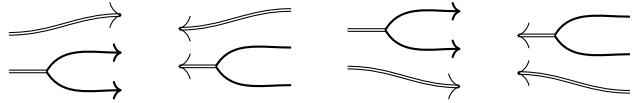
\begin{gather*}
\begin{tikzpicture}[anchorbase,scale=.5]
  \draw [double, ->] (0,1.5) to [out=0,in=180]  (3,2);
  \draw [double] (0,.5) to [out=0,in=180] (1,.5);
  \draw [thick,->] (1,.5) to [out=60,in=180] (3,1);
  \draw [thick,->] (1,.5) to [out=-60,in=180] (3,0);
\end{tikzpicture}
\qquad
\begin{tikzpicture}[anchorbase,scale=.5]
  \draw [double, <-] (0,1.5) to [out=0,in=180]  (3,2);
  \draw [double,<-] (0,.5) to [out=0,in=180] (1,.5);
  \draw [thick] (1,.5) to [out=60,in=180] (3,1);
  \draw [thick] (1,.5) to [out=-60,in=180] (3,0);
\end{tikzpicture}
\qquad
\begin{tikzpicture}[anchorbase,scale=.5]
  \draw [double, ->] (0,.5) to [out=0,in=180]  (3,0);
  \draw [double] (0,1.5) to [out=0,in=180] (1,1.5);
  \draw [thick,->] (1,1.5) to [out=60,in=180] (3,2);
  \draw [thick,->] (1,1.5) to [out=-60,in=180] (3,1);
\end{tikzpicture}
\qquad
\begin{tikzpicture}[anchorbase,scale=.5]
  \draw [double, <-] (0,.5) to [out=0,in=180]  (3,0);
  \draw [double,<-] (0,1.5) to [out=0,in=180] (1,1.5);
  \draw [thick] (1,1.5) to [out=60,in=180] (3,2);
  \draw [thick] (1,1.5) to [out=-60,in=180] (3,1);
\end{tikzpicture}
\end{gather*}
The above moves hold if the following statement is true: two forkslide moves from \eqref{FS5}--\eqref{FS8} and \eqref{FS13}--\eqref{FS16} that are identical but for the orientation of the 2-strand come with opposite constants. This happens to be true in all cases. 

Thus we turn our attention to:
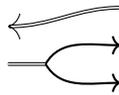
\[
\begin{tikzpicture}[anchorbase,scale=.5]
  \draw [double, <-] (0,1.5) to [out=0,in=180]  (3,2);
  \draw [double] (0,.5) to [out=0,in=180] (1,.5);
  \draw [thick,->] (1,.5) to [out=60,in=180] (3,1);
  \draw [thick,->] (1,.5) to [out=-60,in=180] (3,0);
\end{tikzpicture}
\]
We start again with the case of the $2$-labeled strand passing under the trivalent vertex.
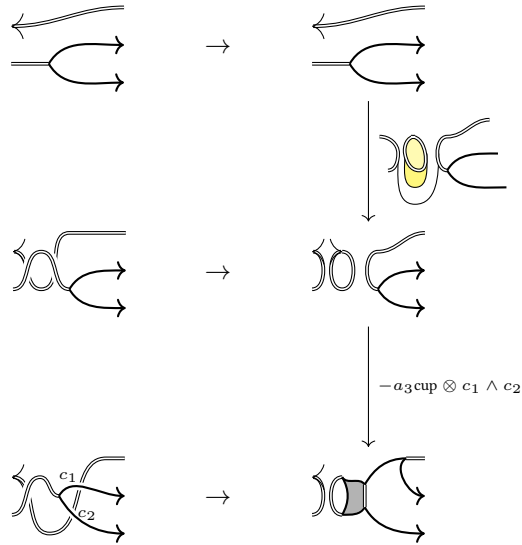
\begin{gather*}
  \begin{tikzpicture}[anchorbase]
    \node (topleft) at (0,0) {
      \begin{tikzpicture}[anchorbase,scale=.5]
        \draw [double, <-] (0,1.5) to [out=0,in=180]  (3,2);
        \draw [double] (0,.5) to [out=0,in=180] (1,.5);
        \draw [thick,->] (1,.5) to [out=60,in=180] (3,1);
        \draw [thick,->] (1,.5) to [out=-60,in=180] (3,0);
      \end{tikzpicture}
    };
    \node at (2,0) {$\rightarrow$};
    \node (topright) at (4,0) {
      \begin{tikzpicture}[anchorbase,scale=.5]
        \draw [double, <-] (0,1.5) to [out=0,in=180]  (3,2);
        \draw [double] (0,.5) to [out=0,in=180] (1,.5);
        \draw [thick,->] (1,.5) to [out=60,in=180] (3,1);
        \draw [thick,->] (1,.5) to [out=-60,in=180] (3,0);
      \end{tikzpicture}
    };
        \node (midleft) at (0,-3) {
      \begin{tikzpicture}[anchorbase,scale=.5]
        \draw [double, <-] (0,1.5) to [out=0,in=180] (.75,.5) to [out=0,in=180] (1.5,2) --  (3,2);
        \draw [white, line width=4] (0,.5) to [out=0,in=180] (.75,1.5) to [out=0,in=180] (1.5,.5);
        \draw [double] (0,.5) to [out=0,in=180] (.75,1.5) to [out=0,in=180] (1.5,.5);
        \draw [thick,->] (1.5,.5) to [out=60,in=180] (3,1);
        \draw [thick,->] (1.5,.5) to [out=-60,in=180] (3,0);
      \end{tikzpicture}
    };
    \node at (2,-3) {$\rightarrow$};
    \node (midright) at (4,-3) {
 \begin{tikzpicture}[anchorbase,scale=.5]
        \draw [double, ->] (0,.5) to [out=0,in=-90] (.3,1) to [out=90,in=0] (0,1.5);
        \draw [double, ->] (.8,1.5) to [out=0,in=90] (1.1,1) to [out=-90,in=0] (.8,.5) to [out=180,in=-90] (.5,1) to [out=90,in=180] (.8,1.5);
        \draw [double] (3,2) to [out=180,in=0] (1.75,1.5) to [out=180,in=180] (1.75,.5);
        \draw [thick,->] (1.75,.5) to [out=60,in=180] (3,1);
        \draw [thick,->] (1.75,.5) to [out=-60,in=180] (3,0);
      \end{tikzpicture}
    };
        \node (bottomleft) at (0,-6) {
      \begin{tikzpicture}[anchorbase,scale=.5]
        \draw [double, <-] (0,1.5) to [out=0,in=180] (1,0) to [out=0,in=180] (2.5,2) --  (3,2);
        \draw [white, line width=4] (0,.5) to [out=0,in=180] (.75,1.5) to [out=0,in=180] (1.25,1);
        \draw [white, line width=3] (1.25,1) to [out=60,in=180] (3,1);
        \draw [white, line width=3] (1.25,1) to [out=-60,in=180] (3,0);
        \draw [double] (0,.5) to [out=0,in=180] (.75,1.5) to [out=0,in=180] (1.25,1);
        \draw [thick,->] (1.25,1) to [out=60,in=180] (3,1);
        \draw [thick,->] (1.25,1) to [out=-60,in=180] (3,0);
        \node at (2,.5) {\tiny $c_2$};
        \node at (1.5,1.5) {\tiny $c_1$};
      \end{tikzpicture}
    };
    \node at (2,-6) {$\rightarrow$};
    \node (bottomright) at (4,-6) {
      \begin{tikzpicture}[anchorbase,scale=.5]
        \fill [black, opacity=.3] (1.4,.75) to [out=-120,in=0] (.8,.5) to [out=60,in=-60] (.8,1.5) to [out=-60,in=120] (1.4,1.3) -- (1.4,.7);
        \draw [double, ->] (0,.5) to [out=0,in=-90] (.3,1) to [out=90,in=0] (0,1.5);
        \draw [double] (.8,1.5) to [out=180,in=90] (.5,1) to [out=-90,in=180] (.8,.5);
        \draw [thick] (.8,.5) to [out=60,in=-60] (.8,1.5);
        \draw [thick] (1.4,.7) to [out=-120,in=0] (.8,.5);
        \draw [thick] (.8,1.5) to [out=-60,in=120] (1.4,1.3);
        \draw [double] (1.4,1.3) -- (1.4,.7);
        \draw [double] (3,2) -- (2.5,2);
        \draw [thick] (2.5,2) to [out=180,in=60] (1.4,1.3);
        \draw [thick,->] (2.5,2) to [out=-150,in=180] (3,1);
        \draw [thick,->] (1.4,.7) to [out=-60,in=180] (3,0);
      \end{tikzpicture}
    };
    \draw [->] (topright) -- (midright) node[midway,right,-] {
      \begin{tikzpicture}[anchorbase, scale=.45]
        \fill [yellow, opacity=.3] (1.1,1) to [out=-80,in=0] (.92,.5) to [out=180,in=-80] (.5,1) to [out=-90,in=180] (.8,0) to [out=0,in=-90] (1.15,1);
        \fill [yellow, opacity=.3] (1.1,1) to [out=100,in=0] (.68,1.5) to [out=180,in=100] (.5,1) to [out=-90,in=180] (.8,0) to [out=0,in=-90] (1.15,1);
        \draw  (.5,1) to [out=-90,in=180] (.8,0) to [out=0,in=-90] (1.15,1);
        \draw (.3,1) to [out=-90,in=180] (.8,-.5) to [out=0,in=-90] (1.5,1);
        \draw [double] (0,.5) to [out=0,in=-90] (.3,1) to [out=90,in=0] (-.25,1.5);
        \draw [double] (.68,1.5) to [out=0,in=100] (1.1,1) to [out=-80,in=0] (.92,.5) to [out=180,in=-80] (.5,1) to [out=100,in=180] (.8,1.5);
        \draw [double] (3,2) to [out=180,in=0] (1.75,1.5) to [out=180,in=180] (1.75,.5);
        \draw [thick] (1.75,.5) to [out=60,in=180] (3.25,1);
        \draw [thick] (1.75,.5) to [out=-60,in=180] (3.5,0);
      \end{tikzpicture}};
    \draw [->] (midright) -- (bottomright) node[midway,right] {\tiny $-a_3\text{cup}\otimes c_1\wedge c_2$};
  \end{tikzpicture}
\end{gather*}
The map induced by the LHS of the move thus equals:

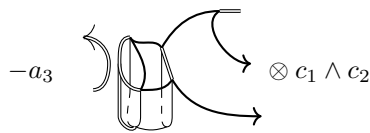
\[
  -a_3\quad
  \begin{tikzpicture}[anchorbase,scale=.7]
    \draw [double] (.5,1) -- (.5,0) to [out=-90,in=-120] (.8,-.2);
    \draw (.9,.5)-- (.9,0) to [out=-90,in=0] (.8,-.2);
    \draw (.7,1.5) -- (.7,.5);
    \draw [dashed] (.7,.5) to [out=-90,in=180] (.8,-.2);
    \draw (1.5,.7) -- (1.5,0) to [out=-90,in=0] (1.4,-.2);
    \draw (1.3,1.3) -- (1.3,.5);
    \draw [dashed] (1.3,.5) -- (1.3,0) to [out=-90,in=180] (1.4,-.2);
    \draw (.8,-.2) to [out=-20,in=-160] (1.4,-.2);
    \draw [double, ->] (0,.5) to [out=0,in=-90] (.3,1) to [out=90,in=0] (-.2,1.5);
    \draw [double] (.7,1.5) to [out=180,in=90] (.5,1) to [out=-90,in=180] (.9,.5);
    \draw [thick] (.9,.5) to [out=60,in=-60] (.7,1.5);
    \draw [thick] (1.5,.7) to [out=-120,in=0] (.9,.5);
    \draw [thick] (.7,1.5) to [out=-60,in=120] (1.3,1.3);
    \draw [double] (1.3,1.3) -- (1.5,.7);
    \draw [double] (2.8,2) -- (2.4,2);
    \draw [thick] (2.4,2) to [out=180,in=60] (1.3,1.3);
    \draw [thick,->] (2.4,2) to [out=-150,in=180] (3,1);
    \draw [thick,->] (1.5,.7) to [out=-60,in=180] (3.2,0);
  \end{tikzpicture}
  \otimes c_1\wedge c_2
\]

On the RHS, one reads:
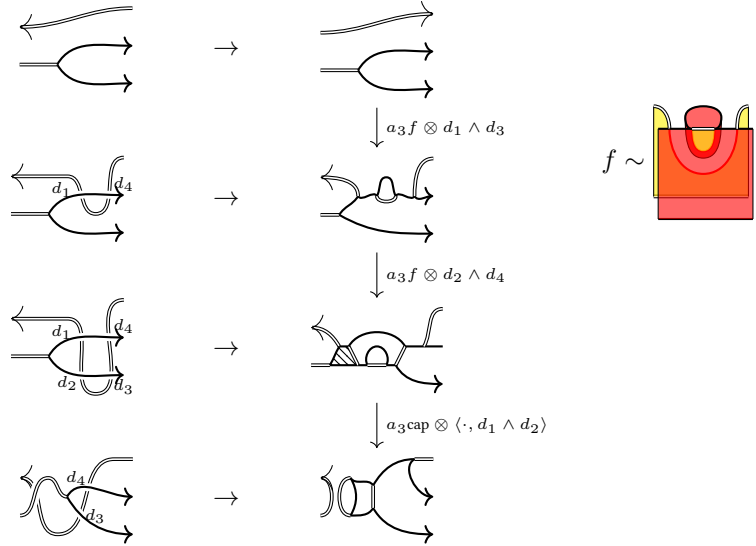
\begin{gather*}
  \begin{tikzpicture}[anchorbase]
    \node (topleft) at (0,0) {
      \begin{tikzpicture}[anchorbase,scale=.5]
        \draw [double, <-] (0,1.5) to [out=0,in=180]  (3,2);
        \draw [double] (0,.5) to [out=0,in=180] (1,.5);
        \draw [thick,->] (1,.5) to [out=60,in=180] (3,1);
        \draw [thick,->] (1,.5) to [out=-60,in=180] (3,0);
      \end{tikzpicture}
    };
    \node at (2,0) {$\rightarrow$};
    \node (topright) at (4,0) {
      \begin{tikzpicture}[anchorbase,scale=.5]
        \draw [double, ->] (0,1.5) to [out=0,in=180]  (3,2);
        \draw [double] (0,.5) to [out=0,in=180] (1,.5);
        \draw [thick,->] (1,.5) to [out=60,in=180] (3,1);
        \draw [thick,->] (1,.5) to [out=-60,in=180] (3,0);
      \end{tikzpicture}
    };
        \node (midleft+) at (0,-2) {
      \begin{tikzpicture}[anchorbase,scale=.5]
        \draw [double, <-] (0,1.5) -- (1.5,1.5) to [out=0,in=180] (2.25,.5) to [out=0,in=180]  (3,2);
        \draw [white, line width=3] (1,.5) to [out=60,in=180] (3,1);
        \draw [white, line width=3] (1,.5) to [out=-60,in=180] (3,0);
        \draw [double] (0,.5) -- (1,.5);
        \draw [thick,->] (1,.5) to [out=60,in=180] (3,1);
        \draw [thick,->] (1,.5) to [out=-60,in=180] (3,0);
        \node at (1.35,1.15) {\tiny $d_1$};
        \node at (3,1.3) {\tiny $d_4$};
      \end{tikzpicture}
    };
    \node at (2,-2) {$\rightarrow$};
    \node (midright+) at (4,-2) {
      \begin{tikzpicture}[anchorbase,scale=.5]
        \draw [double] (0,.5) to (.5,.5);
         \draw [thick,->] (.5,.5) to [out=-30,in=180] (3,0);
         \draw [thick] (.5,.5) to [out=60,in=-120] (1,1);
         \draw [thick] (1,1) to [out=-60,in=120] (1.5,1);
         \draw [double] (1.5,1) to [out=-90,in=-90] (2,1);
         \draw [thick] (1.5,1) to [out=60,in=180] (1.75,1.5) to [out=0,in=120] (2,1);
         \draw [thick] (2,1) to [out=60,in=-120] (2.5,1);
         \draw [thick,->] (2.5,1) to [out=-60,in=180] (3,1);
         \draw [double,->] (1,1) to [out=90,in=0] (0,1.5);
         \draw [double] (2.5,1) to [out=90,in=180] (3,2);
       \end{tikzpicture}
    };
        \node (midleft-) at (0,-4) {
      \begin{tikzpicture}[anchorbase,scale=.5]
        \draw [double, <-] (0,1.5) -- (1.5,1.5) to [out=0,in=180] (2.25,-.5) to [out=0,in=180]  (3,2);
        \draw [white, line width=3] (1,.5) to [out=60,in=180] (3,1);
        \draw [white, line width=3] (1,.5) to [out=-60,in=180] (3,0);
        \draw [double] (0,.5) -- (1,.5);
        \draw [thick,->] (1,.5) to [out=60,in=180] (3,1);
        \draw [thick,->] (1,.5) to [out=-60,in=180] (3,0);
        \node at (1.35,1.15) {\tiny $d_1$};
        \node at (1.5,-.2) {\tiny $d_2$};
        \node at (3,-.3) {\tiny $d_3$};
        \node at (3,1.3) {\tiny $d_4$};
      \end{tikzpicture}
    };
    \node at (2,-4) {$\rightarrow$};
    \node (midright-) at (4,-4) {
      \begin{tikzpicture}[anchorbase,scale=.5]
         \fill[pattern=north west lines] (.5,.5) -- (1.25,.5) to [out=120,in=-60] (1,1) -- (.75,1) to [out=-120,in=60] (.5,.5);
        \draw [double] (0,.5) -- (.5,.5);
        \draw [double,->] (.75,1) to [out=120,in=0] (0,1.5);
        \draw [thick] (.5,.5) to [out=60,in=-120] (.75,1);
        \draw [thick] (.75,1) -- (1,1);
        \draw [thick] (.5,.5) -- (1.25,.5);
        \draw [double] (1,1) -- (1.25,.5);
        \draw [thick] (1.25,.5) -- (1.5,.5);
        \draw [thick] (1.5,.5) to [out=120,in=180] (1.75,1) to [out=0,in=60] (2,.5);
        \draw [double] (1.5,.5) -- (2,.5);
        \draw [thick] (2,.5) -- (2.25,.5);
        \draw [thick] (1,1) to [out=60,in=120] (2.5,1);
        \draw [double] (2.25,.5) -- (2.5,1);
        \draw [thick] (2.5,1) -- (3,1);
        \draw [thick] (3,1) -- (3.5,1);
        \draw [double] (3,1) to [out=60,in=180] (3.5,2);
        \draw [thick,->] (2.25,.5) to [out=-60,in=180] (3.5,0);
      \end{tikzpicture}
    };
        \node (bottomleft) at (0,-6) {
      \begin{tikzpicture}[anchorbase,scale=.5]
        \draw [double, <-] (0,1.5) to [out=0,in=180] (1,0) to [out=0,in=180] (2.5,2) --  (3,2);
        \draw [white, line width=4] (0,.5) to [out=0,in=180] (.75,1.5) to [out=0,in=180] (1.25,1);
        \draw [white, line width=3] (1.25,1) to [out=60,in=180] (3,1);
        \draw [white, line width=3] (1.25,1) to [out=-60,in=180] (3,0);
        \draw [double] (0,.5) to [out=0,in=180] (.75,1.5) to [out=0,in=180] (1.25,1);
        \draw [thick,->] (1.25,1) to [out=60,in=180] (3,1);
        \draw [thick,->] (1.25,1) to [out=-60,in=180] (3,0);
        \node at (2,.5) {\tiny $d_3$};
        \node at (1.5,1.5) {\tiny $d_4$};
      \end{tikzpicture}
    };
    \node at (2,-6) {$\rightarrow$};
    \node (bottomright) at (4,-6) {
      \begin{tikzpicture}[anchorbase,scale=.5]
        \draw [double, ->] (0,.5) to [out=0,in=-90] (.3,1) to [out=90,in=0] (0,1.5);
        \draw [double] (.8,1.5) to [out=180,in=90] (.5,1) to [out=-90,in=180] (.8,.5);
        \draw [thick] (.8,.5) to [out=60,in=-60] (.8,1.5);
        \draw [thick] (1.4,.7) to [out=-120,in=0] (.8,.5);
        \draw [thick] (.8,1.5) to [out=-60,in=120] (1.4,1.3);
        \draw [double] (1.4,1.3) -- (1.4,.7);
        \draw [double] (3,2) -- (2.5,2);
        \draw [thick] (2.5,2) to [out=180,in=60] (1.4,1.3);
        \draw [thick,->] (2.5,2) to [out=-150,in=180] (3,1);
        \draw [thick,->] (1.4,.7) to [out=-60,in=180] (3,0);
      \end{tikzpicture}
    };
    \draw [->] (topright) -- (midright+) node[midway,right] {\tiny $a_3f\otimes d_1\wedge d_3$};
    \draw [->] (midright+) -- (midright-) node[midway,right] {\tiny $a_3f\otimes d_2\wedge d_4$};
    \draw [->] (midright-) -- (bottomright) node[midway,right] {\tiny $a_3\text{cap}\otimes \langle \cdot, d_1\wedge d_2\rangle$};
    \node at (8,-1.5) {$f \sim $ \begin{tikzpicture}[anchorbase,scale=.3]
        \draw [double] (-.2,1) -- (4,1);
        \draw [thick] (0,0) -- (4.2,0);
        \fill [yellow, opacity=.6] (-.2,5) to [out=0,in=90] (.5,4) to [out=-90,in=180] (2,2) to [out=0,in=-90] (3.5,4) to [out=90,in=180] (4,5) -- (4,1) -- (-.2,1) -- (-.2,5);
        \draw (-.2,1) -- (-.2,5);
        \draw (4,1) -- (4,5);
        \fill [red, opacity=.6] (2.8,4.5) to [out=-90,in=0] (2,2.7) to [out=180,in=-90] (1.2,4.5) to [out=-90,in=180] (1.5,4) to [out=-90,in=180] (2,3) to [out=0,in=-90] (2.5,4) to [out=0,in=-90] (2.8,4.5);
        \fill [red, opacity=.6] (2.8,4.5) to [out=-90,in=0] (2,2.7) to [out=180,in=-90] (1.2,4.5) to [out=90,in=180] (2,5) to [out=0,in=90] (2.8,4.5) to [out=-90,in=0] (2.5,4) to [out=0,in=-90] (2.8,4.5);
        \draw (2.8,4.5) to [out=-90,in=0] (2,2.7) to [out=180,in=-90] (1.2,4.5);
        \draw [thick, red] (2.5,4) to [out=-90,in=0] (2,3) to [out=180,in=-90] (1.5,4);
        \fill [yellow, opacity=.6] (1.5,4) -- (2.5,4) to [out=-90,in=0] (2,3) to [out=180,in=-90] (1.5,4);
        \fill [red, opacity=.6] (2.5,4) -- (3.5,4) to [out=-90,in=0] (2,2) to [out=180,in=-90] (.5,4) -- (1.5,4) to [out=-90,in=180] (2,3) to [out=0,in=-90] (2.5,4);
        \fill [red,opacity=.6] (3.5,4) to [out=-90,in=0] (2,2) to [out=180,in=-90] (.5,4) -- (0,4) -- (0,0) -- (4.2,0) -- (4.2,4) -- (3.5,4);
        \draw (0,4) -- (0,0);
        \draw (4.2,4) -- (4.2,0);
        \draw [thick, red] (3.5,4) to [out=-90,in=0] (2,2) to [out=180,in=-90] (.5,4);
        \draw [thick] (0,4) -- (.5,4);
        \draw [double] (.5,4) to [out=90,in=0] (-.2,5);
        \draw [thick] (.5,4) -- (1.5,4);
        \draw [double] (1.5,4) -- (2.5,4);
        \draw [thick] (2.5,4) -- (3.5,4);
        \draw [thick] (2.5,4) to [out=0,in=-90] (2.8,4.5) to [out=90,in=0] (2,5) to [out=180,in=90] (1.2,4.5) to [out=-90,in=180] (1.5,4);
        \draw [thick] (3.5,4) -- (4.2,4);
        \draw [double] (3.5,4) to [out=90,in=180] (4,5);
  \end{tikzpicture}};
  \end{tikzpicture}
\end{gather*}

Composing all steps together, one obtains the very same foam as before. The coefficient is computed as:
\[
  a_3^3\langle d_1\wedge d_3\wedge d_2\wedge d_4,d_1\wedge d_2\rangle=-a_3 d_3\wedge d_4
\]
which matches the expectation after identification of $d_3$ with $c_1$ and $d_4$ with $c_2$.

Again, this computation extends for all configurations.

Then we consider the case where the single strand comes with label $1$. It is enough to check that the coefficients agree after projecting on a non-zero fully deformed configuration. We start with the following check:

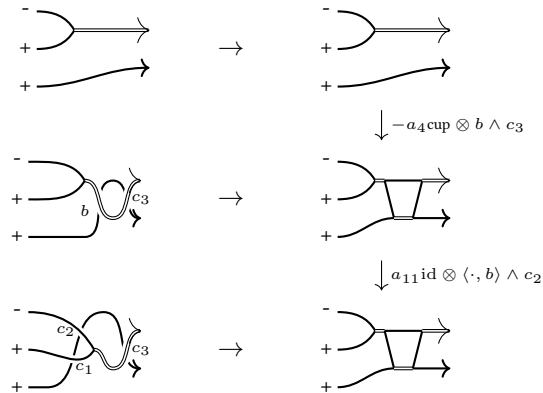
\begin{gather*}
  \begin{tikzpicture}[anchorbase]
    \node (topleft) at (0,0) {
      \begin{tikzpicture}[anchorbase,scale=.5]
        \draw [thick, ->] (0,0) to [out=0,in=180]  (3,.5);
        \draw [thick] (0,2) to [out=0,in=120] (1,1.5);
        \draw [thick] (0,1) to [out=0,in=-120] (1,1.5);
        \draw [double,->] (1,1.5) to [out=0,in=180] (3,1.5);
        \node at (-.3,0) {\small +};
        \node at (-.3,1) {\small +};
        \node at (-.3,2) {\small -};
      \end{tikzpicture}
    };
    \node at (2,0) {$\rightarrow$};
    \node (topright) at (4,0) {
      \begin{tikzpicture}[anchorbase,scale=.5]
        \draw [thick, ->] (0,0) to [out=0,in=180]  (3,.5);
        \draw [thick] (0,2) to [out=0,in=120] (1,1.5);
        \draw [thick] (0,1) to [out=0,in=-120] (1,1.5);
        \draw [double,->] (1,1.5) to [out=0,in=180] (3,1.5);
        \node at (-.3,0) {\small +};
        \node at (-.3,1) {\small +};
        \node at (-.3,2) {\small -};
      \end{tikzpicture}
    };
        \node (midleft) at (0,-2) {
      \begin{tikzpicture}[anchorbase,scale=-.5]
        \draw [thick, <-] (0,1.5) to [out=0,in=180] (.75,.5) to [out=0,in=180] (1.5,2) --  (3,2);
        \draw [white, line width=4] (0,.5) to [out=0,in=180] (.75,1.5) to [out=0,in=180] (1.5,.5);
        \draw [double,<-] (0,.5) to [out=0,in=180] (.75,1.5) to [out=0,in=180] (1.5,.5);
        \draw [thick] (1.5,.5) to [out=60,in=180] (3,1);
        \draw [thick] (1.5,.5) to [out=-60,in=180] (3,0);
        \node at (0,1) {\tiny $c_3$};
        \node at (1.5,1.3) {\tiny $b$};
        \node at (3.3,2) {\small +};
        \node at (3.3,1) {\small +};
        \node at (3.3,0) {\small -};
      \end{tikzpicture}
    };
    \node at (2,-2) {$\rightarrow$};
    \node (midright) at (4,-2) {
      \begin{tikzpicture}[anchorbase,scale=.5]
        \draw [thick] (0,1) to [out=0,in=-120] (1,1.5);
        \draw [thick] (0,2) to [out=0,in=120] (1,1.5);
        \draw [double] (1,1.5) -- (1.25,1.5);
        \draw [thick] (0,0) to [out=0,in=180] (1.5,.5);
        \draw [thick] (1.25,1.5) -- (1.5,.5);
        \draw [double] (1.5,.5) -- (2,.5);
        \draw [thick] (1.25,1.5) -- (2.25,1.5);
        \draw [thick] (2,.5) -- (2.25,1.5);
        \draw [double, ->] (2.25,1.5) -- (3,1.5);
        \draw [thick, ->] (2,.5) -- (3,.5);
        \node at (-.3,0) {\small +};
        \node at (-.3,1) {\small +};
        \node at (-.3,2) {\small -};
      \end{tikzpicture}
    };
        \node (bottomleft) at (0,-4) {
      \begin{tikzpicture}[anchorbase,scale=-.5]
        \draw [thick, <-] (0,1.5) to [out=0,in=180] (1,0) to [out=0,in=180] (2.5,2) --  (3,2);
        \draw [white, line width=4] (0,.5) to [out=0,in=180] (.75,1.5) to [out=0,in=180] (1.25,1);
        \draw [white, line width=3] (1.25,1) to [out=60,in=180] (3,1);
        \draw [white, line width=3] (1.25,1) to [out=-60,in=180] (3,0);
        \draw [double,<-] (0,.5) to [out=0,in=180] (.75,1.5) to [out=0,in=180] (1.25,1);
        \draw [thick] (1.25,1) to [out=60,in=180] (3,1);
        \draw [thick] (1.25,1) to [out=-60,in=180] (3,0);
        \node at (2,.5) {\tiny $c_2$};
        \node at (1.5,1.5) {\tiny $c_1$};
        \node at (0,1) {\tiny $c_3$};
        \node at (3.3,2) {\small +};
        \node at (3.3,1) {\small +};
        \node at (3.3,0) {\small -};
      \end{tikzpicture}
    };
    \node at (2,-4) {$\rightarrow$};
    \node (bottomright) at (4,-4) {
      \begin{tikzpicture}[anchorbase,scale=.5]
        \draw [thick] (0,1) to [out=0,in=-120] (1,1.5);
        \draw [thick] (0,2) to [out=0,in=120] (1,1.5);
        \draw [double] (1,1.5) -- (1.25,1.5);
        \draw [thick] (0,0) to [out=0,in=180] (1.5,.5);
        \draw [thick] (1.25,1.5) -- (1.5,.5);
        \draw [double] (1.5,.5) -- (2,.5);
        \draw [thick] (1.25,1.5) -- (2.25,1.5);
        \draw [thick] (2,.5) -- (2.25,1.5);
        \draw [double, ->] (2.25,1.5) -- (3,1.5);
        \draw [thick, ->] (2,.5) -- (3,.5);
        \node at (-.3,0) {\small +};
        \node at (-.3,1) {\small +};
        \node at (-.3,2) {\small -};
      \end{tikzpicture}
    };
    \draw [->] (topright) -- (midright) node[midway,right] {\tiny $-a_4 \text{cup}\otimes b\wedge c_3$};
    \draw [->] (midright) -- (bottomright) node[midway,right] {\tiny $a_{11}\id\otimes \langle \cdot,b\rangle \wedge c_2$};
  \end{tikzpicture}
\end{gather*}

The above two steps compose into a cup foam with coefficient $-a_4a_{11}c_2\wedge c_3$. On the right-hand side, one gets:
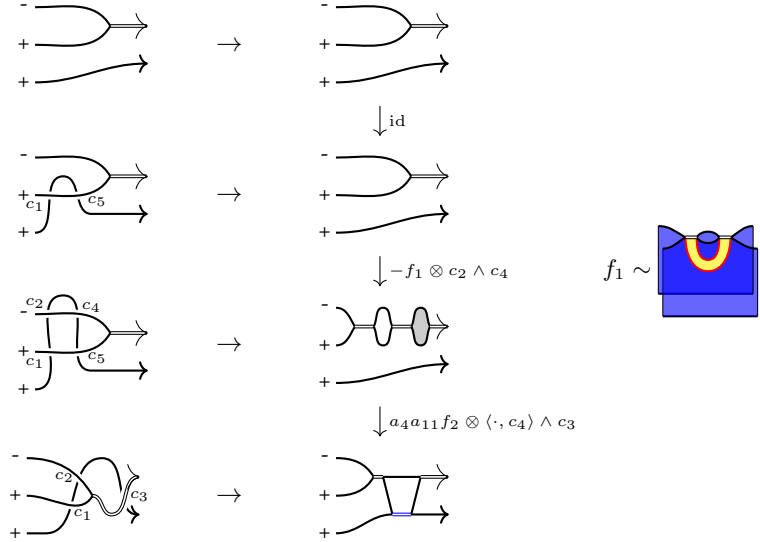
\begin{gather*}
  \begin{tikzpicture}[anchorbase]
    \node (topleft) at (0,0) {
      \begin{tikzpicture}[anchorbase,scale=-.5]
        \draw [thick, <-] (0,1.5) to [out=0,in=180]  (3,2);
        \draw [double,<-] (0,.5) to [out=0,in=180] (1,.5);
        \draw [thick] (1,.5) to [out=60,in=180] (3,1);
        \draw [thick] (1,.5) to [out=-60,in=180] (3,0);
        \node at (3.3,2) {\small +};
        \node at (3.3,1) {\small +};
        \node at (3.3,0) {\small -};
      \end{tikzpicture}
    };
    \node at (2,0) {$\rightarrow$};
    \node (topright) at (4,0) {
      \begin{tikzpicture}[anchorbase,scale=-.5]
        \draw [thick, <-] (0,1.5) to [out=0,in=180]  (3,2);
        \draw [double,<-] (0,.5) to [out=0,in=180] (1,.5);
        \draw [thick] (1,.5) to [out=60,in=180] (3,1);
        \draw [thick] (1,.5) to [out=-60,in=180] (3,0);
        \node at (3.3,2) {\small +};
        \node at (3.3,1) {\small +};
        \node at (3.3,0) {\small -};
      \end{tikzpicture}
    };
        \node (midleft+) at (0,-2) {
      \begin{tikzpicture}[anchorbase,scale=-.5]
        \draw [thick, <-] (0,1.5) -- (1.5,1.5) to [out=0,in=180] (2.25,.5) to [out=0,in=180]  (3,2);
        \draw [white, line width=3] (1,.5) to [out=60,in=180] (3,1);
        \draw [white, line width=3] (1,.5) to [out=-60,in=180] (3,0);
        \draw [double,<-] (0,.5) -- (1,.5);
        \draw [thick] (1,.5) to [out=60,in=180] (3,1);
        \draw [thick] (1,.5) to [out=-60,in=180] (3,0);
        \node at (1.35,1.15) {\tiny $c_5$};
        \node at (3,1.3) {\tiny $c_1$};
        \node at (3.3,2) {\small +};
        \node at (3.3,1) {\small +};
        \node at (3.3,0) {\small -};      \end{tikzpicture}
    };
    \node at (2,-2) {$\rightarrow$};
    \node (midright+) at (4,-2) {
      \begin{tikzpicture}[anchorbase,scale=-.5]
        \draw [thick, <-] (0,1.5) to [out=0,in=180]  (3,2);
        \draw [double,<-] (0,.5) to [out=0,in=180] (1,.5);
        \draw [thick] (1,.5) to [out=60,in=180] (3,1);
        \draw [thick] (1,.5) to [out=-60,in=180] (3,0);
        \node at (3.3,2) {\small +};
        \node at (3.3,1) {\small +};
        \node at (3.3,0) {\small -};
      \end{tikzpicture}
    };
        \node (midleft-) at (0,-4) {
      \begin{tikzpicture}[anchorbase,scale=-.5]
        \draw [thick, <-] (0,1.5) -- (1.5,1.5) to [out=0,in=180] (2.25,-.5) to [out=0,in=180]  (3,2);
        \draw [white, line width=3] (1,.5) to [out=60,in=180] (3,1);
        \draw [white, line width=3] (1,.5) to [out=-60,in=180] (3,0);
        \draw [double,<-] (0,.5) -- (1,.5);
        \draw [thick] (1,.5) to [out=60,in=180] (3,1);
        \draw [thick] (1,.5) to [out=-60,in=180] (3,0);
        \node at (1.35,1.15) {\tiny $c_5$};
        \node at (1.5,-.2) {\tiny $c_4$};
        \node at (3,-.3) {\tiny $c_2$};
        \node at (3,1.3) {\tiny $c_1$};
        \node at (3.3,2) {\small +};
        \node at (3.3,1) {\small +};
        \node at (3.3,0) {\small -};
      \end{tikzpicture}
    };
    \node at (2,-4) {$\rightarrow$};
    \node (midright-) at (4,-4) {
      \begin{tikzpicture}[anchorbase,scale=.5]
        \fill [fill=black, opacity=.2]  (2,1.5) to [out=60,in=180] (2.25,2) to [out=0,in=120] (2.5,1.5) to [out=-120,in=0] (2.25,1) to [out=180,in=-60] (2,1.5);
        \draw [thick,->] (0,0) to [out=0,in=180] (3,.5);
        \draw [thick] (0,1) to [out=0,in=-120] (.5,1.5);
        \draw [thick] (0,2) to [out=0,in=120] (.5,1.5);
        \draw [double] (.5,1.5) -- (1,1.5);
        \draw [thick] (1,1.5) to [out=60,in=180] (1.25,2) to [out=0,in=120] (1.5,1.5);
        \draw [thick] (1,1.5) to [out=-60,in=180] (1.25,1) to [out=0,in=-120] (1.5,1.5);
        \draw [double] (1.5,1.5) -- (2,1.5);
        \draw [thick] (2,1.5) to [out=60,in=180] (2.25,2) to [out=0,in=120] (2.5,1.5);
        \draw [thick] (2,1.5) to [out=-60,in=180] (2.25,1) to [out=0,in=-120] (2.5,1.5);
        \draw [double,->] (2.5,1.5) -- (3,1.5);
        \node at (-.3,0) {\small +};
        \node at (-.3,1) {\small +};
        \node at (-.3,2) {\small -};
      \end{tikzpicture}
    };
        \node (bottomleft) at (0,-6) {
      \begin{tikzpicture}[anchorbase,scale=-.5]
        \draw [thick, <-] (0,1.5) to [out=0,in=180] (1,0) to [out=0,in=180] (2.5,2) --  (3,2);
        \draw [white, line width=4] (0,.5) to [out=0,in=180] (.75,1.5) to [out=0,in=180] (1.25,1);
        \draw [white, line width=3] (1.25,1) to [out=60,in=180] (3,1);
        \draw [white, line width=3] (1.25,1) to [out=-60,in=180] (3,0);
        \draw [double,<-] (0,.5) to [out=0,in=180] (.75,1.5) to [out=0,in=180] (1.25,1);
        \draw [thick] (1.25,1) to [out=60,in=180] (3,1);
        \draw [thick] (1.25,1) to [out=-60,in=180] (3,0);
        \node at (2,.5) {\tiny $c_2$};
        \node at (1.5,1.5) {\tiny $c_1$};
        \node at (0,1) {\tiny $c_3$};
        \node at (3.3,2) {\small +};
        \node at (3.3,1) {\small +};
        \node at (3.3,0) {\small -};
      \end{tikzpicture}
    };
    \node at (2,-6) {$\rightarrow$};
    \node (bottomright) at (4,-6) {
           \begin{tikzpicture}[anchorbase,scale=.5]
        \draw [thick] (0,1) to [out=0,in=-120] (1,1.5);
        \draw [thick] (0,2) to [out=0,in=120] (1,1.5);
        \draw [double] (1,1.5) -- (1.25,1.5);
        \draw [thick] (0,0) to [out=0,in=180] (1.5,.5);
        \draw [thick] (1.25,1.5) -- (1.5,.5);
        \draw [double,blue] (1.5,.5) -- (2,.5);
        \draw [thick] (1.25,1.5) -- (2.25,1.5);
        \draw [thick] (2,.5) -- (2.25,1.5);
        \draw [double, ->] (2.25,1.5) -- (3,1.5);
        \draw [thick, ->] (2,.5) -- (3,.5);
        \node at (-.3,0) {\small +};
        \node at (-.3,1) {\small +};
        \node at (-.3,2) {\small -};
      \end{tikzpicture}
    };
    \draw [->] (topright) -- (midright+) node[midway,right] {\tiny $\id$};
    \draw [->] (midright+) -- (midright-) node[midway,right] {\tiny $-f_1 \otimes c_2\wedge c_4$};
    \draw [->] (midright-) -- (bottomright) node[midway,right] {\tiny $a_4a_{11} f_2\otimes \langle \cdot, c_4\rangle\wedge c_3$};
    \node at (8,-3) {\small $f_1\sim$      \begin{tikzpicture}[anchorbase,scale=-.3]
    \draw [thick] (0,3) -- (4.2,3);
    \draw [thick] (-.2,4) -- (4,4);
    \fill [blue,opacity=.6] (1.5,.5) to [out=60,in=120] (2.5,.5) to [out=90,in=0] (2,1.5) to [out=180,in=90] (1.5,.5);
    \fill [blue, opacity=.6] (-.2,1) to [out=0,in=150] (1,.5) to [out=90,in=180] (2,2) to [out=0,in=90] (3,.5) to [out=30,in=180] (4,1) -- (4,4) -- (-.2,4) -- (-.2,1);
    \draw (-.2,1) -- (-.2,4);
    \draw (4,1) -- (4,4);
    \fill [blue,opacity=.6] (1.5,.5) to [out=-60,in=-120] (2.5,.5) to [out=90,in=0] (2,1.5) to [out=180,in=90] (1.5,.5);
    \fill [yellow,opacity=.6]   (1,.5) -- (1.5,.5) to [out=90,in=180] (2,1.5) to [out=0,in=90] (2.5,.5) -- (3,.5) to [out=90,in=0] (2,2) to [out=180,in=90] (1,.5);
    \fill [blue, opacity=.6] (0,0) to [out=0,in=-150] (1,.5) to [out=90,in=180] (2,2) to [out=0,in=90] (3,.5) to [out=-30,in=180] (4.2,0) -- (4.2,3) -- (0,3) -- (0,0);
    \draw (0,0) -- (0,3);
    \draw (4.2,0) -- (4.2,3);
    \draw [red,thick] (2.5,.5) to [out=90,in=0] (2,1.5) to [out=180,in=90] (1.5,.5);
    \draw [red, thick] (3,.5) to [out=90,in=0] (2,2) to [out=180,in=90] (1,.5);
    \draw [thick] (0,0) to [out=0,in=-150] (1,.5);
    \draw [thick] (-.2,1) to [out=0,in=150] (1,.5);
    \draw [double] (1,.5) -- (1.5,.5);
    \draw [thick] (1.5,.5) to [out=60,in=120] (2.5,.5);
    \draw [thick] (1.5,.5) to [out=-60,in=-120] (2.5,.5);
    \draw [double] (2.5,.5) -- (3,.5);
    \draw [thick] (3,.5) to [out=30,in=180] (4,1);
    \draw [thick] (3,.5) to [out=-30,in=180] (4.2,0);
  \end{tikzpicture}
};
  \end{tikzpicture}
\end{gather*}
Above, $f_2$ is the foam that consists in caping off the digon fill with gray at step 3 and merge two 1-labeled strands so as to create the 2-labeled strand colored in blue at step 4.
Altogether, the composite foam agrees on the nose with the previous one. One then computes as expected:
\[
  -a_4a_{11}\langle c_2\wedge c_4,c_4 \rangle\wedge c_3=-a_4a_{11} c_2\wedge c_3
\]

To deduce other versions of the moves, one can carefully follow the effect on coefficients induced by globally reversing the orientation or changing the over/under information for all crossings, as well as applying a symmetry along a horizontal line. 

So this leaves us with a last configuration, where one only reverses the orientation of the 1-labeled line. One gets for one side of a first version of the move:
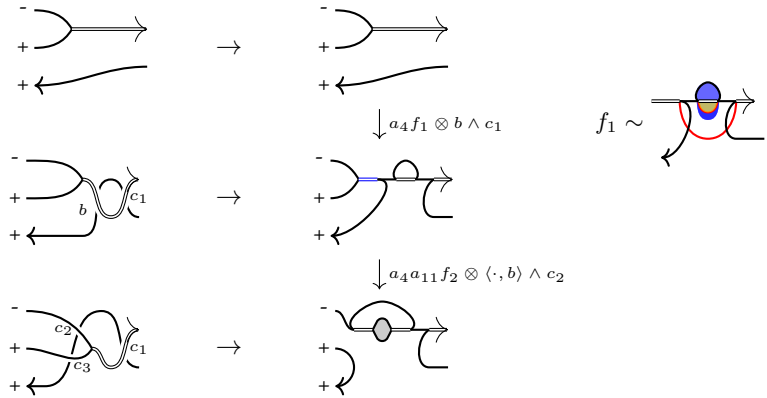
\begin{gather*}
  \begin{tikzpicture}[anchorbase]
    \node (topleft) at (0,0) {
      \begin{tikzpicture}[anchorbase,scale=.5]
        \draw [thick, <-] (0,0) to [out=0,in=180]  (3,.5);
        \draw [thick] (0,2) to [out=0,in=120] (1,1.5);
        \draw [thick] (0,1) to [out=0,in=-120] (1,1.5);
        \draw [double,->] (1,1.5) to [out=0,in=180] (3,1.5);
        \node at (-.3,0) {\small +};
        \node at (-.3,1) {\small +};
        \node at (-.3,2) {\small -};
      \end{tikzpicture}
    };
    \node at (2,0) {$\rightarrow$};
    \node (topright) at (4,0) {
      \begin{tikzpicture}[anchorbase,scale=.5]
        \draw [thick, <-] (0,0) to [out=0,in=180]  (3,.5);
        \draw [thick] (0,2) to [out=0,in=120] (1,1.5);
        \draw [thick] (0,1) to [out=0,in=-120] (1,1.5);
        \draw [double,->] (1,1.5) to [out=0,in=180] (3,1.5);
        \node at (-.3,0) {\small +};
        \node at (-.3,1) {\small +};
        \node at (-.3,2) {\small -};
      \end{tikzpicture}
    };
        \node (midleft) at (0,-2) {
      \begin{tikzpicture}[anchorbase,scale=-.5]
        \draw [thick, ->] (0,1.5) to [out=0,in=180] (.75,.5) to [out=0,in=180] (1.5,2) --  (3,2);
        \draw [white, line width=4] (0,.5) to [out=0,in=180] (.75,1.5) to [out=0,in=180] (1.5,.5);
        \draw [double,<-] (0,.5) to [out=0,in=180] (.75,1.5) to [out=0,in=180] (1.5,.5);
        \draw [thick] (1.5,.5) to [out=60,in=180] (3,1);
        \draw [thick] (1.5,.5) to [out=-60,in=180] (3,0);
        \node at (0,1) {\tiny $c_1$};
        \node at (1.5,1.3) {\tiny $b$};
        \node at (3.3,2) {\small +};
        \node at (3.3,1) {\small +};
        \node at (3.3,0) {\small -};
      \end{tikzpicture}
    };
    \node at (2,-2) {$\rightarrow$};
    \node (midright) at (4,-2) {
      \begin{tikzpicture}[anchorbase,scale=.5]
        \draw [thick] (0,1) to [out=0,in=-120] (.75,1.5);
        \draw [thick] (0,2) to [out=0,in=120] (.75,1.5);
        \draw [double,blue] (.75,1.5) -- (1.25,1.5);
        \draw [thick,->] (1.25,1.5) to [out=0,in=0] (0,0);
        \draw [thick] (1.25,1.5) -- (1.75,1.5);
        \draw [double] (1.75,1.5) -- (2.25,1.5);
        \draw [thick] (2.25,1.5) to [out=30,in=0] (2,2) to [out=180,in=150] (1.75,1.5);
        \draw [thick] (2.25,1.5) -- (2.75,1.5);
        \draw [double,->] (2.75,1.5) -- (3.25,1.5);
        \draw [thick] (3.25,.5) -- (2.75,.5) to [out=180,in=-150] (2.75,1.5);
        \node at (-.3,0) {\small +};
        \node at (-.3,1) {\small +};
        \node at (-.3,2) {\small -};
      \end{tikzpicture}
    };
        \node (bottomleft) at (0,-4) {
      \begin{tikzpicture}[anchorbase,scale=-.5]
        \draw [thick, ->] (0,1.5) to [out=0,in=180] (1,0) to [out=0,in=180] (2.5,2) --  (3,2);
        \draw [white, line width=4] (0,.5) to [out=0,in=180] (.75,1.5) to [out=0,in=180] (1.25,1);
        \draw [white, line width=3] (1.25,1) to [out=60,in=180] (3,1);
        \draw [white, line width=3] (1.25,1) to [out=-60,in=180] (3,0);
        \draw [double,<-] (0,.5) to [out=0,in=180] (.75,1.5) to [out=0,in=180] (1.25,1);
        \draw [thick] (1.25,1) to [out=60,in=180] (3,1);
        \draw [thick] (1.25,1) to [out=-60,in=180] (3,0);
        \node at (2,.5) {\tiny $c_2$};
        \node at (1.5,1.5) {\tiny $c_3$};
        \node at (0,1) {\tiny $c_1$};
        \node at (3.3,2) {\small +};
        \node at (3.3,1) {\small +};
        \node at (3.3,0) {\small -};
      \end{tikzpicture}
    };
    \node at (2,-4) {$\rightarrow$};
    \node (bottomright) at (4,-4) {
      \begin{tikzpicture}[anchorbase,scale=.5]
        \fill [black, opacity=.2] (1,1.5) to [out=60,in=180] (1.25,1.8) to [out=0,in=120] (1.5,1.5) to [out=-120,in=0] (1.25,1.2) to [out=180,in=-60] (1,1.5);
        \draw [thick,->] (0,1) to [out=0,in=90] (.5,.5) to [out=-90,in=0] (0,0);
        \draw [thick] (0,2) to [out=0,in=180] (.5,1.5);
        \draw [double] (.5,1.5) -- (1,1.5);
        \draw [thick] (1,1.5) to [out=60,in=180] (1.25,1.8) to [out=0,in=120] (1.5,1.5);
        \draw [thick] (1,1.5) to [out=-60,in=180] (1.25,1.2) to [out=0,in=-120] (1.5,1.5);
        \draw [double] (1.5,1.5) -- (2,1.5);
        \draw [thick] (2,1.5) to [out=30,in=0] (1.25,2.25) to [out=180,in=150] (.5,1.5);
        \draw [thick] (2,1.5) -- (2.5,1.5);
        \draw [double,->] (2.5,1.5) -- (3,1.5);
        \draw [thick] (3,.5) -- (2.5,.5) to [out=180,in=-150] (2.5,1.5);
        \node at (-.3,0) {\small +};
        \node at (-.3,1) {\small +};
        \node at (-.3,2) {\small -};
      \end{tikzpicture}
    };
    \draw [->] (topright) -- (midright) node[midway,right] {\tiny $a_4 f_1\otimes b\wedge c_1$};
    \draw [->] (midright) -- (bottomright) node[midway,right] {\tiny $a_4a_{11}f_2\otimes \langle \cdot,b\rangle \wedge c_2$};
    \node at (8,-1) {\small $f_1 \sim$
      \begin{tikzpicture}[anchorbase,scale=.5]
        \fill [blue, opacity=.6] (2.3,1.75) to [out=90,in=0] (2,2) to [out=180,in=90] (1.7,1.75) to [out=-90,in=180] (2,1) to [out=0,in=-90] (2.3,1.75);
        \fill [blue, opacity=.6] (2.3,1.75) to [out=-90,in=0] (2.25,1.5) to [out=-90,in=0] (2,1.2) to [out=180,in=-90] (1.75,1.5) to [out=150,in=-90] (1.7,1.75) to [out=-90,in=180] (2,1) to [out=0,in=-90] (2.3,1.75);
        \draw [thick,red] (2.25,1.5) to [out=-90,in=0] (2,1.2) to [out=180,in=-90] (1.75,1.5);
        \draw [thick,red] (2.75,1.5) to [out=-90,in=0] (2,.5) to [out=180,in=-90] (1.25,1.5);
        \fill [yellow, opacity=.6] (2.25,1.5) to [out=-90,in=0] (2,1.2) to [out=180,in=-90] (1.75,1.5);
        \draw [double] (.5,1.5) -- (1.25,1.5);
        \draw [thick,->] (1.25,1.5) to [out=0,in=0] (.75,0);
        \draw [thick] (1.25,1.5) -- (1.75,1.5);
        \draw [double] (1.75,1.5) -- (2.25,1.5);
        \draw [thick] (2.25,1.5) to [out=30,in=0] (2,2) to [out=180,in=150] (1.75,1.5);
        \draw [thick] (2.25,1.5) -- (2.75,1.5);
        \draw [double,->] (2.75,1.5) -- (3.25,1.5);
        \draw [thick] (3.5,.5) -- (2.75,.5) to [out=180,in=-150] (2.75,1.5);
      \end{tikzpicture}
};
  \end{tikzpicture}
\end{gather*}
Above, $f_2$ is the foam that caps off the blue-colored 2-labeled strand before cupping the Gray-filled digon. We record the coefficient:
\[-a_{11}\otimes c_1 \wedge c_2\]

On the right-hand side of the move, one gets:
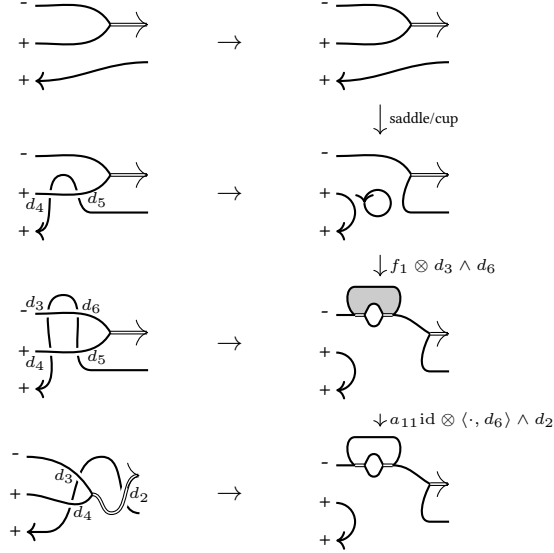
\begin{gather*}
  \begin{tikzpicture}[anchorbase]
    \node (topleft) at (0,0) {
      \begin{tikzpicture}[anchorbase,scale=-.5]
        \draw [thick, ->] (0,1.5) to [out=0,in=180]  (3,2);
        \draw [double,<-] (0,.5) to [out=0,in=180] (1,.5);
        \draw [thick] (1,.5) to [out=60,in=180] (3,1);
        \draw [thick] (1,.5) to [out=-60,in=180] (3,0);
        \node at (3.3,2) {\small +};
        \node at (3.3,1) {\small +};
        \node at (3.3,0) {\small -};
      \end{tikzpicture}
    };
    \node at (2,0) {$\rightarrow$};
    \node (topright) at (4,0) {
      \begin{tikzpicture}[anchorbase,scale=-.5]
        \draw [thick, ->] (0,1.5) to [out=0,in=180]  (3,2);
        \draw [double,<-] (0,.5) to [out=0,in=180] (1,.5);
        \draw [thick] (1,.5) to [out=60,in=180] (3,1);
        \draw [thick] (1,.5) to [out=-60,in=180] (3,0);
        \node at (3.3,2) {\small +};
        \node at (3.3,1) {\small +};
        \node at (3.3,0) {\small -};
      \end{tikzpicture}
    };
        \node (midleft+) at (0,-2) {
      \begin{tikzpicture}[anchorbase,scale=-.5]
        \draw [thick, ->] (0,1.5) -- (1.5,1.5) to [out=0,in=180] (2.25,.5) to [out=0,in=180]  (3,2);
        \draw [white, line width=3] (1,.5) to [out=60,in=180] (3,1);
        \draw [white, line width=3] (1,.5) to [out=-60,in=180] (3,0);
        \draw [double,<-] (0,.5) -- (1,.5);
        \draw [thick] (1,.5) to [out=60,in=180] (3,1);
        \draw [thick] (1,.5) to [out=-60,in=180] (3,0);
        \node at (1.35,1.15) {\tiny $d_5$};
        \node at (3,1.3) {\tiny $d_4$};
        \node at (3.3,2) {\small +};
        \node at (3.3,1) {\small +};
        \node at (3.3,0) {\small -};      \end{tikzpicture}
    };
    \node at (2,-2) {$\rightarrow$};
    \node (midright+) at (4,-2) {
      \begin{tikzpicture}[anchorbase,scale=.5]
        \draw [thick,->] (0,1) to [out=0,in=90] (.5,.5) to [out=-90,in=0] (0,0);
        \draw [thick,->] (.75,.75) arc (-180:180:.35);
        \draw [thick] (0,2) to [out=0,in=120] (2,1.5);
        \draw [thick] (3,.5) -- (2,.5) to [out=180,in=-120] (2,1.5);
        \draw [double,->] (2,1.5) -- (3,1.5);
        \node at (-.3,0) {\small +};
        \node at (-.3,1) {\small +};
        \node at (-.3,2) {\small -};
      \end{tikzpicture}
    };
        \node (midleft-) at (0,-4) {
      \begin{tikzpicture}[anchorbase,scale=-.5]
        \draw [thick, ->] (0,1.5) -- (1.5,1.5) to [out=0,in=180] (2.25,-.5) to [out=0,in=180]  (3,2);
        \draw [white, line width=3] (1,.5) to [out=60,in=180] (3,1);
        \draw [white, line width=3] (1,.5) to [out=-60,in=180] (3,0);
        \draw [double,<-] (0,.5) -- (1,.5);
        \draw [thick] (1,.5) to [out=60,in=180] (3,1);
        \draw [thick] (1,.5) to [out=-60,in=180] (3,0);
        \node at (1.35,1.15) {\tiny $d_5$};
        \node at (1.5,-.2) {\tiny $d_6$};
        \node at (3,-.3) {\tiny $d_3$};
        \node at (3,1.3) {\tiny $d_4$};
        \node at (3.3,2) {\small +};
        \node at (3.3,1) {\small +};
        \node at (3.3,0) {\small -};
      \end{tikzpicture}
    };
    \node at (2,-4) {$\rightarrow$};
    \node (midright-) at (4,-4) {
      \begin{tikzpicture}[anchorbase,scale=.5]
        \fill [black, opacity=.2] (.5,2) -- (.75,2)to [out=60,in=180] (1,2.3) to [out=0,in=120] (1.25,2) -- (1.5,2) to [out=30,in=0] (1.5,2.75) -- (.5,2.75) to [out=180,in=150] (.5,2);
        \draw [thick,->] (0,1) to [out=0,in=90] (.5,.5) to [out=-90,in=0] (0,0);
        \draw [thick] (0,2) -- (.5,2);
        \draw [double] (.5,2) -- (.75,2);
        \draw [thick] (.75,2) to [out=60,in=180] (1,2.3) to [out=0,in=120] (1.25,2);
        \draw [thick] (.75,2) to [out=-60,in=180] (1,1.7) to [out=0,in=-120] (1.25,2);
        \draw [double] (1.25,2) -- (1.5,2);
        \draw [thick] (1.5,2) to [out=30,in=0] (1.5,2.75) -- (.5,2.75) to [out=180,in=150] (.5,2);
        \draw [thick] (1.5,2) to [out=0,in=150] (2.5,1.5);
        \draw [thick] (3,.5) -- (2.5,.5) to [out=180,in=-120] (2.5,1.5);
        \draw [double,->] (2.5,1.5) -- (3,1.5);
        \node at (-.3,0) {\small +};
        \node at (-.3,1) {\small +};
        \node at (-.3,2) {\small -};
      \end{tikzpicture}
    };
        \node (bottomleft) at (0,-6) {
      \begin{tikzpicture}[anchorbase,scale=-.5]
        \draw [thick, ->] (0,1.5) to [out=0,in=180] (1,0) to [out=0,in=180] (2.5,2) --  (3,2);
        \draw [white, line width=4] (0,.5) to [out=0,in=180] (.75,1.5) to [out=0,in=180] (1.25,1);
        \draw [white, line width=3] (1.25,1) to [out=60,in=180] (3,1);
        \draw [white, line width=3] (1.25,1) to [out=-60,in=180] (3,0);
        \draw [double,<-] (0,.5) to [out=0,in=180] (.75,1.5) to [out=0,in=180] (1.25,1);
        \draw [thick] (1.25,1) to [out=60,in=180] (3,1);
        \draw [thick] (1.25,1) to [out=-60,in=180] (3,0);
        \node at (2,.5) {\tiny $d_3$};
        \node at (1.5,1.5) {\tiny $d_4$};
        \node at (0,1) {\tiny $d_2$};
        \node at (3.3,2) {\small +};
        \node at (3.3,1) {\small +};
        \node at (3.3,0) {\small -};
      \end{tikzpicture}
    };
    \node at (2,-6) {$\rightarrow$};
    \node (bottomright) at (4,-6) {
      \begin{tikzpicture}[anchorbase,scale=.5]
        \draw [thick,->] (0,1) to [out=0,in=90] (.5,.5) to [out=-90,in=0] (0,0);
        \draw [thick] (0,2) -- (.5,2);
        \draw [double] (.5,2) -- (.75,2);
        \draw [thick] (.75,2) to [out=60,in=180] (1,2.3) to [out=0,in=120] (1.25,2);
        \draw [thick] (.75,2) to [out=-60,in=180] (1,1.7) to [out=0,in=-120] (1.25,2);
        \draw [double] (1.25,2) -- (1.5,2);
        \draw [thick] (1.5,2) to [out=30,in=0] (1.5,2.75) -- (.5,2.75) to [out=180,in=150] (.5,2);
        \draw [thick] (1.5,2) to [out=0,in=150] (2.5,1.5);
        \draw [thick] (3,.5) -- (2.5,.5) to [out=180,in=-120] (2.5,1.5);
        \draw [double,->] (2.5,1.5) -- (3,1.5);
        \node at (-.3,0) {\small +};
        \node at (-.3,1) {\small +};
        \node at (-.3,2) {\small -};
      \end{tikzpicture}
    };
    \draw [->] (topright) -- (midright+) node[midway,right] {\tiny $\text{saddle/cup}$};
    \draw [->] (midright+) -- (midright-) node[midway,right] {\tiny $f_1 \otimes d_3\wedge d_6$};
    \draw [->] (midright-) -- (bottomright) node[midway,right] {\tiny $a_{11} \id\otimes \langle \cdot, d_6\rangle\wedge d_2$};
  \end{tikzpicture}
\end{gather*}
Above, $f_1$ creates the Gray-colored square. The composite foam equals the one from the LHS on the nose (by an altitude change). As for the coefficient, one reads:
\[
  a_{11}\otimes d_3\wedge d_2=-a_{11}\otimes d_2\wedge d_3
\]
which matches the previous computation since $d_2$ identifies with $c_1$ and $d_3$ with $c_2$.

Symmetric cases (under global orientation reversal, crossing change and mirror image) follow similarly. This completes the check for Move \eqref{StrandAroundVertex}.

\noindent {\bf Sauron:} 
Let us now consider Sauron's move from~\eqref{SlZipSaur}. Let us fix the digon to be as follows:
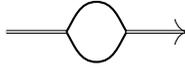
\[
  \begin{tikzpicture}[anchorbase, scale=.4]
    \draw [double] (0,3) -- (2,3);
    \draw [thick] (2,3) to [out=60,in=180] (3,4) to [out=0,in=120] (4,3);
    \draw [thick]  (2,3) to [out=-60,in=180] (3,2) to [out=0,in=-120] (4,3);
    \draw [double,->] (4,3) -- (6,3);
  \end{tikzpicture}
\]

Then the free strand can go over or under, be $1$- or $2$-labeled, and be
oriented up and down. Under the move~\eqref{StrandAroundVertex}, one can reduce
to the cases where it always goes up, leaving us with 4 checks to run. Although
the moves taking place here are not invertible, it appears 
that the space of maps of the correct degree between the first and last
steps is one dimensional. Both versions of the move can only differ by a scalar,
that we compute by going to the fully deformed setting (we compute both
directions: starting from the eye, or finishing with the eye).

Let us start with the following configuration:

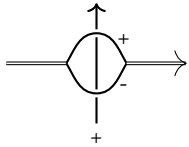
\[
  \begin{tikzpicture}[anchorbase, scale=.4]
    \draw [thick, ->] (3,0) -- (3,4);
    \draw [line width=3,white] (2,2) to [out=60,in=180] (3,3) to [out=0,in=120] (4,2);
    \draw [line width=3,white] (2,2) to [out=-60,in=180] (3,1) to [out=0,in=-120] (4,2);
    \draw [double] (0,2) -- (2,2);
    \draw [thick] (2,2) to [out=60,in=180] (3,3) to [out=0,in=120] (4,2);
    \draw [thick]  (2,2) to [out=-60,in=180] (3,1) to [out=0,in=-120] (4,2);
    \draw [double,->] (4,2) -- (6,2);
    \node at (3,-.5) {\small +};
    \node at (3.9,2.8) {\small +};
    \node at (3.9,1.2) {\small -};
  \end{tikzpicture}
\]

The check amounts to verifying that the coefficients in~\eqref{FS3} that affect the map with the cobordism $F$ equal those of the identity component in~\eqref{FS11}. Similarly, if one lets the strand go under the digon, one has to compare~\eqref{FS1} and~\eqref{FS9}.

Replacing the 1-strand by a 2-strand yields a very similar computation (no need to go down to the fully deformed setting this time), with~\eqref{FS6} and~\eqref{FS14} required to have the same coefficients, and~\eqref{FS5} and~\eqref{FS13} as well.

This concludes the analysis for Sauron's move.

\noindent {\bf Sliding zipper:} We now consider the sliding zipper move from~\eqref{SlZipSaur}. As in Sauron's case, one can reduce to the following situation (with the free strand passing over or under, and being $1$- or $2$-labeled):

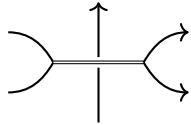
\[
  \begin{tikzpicture}[anchorbase,scale=.4]
    \draw [thick,->] (3,0) -- (3,4);
    \draw [line width=4,white] (1.5,2) -- (4.5,2);
    \draw [thick] (0,1) to [out=0,in=-120] (1.5,2);
    \draw [thick] (0,3) to [out=0,in=120] (1.5,2);
    \draw [double] (1.5,2) -- (4.5,2);
    \draw [thick,->] (4.5,2) to [out=60,in=180] (6,3);
    \draw [thick,->] (4.5,2) to [out=-60,in=180] (6,1);
  \end{tikzpicture}
\]

Again, in the case of a $1$-labeled strand, we go to the fully deformed setting with the following choice:
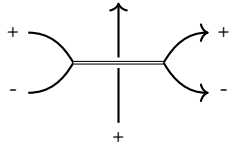
\[
  \begin{tikzpicture}[anchorbase,scale=.4]
    \draw [thick,->] (3,0) -- (3,4);
    \draw [line width=4,white] (1.5,2) -- (4.5,2);
    \draw [thick] (0,1) to [out=0,in=-120] (1.5,2);
    \draw [thick] (0,3) to [out=0,in=120] (1.5,2);
    \draw [double] (1.5,2) -- (4.5,2);
    \draw [thick,->] (4.5,2) to [out=60,in=180] (6,3);
    \draw [thick,->] (4.5,2) to [out=-60,in=180] (6,1);
    \node at (-.5,1) {\small -};
    \node at (-.5,3) {\small +};
    \node at (3,-.5) {\small +};
    \node at (6.5,1) {\small -};
    \node at (6.5,3) {\small +};    
  \end{tikzpicture}
\]

Both sides of the computation yield the same foam, and one wants the coefficients of the components with $F$ in~\eqref{FS3} to be equal to those of the identity components in~\eqref{FS11}. This is true, as it recovers the same condition as for Sauron's move. Similarly if the strand goes under, one should compare~\eqref{FS1} and~\eqref{FS9}.

For the 2-labeled case, one ends up with the same checks as in Sauron's case.

\noindent {\bf Move $C_1$}: We now consider Carter's move~\eqref{C1}. It is made of a cycle around the superposition of two trivalent vertices. Since all maps involved are invertible, it suffices to check from only one preferred starting point. We choose to compare the two maps that go down from the top.

Furthermore, due to the requirement that tangent vectors around a trivalent vertex are parallel, the only configuration is the following one (or the same one with crossing reversed):

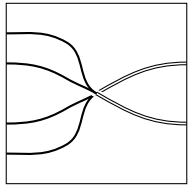
\[

  \]

For each choice of crossing, there are four cases to check, corresponding to two possible orientations for
each of the two trivalent webs.

We argue that the version where the crossing is reversed follows from the one above together with the sliding zipper move from Equation~\eqref{SlZipSaur}. Indeed, considering the following sketched move:
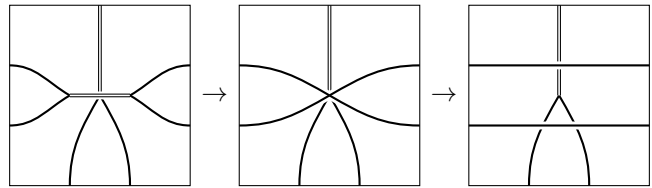
\[
   %
\]
Pushing the upward strand to the right or left uses one or the other instance of the move~\eqref{C1}. Thanks to the sliding zipper move we have already checked, both choices produce the same result.

We thus focus on the picture shown at the beginning of our analysis.  Once again, we will use projective functoriality and a special configuration in fully deformed homology to check the coefficients.

  We first consider:
  \[
   %
\]

On one side one gets:
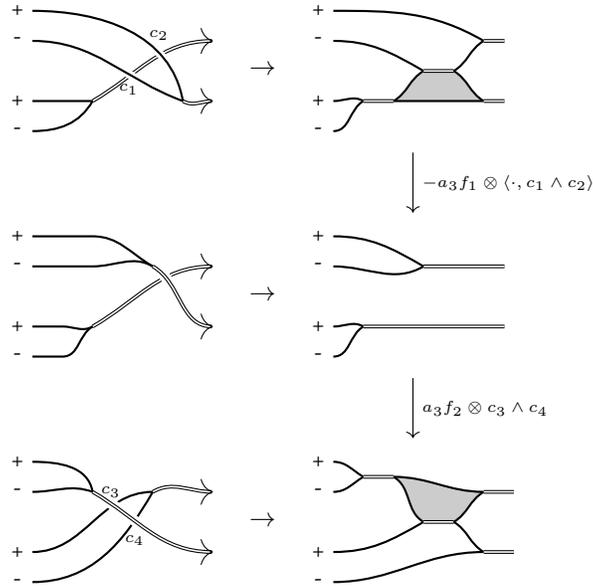
\begin{gather*}
  %
\end{gather*}
Above, $f_1$ and $f_2$ cap off or cup off the shaded area. The coefficient part of the composition is thus $-\langle \cdot, c_1\wedge c_2\rangle \wedge c_3\wedge c_4$.

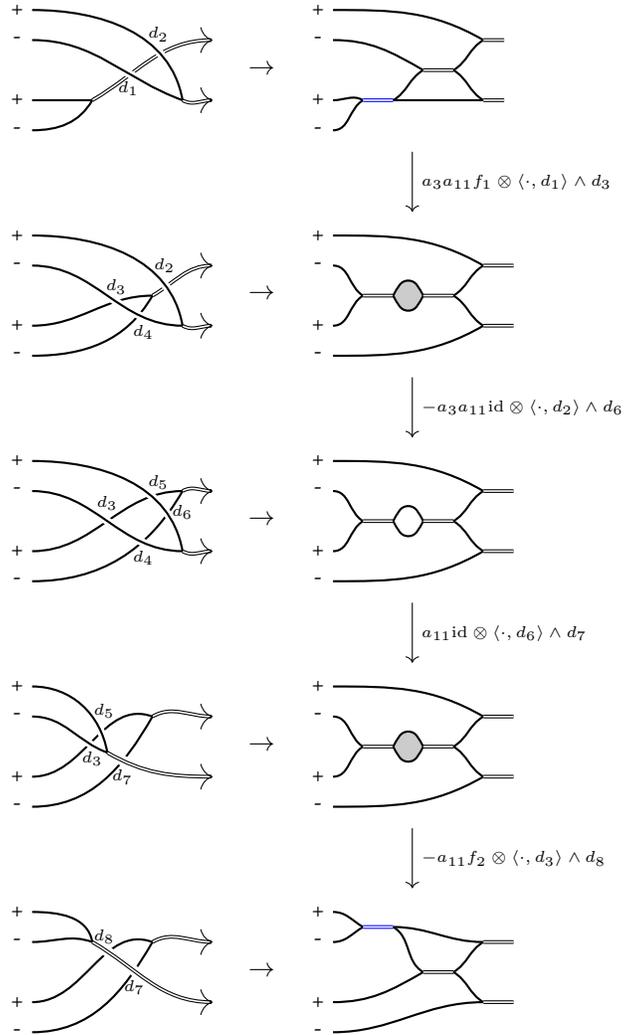
\begin{gather*}
  %
};
\draw [->] (topright) -- (midright+) node[midway,right] {\tiny $a_3a_{11}f_1\otimes \langle \cdot,d_1\rangle\wedge d_3$};
\draw [->] (midright+) -- (midright-) node[midway,right] {\tiny $-a_3a_{11}\id \otimes \langle \cdot,d_2\rangle\wedge d_6$};
\draw [->] (midright-) -- (bottomright) node[midway,right] {\tiny $a_{11}\id \otimes \langle \cdot, d_6\rangle \wedge d_7$};
\draw [->] (bottomright) -- (bottomright-) node[midway,right] {\tiny $-a_{11}f_2 \otimes \langle \cdot, d_3\rangle \wedge d_8$};
    \end{tikzpicture}
\end{gather*}
Above, $f_1$ is the composition of unzipping the blue double strand and cupping off the gray area, and $f_2$ caps off the same gray area before zipping the blue double strand. Altogether, the part of the foam that creates and deletes the digon can be removed, at the expense of a $\frac{1}{2}$ coefficient (see \eqref{deformedbb3}). One then obtains the following (sketched) foam:
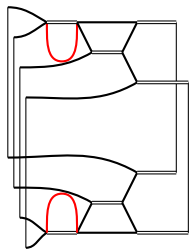
\[
  \begin{tikzpicture}[anchorbase,scale=.4]
    \draw [thick] (.3,0) to [out=0,in=-150] (1,.5);
    \draw [thick] (.1,1) to [out=0,in=150] (1,.5);
    \draw [double] (1,.5) -- (2,.5);
    \draw [thick] (2,.5) -- (2.5,1.5);
    \draw [thick] (-.1,2) to [out=0,in=150] (2.5,1.5);
    \draw [double] (2.5,1.5) -- (3.5,1.5);
    \draw [thick] (3.5,1.5) -- (4,.5);
    \draw [thick] (2,.5) -- (4,.5);
    \draw [thick] (3.5,1.5) -- (4,2.5);
    \draw [thick] (-.3,3) to [out=0,in=150] (4,2.5);
    \draw [double] (4,2.5) -- (5.3,2.5);
    \draw [double] (4,.5) -- (5.7,.5);
    \draw [thick, red] (1,7.5) to [out=-90,in=180] (1.5,6.2) to [out=0,in=-90] (2,7.5); 
    \draw [thick, red] (1,.5) to [out=90,in=180] (1.5,1.8) to [out=0,in=90] (2,.5);
    \draw (.3,0) -- (.3,0+5);
    \draw (.1,1) -- (.1,1+5);
    \draw (-.1,2) -- (-.1,2+5);
    \draw (-.3,3) -- (-.3,3+5);
    \draw (5.3,2.5) -- (5.3,2.5+5);
    \draw (5.7,.5) -- (5.7,.5+5);
    \draw [thick] (-.3,3+5) to [out=0,in=150] (1,2.5+5);
    \draw [thick] (-.1,2+5) to [out=0,in=-150] (1,2.5+5);
    \draw [double] (1,2.5+5) -- (2,2.5+5);
    \draw [thick] (2,2.5+5) -- (2.5,1.5+5);
    \draw [thick] (.1,1+5) to [out=0,in=-150] (2.5,1.5+5);
    \draw [double] (2.5,1.5+5) -- (3.5,1.5+5);
    \draw [thick] (3.5,1.5+5) -- (4,2.5+5);
    \draw [thick] (2,2.5+5) -- (4,2.5+5);
    \draw [thick] (3.5,1.5+5) -- (4,.5+5);
    \draw [thick] (.3,0+5) to [out=0,in=-150] (4,.5+5);
    \draw [double] (4,.5+5) -- (5.7,.5+5);
    \draw [double] (4,2.5+5) -- (5.3,2.5+5);
  \end{tikzpicture}
\]

In the neighborhood of the top and bottom webs, one can apply relation~\eqref{sl2sq2}. One gets two negative signs that cancel, resulting in the same cap/cup foam as in the first part of the move, with an extra 2-labeled disk as on the RHS of Relation~\eqref{deformedsheetedcyl}. Removing it produces a $-2$.

At that point, both foams agree, and the coefficient of the one from the second part of the move reads as:
\[
  \frac{1}{2}(-2)a_3a_{11}(-a_3a_{11})a_{11}(-a_{11})\langle\langle \langle\langle \cdot, d_1\rangle \wedge d_3,d_2\rangle \wedge d_6,d_6\rangle\wedge d_7,d_3\rangle \wedge d_8
\]
This equals:
\[
  -\langle \cdot,d_1\wedge d_2\rangle \wedge d_8\wedge d_7.
\]
This matches the first part of the computation, under the identification of the crossings (note that here $d_7$ corresponds to $c_4$ while $d_8$ corresponds to $c_3$).

A global reversal of orientation causes no change in the foams, except for the orientation of seams, and one can follow the same computation, being careful with all coefficients.

The other two orientations can also be treated together. They appear to be a bit easier as one obtains the same foams on the nose on both sides of the computation.

\noindent {\bf Move $C_3$:} One can use the following isotopy and combine $C_3$ moves with Sauron's move to constrain the orientation of the trivalent vertex, for example to a merge.

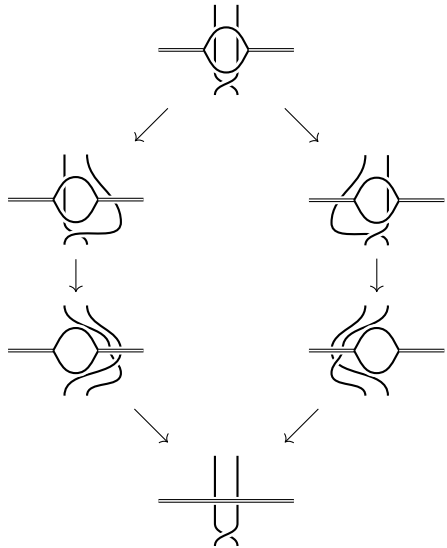
\[

\]

Then using the following isotopy one can restrict the number of orientations to be considered.
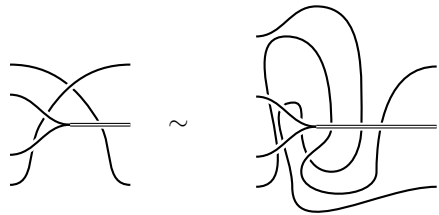
\[
  %
\]

There are thus eight cases to check, which we do by explicit computations.
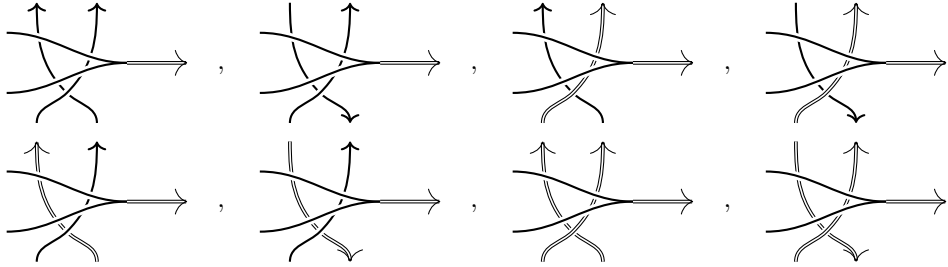
\begin{gather*}
  %
\end{gather*}

We start with the first version, and go to fully deformed homology with the following choice of projectors:
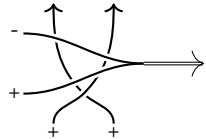
\[
%
\]

The first part of the computation goes as follows:
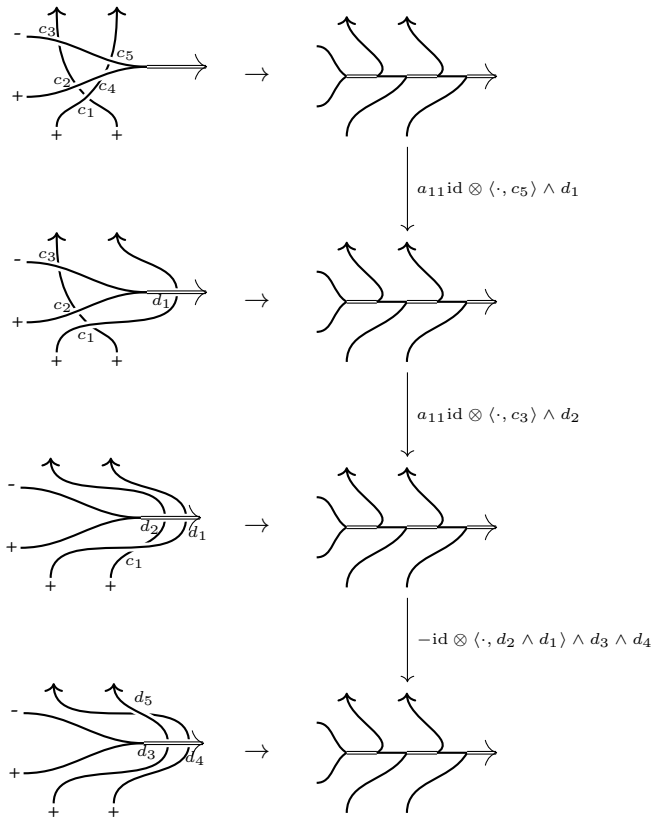
\[
  %
\]
Altogether, one gets: $-\id \otimes \langle \cdot,c_3\wedge c_5\rangle \wedge d_3\wedge d_4$.

The second side goes as follows:

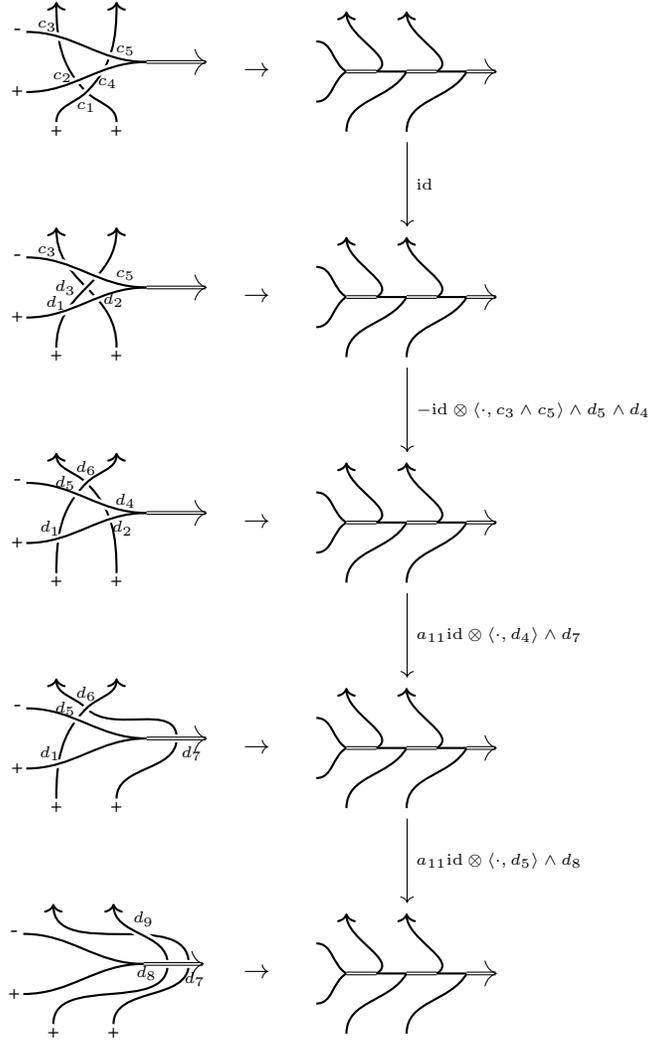
\[
  %
    };
    \draw [->] (topright) -- (midright1) node[midway,right] {\tiny $\id$};
    \draw [->] (midright1) -- (midright2) node[midway,right] {\tiny $-\id\otimes \langle \cdot, c_3\wedge c_5\rangle \wedge d_5\wedge d_4$};
    \draw [->] (midright2) -- (midright3) node[midway,right] {\tiny $a_{11}\id\otimes \langle \cdot, d_4\rangle \wedge d_7$};
    \draw [->] (midright3) -- (bottomright) node[midway,right] {\tiny $a_{11}\id\otimes \langle \cdot, d_5\rangle \wedge d_8$};
  \end{tikzpicture}
\]
Altogether, one gets: $\id \otimes \langle \cdot,c_3\wedge c_5\rangle \wedge d_7\wedge d_8$. This is consistent with the first computation, after identifying $d_7$ with $d_4$ and $d_8$ with $d_3$.

We now look at:
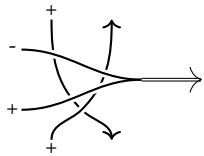
\[
\begin{tikzpicture}[anchorbase,scale=.4]
    \draw [thick,<-] (3,0) to [out=90,in=-45] (2,1) to [out=135,in=-90] (1,4);
    \draw [white, line width=3] (1,0) to [out=90,in=-135] (2,1) to [out=45,in=-90] (3,4);
    \draw [thick,->] (1,0) to [out=90,in=-135] (2,1) to [out=45,in=-90] (3,4);
    \draw [line width=3, white] (0,3) to [out=0,in=180] (4,2);
    \draw [line width=3, white] (0,1) to [out=0,in=180] (4,2);
    \draw [line width=3, white] (4,2) -- (6,2);
    \draw [thick] (0,3) to [out=0,in=180] (4,2);
    \draw [thick] (0,1) to [out=0,in=180] (4,2);
    \draw [double,->] (4,2) -- (6,2);
    \node at (1,-.3) {\small +};
    \node at (1,4.3) {\small +};
    \node at (-.3,1) {\small +};
    \node at (-.3,3) {\small -};
  \end{tikzpicture}
\]

The first part of the computation goes as follows:
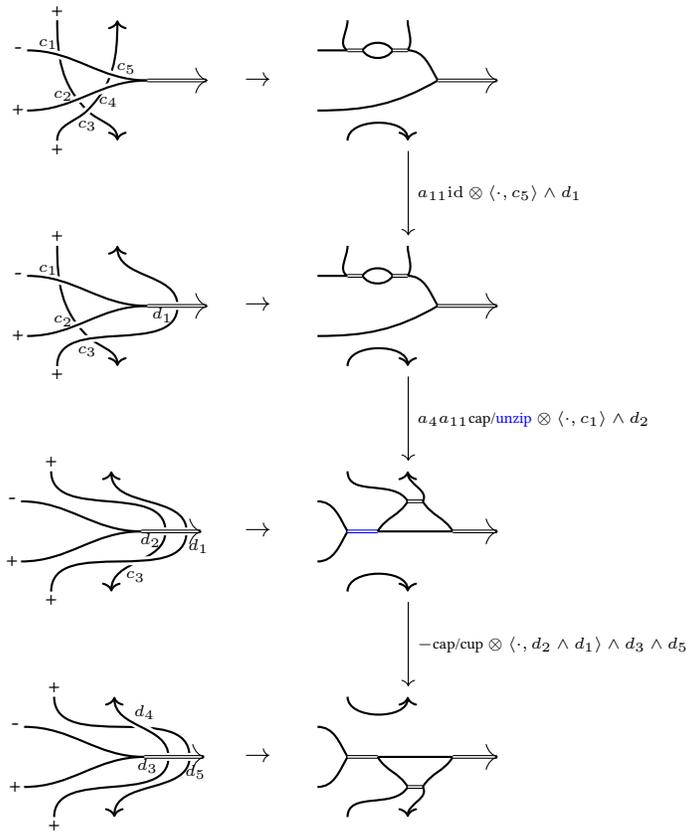
\[
  \begin{tikzpicture}[anchorbase]
    \node (topleft) at (0,0) {
      \begin{tikzpicture}[anchorbase,scale=.4]
        \draw [thick,<-] (3,0) to [out=90,in=-45] (2,1) to [out=135,in=-90] (1,4);
        \draw [white, line width=3] (1,0) to [out=90,in=-135] (2,1) to [out=45,in=-90] (3,4);
        \draw [thick,->] (1,0) to [out=90,in=-135] (2,1) to [out=45,in=-90] (3,4);
        \draw [line width=3, white] (0,3) to [out=0,in=180] (4,2);
        \draw [line width=3, white] (0,1) to [out=0,in=180] (4,2);
        \draw [line width=3, white] (4,2) -- (6,2);
        \draw [thick] (0,3) to [out=0,in=180] (4,2);
        \draw [thick] (0,1) to [out=0,in=180] (4,2);
        \draw [double,->] (4,2) -- (6,2);
        \node at (1,-.3) {\small +};
        \node at (1,4.3) {\small +};
        \node at (-.3,1) {\small +};
        \node at (-.3,3) {\small -};
        \node at (2,.5) {\tiny $c_3$};
        \node at (1.2,1.5) {\tiny $c_2$};
        \node at (.7,3.2) {\tiny $c_1$};
        \node at (2.7,1.3) {\tiny $c_4$};
        \node at (3.3,2.4) {\tiny $c_5$};
      \end{tikzpicture}
    };
    \node at (2,0) {$\rightarrow$};
    \node (topright) at (4,0) {
      \begin{tikzpicture}[anchorbase,scale=.4]
        \draw [thick] (0,3) to [out=0,in=180] (1,3);
        \draw [thick] (1,4) to [out=-90,in=150] (1,3);
        \draw [double] (1,3) -- (1.5,3);
        \draw [thick] (1.5,3) to [out=60,in=120] (2.5,3);
        \draw [thick] (1.5,3) to [out=-60,in=-120] (2.5,3);
        \draw [double] (2.5,3) -- (3,3);
        \draw [thick] (3,3) to [out=60,in=-90] (3,4);
        \draw [thick] (3,3) to [out=0,in=120] (4,2);
        \draw [thick] (0,1) to [out=0,in=-150] (4,2);
        \draw [double,->] (4,2) -- (6,2);
        \draw [thick,->] (1,0) to [out=90,in=90] (3,0);
      \end{tikzpicture}
      };
    \node at (0,-3) (midleft1) {
      \begin{tikzpicture}[anchorbase,scale=.4]
        \draw [thick,<-] (3,0) to [out=90,in=-45] (2,1) to [out=135,in=-90] (1,4);
        \draw [white, line width=3](1,0) to [out=90,in=-90] (5,2) to [out=90,in=-90] (3,4);
        \draw [thick,->] (1,0) to [out=90,in=-90] (5,2) to [out=90,in=-90] (3,4);
        \draw [line width=3, white] (0,3) to [out=0,in=180] (4,2);
        \draw [line width=3, white] (0,1) to [out=0,in=180] (4,2);
        \draw [line width=3, white] (4,2) -- (6,2);
        \draw [thick] (0,3) to [out=0,in=180] (4,2);
        \draw [thick] (0,1) to [out=0,in=180] (4,2);
        \draw [double,->] (4,2) -- (6,2);
        \node at (1,-.3) {\small +};
        \node at (1,4.3) {\small +};
        \node at (-.3,1) {\small +};
        \node at (-.3,3) {\small -};
        \node at (2,.5) {\tiny $c_3$};
        \node at (1.2,1.5) {\tiny $c_2$};
        \node at (.7,3.2) {\tiny $c_1$};
        \node at (4.5,1.7) {\tiny $d_1$};
             \end{tikzpicture}
      };
    \node at (2,-3) {$\rightarrow$};
    \node (midright1) at (4,-3) {
      \begin{tikzpicture}[anchorbase,scale=.4]
        \draw [thick] (0,3) to [out=0,in=180] (1,3);
        \draw [thick] (1,4) to [out=-90,in=150] (1,3);
        \draw [double] (1,3) -- (1.5,3);
        \draw [thick] (1.5,3) to [out=60,in=120] (2.5,3);
        \draw [thick] (1.5,3) to [out=-60,in=-120] (2.5,3);
        \draw [double] (2.5,3) -- (3,3);
        \draw [thick] (3,3) to [out=60,in=-90] (3,4);
        \draw [thick] (3,3) to [out=0,in=120] (4,2);
        \draw [thick] (0,1) to [out=0,in=-150] (4,2);
        \draw [double,->] (4,2) -- (6,2);
        \draw [thick,->] (1,0) to [out=90,in=90] (3,0);
      \end{tikzpicture}
          };
    \node at (0,-6) (midleft2) {
      \begin{tikzpicture}[anchorbase,scale=.4]
        \draw [thick,<-] (3,0) to [out=90,in=-90] (4.8,2) to [out=90,in=-90] (1,4);
        \draw [white, line width=3] (1,0) to [out=90,in=-90] (5.5,2) to [out=90,in=-90] (3,4);
        \draw [thick,->] (1,0) to [out=90,in=-90] (5.5,2) to [out=90,in=-90] (3,4);
        \draw [line width=3, white] (0,3) to [out=0,in=180] (4,2);
        \draw [line width=3, white] (0,1) to [out=0,in=180] (4,2);
        \draw [line width=3, white] (4,2) -- (6,2);
        \draw [thick] (0,3) to [out=0,in=180] (4,2);
        \draw [thick] (0,1) to [out=0,in=180] (4,2);
        \draw [double,->] (4,2) -- (6,2);
        \node at (1,-.3) {\small +};
        \node at (1,4.3) {\small +};
        \node at (-.3,1) {\small +};
        \node at (-.3,3) {\small -};
        \node at (3.8,.5) {\tiny $c_3$};
        \node at (4.3,1.7) {\tiny $d_2$};
        \node at (5.9,1.5) {\tiny $d_1$};
             \end{tikzpicture}
    };
    \node at (2,-6) {$\rightarrow$};
    \node (midright2) at (4,-6) {
      \begin{tikzpicture}[anchorbase,scale=.4]
        \draw [thick] (0,1) to [out=0,in=-120] (1,2);
        \draw [thick]  (0,3) to [out=0,in=120] (1,2);
        \draw [double,blue] (1,2) -- (2,2);
        \draw [thick] (2,2) to [out=60,in=-120] (3,3);
        \draw [thick] (1,4) to [out=-90,in=120] (3,3);
        \draw [double] (3,3) -- (3.5,3);
        \draw [thick,->] (3.5,3) to [out=60,in=-90] (3,4);
        \draw [thick] (3.5,3) to [out=-60,in=120] (4.5,2);
        \draw [thick] (2,2) -- (4.5,2);
        \draw [double,->] (4.5,2) -- (6,2);
        \draw [thick,->] (1,0) to [out=90,in=90] (3,0);
      \end{tikzpicture}
      };
   \node at (0,-9) (bottomleft) {
      \begin{tikzpicture}[anchorbase,scale=.4]
        \draw [thick,<-] (3,0) to [out=90,in=-90] (5.5,2) to [out=90,in=-90] (1,4);
        \draw [white, line width=3] (1,0) to [out=90,in=-90] (4.8,2) to [out=90,in=-90] (3,4);
        \draw [thick,->] (1,0) to [out=90,in=-90] (4.8,2) to [out=90,in=-90] (3,4);
        \draw [line width=3, white] (0,3) to [out=0,in=180] (4,2);
        \draw [line width=3, white] (0,1) to [out=0,in=180] (4,2);
        \draw [line width=3, white] (4,2) -- (6,2);
        \draw [thick] (0,3) to [out=0,in=180] (4,2);
        \draw [thick] (0,1) to [out=0,in=180] (4,2);
        \draw [double,->] (4,2) -- (6,2);
        \node at (1,-.3) {\small +};
        \node at (1,4.3) {\small +};
        \node at (-.3,1) {\small +};
        \node at (-.3,3) {\small -};
        \node at (4.1,1.7) {\tiny $d_3$};
        \node at (5.7,1.5) {\tiny $d_5$};
        \node at (4,3.5) {\tiny $d_4$};
      \end{tikzpicture}
    };
    \node at (2,-9) {$\rightarrow$};
    \node (bottomright) at (4,-9) {
      \begin{tikzpicture}[anchorbase,scale=.4]
        \draw [thick] (0,1) to [out=0,in=-120] (1,2);
        \draw [thick]  (0,3) to [out=0,in=120] (1,2);
        \draw [double] (1,2) -- (2,2);
        \draw [thick] (2,2) to [out=-60,in=120] (3,1);
        \draw [thick] (1,0) to [out=90,in=-120] (3,1);
        \draw [double] (3,1) -- (3.5,1);
        \draw [thick,->] (3.5,1) to [out=-60,in=90] (3,0);
        \draw [thick] (3.5,1) to [out=60,in=-120] (4.5,2);
        \draw [thick] (2,2) -- (4.5,2);
        \draw [double,->] (4.5,2) -- (6,2);
        \draw [thick,->] (1,4) to [out=-90,in=-90] (3,4);
      \end{tikzpicture}
    };
    \draw [->] (topright) -- (midright1) node[midway,right] {\tiny $a_{11}\id\otimes \langle \cdot,c_5\rangle\wedge d_1$};
    \draw [->] (midright1) -- (midright2) node[midway,right] {\tiny $a_4a_{11}\text{cap/\textcolor{blue}{unzip}}\otimes \langle \cdot, c_1\rangle \wedge d_2$};
      \draw [->] (midright2) -- (bottomright) node[midway,right] {\tiny $-\text{cap/cup}\otimes \langle \cdot, d_2\wedge d_1\rangle \wedge d_3\wedge d_5$};
  \end{tikzpicture}
\]
Altogether, one gets as a coefficient of the foam : $-a_4 \langle \cdot,c_1\wedge c_5\rangle \wedge d_3\wedge d_5$.

The second side goes as follows:

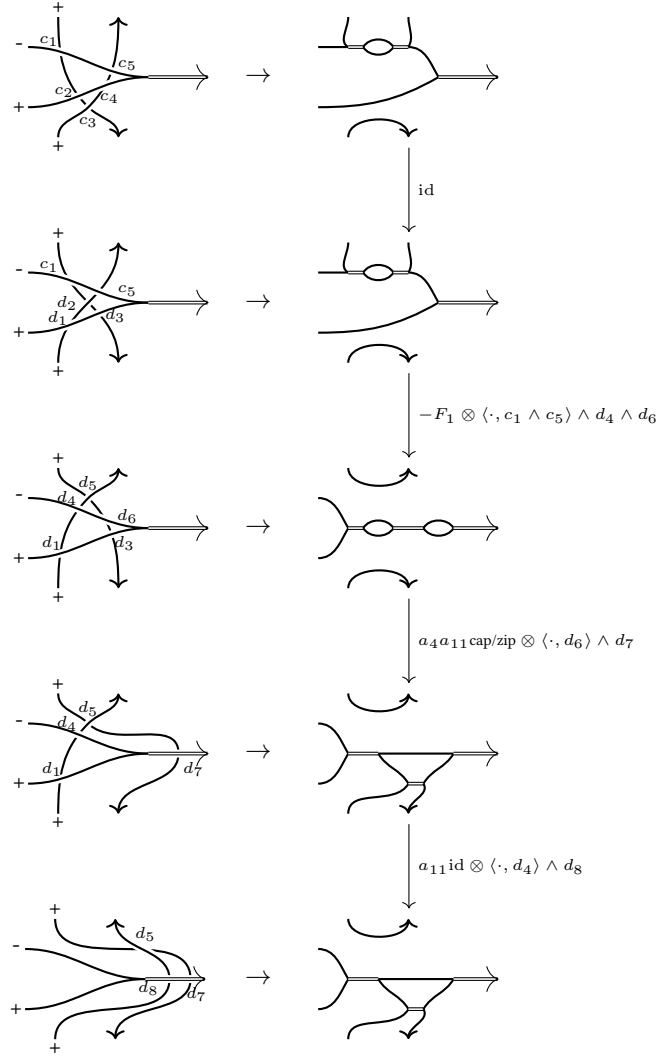
\[
  \begin{tikzpicture}[anchorbase]
    \node (topleft) at (0,0) {
      \begin{tikzpicture}[anchorbase,scale=.4]
        \draw [thick,<-] (3,0) to [out=90,in=-45] (2,1) to [out=135,in=-90] (1,4);
        \draw [white, line width=3] (1,0) to [out=90,in=-135] (2,1) to [out=45,in=-90] (3,4);
        \draw [thick,->] (1,0) to [out=90,in=-135] (2,1) to [out=45,in=-90] (3,4);
        \draw [line width=3, white] (0,3) to [out=0,in=180] (4,2);
        \draw [line width=3, white] (0,1) to [out=0,in=180] (4,2);
        \draw [line width=3, white] (4,2) -- (6,2);
        \draw [thick] (0,3) to [out=0,in=180] (4,2);
        \draw [thick] (0,1) to [out=0,in=180] (4,2);
        \draw [double,->] (4,2) -- (6,2);
        \node at (1,-.3) {\small +};
        \node at (1,4.3) {\small +};
        \node at (-.3,1) {\small +};
        \node at (-.3,3) {\small -};
        \node at (2,.5) {\tiny $c_3$};
        \node at (1.2,1.5) {\tiny $c_2$};
        \node at (.7,3.2) {\tiny $c_1$};
        \node at (2.7,1.3) {\tiny $c_4$};
        \node at (3.3,2.4) {\tiny $c_5$};
      \end{tikzpicture}
    };
    \node at (2,0) {$\rightarrow$};
    \node (topright) at (4,0) {
       \begin{tikzpicture}[anchorbase,scale=.4]
        \draw [thick] (0,3) to [out=0,in=180] (1,3);
        \draw [thick] (1,4) to [out=-90,in=150] (1,3);
        \draw [double] (1,3) -- (1.5,3);
        \draw [thick] (1.5,3) to [out=60,in=120] (2.5,3);
        \draw [thick] (1.5,3) to [out=-60,in=-120] (2.5,3);
        \draw [double] (2.5,3) -- (3,3);
        \draw [thick] (3,3) to [out=60,in=-90] (3,4);
        \draw [thick] (3,3) to [out=0,in=120] (4,2);
        \draw [thick] (0,1) to [out=0,in=-150] (4,2);
        \draw [double,->] (4,2) -- (6,2);
        \draw [thick,->] (1,0) to [out=90,in=90] (3,0);
      \end{tikzpicture}
      };
    \node at (0,-3) (midleft1) {
      \begin{tikzpicture}[anchorbase,scale=.4]
        \draw [thick,<-] (3,0) to [out=90,in=-45] (2,2) to [out=135,in=-90] (1,4);
        \draw [white, line width=3] (1,0) to [out=90,in=-135] (2,2) to [out=45,in=-90] (3,4);
        \draw [thick,->] (1,0) to [out=90,in=-135] (2,2) to [out=45,in=-90] (3,4);
        \draw [line width=3, white] (0,3) to [out=0,in=180] (4,2);
        \draw [line width=3, white] (0,1) to [out=0,in=180] (4,2);
        \draw [line width=3, white] (4,2) -- (6,2);
        \draw [thick] (0,3) to [out=0,in=180] (4,2);
        \draw [thick] (0,1) to [out=0,in=180] (4,2);
        \draw [double,->] (4,2) -- (6,2);
        \node at (1,-.3) {\small +};
        \node at (1,4.3) {\small +};
        \node at (-.3,1) {\small +};
        \node at (-.3,3) {\small -};
        \node at (1,1.5) {\tiny $d_1$};
        \node at (2.9,1.6) {\tiny $d_3$};
        \node at (.7,3.2) {\tiny $c_1$};
        \node at (1.3,2) {\tiny $d_2$};
        \node at (3.3,2.4) {\tiny $c_5$};
      \end{tikzpicture}}
    ;
    \node at (2,-3) {$\rightarrow$};
    \node (midright1) at (4,-3) {
       \begin{tikzpicture}[anchorbase,scale=.4]
        \draw [thick] (0,3) to [out=0,in=180] (1,3);
        \draw [thick] (1,4) to [out=-90,in=150] (1,3);
        \draw [double] (1,3) -- (1.5,3);
        \draw [thick] (1.5,3) to [out=60,in=120] (2.5,3);
        \draw [thick] (1.5,3) to [out=-60,in=-120] (2.5,3);
        \draw [double] (2.5,3) -- (3,3);
        \draw [thick] (3,3) to [out=60,in=-90] (3,4);
        \draw [thick] (3,3) to [out=0,in=120] (4,2);
        \draw [thick] (0,1) to [out=0,in=-150] (4,2);
        \draw [double,->] (4,2) -- (6,2);
        \draw [thick,->] (1,0) to [out=90,in=90] (3,0);
      \end{tikzpicture}
      };
    \node at (0,-6) (midleft2) {
      \begin{tikzpicture}[anchorbase,scale=.4]
        \draw [thick,<-] (3,0) to [out=90,in=-45] (2,3) to [out=135,in=-90] (1,4);
        \draw [white, line width=3] (1,0) to [out=90,in=-135] (2,3) to [out=45,in=-90] (3,4);
        \draw [thick,->] (1,0) to [out=90,in=-135] (2,3) to [out=45,in=-90] (3,4);
        \draw [line width=3, white] (0,3) to [out=0,in=180] (4,2);
        \draw [line width=3, white] (0,1) to [out=0,in=180] (4,2);
        \draw [line width=3, white] (4,2) -- (6,2);
        \draw [thick] (0,3) to [out=0,in=180] (4,2);
        \draw [thick] (0,1) to [out=0,in=180] (4,2);
        \draw [double,->] (4,2) -- (6,2);
        \node at (1,-.3) {\small +};
        \node at (1,4.3) {\small +};
        \node at (-.3,1) {\small +};
        \node at (-.3,3) {\small -};
        \node at (.8,1.5) {\tiny $d_1$};
        \node at (3.2,1.5) {\tiny $d_3$};
        \node at (3.3,2.4) {\tiny $d_6$};
        \node at (1.3,3) {\tiny $d_4$};
        \node at (2,3.5) {\tiny $d_5$};
      \end{tikzpicture}
    };
    \node at (2,-6) {$\rightarrow$};
    \node (midright2) at (4,-6) {
      \begin{tikzpicture}[anchorbase,scale=.4]
        \draw [thick] (0,1) to [out=0,in=-120] (1,2);
        \draw [thick] (0,3) to [out=0,in=120] (1,2);
        \draw [double] (1,2) -- (1.5,2);
        \draw [thick] (1.5,2) to [out=60,in=120] (2.5,2);
        \draw [thick] (1.5,2) to [out=-60,in=-120] (2.5,2);
        \draw [double] (2.5,2) -- (3.5,2);
        \draw [thick] (3.5,2) to [out=60,in=120] (4.5,2);
        \draw [thick] (3.5,2) to [out=-60,in=-120] (4.5,2);
        \draw [double,->] (4.5,2) -- (6,2);
        \draw [thick,->] (1,4) to [out=-90,in=-90] (3,4);
        \draw [thick,->] (1,0) to [out=90,in=90] (3,0);
      \end{tikzpicture}
      };
    \node at (0,-9) (midleft3) {
      \begin{tikzpicture}[anchorbase,scale=.4]
        \draw [thick,<-] (3,0) to [out=90,in=-90] (5,2) to [out=90,in=-45] (2,3) to [out=135,in=-90] (1,4);
        \draw [white, line width=3] (1,0) to [out=90,in=-135] (2,3) to [out=45,in=-90] (3,4);
        \draw [thick,->] (1,0) to [out=90,in=-135] (2,3) to [out=45,in=-90] (3,4);
        \draw [line width=3, white] (0,3) to [out=0,in=180] (4,2);
        \draw [line width=3, white] (0,1) to [out=0,in=180] (4,2);
        \draw [line width=3, white] (4,2) -- (6,2);
        \draw [thick] (0,3) to [out=0,in=180] (4,2);
        \draw [thick] (0,1) to [out=0,in=180] (4,2);
        \draw [double,->] (4,2) -- (6,2);
        \node at (1,-.3) {\small +};
        \node at (1,4.3) {\small +};
        \node at (-.3,1) {\small +};
        \node at (-.3,3) {\small -};
        \node at (.8,1.5) {\tiny $d_1$};
        \node at (5.5,1.5) {\tiny $d_7$};
        \node at (1.3,3) {\tiny $d_4$};
        \node at (2,3.5) {\tiny $d_5$};
      \end{tikzpicture}
    };
    \node at (2,-9) {$\rightarrow$};
    \node (midright3) at (4,-9) {
  \begin{tikzpicture}[anchorbase,scale=.4]
        \draw [thick] (0,1) to [out=0,in=-120] (1,2);
        \draw [thick]  (0,3) to [out=0,in=120] (1,2);
        \draw [double] (1,2) -- (2,2);
        \draw [thick] (2,2) to [out=-60,in=120] (3,1);
        \draw [thick] (1,0) to [out=90,in=-120] (3,1);
        \draw [double] (3,1) -- (3.5,1);
        \draw [thick,->] (3.5,1) to [out=-60,in=90] (3,0);
        \draw [thick] (3.5,1) to [out=60,in=-120] (4.5,2);
        \draw [thick] (2,2) -- (4.5,2);
        \draw [double,->] (4.5,2) -- (6,2);
        \draw [thick,->] (1,4) to [out=-90,in=-90] (3,4);
      \end{tikzpicture}    
      };
   \node at (0,-12) (bottomleft) {
      \begin{tikzpicture}[anchorbase,scale=.4]
        \draw [thick,<-] (3,0) to [out=90,in=-90] (5.5,2) to [out=90,in=-90] (1,4);
        \draw [white, line width=3] (1,0) to [out=90,in=-90] (4.8,2) to [out=90,in=-90] (3,4);
        \draw [thick,->] (1,0) to [out=90,in=-90] (4.8,2) to [out=90,in=-90] (3,4);
        \draw [line width=3, white] (0,3) to [out=0,in=180] (4,2);
        \draw [line width=3, white] (0,1) to [out=0,in=180] (4,2);
        \draw [line width=3, white] (4,2) -- (6,2);
        \draw [thick] (0,3) to [out=0,in=180] (4,2);
        \draw [thick] (0,1) to [out=0,in=180] (4,2);
        \draw [double,->] (4,2) -- (6,2);
        \node at (1,-.3) {\small +};
        \node at (1,4.3) {\small +};
        \node at (-.3,1) {\small +};
        \node at (-.3,3) {\small -};
        \node at (4.1,1.7) {\tiny $d_8$};
        \node at (5.7,1.5) {\tiny $d_7$};
        \node at (4,3.5) {\tiny $d_5$};
      \end{tikzpicture}
    };
    \node at (2,-12) {$\rightarrow$};
    \node (bottomright) at (4,-12) {
  \begin{tikzpicture}[anchorbase,scale=.4]
        \draw [thick] (0,1) to [out=0,in=-120] (1,2);
        \draw [thick]  (0,3) to [out=0,in=120] (1,2);
        \draw [double] (1,2) -- (2,2);
        \draw [thick] (2,2) to [out=-60,in=120] (3,1);
        \draw [thick] (1,0) to [out=90,in=-120] (3,1);
        \draw [double] (3,1) -- (3.5,1);
        \draw [thick,->] (3.5,1) to [out=-60,in=90] (3,0);
        \draw [thick] (3.5,1) to [out=60,in=-120] (4.5,2);
        \draw [thick] (2,2) -- (4.5,2);
        \draw [double,->] (4.5,2) -- (6,2);
        \draw [thick,->] (1,4) to [out=-90,in=-90] (3,4);
      \end{tikzpicture}
    };
    \draw [->] (topright) -- (midright1) node[midway,right] {\tiny $\id$};
    \draw [->] (midright1) -- (midright2) node[midway,right] {\tiny $-F_1\otimes \langle \cdot, c_1\wedge c_5\rangle \wedge d_4\wedge d_6$};
    \draw [->] (midright2) -- (midright3) node[midway,right] {\tiny $a_4a_{11}\text{cap/zip}\otimes \langle \cdot, d_6\rangle \wedge d_7$};
    \draw [->] (midright3) -- (bottomright) node[midway,right] {\tiny $a_{11}\id\otimes \langle \cdot, d_4\rangle \wedge d_8$};
  \end{tikzpicture}
\]
Altogether, one gets the same foam on the nose, and as a coefficient $\langle \cdot,c_1\wedge c_5\rangle \wedge d_7\wedge d_8$. This is consistent with the first computation, after identifying $d_7$ with $d_5$ and $d_8$ with $d_3$.

Checks for other labelings are run following the same strategy. We omit the details.

\noindent {\bf Move $C_4$:} The checks for the move~\eqref{C4} are parallel to the ones for~\eqref{C3}. The foams and the crossings numbering can be kept, and one simply looks for differences in the coefficients. We omit the details.

\noindent {\bf Move $C_2$:} 
We now consider the move~\eqref{C2}. It corresponds to the following non-generic situations:
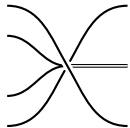
\[
  \begin{tikzpicture}[anchorbase,scale=.4]
    \draw [thick] (0,0) to [out=0,in=-120] (2,2) to [out=60,in=180] (4,4);
    \draw [line width=3, white] (0,3) to [out=0,in=180] (2,2);
    \draw [line width=3, white] (0,1) to [out=0,in=180] (2,2);
    \draw [line width=3, white] (2,2) -- (4,2);
    \draw [thick] (0,3) to [out=0,in=180] (2,2);
    \draw [thick] (0,1) to [out=0,in=180] (2,2);
    \draw [double] (2,2) -- (4,2);
    \draw [line width=3, white] (0,4) to [out=0,in=120] (2,2) to [out=-60,in=180] (4,0);
    \draw [thick] (0,4) to [out=0,in=120] (2,2) to [out=-60,in=180] (4,0);
  \end{tikzpicture}
\]
where the free strands are at the top and at the bottom respectively, and can be
$1$ or $2$ colored.

By looking at the following singular situation:
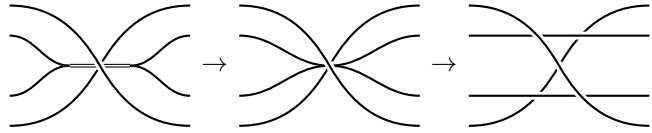
\[
  \begin{tikzpicture}[anchorbase,scale=.4]
    \draw [thick] (0,0) to [out=0,in=-120] (3,2) to [out=60,in=180] (6,4);
    \draw [line width=3, white] (0,3) to [out=0,in=180] (2,2);
    \draw [line width=3, white] (0,1) to [out=0,in=180] (2,2);
    \draw [line width=3, white] (2,2) -- (4,2);
    \draw [line width=3, white] (4,2) to [out=0,in=180] (6,3);
    \draw [line width=3, white] (4,2) to [out=0,in=180] (6,1);
    \draw [thick] (0,3) to [out=0,in=180] (2,2);
    \draw [thick] (0,1) to [out=0,in=180] (2,2);
    \draw [double] (2,2) -- (4,2);
    \draw [thick] (4,2) to [out=0,in=180] (6,3);
    \draw [thick] (4,2) to [out=0,in=180] (6,1);
    \draw [line width=3, white] (0,4) to [out=0,in=120] (3,2) to [out=-60,in=180] (6,0);
    \draw [thick] (0,4) to [out=0,in=120] (3,2) to [out=-60,in=180] (6,0);
  \end{tikzpicture}
\rightarrow
  \begin{tikzpicture}[anchorbase,scale=.4]
    \draw [thick] (0,0) to [out=0,in=-120] (3,2) to [out=60,in=180] (6,4);
    \draw [line width=3, white] (0,3) to [out=0,in=180] (2,2);
    \draw [line width=3, white] (0,1) to [out=0,in=180] (2,2);
    \draw [line width=3, white] (3,2) to [out=0,in=180] (6,3);
    \draw [line width=3, white] (3,2) to [out=0,in=180] (6,1);
    \draw [thick] (0,3) to [out=0,in=180] (3,2);
    \draw [thick] (0,1) to [out=0,in=180] (3,2);
    \draw [thick] (3,2) to [out=0,in=180] (6,3);
    \draw [thick] (3,2) to [out=0,in=180] (6,1);
    \draw [line width=3, white] (0,4) to [out=0,in=120] (3,2) to [out=-60,in=180] (6,0);
    \draw [thick] (0,4) to [out=0,in=120] (3,2) to [out=-60,in=180] (6,0);
  \end{tikzpicture}
\rightarrow
  \begin{tikzpicture}[anchorbase,scale=.4]
    \draw [thick] (0,0) to [out=0,in=-120] (3,2) to [out=60,in=180] (6,4);
    \draw [line width=3, white] (0,3) -- (6,3);
    \draw [line width=3, white] (0,1) -- (6,1);
    \draw [thick] (0,3) -- (6,3);
    \draw [thick] (0,1) -- (6,1);
    \draw [line width=3, white] (0,4) to [out=0,in=120] (3,2) to [out=-60,in=180] (6,0);
    \draw [thick] (0,4) to [out=0,in=120] (3,2) to [out=-60,in=180] (6,0);
  \end{tikzpicture}
\]
one can pass from the merge to the split case, using~\ref{SlZipSaur}. It thus suffices to prove one of the two cases.

Furthermore, one can also constrain the number of orientations to be checked using the following isotopy:
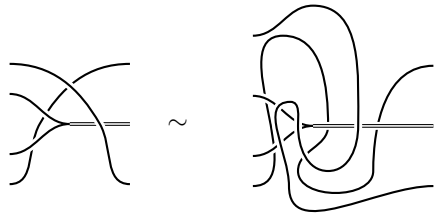
\[
  \begin{tikzpicture}[anchorbase,scale=.4]
    \draw [thick] (0,0) to [out=0,in=-120] (1,2) to [out=60,in=180] (4,4);
    \draw [line width=3, white] (0,3) to [out=0,in=180] (2,2);
    \draw [line width=3, white] (0,1) to [out=0,in=180] (2,2);
    \draw [line width=3, white] (2,2) -- (4,2);
    \draw [thick] (0,3) to [out=0,in=180] (2,2);
    \draw [thick] (0,1) to [out=0,in=180] (2,2);
    \draw [double] (2,2) -- (4,2);
    \draw [line width=3, white] (0,4) to [out=0,in=120] (3,2) to [out=-60,in=180] (4,0);
    \draw [thick] (0,4) to [out=0,in=120] (3,2) to [out=-60,in=180] (4,0);
  \end{tikzpicture}
  \quad \sim \quad
  \begin{tikzpicture}[anchorbase,scale=.4]
    \draw [thick] (0,0) to [out=0,in=180] (1,5) to [out=0,in=90] (2.5,2) to [out=-90,in=90] (1.5,.5) to [out=-90,in=-90] (4,.5) to [out=90,in=180] (6,4);
    \draw [line width=3, white] (0,3) to [out=0,in=180] (2,2);
    \draw [line width=3, white] (0,1) to [out=0,in=180] (2,2);
    \draw [line width=3, white] (2,2) -- (6,2);
    \draw [thick] (0,3) to [out=0,in=180] (2,2);
    \draw [thick] (0,1) to [out=0,in=180] (2,2);
    \draw [double] (2,2) -- (6,2);
    \draw [line width=3, white] (0,5) to [out=0,in=180] (2,6) to [out=0,in=90] (3.5,2) to [out=-90,in=0] (2.5,.5) to [out=180,in=-90] (1.5,2) to [out=90,in=0] (1.25,2.8) to [out=180,in=90] (1.2,0) to [out=-90,in=180] (6,0);
    \draw [thick] (0,5) to [out=0,in=180] (2,6) to [out=0,in=90] (3.5,2) to [out=-90,in=0] (2.5,.5) to [out=180,in=-90] (1.5,2) to [out=90,in=0] (1.25,2.8) to [out=180,in=90] (1.2,0) to [out=-90,in=180] (6,0);
  \end{tikzpicture}
\]

In the $1$-labeled case, we choose to run checks in the following two cases (we indicate which generators we choose when going down to the fully deformed projection):
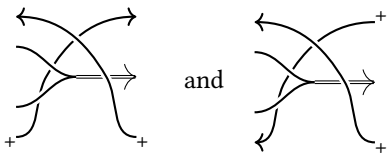
\[
  \begin{tikzpicture}[anchorbase,scale=.4]
    \draw [thick,->] (0,0) to [out=0,in=-120] (1,2) to [out=60,in=180] (4,4);
    \draw [line width=3, white] (0,3) to [out=0,in=180] (2,2);
    \draw [line width=3, white] (0,1) to [out=0,in=180] (2,2);
    \draw [line width=3, white] (2,2) -- (4,2);
    \draw [thick] (0,3) to [out=0,in=180] (2,2);
    \draw [thick] (0,1) to [out=0,in=180] (2,2);
    \draw [double,->] (2,2) -- (4,2);
    \draw [line width=3, white] (0,4) to [out=0,in=120] (3,2) to [out=-60,in=180] (4,0);
    \draw [thick,<-] (0,4) to [out=0,in=120] (3,2) to [out=-60,in=180] (4,0);
    \node at (-.2,-.2) {\small +};
    \node at (4.2,-.2) {\small +};
  \end{tikzpicture}
  \quad
  \text{and}
  \quad
    \begin{tikzpicture}[anchorbase,scale=.4]
    \draw [thick,<-] (0,0) to [out=0,in=-120] (1,2) to [out=60,in=180] (4,4);
    \draw [line width=3, white] (0,3) to [out=0,in=180] (2,2);
    \draw [line width=3, white] (0,1) to [out=0,in=180] (2,2);
    \draw [line width=3, white] (2,2) -- (4,2);
    \draw [thick] (0,3) to [out=0,in=180] (2,2);
    \draw [thick] (0,1) to [out=0,in=180] (2,2);
    \draw [double,->] (2,2) -- (4,2);
    \draw [line width=3, white] (0,4) to [out=0,in=120] (3,2) to [out=-60,in=180] (4,0);
    \draw [thick,<-] (0,4) to [out=0,in=120] (3,2) to [out=-60,in=180] (4,0);
    \node at (4.2,4.2) {\small +};
    \node at (4.2,-.2) {\small +};
  \end{tikzpicture}
\]

The checks in the case where both free strands are $2$-labeled are very easy. In the case where only one strand is $2$-labeled, the proof is very close to the ones for moves~\eqref{C3} and~\eqref{C4}. Indeed, they are related to those moves by crossing changes, and changing the over/under information in a crossing involving a $2$-labeled strand is easy to track through computations. One can thus go through the proofs for moves~\eqref{C3} and~\eqref{C4} and track changes explicitly.

\noindent {\bf R1 move passing through a trivalent vertex~\eqref{R1Triv}:} 

Here we will consider a reordered version of the move, by first letting the twist pass through the trivalent vertex, and then the curl. Both versions are equivalent thanks to the move~\eqref{TwTrivStrand} that will be proved independently just next.

We will pack a good deal of the technical computation into Lemma~\ref{lem:curlthrufork}. Before stating it, let us observe the following phenomenon, which easily follows from changing the curl into a twist and using $\rm{MM}_8$.

\begin{lemma} \label{lem:curloverstrand}
  In fully deformed homology, moving a $1$-labeled curl over a $1$-labeled strand induces:
\begin{equation}
  \begin{tikzpicture}[anchorbase]
    \node at (-4,0) {
      \begin{tikzpicture}[anchorbase,scale=.5]
        \draw [thick] (2,0) -- (0,2);
        \draw [thick, ->] (1.1,1.9) to [out=-135,in=-135] (1.5,1.5) -- (2,2);
        \draw [ultra thick, white] (0,0) -- (1,1) -- (1.5,1.5) to [out=45,in=-45] (1.1,1.9);
        \draw [thick] (0,0) -- (1,1) -- (1.5,1.5) to [out=45,in=45] (1.1,1.9);
        \node at (-.3,-.3) {\tiny $\pm$};
      \end{tikzpicture}
    };
    \node at (-2,0) {$\rightarrow$};
    \node (A) at (0,0) {
      \begin{tikzpicture}[anchorbase,scale=.5]
        \draw [thick] (2,0) -- (0,2);
        \draw [ultra thick, white] (0,0) -- (2,2);
        \draw [thick,->] (0,0) -- (2,2);
        \draw [thick] (1.25,1.75) circle (.2);
        \node at (-.3,-.3) {\tiny $\pm$};
      \end{tikzpicture}};
    \node at (-4,-4) {
      \begin{tikzpicture}[anchorbase,scale=.5]
        \draw [thick] (2,0) -- (0,2);
        \draw [ultra thick, white] (.1,.9) to [out=-135,in=-135] (.5,.5) -- (2,2);
        \draw [thick, ->] (.1,.9) to [out=-135,in=-135] (.5,.5) -- (2,2);
        \draw [ultra thick, white] (0,0) --  (.5,.5) to [out=45,in=-45] (.1,.9);
        \draw [thick] (0,0) -- (.5,.5) to [out=45,in=45] (.1,.9);
        \node at (-.3,-.3) {\tiny $\pm$};
      \end{tikzpicture}
    };
    \node at (-2,-4) {$\rightarrow$};
        \node (E) at (0,-4) {
      \begin{tikzpicture}[anchorbase,scale=.5]
        \draw [thick] (2,0) -- (0,2);
        \draw [ultra thick, white] (0,0) -- (2,2);
        \draw [thick,->] (0,0) -- (2,2);
        \draw [thick] (.25,.75) circle (.2);
        \node at (-.3,-.3) {\tiny $\pm$};
      \end{tikzpicture}
    };
    \draw [green,->] (A) edge node [right] {\tiny $\pm 2$cap/cup} (E);
  \end{tikzpicture}
\end{equation}
\end{lemma}

The next lemma aims at describing the effect of passing a curl through a trivalent vertex. The two sides of the move can be summarized as follows:
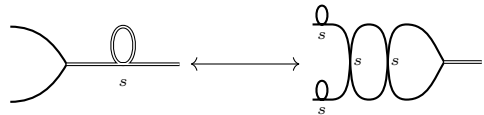
\[
\begin{tikzpicture}[anchorbase]
  \node at (0,0) (L) {
    \begin{tikzpicture}[scale=.5]
      \draw [thick] (0,1) to [out=0,in=120] (1.5,0);
      \draw [thick] (0,-1) to [out=0,in=-120] (1.5,0);
      \draw [double] (1.5,0) -- (3,0) to [out=0,in=0] (3,1) to [out=180,in=180] (3,0) -- (4.5,0);
      \node at (3,-.5) {\tiny $s$};
    \end{tikzpicture}
  };
  \node at (4,0) (R) {
    \begin{tikzpicture}[scale=.5]
      \draw [thick] (0,1)-- (.25,1) to [out=0,in=0] (.25,1.5) to [out=180,in=180] (.25,1) --  (.5,1)  to [out=0,in=180] (1.5,-1) to [out=0,in=180] (2.5,1) to [out=0,in=120] (3.5,0);
      \draw [thick] (0,-1)-- (.25,-1) to [out=0,in=0] (.25,-.5) to [out=180,in=180] (.25,-1) --  (.5,-1)  to [out=0,in=180] (1.5,1) to [out=0,in=180] (2.5,-1) to [out=0,in=-120] (3.5,0);
      \draw [double] (3.5,0) -- (4.5,0);
      \node at (.25,.7) {\tiny $s$};
      \node at (.25,-1.3) {\tiny $s$};
      \node at (1.2,0) {\tiny $s$};
      \node at (2.2,0) {\tiny $s$};
    \end{tikzpicture}
  };
  \draw [<->] (L) -- (R);
\end{tikzpicture}
\]

In the above picture we use the following conventions: $s=\pm$, and depending on the sign the crossings are resolved positively or negatively (with respect to any orientation, as the result in independent on the choice of orientation). One could also consider curls on the other side of the strands. Altogether, one gets eight versions of the move ($s=\pm$, curls above or below, and orientation running left (we take a variable $\eta=-1$) or right ($\eta=1$)). The two crossings are labeled $c_1$ and $c_2$, in the order prescribed by the orientation of the strands.

We denote by $F$ the following typical foam in fully deformed homology (to be adapted to each of the situations described above):

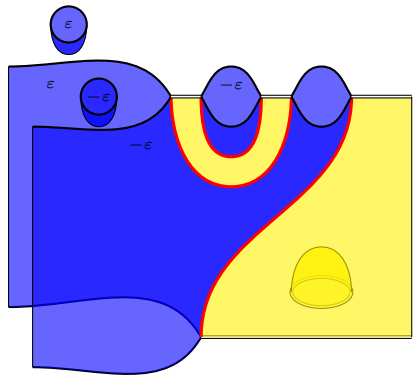
\[

\]

The mirror image of the same foam (used in the reverse direction of the move) will be denoted $\bar{F}$.

\begin{lemma} \label{lem:curlthrufork} With the above notations, the move
  induces the map: 
  \[
  2\eta\epsilon s a_3 F\otimes \cdot \wedge c_1\wedge c_2
  \]
  with inverse
  \[
  -2\eta\epsilon s a_3 \bar{F}\otimes \langle \cdot, c_1\wedge c_2\rangle
  \]
\end{lemma}

\begin{proof}
  We fully detail the following case.

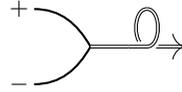
\[
    %
  
\]

The computation first goes through the following steps (below, crossings are numbered, corresponding to crossing variables $c_{i}$'s):
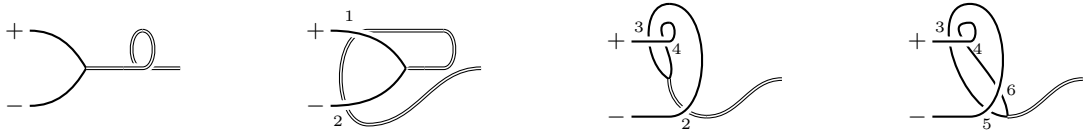
\[
  %
\]

Composing the corresponding elementary moves~\eqref{FS7} and~\eqref{FS1} (twice), one gets:
\[
  -a_{17} F \otimes c_3\wedge c_6
\]
where $F$ is the foam (we haven't colored all facets in order to ease the readability of the picture):
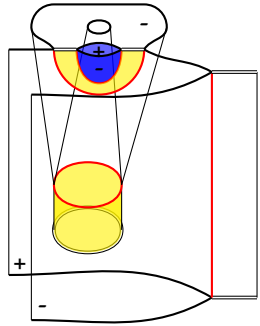
\[
  %
\]

Next, one slides a curl over a strand as follows:

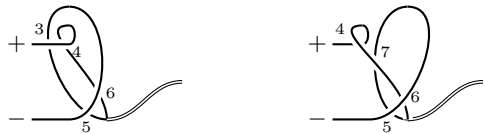
\[
  %
\]

This is controlled by Lemma~\ref{lem:curloverstrand}, and has the effect of caping off the inner circle, and cupping it back outside.

Then an $R_{III}$ move is performed.

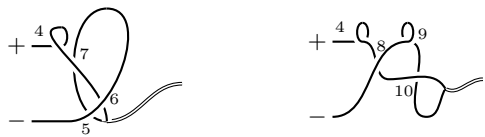
\[
  %
\]
We invoke~\eqref{R31}, and get at this point: $-2a_{17}F'\otimes c_8\wedge c_{10}$, where $F'$ is the result of post-composing $F$ with:

\[
  %
\]

One can apply relations~\eqref{deformedbb3} and \eqref{deformedsl2NH_enh} to reduce to: $2a_{17}F''\otimes c_8\wedge c_{10}$, with $F''$:

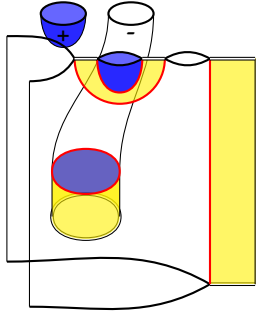
\[
  %
\]

Finally, we again slide a curl across a vertex, resulting in another cap/cup morphism. One uses relation~\eqref{deformedbb3} to find that the move induces $2a_{17}F'''\otimes c_8\wedge c_{10}$ as expected, with $F'''$:
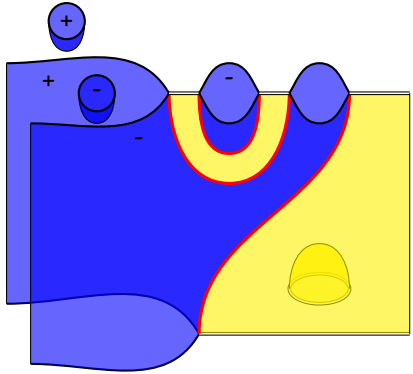
\[
%
\]

This matches the expected result, with: $\epsilon=+$ and $s=-$ and $\eta=+$.

This finishes the proof of Lemma~\ref{lem:curlthrufork}.
\end{proof}

We now return to the study of the move~\eqref{R1Triv}. We make choices for orientation and crossing signs, and decompose the complicated side into the following elementary moves:

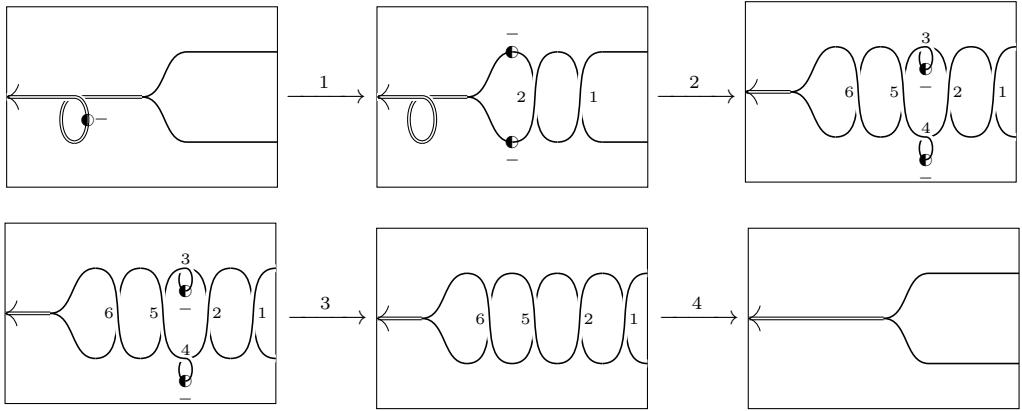
\[
    %
\]

Furthermore we study the effect in fully deformed homology when the upper strand from the trivalent vertex is labeled $+$.
\begin{itemize}
\item Step $1$ is governed by~\ref{TT3}, in the upward direction, giving $4a_3\zeta \rm{cup}^2\otimes c_1\wedge c_2$;
\item Step $2$ is governed by Lemma~\ref{lem:curlthrufork} with $\epsilon=-$, $s=+$, $\eta=+$, giving $-2a_3\otimes c_5\wedge c_6$ with foam given by the lemma;
\item Step $3$ consists of two moves~\ref{R11} performed together. The caps will close off cups from the previous step, each one labeled with a different idempotent, thus producing $\frac{1}{2}\frac{-1}{2}=\frac{-1}{4}$ thanks to Relation~\ref{deformedsphere};
\item Step $4$ consists of two moves~\ref{R22}, with effect $\langle \cdot,c_1\wedge c_6\wedge c_2 \wedge c_5\rangle$ on the exterior part, and the following foam:
  
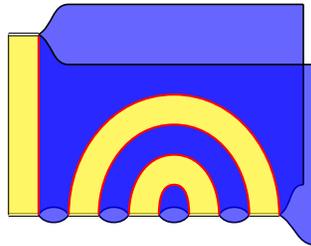
\[
  \begin{tikzpicture}[anchorbase,scale=.4]
    \draw [double] (0,0) -- (1,0);
    \draw [semithick] (1,0) to [out=60,in=120] (2,0);
    \draw [semithick] (1,0) to [out=-60,in=-120] (2,0);
    \draw [double] (2,0) -- (3,0);
    \draw [semithick] (3,0) to [out=60,in=120] (4,0);
    \draw [semithick] (3,0) to [out=-60,in=-120] (4,0);
    \draw [double] (4,0) -- (5,0);
    \draw [semithick] (5,0) to [out=60,in=120] (6,0);
    \draw [semithick] (5,0) to [out=-60,in=-120] (6,0);
    \draw [double] (6,0) -- (7,0);
    \draw [semithick] (7,0) to [out=60,in=120] (8,0);
    \draw [semithick] (7,0) to [out=-60,in=-120] (8,0);
    \draw [double] (8,0) -- (9,0);
    \draw [semithick] (9,0) to [out=30,in=180] (9.8,1);
    \draw [semithick] (9,0) to [out=-30,in=180] (10.2,-1);
    \fill [opacity=.6,yellow] (0,0) -- (0,6) -- (1,6) -- (1,0);
    \draw (0,0) -- (0,6);
    \fill [blue, opacity=.6] (9.8,1) -- (9.8,7) -- (2,7) to [out=180,in=30] (1,6) -- (1,0) to [out=60,in=120] (2,0) to [out=90,in=180] (5.5,4) to [out=0,in=90] (9,0) to [out=30,in=180] (9.8,1);
    \draw (9.8,1) -- (9.8,7);
    \fill [blue, opacity=.6] (10.2,-1) -- (10.2,5) -- (2,5) to [out=180,in=-30] (1,6) -- (1,0) to [out=-60,in=-120] (2,0) to [out=90,in=180] (5.5,4) to [out=0,in=90] (9,0) to [out=-30,in=180] (10.2,-1);
    \draw (9.8,1) -- (9.8,7);
    \draw (10.2,-1) -- (10.2,5);
    \fill [yellow, opacity=.6] (2,0) to [out=90,in=180] (5.5,4) to [out=0,in=90] (9,0) -- (8,0) to [out=90,in=0] (5.5,3) to [out=180,in=90] (3,0) -- (2,0);
    \fill [blue, opacity=.6] (4,0) to [out=120,in=60] (3,0) to [out=90,in=180] (5.5,3) to [out=0,in=90] (8,0) to [out=120,in=60] (7,0) to [out=90,in=0] (5.5,2) to [out=180,in=90] (4,0);
    \fill [blue, opacity=.6] (4,0) to [out=-120,in=-60] (3,0) to [out=90,in=180] (5.5,3) to [out=0,in=90] (8,0) to [out=-120,in=-60] (7,0) to [out=90,in=0] (5.5,2) to [out=180,in=90] (4,0);
    \fill [yellow, opacity=.6] (4,0) to [out=90,in=180] (5.5,2) to [out=0,in=90] (7,0) -- (6,0) to [out=90,in=0] (5.5,1) to [out=180,in=90] (5,0) -- (4,0);
    \fill [blue, opacity=.6] (5,0) to [out=60,in=120] (6,0) to [out=90,in=0] (5.5,1) to [out=180,in=90] (5,0);
    \fill [blue, opacity=.6] (5,0) to [out=-60,in=-120] (6,0) to [out=90,in=0] (5.5,1) to [out=180,in=90] (5,0);
    \draw [red, thick] (5,0) to [out=90,in=180] (5.5,1) to [out=0,in=90] (6,0);
    \draw [red, thick] (4,0) to [out=90,in=180] (5.5,2) to [out=0,in=90] (7,0);
    \draw [red, thick] (3,0) to [out=90,in=180] (5.5,3) to [out=0,in=90] (8,0);
    \draw [red, thick] (2,0) to [out=90,in=180] (5.5,4) to [out=0,in=90] (9,0);
    \draw [red, thick] (1,0) -- (1,6);
    \draw [double] (0,6) -- (1,6);
    \draw [semithick] (1,6) to [out=30,in=180] (2,7) -- (9.8,7);
    \draw [semithick] (1,6) to [out=-30,in=180] (2,5) -- (10.2,5);
  \end{tikzpicture}
\]
\end{itemize}

Assembling the pieces together, one gets $2\zeta$ times the following foam :
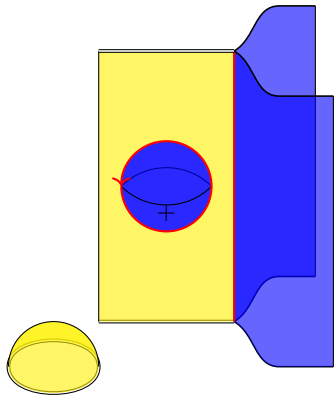
\[
  \begin{tikzpicture}[anchorbase,scale=.6]
    \draw [double] (0,0) -- (3,0);
    \draw [semithick] (3,0) to [out=30,in=180] (4,1) -- (4.8,1);
    \draw [semithick] (3,0) to [out=-30,in=180] (4,-1) -- (5.2,-1);
    \fill [opacity=.6,yellow] (0,0) -- (0,6) -- (3,6) -- (3,0) -- (0,0) -- (.5,3)  arc(-180:180:1);
    \draw (0,0) -- (0,6);
    \fill [opacity=.6,blue] (.5,3)  arc(-180:180:1);
    \draw (.5,3) to [out=45,in=135] (2.5,3);
    \fill [opacity=.6,blue] (.5,3)  arc(-180:180:1);
    \draw (.5,3) to [out=-45,in=-135] (2.5,3);
    \fill [blue, opacity=.6] (4.8,1) -- (4.8,7) -- (4,7) to [out=180,in=30] (3,6) -- (3,0) to [out=30,in=180] (4,1) -- (4.8,1);
    \draw (4.8,1) -- (4.8,7);
    \fill [blue, opacity=.6] (5.2,-1) -- (5.2,5) -- (4,5) to [out=180,in=-30] (3,6) -- (3,0) to [out=-30,in=180] (4,-1) -- (5.2,-1);
    \draw (5.2,-1) -- (5.2,5);
    \draw [red, thick] (3,0) -- (3,6);
    \draw [red,thick,->] (.5,3)  arc(-180:180:1);
    \draw [double] (0,6) -- (3,6);
    \draw [semithick] (3,6) to [out=30,in=180] (4,7) -- (4.8,7);
    \draw [semithick] (3,6) to [out=-30,in=180] (4,5) -- (5.2,5);
    \node at (1.5,2.4) {\small \bf $+$};
    \draw [double] (-2,-1) to [out=90,in=90] (0,-1) to [out=-90,in=-90] (-2,-1);
    \fill [yellow, opacity=.6] (-2,-1) to [out=90,in=90] (0,-1) to [out=90,in=0] (-1,0) to [out=180,in=90] (-2,-1);
    \fill [yellow, opacity=.6] (-2,-1) to [out=-90,in=-90] (0,-1) to [out=90,in=0] (-1,0) to [out=180,in=90] (-2,-1);
    \draw (-2,-1) to [out=90,in=180] (-1,0) to [out=0,in=90] (0,-1); 
  \end{tikzpicture}
\]

The blister has the $+$ decoration on the front side, and thanks to Relation~\eqref{deformedbb3} it produces a coefficient $\frac{1}{2}$.

Altogether, one thus gets $\zeta$ times the same foam as above with the blister removed, which is exactly what one reads from the other side of the move, that simply consists of the move~\eqref{R13}.

All other cases follow under a very similar process. Indeed, reversing the orientation will introduce a sign at step $2$ and change the sign from the blister, changing the side of the curl is of no effect, and the main effect in changing all crossings and twists signs consists in the fact that at step $1$ one looses the coefficient $4$, but it appears back again at step $3$ (see move~\eqref{R12}).

\noindent {\bf Trivalent twist and strand~\eqref{TwTrivStrand}:} 

There are a priori $16$ versions of this move to check: the orientation on the fork can run in one or the other direction, the strand can pass over or under the fork, with any orientation, and the twist can be positive or negative. We argue that it is enough to check only the case where the fork is a split, the orientation of the strand runs leftwards, and the sign of the twist is positive, leaving us with two cases (strand passing over or under).

Indeed, using the movie move~\eqref{StrandAroundVertex} one can revert the orientation of the strand. Using the (yet unproved but the proof will be independent) move~\eqref{TwMorseTriv}, one can revert the sign of the twist. Finally, the move~\eqref{eq:TwZip} (also yet to prove but again with independent proof) allows to pass from a split to a merge trivalent vertex.

Let us first focus on the following case (the other one follows the very same lines):

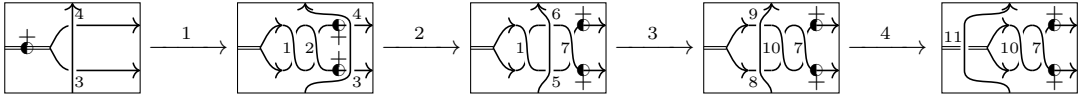
\[

    \]

    We go down to fully deformed homology, coloring the free strand with a $+$, and the split vertex with a $+$ for the bottom strand and a $-$ for the top one.

    \begin{itemize}
    \item Step $1$ is controlled by~\eqref{TT2}, giving $a_3\zeta\rm{cup}^2\otimes c_1\wedge c_2$;
    \item Step $2$ is controlled by~\eqref{R327} with $(+,-,+)$, inducing $-\langle \cdot,c_2\wedge c_4\rangle \wedge c_5\wedge c_7$ times the following foam:

  \[
  %
\]
\item Step $3$ uses again~\eqref{R327} and produces $-\id \otimes \langle\cdot,c_1\wedge c_5\rangle \wedge c_9\wedge c_{10}$;
\item Step $4$ relies on~\eqref{FS9} and produces $a_{11}\otimes \langle \cdot,c_9\rangle \wedge c_{11}$ times the  foam the left part of which is (elsewhere it is just an identity):
    \[
      %
\]
    \end{itemize}

    Altogether, one reads $a_3a_{11}\zeta\langle \cdot,c_4\rangle \wedge c_7\wedge c_{10}\wedge c_{11}$ times a foam that contains a blister on the left, evaluating to $\frac{-1}{2}$ using~\eqref{deformedbb3}, and otherwise is obtained by composition of the following steps (the gray ellipses in the last two pictures highlight the locus where the last zipped cup takes place):
    \[
    %
\]

One can then perform Relation~\eqref{deformedsl2NH_enh} just under the leftmost of the two digons at the top, at the cost of a $-2$ coefficient that will balance the coefficient that arose from the blister. As a result, both digons are created by cups. This will make the identification with the other side easier.

Indeed consider now:
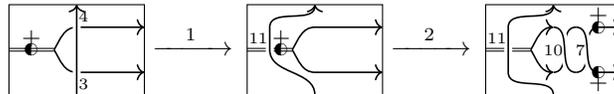
\[
  %
    \]

    The first step is again Move~\ref{FS9}, and the second one is~\eqref{TT2} again. They assemble to give the very same result as above.

    \noindent {\bf Cup and twist~\eqref{eq:TwCup}:} 
    There are four cases to check, which are very similar. We consider:
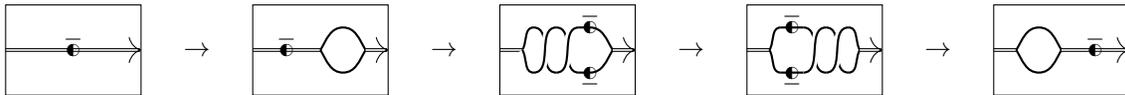
\[
  %
\]
we want to see that the map equals the one where the digon is opened right away to the left of the twist. The $\Hom$-space of degree zero between both sides of the above picture is $1$-dimensional Thus we will prove that this map induces an identity map on any of the two  non-zero surviving generators in the deformed space. Let us assume that we have a minus on the upper strand of the digon.

The second step is controlled by the generator~\eqref{TT4}, going upwards, the third one induces an identity, while the fourth one is~\eqref{TT3}. Altogether, one gets $-4$ times a cup together with two blisters that will give opposite signs through Equation~\eqref{deformedbb3}, thus producing simply a cup. In reversed direction (with caps), the same holds true for the same reasons. The case of a positive twist is similar.

\noindent {\bf Zip and twist~\eqref{eq:TwZip}:}

The $\Hom$ space between the bottom and the top of the move is one dimensional, so it suffices to check that the equality holds in deformed homology for any of the surviving generators. Here \eqref{TT1} and \eqref{TT2} appear on one or the other side, using both times the upward version (or the downward if reversing the time frame) (and if one reverses the twist sign then \eqref{TT3} and \eqref{TT4} are at play). One gets on the nose the same foam, and coefficients agree. 

\noindent {\bf Morse point on twist and trivalent vertex~\eqref{TwMorseTriv}:}

We start with the following version of the move:
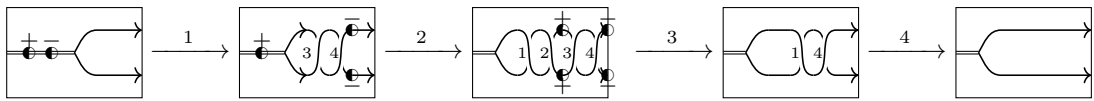
\[
  \begin{tikzpicture}[anchorbase,scale=.3]
    \draw [double] (0,2) -- (3,2);
    \draw [semithick,->] (3,2) to [out=60,in=180] (4,3) -- (6,3);
    \draw [semithick,->] (3,2) to [out=-60,in=180] (4,1) -- (6,1);
    \node at (1,2) {\tiny $\LEFTcircle$};
    \node at (1,2.5) {\small $+$};
    \node at (2,2) {\tiny $\LEFTcircle$};
    \node at (2,2.5) {\small $-$};
    \draw (0,0) rectangle (6,4);
  \end{tikzpicture}
 \xrightarrow{\quad 1\quad}
  \begin{tikzpicture}[anchorbase,scale=.3]
    \draw [double] (0,2) -- (2,2);
    \draw [semithick,->] (2,2) to [out=60,in=180] (3,3);
    \draw [semithick,->] (2,2) to [out=-60,in=180] (3,1);
    \node at (1,2) {\tiny $\LEFTcircle$};
    \node at (1,2.5) {\small $+$};
    \draw [white, line width=3] (3,3) to [out=0,in=180] (4,1);
    \draw [semithick] (3,3) to [out=0,in=180] (4,1);
    \draw [white, line width=3] (3,1) to [out=0,in=180] (4,3);
    \draw [semithick] (3,1) to [out=0,in=180] (4,3);
    \draw [white, line width=3] (4,3) to [out=0,in=180] (5,1);
    \draw [semithick] (4,3) to [out=0,in=180] (5,1);
    \draw [white, line width=3] (4,1) to [out=0,in=180] (5,3);
    \draw [semithick] (4,1) to [out=0,in=180] (5,3);
    \draw [semithick,->] (5,3) -- (6,3);
    \draw [semithick,->] (5,1) -- (6,1);
    \node at (5,3) {\tiny $\LEFTcircle$};
    \node at (5,3.5) {\small $-$};
    \node at (5,1) {\tiny $\LEFTcircle$};
    \node at (5,.5) {\small $-$};
    \draw (0,0) rectangle (6,4);
    \node at (3,2) {\tiny $3$};
    \node at (4.2,2) {\tiny $4$};
      \end{tikzpicture}
 \xrightarrow{\quad 2\quad}
  \begin{tikzpicture}[anchorbase,scale=.3]
    \draw [double] (0,2) -- (1,2);
    \draw [semithick] (1,2) to [out=60,in=180] (2,3);
    \draw [semithick] (1,2) to [out=-60,in=180] (2,1);
    \draw [white, line width=3] (2,1) to [out=0,in=180] (3,3);
    \draw [semithick] (2,1) to [out=0,in=180] (3,3);
    \draw [white, line width=3] (2,3) to [out=0,in=180] (3,1);
    \draw [semithick] (2,3) to [out=0,in=180] (3,1);
    \draw [white, line width=3] (3,1) to [out=0,in=180] (4,3);
    \draw [semithick] (3,1) to [out=0,in=180] (4,3);
    \draw [white, line width=3] (3,3) to [out=0,in=180] (4,1);
    \draw [semithick] (3,3) to [out=0,in=180] (4,1);
    \draw [white, line width=3] (4,3) to [out=0,in=180] (5,1);
    \draw [semithick] (4,3) to [out=0,in=180] (5,1);
    \draw [white, line width=3] (4,1) to [out=0,in=180] (5,3);
    \draw [semithick] (4,1) to [out=0,in=180] (5,3);
    \draw [white, line width=3] (5,3) to [out=0,in=180] (6,1);
    \draw [semithick,->] (5,3) to [out=0,in=180] (6,1);
    \draw [white, line width=3] (5,1) to [out=0,in=180] (6,3);
    \draw [semithick,->] (5,1) to [out=0,in=180] (6,3);
    \node at (4,3) {\tiny $\LEFTcircle$};
    \node at (4,3.5) {\small $+$};
    \node at (4,1) {\tiny $\LEFTcircle$};
    \node at (4,.5) {\small $+$};
    \node at (6,3) {\tiny $\LEFTcircle$};
    \node at (6,3.5) {\small $-$};
    \node at (6,1) {\tiny $\LEFTcircle$};
    \node at (6,.5) {\small $-$};
    \draw (0,0) rectangle (6,4);
    \node at (2.2,2) {\tiny $1$};
    \node at (3.2,2) {\tiny $2$};
    \node at (4.2,2) {\tiny $3$};
    \node at (5.2,2) {\tiny $4$};
      \end{tikzpicture}
 \xrightarrow{\quad 3\quad}
  \begin{tikzpicture}[anchorbase,scale=.3]
    \draw [double] (0,2) -- (1,2);
    \draw [semithick] (1,2) to [out=60,in=180] (2,3);
    \draw [semithick] (1,2) to [out=-60,in=180] (2,1);
    \draw [semithick] (2,3) -- (3,3);
    \draw [semithick] (2,1) -- (3,1);
    \draw [white, line width=3] (3,1) to [out=0,in=180] (4,3);
    \draw [semithick] (3,1) to [out=0,in=180] (4,3);
    \draw [white, line width=3] (3,3) to [out=0,in=180] (4,1);
    \draw [semithick] (3,3) to [out=0,in=180] (4,1);
    \draw [white, line width=3] (4,3) to [out=0,in=180] (5,1);
    \draw [semithick,->] (4,3) to [out=0,in=180] (5,1) -- (6,1);
    \draw [white, line width=3] (4,1) to [out=0,in=180] (5,3);
    \draw [semithick,->] (4,1) to [out=0,in=180] (5,3) -- (6,3);
    \draw (0,0) rectangle (6,4);
    \node at (3.2,2) {\tiny $1$};
    \node at (4.2,2) {\tiny $4$};
      \end{tikzpicture}
 \xrightarrow{\quad 4\quad}
  \begin{tikzpicture}[anchorbase,scale=.3]
    \draw [double] (0,2) -- (1,2);
    \draw [semithick,->] (1,2) to [out=60,in=180] (2,3) -- (6,3);
    \draw [semithick,->] (1,2) to [out=-60,in=180] (2,1) -- (6,1);
    \draw (0,0) rectangle (6,4);
  \end{tikzpicture}
\]

\begin{itemize}
\item Step $1$ is~\eqref{TT4} in the upward direction, producing $4a_3\zeta\rm{cup}^2\otimes c_3\wedge c_4$;
\item Step $2$ is~\eqref{TT2} in the upward direction, producing $-a_3\zeta\rm{cup}^2\otimes c_1\wedge c_2$;
\item Step $3$ is~\eqref{R21}, producing $\langle \cdot, c_2\wedge c_3\rangle$ with the following foam placed over the middle digons:
  
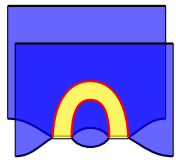
\[
  \begin{tikzpicture}[anchorbase,scale=.5]
    \draw [thick] (0,0) to [out=0,in=-150] (1,.5);
    \draw [thick] (-.2,1) to [out=0,in=150] (1,.5);
    \draw [double] (1,.5) -- (1.5,.5);
    \draw [thick] (1.5,.5) to [out=60,in=120] (2.5,.5);
    \draw [thick] (1.5,.5) to [out=-60,in=-120] (2.5,.5);
    \draw [double] (2.5,.5) -- (3,.5);
    \draw [thick] (3,.5) to [out=30,in=180] (4,1);
    \draw [thick] (3,.5) to [out=-30,in=180] (4.2,0);
    \fill [blue,opacity=.6] (1.5,.5) to [out=60,in=120] (2.5,.5) to [out=90,in=0] (2,1.5) to [out=180,in=90] (1.5,.5);
    \fill [blue, opacity=.6] (-.2,1) to [out=0,in=150] (1,.5) to [out=90,in=180] (2,2) to [out=0,in=90] (3,.5) to [out=30,in=180] (4,1) -- (4,4) -- (-.2,4) -- (-.2,1);
    \draw (-.2,1) -- (-.2,4);
    \draw (4,1) -- (4,4);
    \fill [blue,opacity=.6] (1.5,.5) to [out=-60,in=-120] (2.5,.5) to [out=90,in=0] (2,1.5) to [out=180,in=90] (1.5,.5);
    \fill [yellow,opacity=.6]   (1,.5) -- (1.5,.5) to [out=90,in=180] (2,1.5) to [out=0,in=90] (2.5,.5) -- (3,.5) to [out=90,in=0] (2,2) to [out=180,in=90] (1,.5);
    \fill [blue, opacity=.6] (0,0) to [out=0,in=-150] (1,.5) to [out=90,in=180] (2,2) to [out=0,in=90] (3,.5) to [out=-30,in=180] (4.2,0) -- (4.2,3) -- (0,3) -- (0,0);
    \draw (0,0) -- (0,3);
    \draw (4.2,0) -- (4.2,3);
    \draw [red,thick] (2.5,.5) to [out=90,in=0] (2,1.5) to [out=180,in=90] (1.5,.5);
    \draw [red, thick] (3,.5) to [out=90,in=0] (2,2) to [out=180,in=90] (1,.5);
    \draw [thick] (0,3) -- (4.2,3);
    \draw [thick] (-.2,4) -- (4,4);
  \end{tikzpicture}
\]
\item Step $4$ is again~\eqref{R21} with $\langle \cdot,c_1\wedge c_4\rangle$ and same local foam as above.
\end{itemize}

Altogether, one gets a coefficient $-4\langle c_3\wedge c_4\wedge c_1\wedge c_2,c_2\wedge c_3\wedge c_1\wedge c_4\rangle=-4$ and a foam the sketch of which is as follows:
\[
  \begin{tikzpicture}[anchorbase,scale=.6]
    \draw [double] (0,-2) -- (9,-2);
    \draw [semithick] (9,-2) to [out=60,in=180] (9.8,-1.5);
    \draw [semithick] (9,-2) to [out=-60,in=180] (10.2,-2.5);
    \draw [thick, red] (1,5) -- (1,0) to [out=-90,in=180] (1.5,-1) to [out=0,in=-90] (2,0) to [out=90,in=180] (5.5,4) to [out=0,in=90] (9,0) -- (9,-2);
    \draw [thick, red] (3,0) to [out=90,in=180] (5.5,3) to [out=0,in=90] (8,0) to [out=-90,in=0] (7.5,-1) to [out=180,in=-90] (7,0) to [out=90,in=0] (5.5,2) to [out=180,in=90] (4,0) to [out=-90,in=0] (3.5,-1) to [out=180,in=-90] (3,0);
    \draw [thick,red] (5,0) to [out=90,in=180] (5.5,1) to [out=0,in=90] (6,0) to [out=-90,in=0] (5.5,-1) to [out=180,in=-90] (5,0);
    \draw [double] (0,0) -- (1,0);
    \draw [semithick] (1,0) to [out=60,in=120] (2,0);
    \draw [semithick] (1,0) to [out=-60,in=-120] (2,0);
    \draw [double] (2,0) -- (3,0);
    \draw [semithick] (3,0) to [out=60,in=120] (4,0);
    \draw [semithick] (3,0) to [out=-60,in=-120] (4,0);
    \draw [double] (4,0) -- (5,0);
    \draw [semithick] (5,0) to [out=60,in=120] (6,0);
    \draw [semithick] (5,0) to [out=-60,in=-120] (6,0);
    \draw [double] (6,0) -- (7,0);
    \draw [semithick] (7,0) to [out=60,in=120] (8,0);
    \draw [semithick] (7,0) to [out=-60,in=-120] (8,0);
    \draw [double] (8,0) -- (9,0);
    \draw [semithick] (9,0) to [out=60,in=180] (9.8,.5);
    \draw [semithick] (9,0) to [out=-60,in=180] (10.2,-.5);
    \draw (0,-2) -- (0,5);
    \draw (9.8,-1.5) -- (9.8,5.5);
    \draw (10.2,-2.5) -- (10.2,4.5);
    \draw [double] (0,5) -- (1,5);
    \draw [semithick] (1,5) to [out=60,in=180] (2,5.5) -- (9.8,5.5);
    \draw [semithick] (1,5) to [out=-60,in=180] (2,4.5) -- (10.2,4.5);
  \end{tikzpicture}
\]

One sees two blisters, with same seam orientation but opposite idempotent labeling, thus producing a global coefficient $\frac{-1}{4}$ when removed, thanks to Equation~\eqref{deformedbb3}. This coefficient balances the $-4$ coefficient that we computed before, and altogether one simply gets an identity foam, as expected.

Looking at other configurations, one can observe that exchanging the order of the twist yields the same global computation. Changing the orientation of the trivalent vertex involves Moves~\eqref{TT1} and~\eqref{TT3} in place of~\eqref{TT2} and~\eqref{TT4} and yields the same result.

This concludes the proof, as all moves have now been checked.

\end{proof}

\section{Positivity of skein modules}

\subsection{Linear complexes}

Much of the material we will use throughout this section has been studied in
detail in joint work with P. Wedrich~\cite{QW_SkeinCat}, especially in Section 3, to which we refer for
details. We will only briefly recall the main facts, and build upon them.

Recall that the category of foams over the surface $\Su$, denoted $\Sfoam$, is
equivalent to a reduced category $\Sfoamred$, obtained by simplifying webs until
they contain no inessential 1-labeled components. A crucial result is
\cite[Corollary 3.18]{QW_SkeinCat}, stating that, provided the base surface $\Su$
is not $S^2$, the foam categories are positively graded. This
surface-specific grading is just the quantum grading on the reduced category.

This allows for projection functors $\Sfoam\rightarrow \Sfoam_0$ to the degree-zero subcategory, that are especially useful when one goes to the homotopy category.

The next question is to describe simple objects, that is those having a $1$-dimensional ring of endomorphisms of degree zero. To do so, recall the Jones-Wenzl projectors for $\glnn{2}$, akin to the ones for $\slnn{2}$~\cite{Wen}:
\begin{definition}
  The $\glnn{2}$ Jones--Wenzl projectors $P_m\in \Webq\otimes \Q(q)$ are defined by $P_1=\id_1$ and then:
 \[\begin{tikzpicture}[anchorbase, scale=.3]
\fill[black,opacity=.2] (0,1) rectangle (3,3);
\draw[thick] (0,1) rectangle (3,3);
\draw [very thick] (.5,0) to (.5,1);
\draw [thick, dotted] (.7,0.5) to (1.3,.5);
\draw [very thick] (1.5,0) to (1.5,1);
\draw [very thick] (2.5,0) to (2.5,1);
\draw [very thick,->] (.5,3) to (.5,4);
\draw [thick, dotted] (.7,3.25) to (1.3,3.25);
\draw [very thick,->] (1.5,3) to (1.5,4);
\draw [very thick,->] (2.5,3) to (2.5,4);
\node at (1.5,1.9) {$P_{m+1}$};
\end{tikzpicture}
   \;\; := \;\;
   \begin{tikzpicture}[anchorbase, scale=.3]
\fill[black,opacity=.2] (0,1) rectangle (2,3);
\draw[thick] (0,1) rectangle (2,3);
\draw [very thick] (.5,0) to (.5,1);
\draw [thick, dotted] (.7,0.5) to (1.3,.5);
\draw [very thick] (1.5,0) to (1.5,1);
\draw [very thick,->] (2.5,0) to (2.5,4);
\draw [very thick,->] (.5,3) to (.5,4);
\draw [thick, dotted] (.7,3.25) to (1.3,3.25);
\draw [very thick,->] (1.5,3) to (1.5,4);
\node at (1,1.9) {$P_{m}$};
\end{tikzpicture}
    \;-\;\frac{m}{m+1}\;
       \begin{tikzpicture}[anchorbase, scale=.3]
\fill[black,opacity=.2] (0,.5) rectangle (2,1.5);
\draw[thick] (0,.5) rectangle (2,1.5);
\fill[black,opacity=.2] (0,2.5) rectangle (2,3.5);
\draw[thick] (0,2.5) rectangle (2,3.5);
\draw [very thick] (.5,0) to (.5,.5);
\draw [very thick] (1.5,0) to (1.5,.5);
\draw [double] (2,1.75) to (2,2.25);
\draw [very thick] (2.5,0) to (2.5,1.5)to [out=90,in=90] (1.5,1.5);
\draw [very thick,->] (1.5,2.5) to [out=270,in=270] (2.5,2.5) to (2.5,4); 
\draw [very thick] (.5,1.5) to (.5,2.5);
\draw [thick, dotted] (.7,2) to (1.3,2);
\draw [very thick,->] (.5,3.5) to (.5,4);
\draw [very thick,->] (1.5,3.5) to (1.5,4);
\node at (1,.95) {\tiny$P_{m}$};
\node at (1,2.95) {\tiny$P_{m}$};
\end{tikzpicture}   
\]
\end{definition}

By sending these projectors $P_m$ to rotation foams $P_m\times \Ss^1$, one
obtains a family of idempotent foams. We use these foams to define Jones-Wenzl
foams $L^F_S$ attached to an integer lamination of the surface in
\cite[Definition 3.28]{QW_SkeinCat}. Those are categorical analogs of the basis
elements $L_S\in \Su\basisJW$. To get the full basis, one has to further
superpose those foams with identities on 2-labeled curves, see \cite[Definition
3.30]{QW_SkeinCat}. Lemma 3.34 then identifies this basis with the decategorified one,
and Corollary 3.35 shows that identity foams over webs decompose as sums of such
Jones-Wenzl idempotent foams. Finally, Corollary 3.37 states that these
projectors describe all simple objects of $\Kar(\Sfoam_0)$, that form the set
$\Su\basisJW^F$.

One could recast the non-negative grading construction by introducing a t-structure on the homotopy category $\HC(\Sfoamred)$. Since we do not explore the outcomes of this $t$-structure in the present paper, we do not fully detail it. All we will care about is its heart, that we describe in terms of \emph{linear complexes} as follows.

\begin{definition}
 In $\HC(\Sfoamred)$, we call a complex \emph{linear} if it is homotopy equivalent to one whose terms have balanced homological and quantum shifts. In other words and following the notations from the previous section, their shifts write $t^lq^{-l}$.
\end{definition}

More generally, $\HC(\Sfoamred)^{\geq 0}$ has objects, complexes with shifts of
the kind $t^{l}q^{-l-r}$ with $r\geq 0$. Similarly $\HC(\Sfoamred)^{\leq 0}$ has
shifts $t^lq^{-l+r}$, $r\geq 0$. Because of the non-negative grading on
morphisms, one can check that there are no morphisms between objects in the positive part and objects in the non-positive part of the t-structure\footnote{This is not the usual definition: one usually requires that there are no positive degree morphisms... We stick to the usual Khovanov convention that dots are of positive degree, and thus we twist the definition of a t-structure.}.

The following lemma is an easy but crucial observation regarding positivity statements.

\begin{lemma}
  Linear complexes (in other words, the heart of the t-structure) decategorify to skein elements that are positive in $(-q)$ in the $JW$ basis.
\end{lemma}

\begin{proof}
Any complex can be made homotopy equivalent to a reduced one with objects given by simple objects. Simples in $\HC(\Sfoamred)$ are described by $\Su\basisJW^F$. The objects thus decategorify to elements of the Jones-Wenzl basis. The positivity property in $-q$ follows since homological and grading shifts balance for linear complexes.
\end{proof}
This is exactly what is needed for the Jones-Wenzl version of the Fock-Goncharov-Thurston positivity conjecture, as explained in \cite[Section 2.1]{QW_SkeinCat}.

Recall now from \cite[Proposition 5.3]{QW_SkeinCat} (which is now fully proven
thanks to Theorem \ref{thm:functoriality}) that we have a bifunctor:
\[
\ast \;\colon\; \Sfoam\times \Sfoam \rightarrow \HC(\Sfoam)
\]

What we would like to do is consider linear complexes $A$ and $B$ in $\HC(\Sfoamred)$, define $A\ast B$ and prove that it remains linear. It is enough to prove it for products of elements of $\Su\basisJW^F$ (in which case we will say that $\Su\basisJW^F$ is linear). This would imply the positivity statement for the Jones-Wenzl basis. Indeed, the skein relations at the categorical level that are used to define the $\ast$ operation decategorify on the nose to the skein product, so if $A\ast B$ is linear, then $[A]\ast [B]=[A\ast B]$ is positive. 

Here is a difficulty in this approach: in $A$ and $B$, the differentials are defined up to homotopy. One can form $A\ast B$ by superposing objects and differentials maps, using the bifunctor previously defined. However, the differentials in the resulting complex will only be well-defined \emph{up to homotopy}. This might indicate that the right setup would be some kind of $A_{\infty}$-structure. As we do not need a full understanding of the theoretical properties of the functor, we will run a proof where all we care about is some choice for a realization of $A\ast B$ (specifically, one starts with a particular choice of realization for $A$ and $B$ and produces a complex $A\ast B$, even if this is not canonical). In order to highlight the difference, we change notations as follows.

\begin{definition}
  Given $A$ and $B$ representatives of complexes in $\HC(\Kar(\Sfoam))$ (with chosen isotopy representatives for curves, idempotents and differentials), we define $A\# B$ as the complex obtained by superposition of objects and morphisms and application of Khovanov's functor.
\end{definition}

Functoriality of Khovanov homology ensures that given other representatives $A'$ and $B'$, the resulting complexes $A\# B$ and $A'\# B'$ are homotopy equivalent.

Before going on with the positivity proof, let us state the following technical result.

\begin{lemma} \label{lem:dirSum}
  The bifunctor $\ast \;\colon\; \Sfoam\times \Sfoam \rightarrow \HC(\Sfoam)$ preserves direct sum decomposition.
\end{lemma}
\begin{proof}
  Consider a decomposition of the identity on some web W, $\id_W=p_1+p_2$ with
  $p_i^2=p_i$ and $p_1p_2=p_2p_1=0$. It is enough to prove the result for
  superposition with the identity over some web $W'$. Then, $id_W\ast
  id_{W'}=id_{W\ast W'}$ by definition of Khovanov's functor, and $\id_{W\ast
  W'}=p_1\ast \id_W' + p_2\ast \id_W'$. The result follows if the idempotent
  relations on $p_1$ and $p_2$ still hold even in presence of $\id_{W'}$. This
  is true, as they are obtained by successions of foam relations that happen in
  balls and these balls can be moved away from $\id_{W'}$ using the
  functoriality result. 
Let us for example focus on the equality $p_1^2=p_1$. Superposing $\id_{W'}$, the functoriality result tells us that $(p_1\ast \id_{W'})^2=p_1\ast \id_{W'}+h$ with $h$ an homotopy of $W\ast W'$. Since in the homotopy category homotopies are sent to zero, the idempotent relation does hold. We proceed similarly with the other relations to prove that $(p_1\ast \id_{W'})+(p_2\ast \id_{W'})$ is an idempotent decomposition. 

\end{proof}

\subsection{Positivity}

\begin{lemma}\label{lem:reduc} If for all pairs $(\gamma_1,\gamma_2)$ of
$1$-labeled multicurves, $\gamma_1\ast \gamma_2$ is linear, then $\Su\basisJW^F$
is linear.
\end{lemma}
\begin{proof}
  Assume that $L_S^{F}\ast \bigwedge^{(x)}$ and ${L'}_S^{F}\ast \bigwedge^{(x')}$ are two basis elements in $\Su\basisJW^F$. Consider the ($1$-labeled) multicurves $\gamma_1$ associated to the integer lamination $L$ and $\gamma_2$ associated to the integer lamination $L'$.

  First notice that $(L_S^{F}\ast \bigwedge^{x})\# ({L'}_S^{F}\ast \bigwedge^{(x')})$ is linear if and only if $(L_S^{F})\# ({L'}_S^{F})$ is linear. Indeed, one can pass between the two versions via inverting $(1,2)$ or $(2,2)$ labeled crossings and superposing $2$-labeled curves (what we call a $\ast \bigwedge$ operation in~\cite{QW_SkeinCat}), both of which are invertible processes and preserve the property of being linear (see~\eqref{eq:thickcrossing-2}).

  Now, in the Karoubi envelope, $\gamma_1=L_S^{F}\oplus \cdots$ and $\gamma_2={L'}_S^{F}$. Thus when taking the product and thanks to Lemma~\ref{lem:dirSum}, $L_S^{F}\# {L'}_S^{F}$ is a summand. Since summands of a linear complex are linear, the result follows.
\end{proof}

The previous lemma allows to drastically simplify the proof, using a phenomenon
that is inaccessible at the decategorified level. Indeed, we are using here that a summand of a linear object is still linear. The corresponding statement at the decategorified level would state that if $A+B$ is positive, then $A$ and $B$ are positive, which is obviously wrong.

We now come to the main positivity theorem.

\begin{theorem}\label{thm:positivity}
  If $\Su$ is not a torus, then $\Su\basisJW$ is linear. At the decategorified level, this implies that $\Su\basisJW$ is positive.
\end{theorem}

\begin{proof}
  It suffices to prove the categorical statement. Using Lemma~\ref{lem:reduc}, one can focus on multicurves only. Consider two multicurves $\gamma_1$ and $\gamma_2$.
 
Notice that for $3$ multicurves $\gamma$, $\gamma'$ and $\gamma''$, we have an homotopy equivalence: $\left(\gamma\ast \gamma'\right) \# \gamma''\simeq \gamma \ast \gamma'\ast \gamma''\simeq \gamma\# (\gamma'\ast \gamma'')$ (the middle term takes an extension of $\ast$ intro a tri-functor, just by taking a 3-stage stacking process). It is enough to prove that any of the three complexes is linear to know the result for all of them. Thus by writing $\gamma_1$ as a product of simple closed curves, it suffices by induction to prove positivity for $\gamma_1$ a simple closed curve.

  Assume that $\gamma_1\ast \gamma_2$ is not linear, and that $\gamma_2$ is minimal for this property with respect to the absolute intersection number $|\gamma_1\cdot \gamma_2|=c$ (in particular the two curves are in minimal position). The case where $c=0$ cannot happen, as then the product creates no crossing, and $c=1$ is also excluded as no trivial circles are created with only one crossing (only one connected component of $\gamma_2$ is involved, and both resolutions are a non-trivial, connected curve). Indeed, notice that only the creation of inessential circles can make the product of two curves non-linear, since the smoothing process from Equation~\eqref{eq:thickcrossing-2} locally preserves linearity.

  Thus $c\geq 2$.

  The picture below illustrates several steps of the following proof.

  
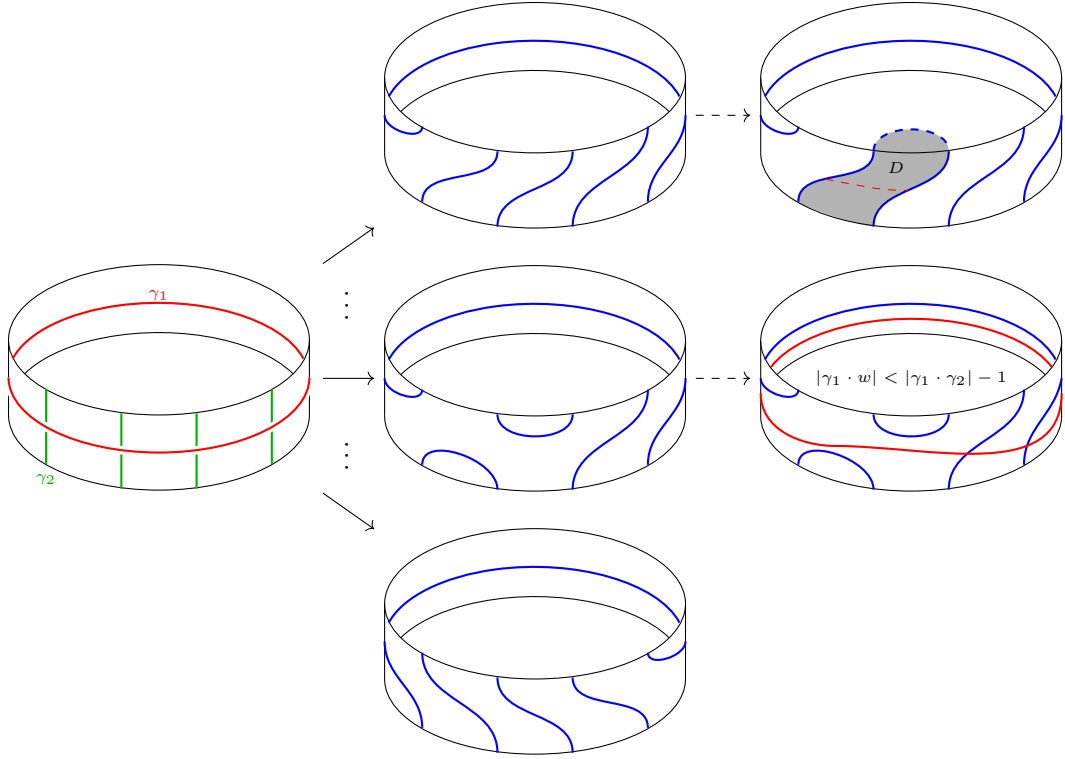
\[
\begin{tikzpicture}[anchorbase]
  \node (base) at (0,0) {
    \begin{tikzpicture}[anchorbase]
      \draw (1.78,-.45) arc (27:153:2 and 1);
      \draw (2,-1) arc (0:-180: 2 and 1);
      \draw (0,0) ellipse (2 and 1);
      \draw (2,0) -- (2,-1);
      \draw (-2,0) -- (-2,-1);
      \draw [green, thick] (-1.5, -1.65) -- (-1.5,-.65);
      \draw [green, thick] (-.5, -1.98) -- (-.5,-.98);
      \draw [green, thick] (.5, -1.98) -- (.5,-.98);
      \draw [green, thick] (1.5, -1.65) -- (1.5,-.65);
      \node [green] at (-1.5,-1.85) {\tiny $\gamma_2$};
      \draw [white, double=red, very thick, double distance=.8pt] (2,-.5) arc (0:-180:2 and 1);
      \draw [red, thick] (1.92,-.25) arc (15:165:2 and 1);
      \node [red] at (0,.6) {\tiny $\gamma_1$};
    \end{tikzpicture}
  };
  \node (extreme1) at (5,3.5) {
    \begin{tikzpicture}[anchorbase]
      \draw (1.78,-.45) arc (27:153:2 and 1);
      \draw (2,-1) arc (0:-180: 2 and 1);
      \draw (0,0) ellipse (2 and 1);
      \draw (2,0) -- (2,-1);
      \draw (-2,0) -- (-2,-1);
      \draw [blue, thick] (2,-.5) to [out=-90,in=90] (1.5,-1.65);
      \draw [blue, thick] (1.5,-.65) to [out=-90,in=90] (.5,-1.98);
      \draw [blue, thick] (.5,-.98) to [out=-90,in=90] (-.5,-1.98);
      \draw [blue, thick] (-.5,-.98) to [out=-90,in=90] (-1.5,-1.65);
      \draw [blue, thick] (-1.5,-.65) to [out=-90,in=-90] (-2,-.5);
      \draw [blue, thick] (1.92,-.25) arc (15:165:2 and 1);
        \end{tikzpicture}
  };  
  \node (extreme1disk) at (10,3.5) {
    \begin{tikzpicture}[anchorbase]
      \draw (1.78,-.45) arc (27:153:2 and 1);
      \draw (2,-1) arc (0:-180: 2 and 1);
      \draw (0,0) ellipse (2 and 1);
      \draw (2,0) -- (2,-1);
      \draw (-2,0) -- (-2,-1);
      \draw [blue, thick] (2,-.5) to [out=-90,in=90] (1.5,-1.65);
      \draw [blue, thick] (1.5,-.65) to [out=-90,in=90] (.5,-1.98);
      \draw [blue, thick] (.5,-.98) to [out=-90,in=90] (-.5,-1.98);
      \draw [blue, thick] (-.5,-.98) to [out=-90,in=90] (-1.5,-1.65);
      \draw [blue, thick] (-1.5,-.65) to [out=-90,in=-90] (-2,-.5);
      \draw [blue, thick] (1.92,-.25) arc (15:165:2 and 1);
      \fill [black,opacity=.3] (-1.5,-1.65) to [out=90,in=-90] (-.5,-.98) to [out=90,in=90] (.5,-.98) to [out=-90,in=90] (-.5,-1.98) to [out=175,in=160] (-1.5,-1.65);
      \draw [red,dashed] (-1.14,-1.35) to [out=-10,in=180] (-.05,-1.5);
      \node at (-.2,-1.2) {\tiny $D$};
      \draw [thick, blue, dashed] (-.5,-.98) to [out=90,in=90] (.5,-.98);
        \end{tikzpicture}
  };  
  \node (mixed) at (5,0) {
    \begin{tikzpicture}[anchorbase]
      \draw (1.78,-.45) arc (27:153:2 and 1);
      \draw (2,-1) arc (0:-180: 2 and 1);
      \draw (0,0) ellipse (2 and 1);
      \draw (2,0) -- (2,-1);
      \draw (-2,0) -- (-2,-1);
      \draw [blue, thick] (2,-.5) to [out=-90,in=90] (1.5,-1.65);
      \draw [blue, thick] (1.5,-.65) to [out=-90,in=90] (.5,-1.98);
      \draw [blue, thick] (.5,-.98) to [out=-90,in=-90] (-.5,-.98);
      \draw [blue, thick] (-.5,-1.98) to [out=90,in=90] (-1.5,-1.65);
      \draw [blue, thick] (-1.5,-.65) to [out=-90,in=-90] (-2,-.5);
      \draw [blue, thick] (1.92,-.25) arc (15:165:2 and 1);
        \end{tikzpicture}
  };  
  \node (extreme2) at (5,-3.5) {
    \begin{tikzpicture}[anchorbase]
      \draw (1.78,-.45) arc (27:153:2 and 1);
      \draw (2,-1) arc (0:-180: 2 and 1);
      \draw (0,0) ellipse (2 and 1);
      \draw (2,0) -- (2,-1);
      \draw (-2,0) -- (-2,-1);
      \draw [blue, thick] (2,-.5) to [out=-90,in=-90] (1.5,-.65);
      \draw [blue, thick] (1.5,-1.65) to [out=90,in=-90] (.5,-.98);
      \draw [blue, thick] (.5,-1.98) to [out=90,in=-90] (-.5,-.98);
      \draw [blue, thick] (-.5,-1.98) to [out=90,in=-90] (-1.5,-.65);
      \draw [blue, thick] (-1.5,-1.65) to [out=90,in=-90] (-2,-.5);
      \draw [blue, thick] (1.92,-.25) arc (15:165:2 and 1);
        \end{tikzpicture}
  };
  \node (mixedinter) at (10,0) {
    \begin{tikzpicture}[anchorbase]
      \draw (1.78,-.45) arc (27:153:2 and 1);
      \draw (2,-1) arc (0:-180: 2 and 1);
      \draw (0,0) ellipse (2 and 1);
      \draw (2,0) -- (2,-1);
      \draw (-2,0) -- (-2,-1);
      \draw [blue, thick] (2,-.5) to [out=-100,in=90] (1.5,-1.65);
      \draw [blue, thick] (1.5,-.65) to [out=-90,in=90] (.5,-1.98);
      \draw [blue, thick] (.5,-.98) to [out=-90,in=-90] (-.5,-.98);
      \draw [blue, thick] (-.5,-1.98) to [out=90,in=90] (-1.5,-1.65);
      \draw [blue, thick] (-1.5,-.65) to [out=-90,in=-90] (-2,-.5);
      \draw [blue, thick] (1.92,-.25) arc (15:165:2 and 1);
      \draw [red, thick] (2,-.7) to [out=-90,in=0] (-1,-1.4) to [out=180,in=-90] (-2,-.7);
      \draw [red, thick] (1.87,-.35) arc (21:159:2 and 1);
      \node at (0,-.5) {\tiny $|\gamma_1\cdot w|< |\gamma_1\cdot \gamma_2|-1$};
    \end{tikzpicture}
  };
  \draw [->] (base) -- (extreme1);
  \draw [->] (base) -- (mixed);
  \draw [->] (base) -- (extreme2);
  \draw [dashed,->] (extreme1) -- (extreme1disk);
  \draw [dashed,->] (mixed) -- (mixedinter);
  \node [rotate=90] at (2.5,1) {$\dots$};
  \node [rotate=90] at (2.5,-1) {$\dots$};
\end{tikzpicture}
\]

  Looking at a tubular neighborhood of $\gamma_1$, one sees a succession of crossings, as illustrated above on the left.

  We claim that the two extreme smoothings (illustrated at top and bottom of the
  previous picture), where all crossings are resolved in a parallel way, are
  linear. Indeed, consider an all-parallel resolution (the same argument covers
  both resolutions). We will prove that this resolution cannot contain a circle
  bounding a disk. Assume such a disk $D$ exists. Then $\gamma_1\cap
  \mathring{D}$ draws arcs on D. Consider an outermost arc. (By outermost, we
  mean that this arc cuts $D$ into two connected components, one of which
  contains no other arc.) Then this disk creates an isotopy of $\gamma_2$ that
  reduces the intersection number between $\gamma_1$ and $\gamma_2$ (see the top
  right part of the picture), contradicting the fact that the curves are in
  minimal position.

  Consider a multicurve $w$ (possibly superposed with a 2-labeled part) that
  comes from a not-all-parallel smoothing and that breaks linearity:
  $t^lq^{-l-r}w_{JW}$ appears as an object in a reduced complex for
  $\gamma_1\ast \gamma_2$ with $r\neq 0$. We claim that $|\gamma_1\cdot w|<c-1$:
  indeed, mixed smoothings contain turnbacks, as illustrated on the right part
  of the picture.

  If $r>0$, assume that $w$ maximizes $r$ and is rightmost (in homological grading) amongst terms that share the same value of $r$. Then:
  \[
  \Hom(t^lq^{-l-r}w,\gamma_1\ast\gamma_2)\neq \{0\}
\]
Indeed, the Jones-Wenzl projector composed with the inclusion map is a chain
map, as the differential on $t^lq^{-l-r}w_{JW}$ has to be zero: since $w_{JW}$
is semi-simple in degree zero, a degree zero map would be a multiple of the
identity and would yield a reduction of the complex, and no higher degree maps
can exist as that would necessitate a target with higher degree shift. Then
since differentials are of positive degree, this map cannot be nullhomotopic.

  We will use a duality result, stating that, in the category $\Sfoam$, one has, given three multicurves $A$, $B$ and $C$:
  \[
\HOM(A,B\ast C) \simeq \HOM(B^{\ast}\ast A,C)
\]
Above $B^{\ast}$ stands for $B$ with reverse orientation. To see this, one can track the usual duality process. Taking $\phi$ in $\HOM(A,B\ast C)$, one forms:
\[
  (c_{B^{\ast}\ast B}\ast \id_C)\circ (\id_{B^{\ast}}\ast \phi)
\]
Above $c_{B^{\ast} B}$ stands for
the coevaluation, in other terms, for the cakepan cobordism on $B$. This process
is well-defined, since an isotopy change in $\phi$ induces a homotopy between
the resulting morphisms. The process in the other direction is similar, and it
is easy to see that the compositions yield identities up to isotopies.

 Back to our problem, we get that:
  \[
  \Hom(t^lq^{-l-r}\gamma_1^{\ast}\ast w, \gamma_2)\neq \{0\}
\]
This implies the following:
\[
  \Hom(\gamma_1^{\ast}\ast w,t^{-l}q^{l+r}\gamma_2)\neq \{0\}.
\]
  Since $r>0$ and the category is non-negatively graded, this implies that
  $\gamma_1^\ast\ast w$ is not linear, which contradicts the assumption on
  the minimality of $c$.

  If $r<0$, we focus on a leftmost $w$ and we run the same argument with:
  \[
  \Hom(\gamma_1\ast\gamma_2,t^lq^{-l-r}w)\neq \{0\}.
\]

  \end{proof}


\appendix
\section{Elementary cobordisms} \label{app:maps}

Throughout this section, the following parameters will appear: $a_3$, $a_4$, $a_{11}$, $\zeta$. They take value in $\{\pm 1\}$. The numbering comes from the determination of the moves below: we started with independent scalars for each move, and then throughout the invariance proof had to impose relations. Notice that because we wanted to ensure compatibility with classical invariance proofs, we started from scalar choices for classical $\rm{R}_{I}$ and $\rm{R}_{II}$ moves. It is possible that in full generality one could have extra parameters controlling the maps assigned to generating moves.

While describing these maps, we will only emphasize those foams that are not the obvious ones, and draw a version of them on the side. To save space, we abuse notation and do not differentiate between a foam and its upside-down version, hence the use of the same letter in the downward and upward arrows.

\subsection{\texorpdfstring{$\rm{R}_1$}{R1}}

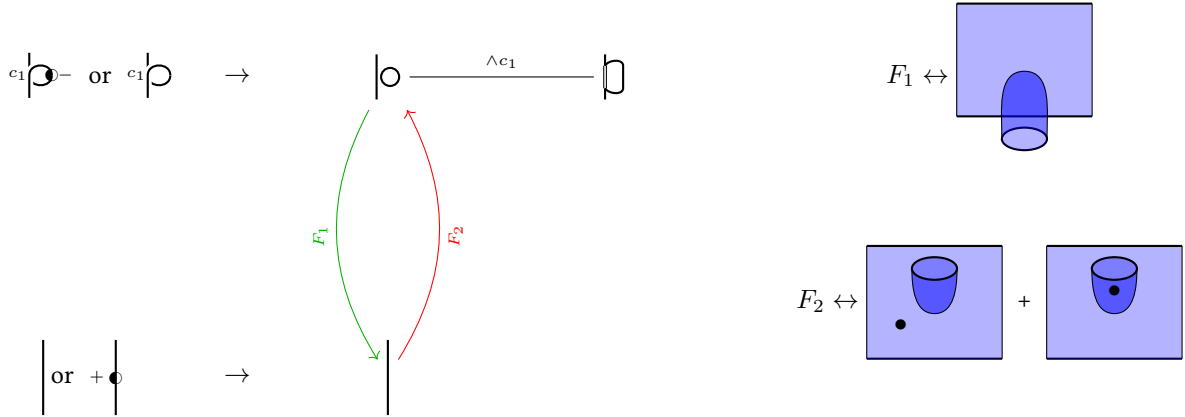
\begin{equation}
  \label{R11}

\end{equation}

In the next versions of the $R_I$ move, we only choose one of the two options for the twist. The other one is similar.

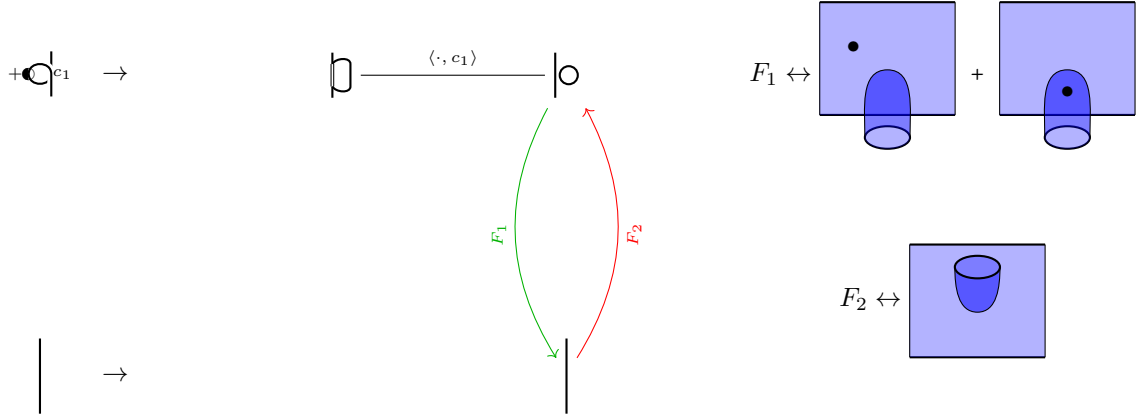
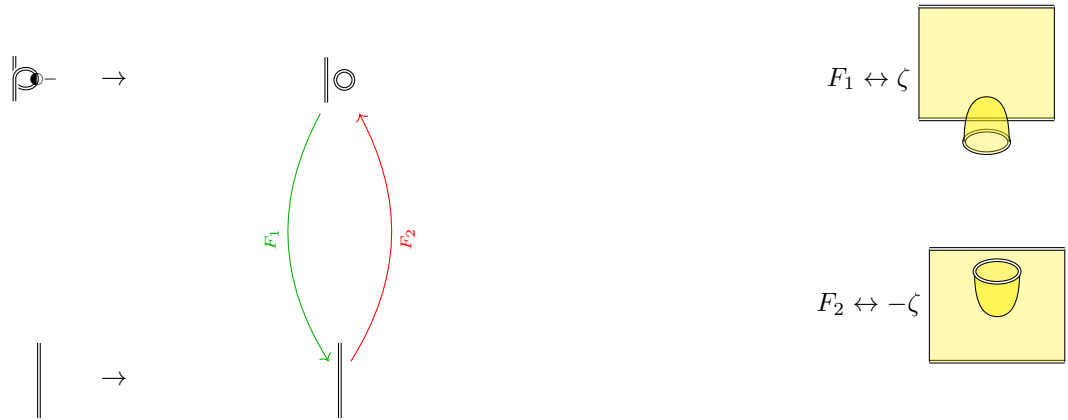
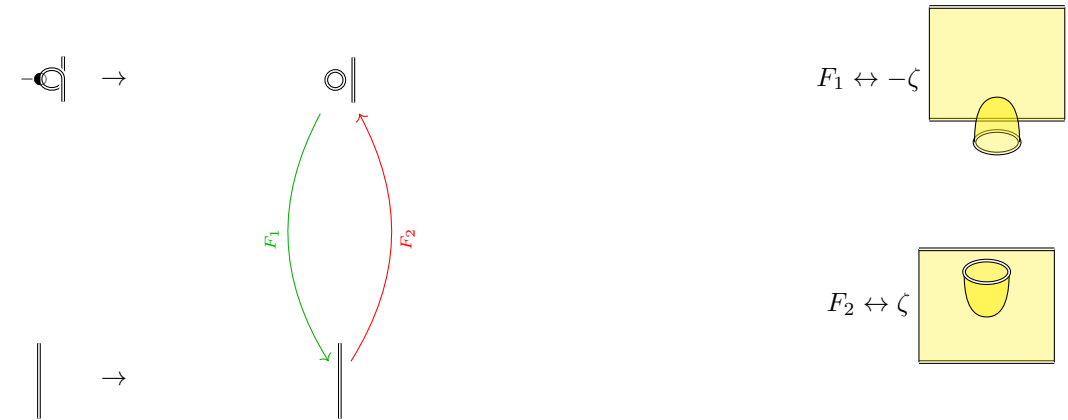
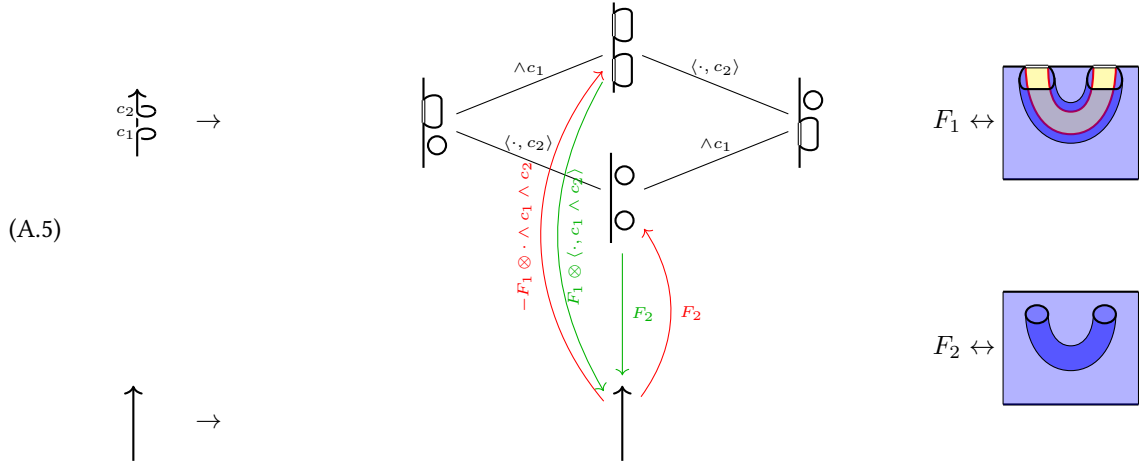
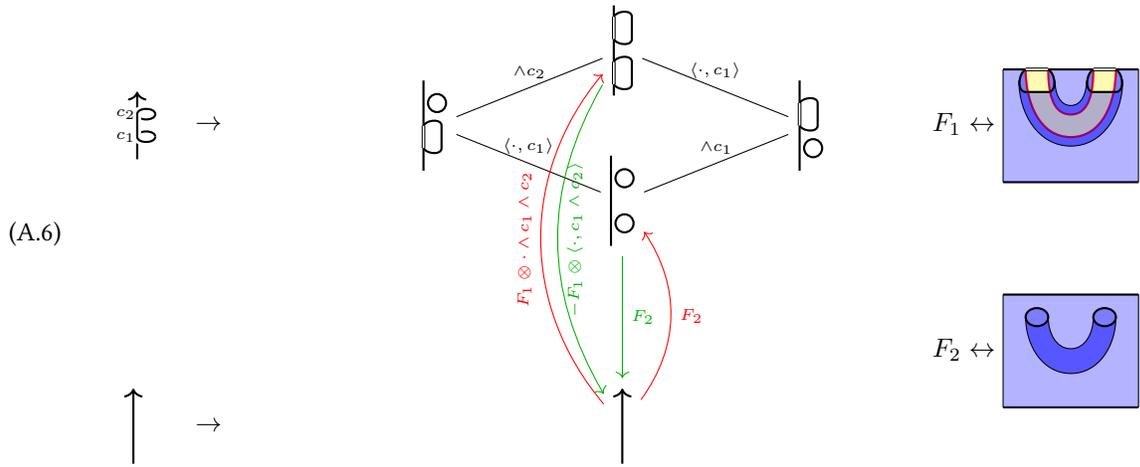
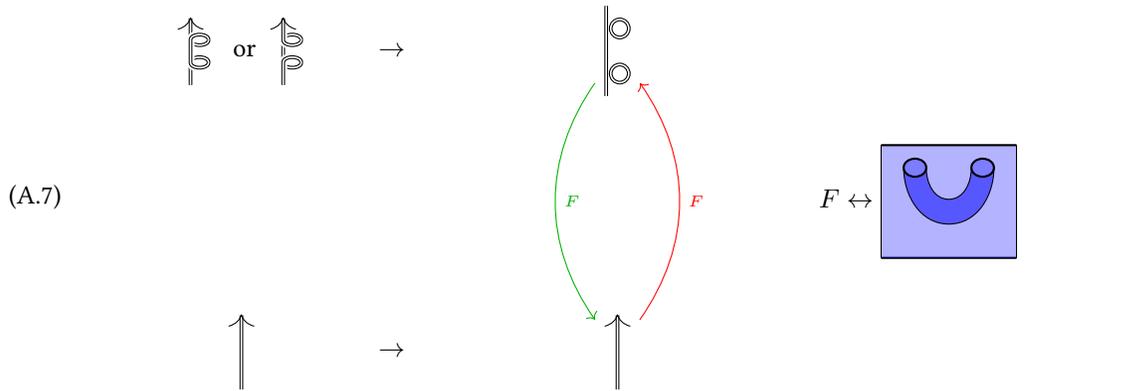
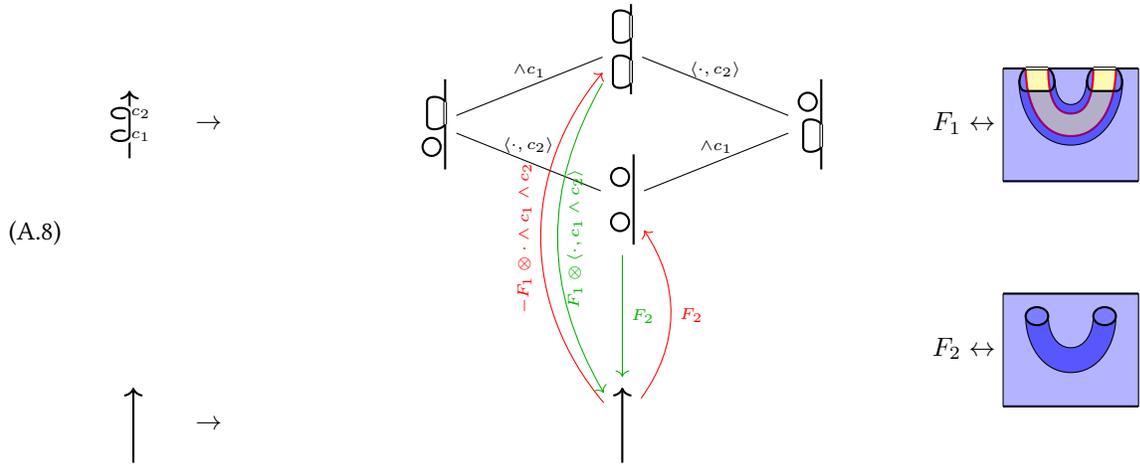
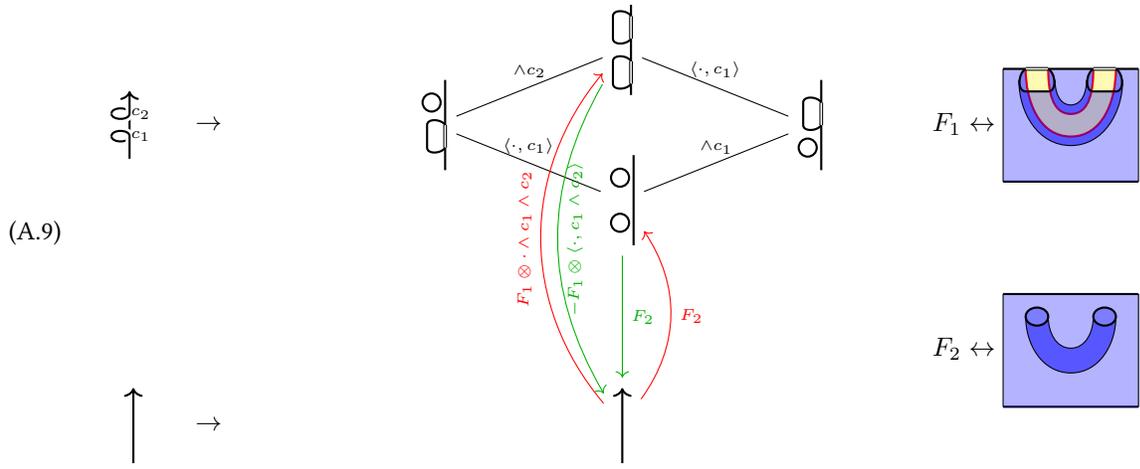
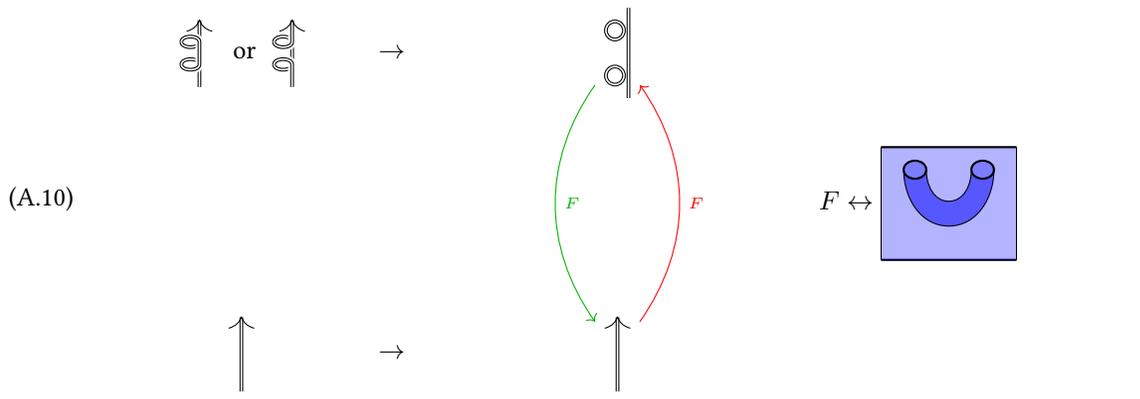
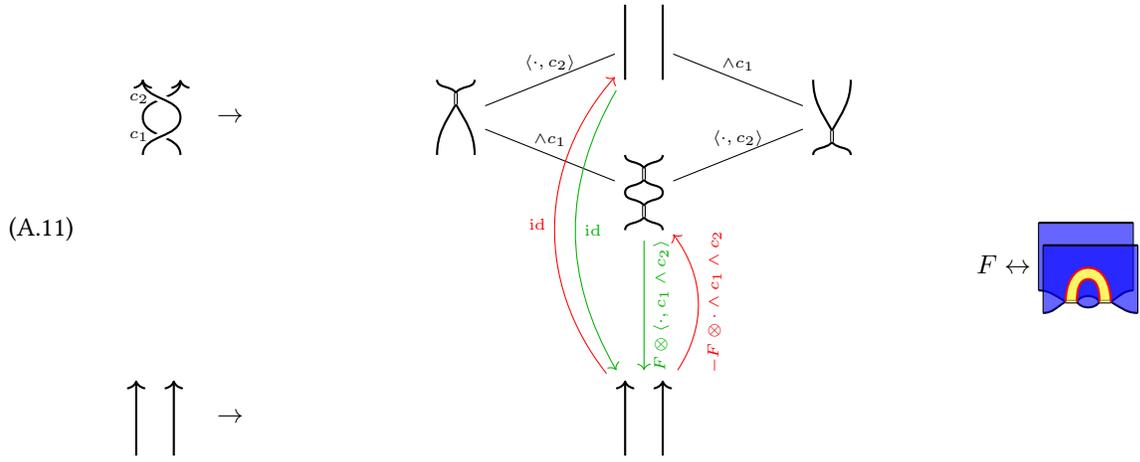
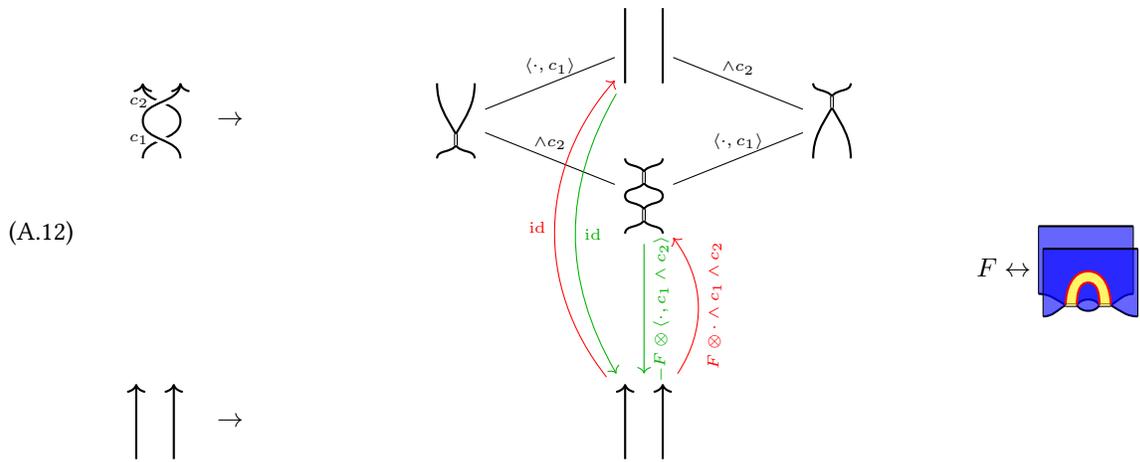
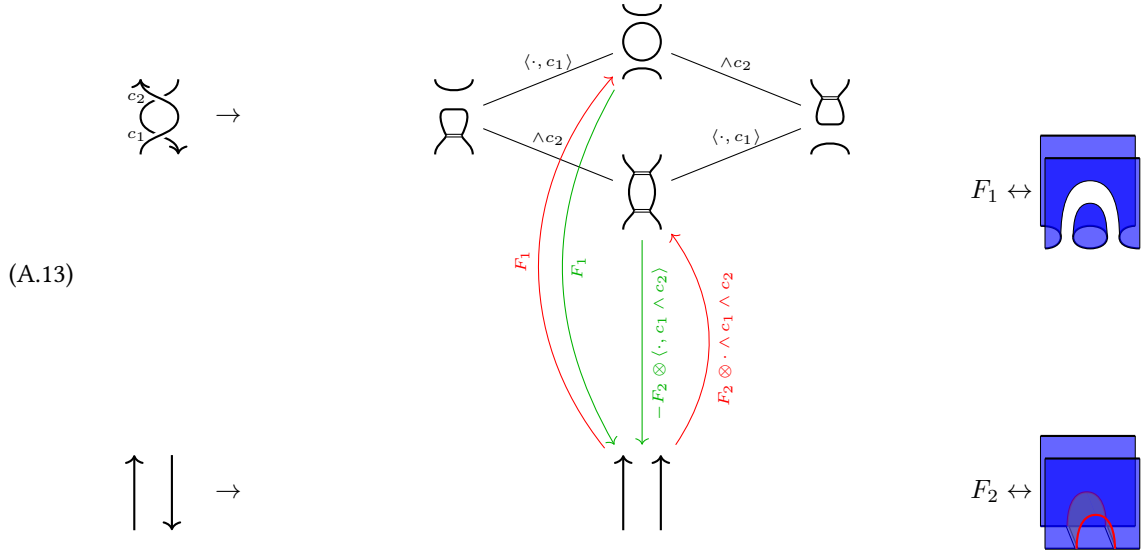
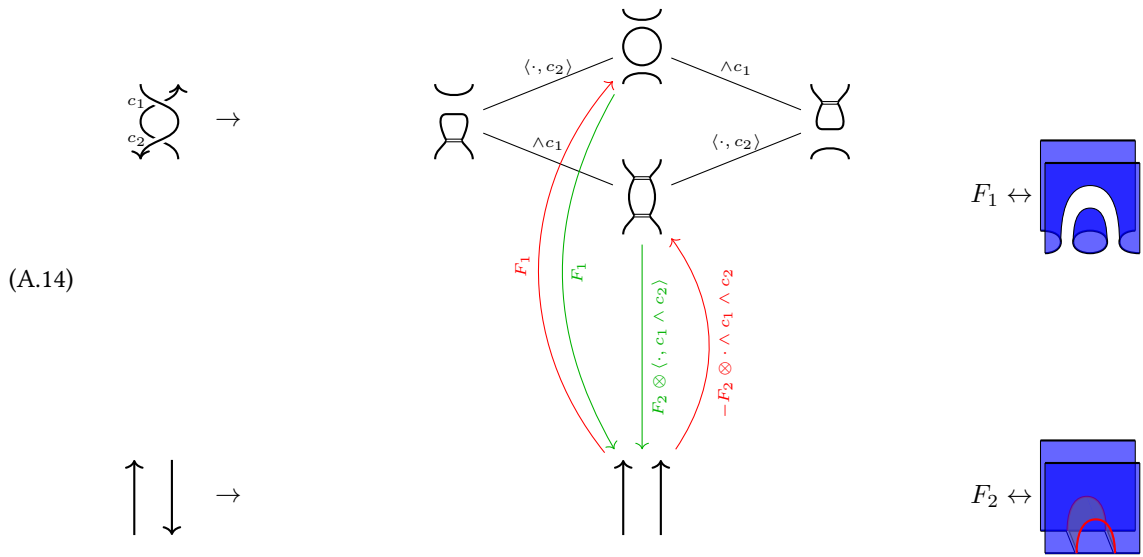
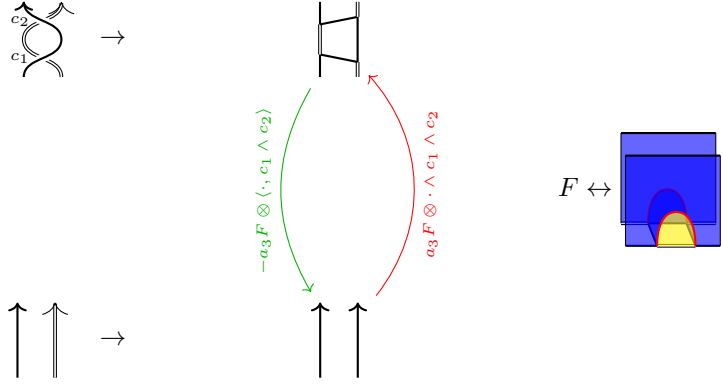
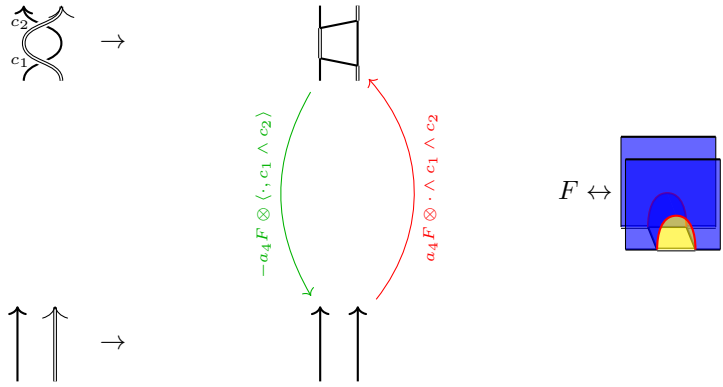
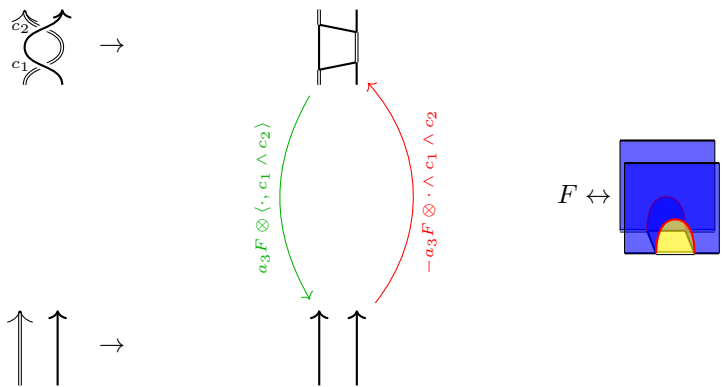
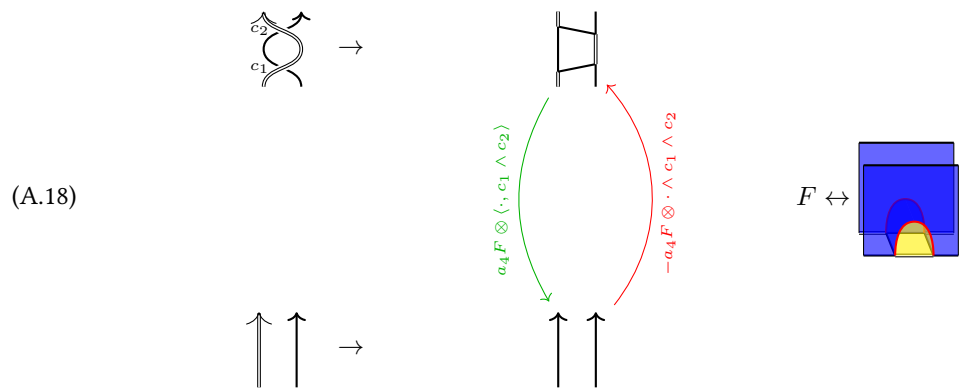
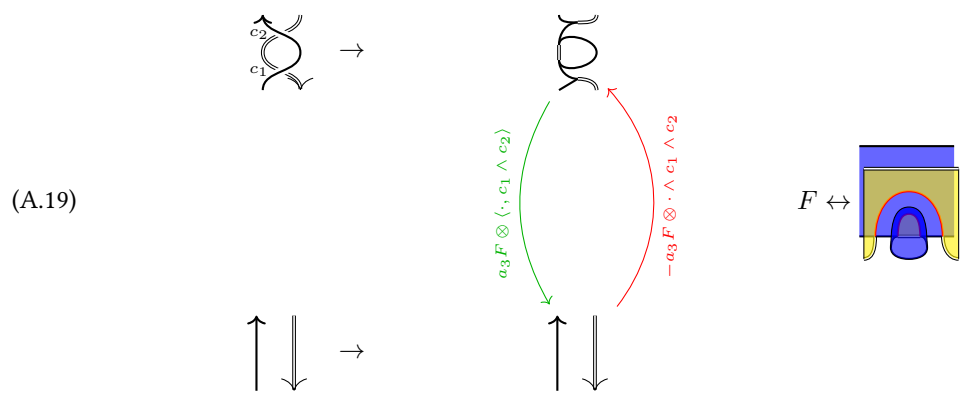
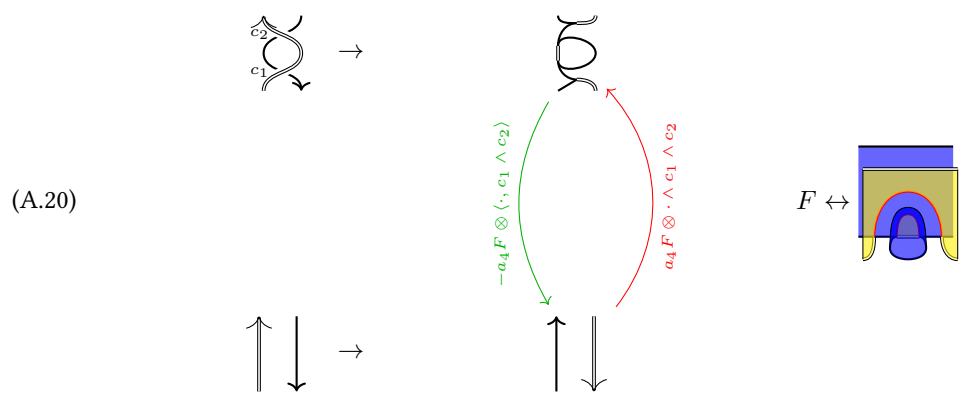
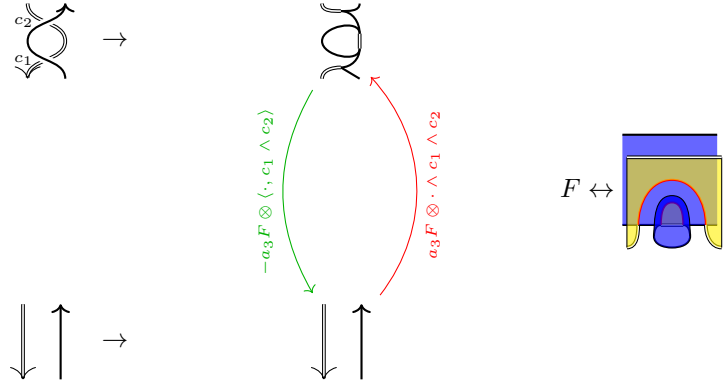
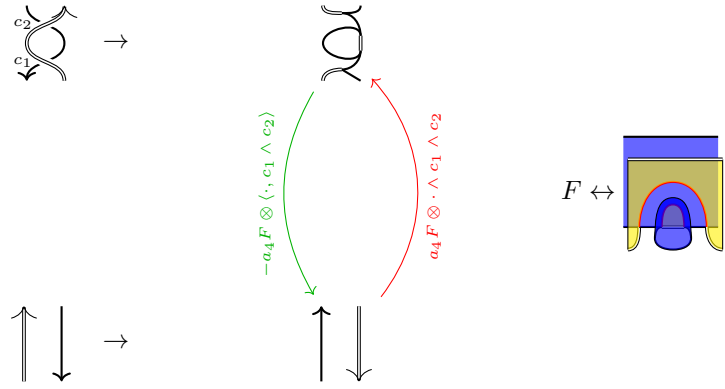
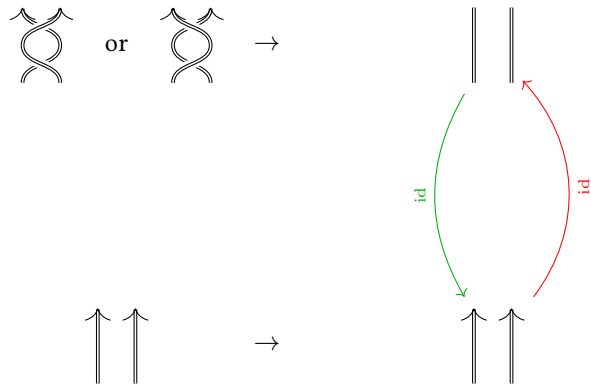
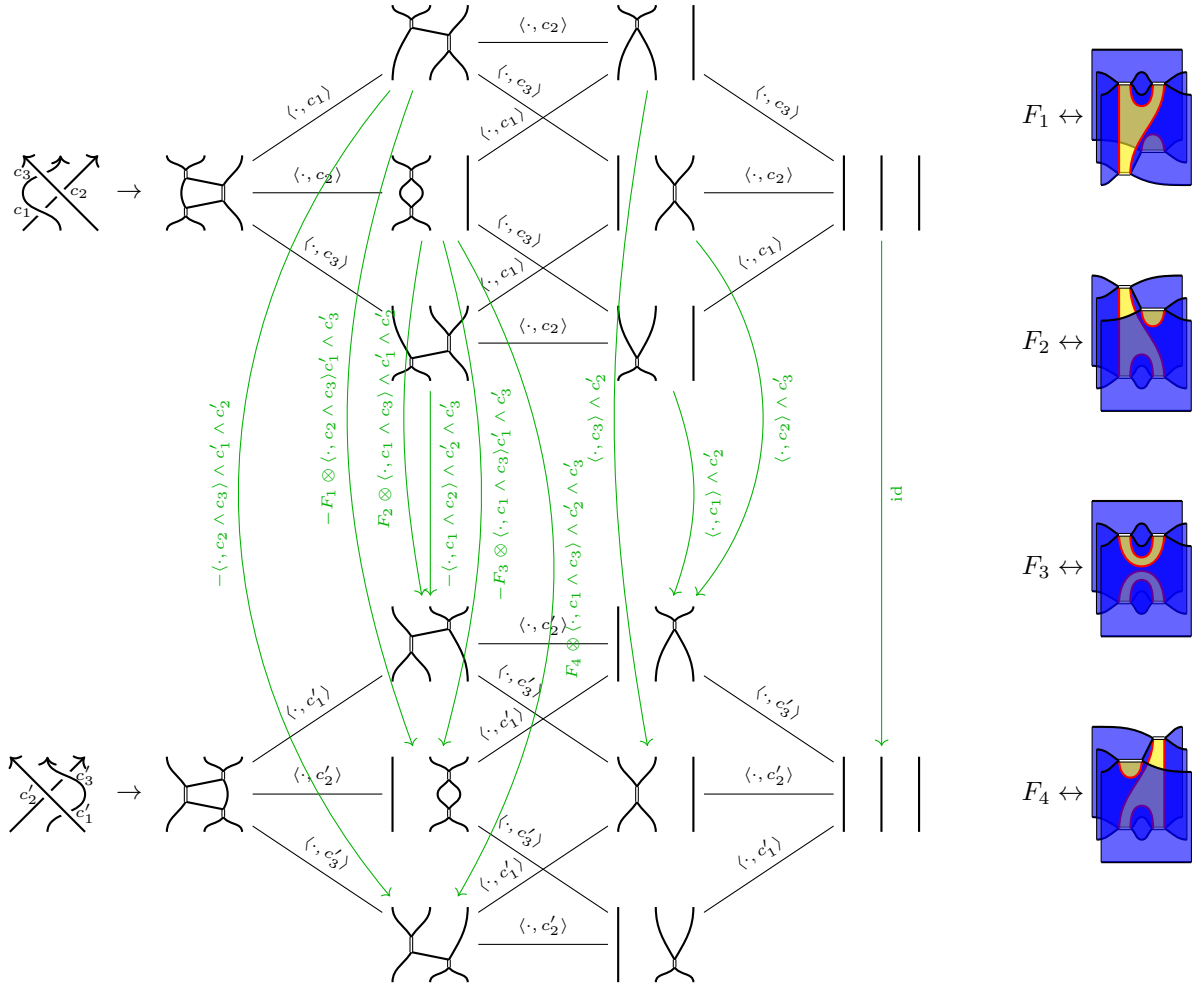
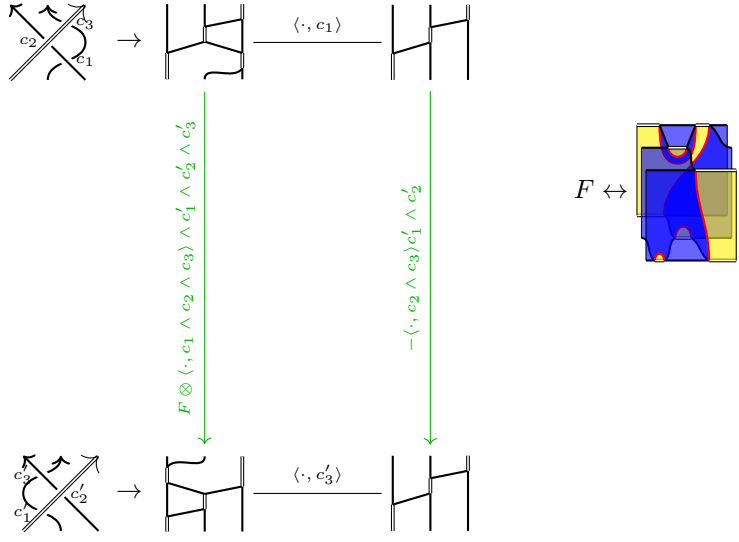
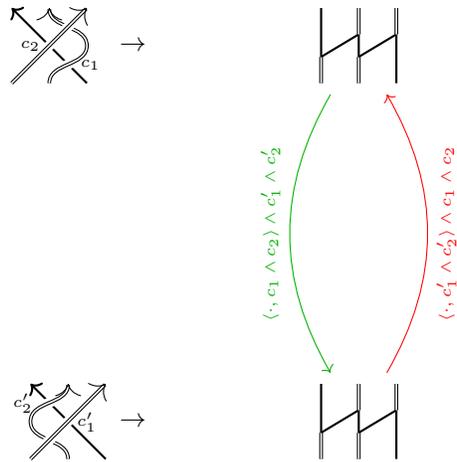
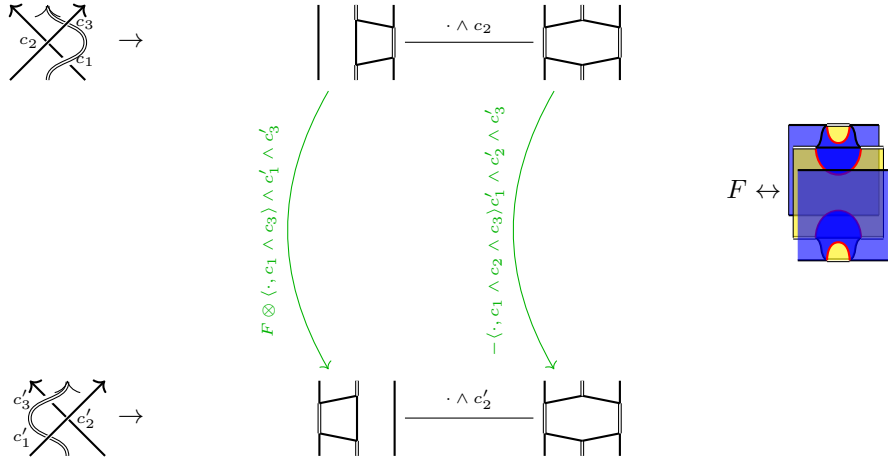
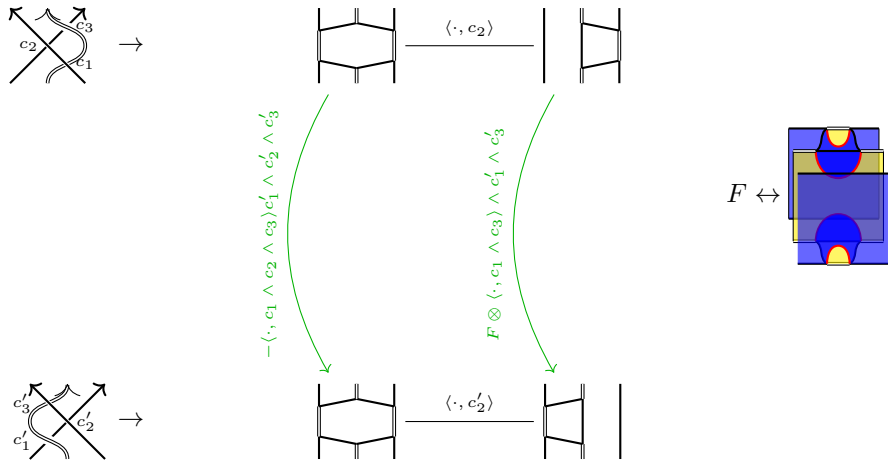
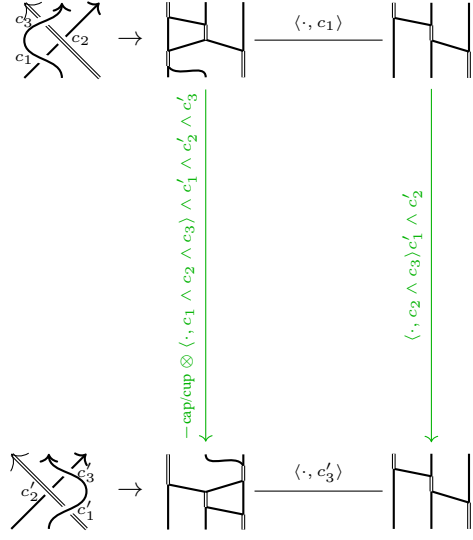
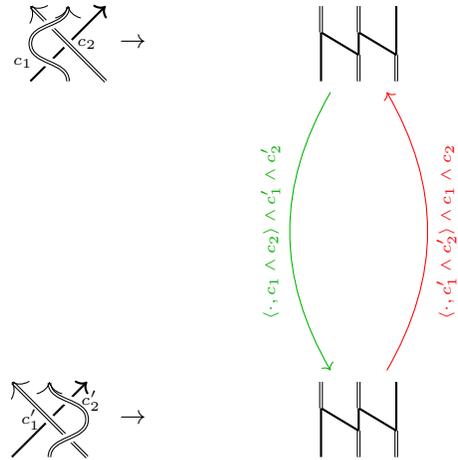
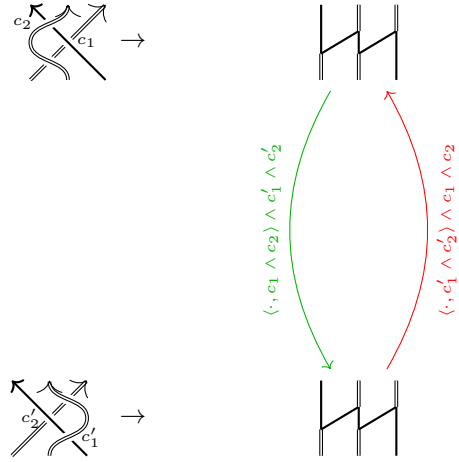
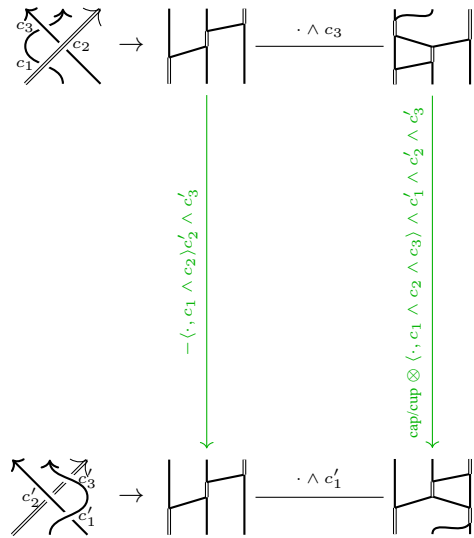
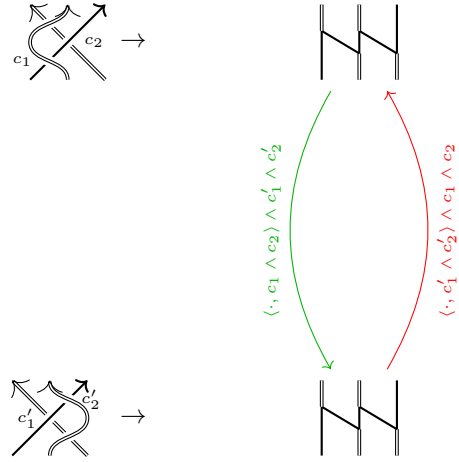
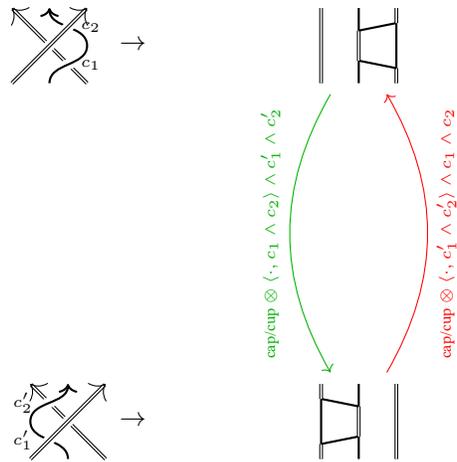
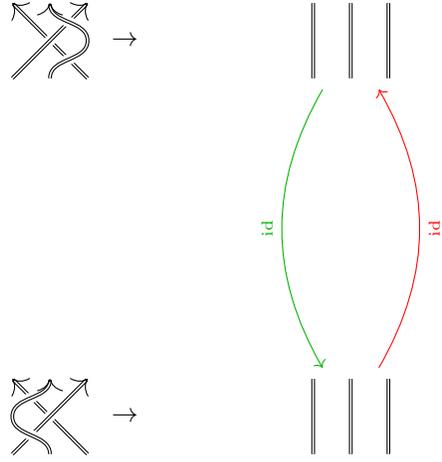
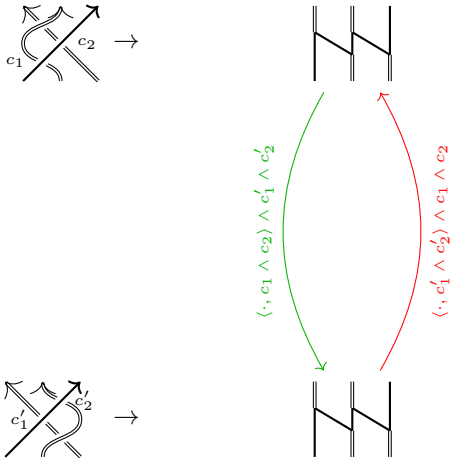
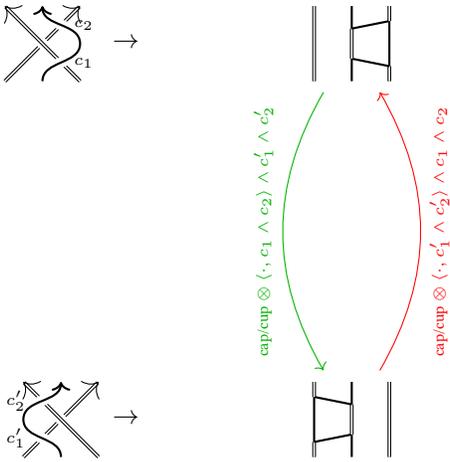
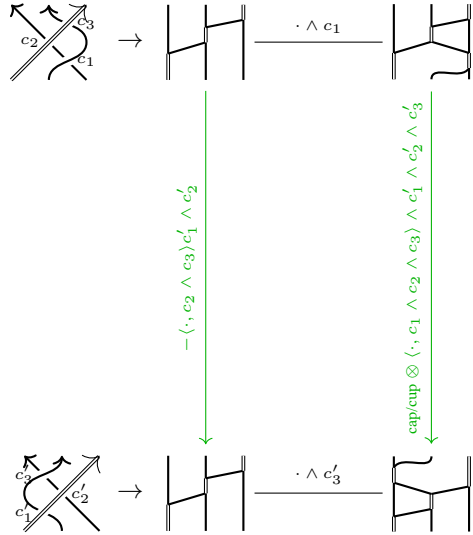
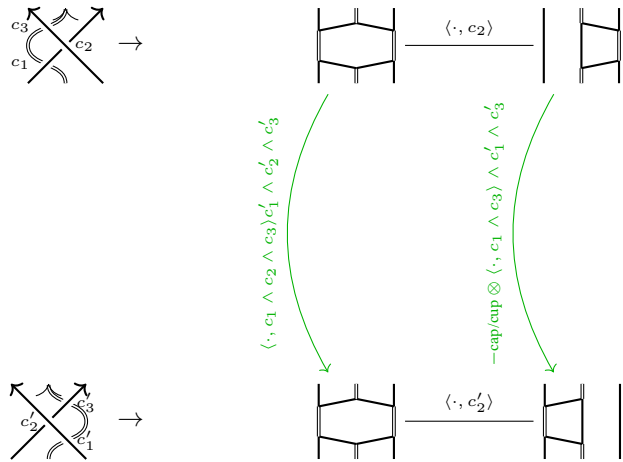
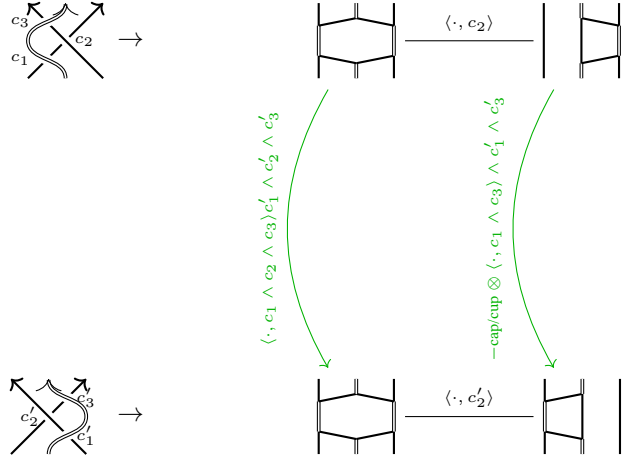
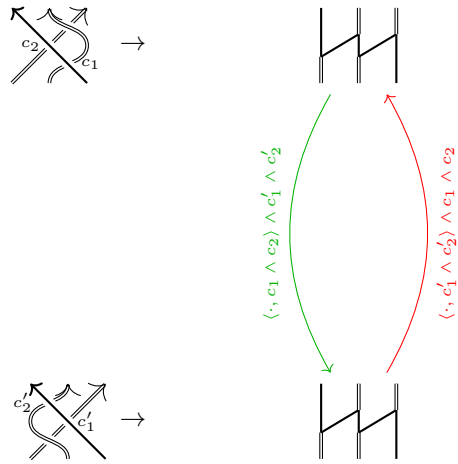
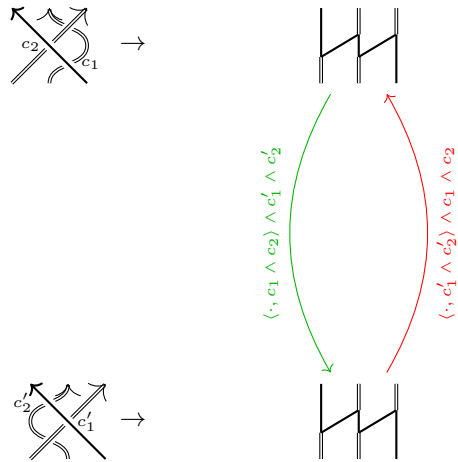
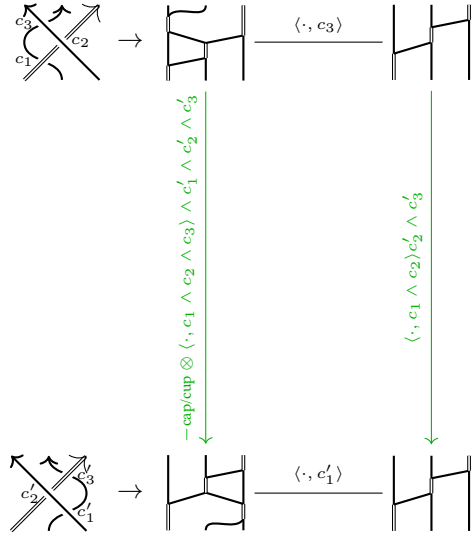
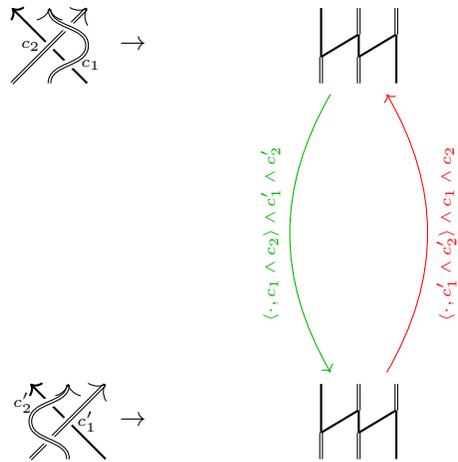
\begin{equation}
  \label{R12}
  %
\end{equation}

\subsection{\texorpdfstring{$\rm{R}'_1$}{R'1}}

\begin{equation}
  \label{R'11}
  %
\end{equation}

\subsection{\texorpdfstring{$\rm{R}_2$}{R2}}

\begin{equation}
  \label{R21}
  %
\end{equation}

\subsection{\texorpdfstring{$\rm{R}_3$}{R3}}

\begin{equation}
  %
    };
     \draw (A) -- (B) node[midway,above,rotate=30] {\tiny $\langle \cdot,c_1\rangle$};
    \draw (A) -- (C) node[midway,above] {\tiny $\langle \cdot,c_2\rangle$};
    \draw (A) -- (D) node[midway,above,rotate=-30] {\tiny $\langle \cdot,c_3\rangle$};
    \draw (B) -- (E) node[midway,above] {\tiny $\langle \cdot,c_2\rangle$};
    \draw (B) -- (F) node[near start,above,rotate=-30] {\tiny $\langle \cdot,c_3\rangle$};
    \draw (C) -- (E) node[near start,above,rotate=30] {\tiny $\langle \cdot,c_1\rangle$};
    \draw (C) -- (G) node[near start,above,rotate=-30] {\tiny $\langle \cdot,c_3\rangle$};
    \draw (D) -- (F) node[near start,above,rotate=30] {\tiny $\langle \cdot,c_1\rangle$};
    \draw (D) -- (G) node[midway,above] {\tiny $\langle \cdot,c_2\rangle$};
    \draw (E) -- (H) node[midway,above,rotate=-30] {\tiny $\langle \cdot,c_3\rangle$};
    \draw (F) -- (H) node[midway,above] {\tiny $\langle \cdot,c_2\rangle$};
    \draw (G) -- (H) node[midway,above,rotate=30] {\tiny $\langle \cdot,c_1\rangle$};
    \node at (-4,-8) {
      \begin{tikzpicture}[anchorbase,scale=.5]
        \draw [line width=4, white] (0,0) -- (2,2);
        \draw [thick, ->] (0,0) -- (2,2);
        \draw [line width=4, white] (1,0) to [out=90,in=-90] (2,1) to [out=90,in=-90] (1,2);
        \draw [thick, ->] (1,0) to [out=90,in=-90] (2,1) to [out=90,in=-90] (1,2); 
        \draw [line width=4, white] (2,0) -- (0,2);
        \draw [thick, ->] (2,0) -- (0,2);
        \node at (2,.5) {\tiny $c'_1$};
        \node at (.5,1) {\tiny $c'_2$};
        \node at (2,1.5) {\tiny $c'_3$};
      \end{tikzpicture}
    };
    \node at (-3,-8) {$\rightarrow$};
        \node (A') at (-2,-8) {
      \begin{tikzpicture}[anchorbase,scale=.5, xscale=-1]
        \draw [thick] (0,0) to [out=90,in=-120] (.5,.4);
        \draw [thick] (1,0) to [out=90,in=-60] (.5,.4);
        \draw [double] (.5,.4) -- (.5,.6);
        \draw [thick] (.5,.6) -- (1.5,.8);
        \draw [thick] (2,0) to [out=90,in=-60] (1.5,.8);
        \draw [double] (1.5,.8) -- (1.5,1.2);
        \draw [thick] (1.5,1.2) -- (.5,1.4);
        \draw [thick] (.5,.6) to [out=120,in=-120] (.5,1.4);
        \draw [double] (.5,1.4) -- (.5,1.6);
        \draw [thick] (.5,1.6) to [out=120,in=-90] (0,2);
        \draw [thick] (.5,1.6) to [out=60,in=-90] (1,2);
        \draw [thick] (1.5,1.2) to [out=60,in=-90] (2,2);
      \end{tikzpicture}
    };
        \node (B') at (1,-6) {
      \begin{tikzpicture}[anchorbase,scale=.5, xscale=-1]
        \draw [thick] (0,0) to [out=90,in=-120] (.5,1.4);
        \draw [thick] (1,0) to [out=90,in=-120] (1.5,.8);
        \draw [thick] (2,0) to [out=90,in=-60] (1.5,.8);
        \draw [double] (1.5,.8) -- (1.5,1.2);
        \draw [thick] (1.5,1.2) -- (.5,1.4);
         \draw [double] (.5,1.4) -- (.5,1.6);
        \draw [thick] (.5,1.6) to [out=120,in=-90] (0,2);
        \draw [thick] (.5,1.6) to [out=60,in=-90] (1,2);
        \draw [thick] (1.5,1.2) to [out=60,in=-90] (2,2);
      \end{tikzpicture}
    };
        \node (C') at (1,-8) {
      \begin{tikzpicture}[anchorbase,scale=.5, xscale=-1]
        \draw [thick] (0,0) to [out=90,in=-120] (.5,.4);
        \draw [thick] (1,0) to [out=90,in=-60] (.5,.4);
        \draw [double] (.5,.4) -- (.5,.6);
        \draw [thick] (.5,.6) to [out=60,in=-90] (.8,1) to [out=90,in=-60] (.5,1.4);
        \draw [thick] (2,0) -- (2,2);
        \draw [thick] (.5,.6) to [out=120,in=-90] (.2,1) to [out=90,in=-120] (.5,1.4);
        \draw [double] (.5,1.4) -- (.5,1.6);
        \draw [thick] (.5,1.6) to [out=120,in=-90] (0,2);
        \draw [thick] (.5,1.6) to [out=60,in=-90] (1,2);
      \end{tikzpicture}
    };
        \node (D') at (1,-10) {
      \begin{tikzpicture}[anchorbase,scale=.5, xscale=-1]
        \draw [thick] (0,0) to [out=90,in=-120] (.5,.4);
        \draw [thick] (1,0) to [out=90,in=-60] (.5,.4);
        \draw [double] (.5,.4) -- (.5,.6);
        \draw [thick] (.5,.6) -- (1.5,.8);
        \draw [thick] (2,0) to [out=90,in=-60] (1.5,.8);
        \draw [double] (1.5,.8) -- (1.5,1.2);
        \draw [thick] (1.5,1.2) to [out=120,in=-90] (1,2);
        \draw [thick] (.5,.6) to [out=120,in=-90] (0,2);
        \draw [thick] (1.5,1.2) to [out=60,in=-90] (2,2);
      \end{tikzpicture}
    };
    \node (E') at (4,-6) {
      \begin{tikzpicture}[anchorbase,scale=.5, xscale=-1]
        \draw [thick] (0,0) to [out=90,in=-120] (.5,1.4);
        \draw [thick] (1,0) to [out=90,in=-60] (.5,1.4);
        \draw [thick] (2,0) -- (2,2);
        \draw [double] (.5,1.4) -- (.5,1.6);
        \draw [thick] (.5,1.6) to [out=120,in=-90] (0,2);
        \draw [thick] (.5,1.6) to [out=60,in=-90] (1,2);
      \end{tikzpicture}
    };
    \node (F') at (4,-8) {
      \begin{tikzpicture}[anchorbase,scale=.5, xscale=-1]
        \draw [thick] (0,0) -- (0,2);
        \draw [thick] (1,0) to [out=90,in=-120] (1.5,.8);
        \draw [thick] (2,0) to [out=90,in=-60] (1.5,.8);
        \draw [double] (1.5,.8) -- (1.5,1.2);
        \draw [thick] (1.5,1.2) to [out=120,in=-90] (1,2);
        \draw [thick] (1.5,1.2) to [out=60,in=-90] (2,2);
      \end{tikzpicture}
    };
    \node (G') at (4,-10) {
      \begin{tikzpicture}[anchorbase,scale=.5, xscale=-1]
        \draw [thick] (0,0) to [out=90,in=-120] (.5,.4);
        \draw [thick] (1,0) to [out=90,in=-60] (.5,.4);
        \draw [thick] (2,0) -- (2,2);
        \draw [double] (.5,.4) -- (.5,.6);
        \draw [thick] (.5,.6) to [out=120,in=-90] (0,2);
        \draw [thick] (.5,.6) to [out=60,in=-90] (1,2);
      \end{tikzpicture}
    };
    \node (H') at (7,-8) {
      \begin{tikzpicture}[anchorbase,scale=.5, xscale=-1]
        \draw [thick] (0,0) -- (0,2);
        \draw [thick] (1,0) -- (1,2);
        \draw [thick] (2,0) -- (2,2);
      \end{tikzpicture}
    };
     \draw (A') -- (B') node[midway,above,rotate=30] {\tiny $\langle \cdot,c'_1\rangle$};
    \draw (A') -- (C') node[midway,above] {\tiny $\langle \cdot,c'_2\rangle$};
    \draw (A') -- (D') node[midway,above,rotate=-30] {\tiny $\langle \cdot,c'_3\rangle$};
    \draw (B') -- (E') node[midway,above] {\tiny $\langle \cdot,c'_2\rangle$};
    \draw (B') -- (F') node[near start,above,rotate=-30] {\tiny $\langle \cdot,c'_3\rangle$};
    \draw (C') -- (E') node[near start,above,rotate=30] {\tiny $\langle \cdot,c'_1\rangle$};
    \draw (C') -- (G') node[near start,above,rotate=-30] {\tiny $\langle \cdot,c'_3\rangle$};
    \draw (D') -- (F') node[near start,above,rotate=30] {\tiny $\langle \cdot,c'_1\rangle$};
    \draw (D') -- (G') node[midway,above] {\tiny $\langle \cdot,c'_2\rangle$};
    \draw (E') -- (H') node[midway,above,rotate=-30] {\tiny $\langle \cdot,c'_3\rangle$};
    \draw (F') -- (H') node[midway,above] {\tiny $\langle \cdot,c'_2\rangle$};
    \draw (G') -- (H') node[midway,above,rotate=30] {\tiny $\langle \cdot,c'_1\rangle$};
    \draw [green,->] (B) edge[bend right=40] node [above,rotate=90] {\tiny $-\langle \cdot,c_2\wedge c_3\rangle \wedge c'_1\wedge c'_2$} (D');
    \draw [green,->] (B) edge[bend right=20] node [above,rotate=90] {\tiny $-F_1\otimes \langle \cdot,c_2\wedge c_3\rangle c'_1\wedge c'_3$} (C');
    \draw [green,->] (C) edge[bend right=10] node [above,rotate=90] {\tiny $F_2\otimes \langle \cdot,c_1\wedge c_3\rangle \wedge c'_1\wedge c'_2$} (B');
    \draw [green,->] (C) edge[bend left=15] node [below,rotate=90] {\tiny $-F_3\otimes\langle \cdot,c_1\wedge c_3\rangle  c'_1\wedge c'_3$} (C');
    \draw [green,->] (C) edge[bend left=30] node [below,rotate=90] {\tiny $F_4\otimes \langle \cdot,c_1\wedge c_3\rangle \wedge c'_2\wedge c'_3$} (D');
    \draw [green,->] (D) edge node [below,rotate=90] {\tiny $-\langle \cdot,c_1\wedge c_2\rangle \wedge c'_2\wedge c'_3 $} (B');
    \draw [green,->] (E) edge[bend right=10] node [above,rotate=90] {\tiny $\langle \cdot,c_3\rangle \wedge c'_2$} (F');
    \draw [green,->] (F) edge[bend left=40] node [below,rotate=90] {\tiny $\langle \cdot,c_2\rangle \wedge c'_3$} (E');
    \draw [green,->] (G) edge[bend left=20] node [below,rotate=90] {\tiny $\langle \cdot,c_1\rangle \wedge c'_2$} (E');
    \draw [green,->] (H) edge node [below,rotate=90] {\tiny $\id $} (H');
        \node at (10,1) {$F_1\leftrightarrow$  
      \begin{tikzpicture}[anchorbase,scale=.3]
        \draw [thick] (.2,0) to [out=0,in=-150] (1,.5);
        \draw [thick] (0,1) to [out=0,in=150] (1,.5);
        \draw [double] (1,.5) -- (1.5,.5);
        \draw [thick] (1.5,.5) -- (2,1.5);
        \draw [thick] (-.2,2) to [out=0,in=150] (2,1.5);
        \draw [double] (2,1.5) -- (3,1.5);
        \draw [thick] (1.5,.5) to [out=-30,in=180] (4.2,0);
        \draw [thick] (3,1.5) to [out=-30,in=180] (4,1);
        \draw [thick] (3,1.5) to [out=30,in=180] (3.8,2);
        \fill [blue, opacity=.6] (-.2,2) to [out=0,in=150] (2,1.5) to [out=90,in=180] (2.5,2.25) to [out=0,in=90] (3,1.5) to [out=30,in=180] (3.8,2) -- (3.8,6) -- (-.2,6) -- (-.2,2);
        \fill [yellow, opacity=.6] (2,1.5) to [out=90,in=180] (2.5,2.25) to [out=0,in=90] (3,1.5);
        \draw (-.2,2) -- (-.2,6);
        \draw (3.8,2) -- (3.8,6);
        \fill [blue, opacity=.6] (0,1) to [out=0,in=150] (1,.5) -- (1,4.5) to [out=150,in=0] (0,5) -- (0,1);
        \draw (0,1) -- (0,5);
        \fill [blue,opacity=.6] (1.5,.5) -- (2,1.5) to [out=90,in=180] (2.5,2.25) to [out=0,in=90] (3,1.5) to [out=-30,in=180] (4,1) -- (4,5) to [out=180,in=30] (3,4.5) to [out=-90,in=90] (1.5,.5);
        \draw (4,1) -- (4,5);
        \draw [thick, red] (2,1.5) to [out=90,in=180] (2.5,2.25) to [out=0,in=90] (3,1.5);
        \fill [blue, opacity=.6] (1.5,4.5) to [out=30,in=180] (2,5) to [out=0,in=150] (2.5,4.5) to [out=-90,in=0] (2,3.5) to [out=180,in=-90] (1.5,4.5);
        \fill [blue, opacity=.6] (1.5,4.5) to [out=-30,in=180] (2,4) to [out=0,in=-150] (2.5,4.5) to [out=-90,in=0] (2,3.5) to [out=180,in=-90] (1.5,4.5);
        \fill [yellow, opacity=.6]  (2.5,4.5) to [out=-90,in=0] (2,3.5) to [out=180,in=-90] (1.5,4.5) -- (1,4.5) to [out=-90,in=90] (1,.5) -- (1.5,.5) to [out=90,in=-90] (3,4.5) -- (2.5,4.5);
        \fill [blue, opacity=.6] (.2,0) to [out=0,in=-150] (1,.5) -- (1,4.5) to [out=-150,in=0] (.2,4);
        \fill [blue, opacity=.6] (1.5,.5) to [out=90,in=-90] (3,4.5) to [out=-30,in=180] (4.2,4) -- (4.2,0) to [out=180,in=-30] (1.5,.5);
        \draw (.2,0) -- (.2,4);
        \draw (4.2,0) -- (4.2,4);
        \draw [thick, red] (2.5,4.5) to [out=-90,in=0] (2,3.5) to [out=180,in=-90] (1.5,4.5);
        \draw [thick, red] (1,.5) to [out=90,in=-90] (1,4.5);
        \draw [thick, red] (1.5,.5) to [out=90,in=-90] (3,4.5);
        \draw [thick] (-.2,6) -- (3.8,6);
        \draw [thick] (0,5) to [out=0,in=150] (1,4.5);
        \draw [thick] (.2,4) to [out=0,in=-150] (1,4.5);
        \draw [double] (1,4.5) -- (1.5,4.5);
        \draw [thick] (1.5,4.5) to [out=30,in=180] (2,5) to [out=0,in=150] (2.5,4.5);
        \draw [thick] (1.5,4.5) to [out=-30,in=180] (2,4) to [out=0,in=-150] (2.5,4.5);
        \draw [double] (2.5,4.5) -- (3,4.5);
        \draw [thick] (3,4.5) to [out=30,in=180] (4,5);
        \draw [thick] (3,4.5) to [out=-30,in=180] (4.2,4);
  \end{tikzpicture}
  };
        \node at (10,-2) {$F_2\leftrightarrow$  
      \begin{tikzpicture}[anchorbase,scale=.3,yscale=-1]
        \draw [thick] (.2,6) -- (4.2,6);
        \draw [thick] (0,5) to [out=0,in=150] (1,4.5);
        \draw [thick] (-.2,4) to [out=0,in=-150] (1,4.5);
        \draw [double] (1,4.5) -- (1.5,4.5);
        \draw [thick] (1.5,4.5) to [out=30,in=180] (2,5) to [out=0,in=150] (2.5,4.5);
        \draw [thick] (1.5,4.5) to [out=-30,in=180] (2,4) to [out=0,in=-150] (2.5,4.5);
        \draw [double] (2.5,4.5) -- (3,4.5);
        \draw [thick] (3,4.5) to [out=30,in=180] (4,5);
        \draw [thick] (3,4.5) to [out=-30,in=180] (3.8,4);
        \fill [yellow, opacity=.6]  (2.5,4.5) to [out=-90,in=0] (2,3.5) to [out=180,in=-90] (1.5,4.5) -- (1,4.5) to [out=-90,in=90] (1,.5) -- (1.5,.5) to [out=90,in=-90] (3,4.5) -- (2.5,4.5);
        \fill [blue, opacity=.6] (-.2,0) to [out=0,in=-150] (1,.5) -- (1,4.5) to [out=-150,in=0] (-.2,4);
        \fill [blue, opacity=.6] (1.5,.5) to [out=90,in=-90] (3,4.5) to [out=-30,in=180] (3.8,4) -- (3.8,0) to [out=180,in=-30] (1.5,.5);
        \draw (-.2,0) -- (-.2,4);
        \draw (3.8,0) -- (3.8,4);
        \fill [blue, opacity=.6] (1.5,4.5) to [out=30,in=180] (2,5) to [out=0,in=150] (2.5,4.5) to [out=-90,in=0] (2,3.5) to [out=180,in=-90] (1.5,4.5);
        \fill [blue, opacity=.6] (1.5,4.5) to [out=-30,in=180] (2,4) to [out=0,in=-150] (2.5,4.5) to [out=-90,in=0] (2,3.5) to [out=180,in=-90] (1.5,4.5);
        \fill [blue, opacity=.6] (0,1) to [out=0,in=150] (1,.5) -- (1,4.5) to [out=150,in=0] (0,5) -- (0,1);
        \draw (0,1) -- (0,5);
        \fill [blue,opacity=.6] (1.5,.5) -- (2,1.5) to [out=90,in=180] (2.5,2.25) to [out=0,in=90] (3,1.5) to [out=-30,in=180] (4,1) -- (4,5) to [out=180,in=30] (3,4.5) to [out=-90,in=90] (1.5,.5);
        \draw (4,1) -- (4,5);
         \draw [thick, red] (1,.5) to [out=90,in=-90] (1,4.5);
        \draw [thick, red] (1.5,.5) to [out=90,in=-90] (3,4.5);
        \draw [thick, red] (2.5,4.5) to [out=-90,in=0] (2,3.5) to [out=180,in=-90] (1.5,4.5);
        \fill [blue, opacity=.6] (.2,2) to [out=0,in=150] (2,1.5) to [out=90,in=180] (2.5,2.25) to [out=0,in=90] (3,1.5) to [out=30,in=180] (4.2,2) -- (4.2,6) -- (.2,6) -- (.2,2);
        \fill [yellow, opacity=.6] (2,1.5) to [out=90,in=180] (2.5,2.25) to [out=0,in=90] (3,1.5);
        \draw [thick, red] (2,1.5) to [out=90,in=180] (2.5,2.25) to [out=0,in=90] (3,1.5);
        \draw (.2,2) -- (.2,6);
        \draw (4.2,2) -- (4.2,6);
        \draw [thick] (-.2,0) to [out=0,in=-150] (1,.5);
        \draw [thick] (0,1) to [out=0,in=150] (1,.5);
        \draw [double] (1,.5) -- (1.5,.5);
        \draw [thick] (1.5,.5) -- (2,1.5);
        \draw [thick] (.2,2) to [out=0,in=150] (2,1.5);
        \draw [double] (2,1.5) -- (3,1.5);
        \draw [thick] (1.5,.5) to [out=-30,in=180] (3.8,0);
        \draw [thick] (3,1.5) to [out=-30,in=180] (4,1);
        \draw [thick] (3,1.5) to [out=30,in=180] (4.2,2);
  \end{tikzpicture}
  };
        \node at (10,-5) {$F_3\leftrightarrow$  
      \begin{tikzpicture}[anchorbase,scale=.3]
        \draw [thick] (.2,0) -- (4.2,0);
        \draw [thick] (0,1) to [out=0,in=-150] (1,1.5);
        \draw [thick] (-.2,2) to [out=0,in=150] (1,1.5);
        \draw [double] (1,1.5) -- (1.5,1.5);
        \draw [thick] (1.5,1.5) to [out=-30,in=180] (2,1) to [out=0,in=-150] (2.5,1.5);
        \draw [thick] (1.5,1.5) to [out=30,in=180] (2,2) to [out=0,in=150] (2.5,1.5);
        \draw [double] (2.5,1.5) -- (3,1.5);
        \draw [thick] (3,1.5) to [out=-30,in=180] (4,1);
        \draw [thick] (3,1.5) to [out=30,in=180] (3.8,2);
        \fill [blue, opacity=.6] (-.2,2) to [out=0,in=150] (1,1.5) to [out=90,in=180] (2,2.9) to [out=0,in=90] (3,1.5) to [out=30,in=180] (3.8,2) -- (3.8,6) -- (-.2,6) -- (-.2,2);
        \fill [yellow, opacity=.6] (1,1.5) to [out=90,in=180] (2,2.9) to [out=0,in=90] (3,1.5) -- (2.5,1.5) to [out=90,in=0] (2,2.5) to [out=180,in=90] (1.5,1.5) -- (1,1.5);
        \fill [blue, opacity=.6] (2.5,1.5) to [out=90,in=0] (2,2.5) to [out=180,in=90] (1.5,1.5) (1.5,1.5) to [out=30,in=180] (2,1) to [out=0,in=150] (2.5,1.5);
        \fill [blue, opacity=.6] (2.5,1.5) to [out=90,in=0] (2,2.5) to [out=180,in=90] (1.5,1.5) (1.5,1.5) to [out=-30,in=180] (2,1) to [out=0,in=-150] (2.5,1.5);
        \draw (-.2,2) -- (-.2,6);
        \draw (3.8,2) -- (3.8,6);
        \fill [blue, opacity=.6] (0,1) to [out=0,in=-150] (1,1.5) to [out=90,in=180] (2,2.9) to [out=0,in=90] (3,1.5) to [out=-30,in=180] (4,1) -- (4,5) to [out=-180,in=30] (3,4.5) to [out=-90,in=0] (2,3.1) to [out=180,in=-90] (1,4.5) to [out=150,in=0] (0,5) -- (0,1);
        \draw [red, thick] (1,1.5) to [out=90,in=180] (2,2.9) to [out=0,in=90] (3,1.5);
        \draw [red, thick] (2.5,1.5) to [out=90,in=0] (2,2.5) to [out=180,in=90] (1.5,1.5);
        \draw (0,1) -- (0,5);
        \draw (4,1) -- (4,5);
        \fill [blue, opacity=.6] (1.5,4.5) to [out=30,in=180] (2,5) to [out=0,in=150] (2.5,4.5) to [out=-90,in=0] (2,3.5) to [out=180,in=-90] (1.5,4.5);
        \fill [blue, opacity=.6] (1.5,4.5) to [out=-30,in=180] (2,4) to [out=0,in=-150] (2.5,4.5) to [out=-90,in=0] (2,3.5) to [out=180,in=-90] (1.5,4.5);
        \fill [yellow, opacity=.6]  (2.5,4.5) to [out=-90,in=0] (2,3.5) to [out=180,in=-90] (1.5,4.5) -- (1,4.5) to [out=-90,in=180] (2,3.1) to [out=0,in=-90] (3,4.5) -- (2.5,4.5);
       \fill [blue, opacity=.6] (.2,4) to [out=0,in=-150]  (1,4.5) to [out=-90,in=180] (2,3.1) to [out=0,in=-90] (3,4.5) to [out=-30,in=180] (4.2,4) -- (4.2,0) -- (.2,0) -- (.2,4);
        \draw [thick, red] (2.5,4.5) to [out=-90,in=0] (2,3.5) to [out=180,in=-90] (1.5,4.5);
        \draw [thick, red] (1,4.5) to [out=-90,in=180] (2,3.1) to [out=0,in=-90] (3,4.5);
        \draw (.2,0) -- (.2,4);
        \draw (4.2,0) -- (4.2,4);
        \draw [thick] (-.2,6) -- (3.8,6);
        \draw [thick] (0,5) to [out=0,in=150] (1,4.5);
        \draw [thick] (.2,4) to [out=0,in=-150] (1,4.5);
        \draw [double] (1,4.5) -- (1.5,4.5);
        \draw [thick] (1.5,4.5) to [out=30,in=180] (2,5) to [out=0,in=150] (2.5,4.5);
        \draw [thick] (1.5,4.5) to [out=-30,in=180] (2,4) to [out=0,in=-150] (2.5,4.5);
        \draw [double] (2.5,4.5) -- (3,4.5);
        \draw [thick] (3,4.5) to [out=30,in=180] (4,5);
        \draw [thick] (3,4.5) to [out=-30,in=180] (4.2,4);
  \end{tikzpicture}
  };
        \node at (10,-8) {$F_4\leftrightarrow$  
      \begin{tikzpicture}[anchorbase,scale=.3,rotate=180]
        \draw [thick] (-.2,6) -- (3.8,6);
        \draw [thick] (0,5) to [out=0,in=150] (1,4.5);
        \draw [thick] (.2,4) to [out=0,in=-150] (1,4.5);
        \draw [double] (1,4.5) -- (1.5,4.5);
        \draw [thick] (1.5,4.5) to [out=30,in=180] (2,5) to [out=0,in=150] (2.5,4.5);
        \draw [thick] (1.5,4.5) to [out=-30,in=180] (2,4) to [out=0,in=-150] (2.5,4.5);
        \draw [double] (2.5,4.5) -- (3,4.5);
        \draw [thick] (3,4.5) to [out=30,in=180] (4,5);
        \draw [thick] (3,4.5) to [out=-30,in=180] (4.2,4);
        \fill [yellow, opacity=.6]  (2.5,4.5) to [out=-90,in=0] (2,3.5) to [out=180,in=-90] (1.5,4.5) -- (1,4.5) to [out=-90,in=90] (1,.5) -- (1.5,.5) to [out=90,in=-90] (3,4.5) -- (2.5,4.5);
        \fill [blue, opacity=.6] (.2,0) to [out=0,in=-150] (1,.5) -- (1,4.5) to [out=-150,in=0] (.2,4);
        \fill [blue, opacity=.6] (1.5,.5) to [out=90,in=-90] (3,4.5) to [out=-30,in=180] (4.2,4) -- (4.2,0) to [out=180,in=-30] (1.5,.5);
        \draw (.2,0) -- (.2,4);
        \draw (4.2,0) -- (4.2,4);
        \fill [blue, opacity=.6] (1.5,4.5) to [out=30,in=180] (2,5) to [out=0,in=150] (2.5,4.5) to [out=-90,in=0] (2,3.5) to [out=180,in=-90] (1.5,4.5);
        \fill [blue, opacity=.6] (1.5,4.5) to [out=-30,in=180] (2,4) to [out=0,in=-150] (2.5,4.5) to [out=-90,in=0] (2,3.5) to [out=180,in=-90] (1.5,4.5);
        \draw [thick, red] (2.5,4.5) to [out=-90,in=0] (2,3.5) to [out=180,in=-90] (1.5,4.5);
        \fill [blue, opacity=.6] (0,1) to [out=0,in=150] (1,.5) -- (1,4.5) to [out=150,in=0] (0,5) -- (0,1);
        \draw (0,1) -- (0,5);
        \fill [blue,opacity=.6] (1.5,.5) -- (2,1.5) to [out=90,in=180] (2.5,2.25) to [out=0,in=90] (3,1.5) to [out=-30,in=180] (4,1) -- (4,5) to [out=180,in=30] (3,4.5) to [out=-90,in=90] (1.5,.5);
        \draw (4,1) -- (4,5);
        \draw [thick, red] (1,.5) to [out=90,in=-90] (1,4.5);
        \draw [thick, red] (1.5,.5) to [out=90,in=-90] (3,4.5);
        \fill [blue, opacity=.6] (-.2,2) to [out=0,in=150] (2,1.5) to [out=90,in=180] (2.5,2.25) to [out=0,in=90] (3,1.5) to [out=30,in=180] (3.8,2) -- (3.8,6) -- (-.2,6) -- (-.2,2);
        \fill [yellow, opacity=.6] (2,1.5) to [out=90,in=180] (2.5,2.25) to [out=0,in=90] (3,1.5);
        \draw (-.2,2) -- (-.2,6);
        \draw (3.8,2) -- (3.8,6);
       \draw [thick, red] (2,1.5) to [out=90,in=180] (2.5,2.25) to [out=0,in=90] (3,1.5);
        \draw [thick] (.2,0) to [out=0,in=-150] (1,.5);
        \draw [thick] (0,1) to [out=0,in=150] (1,.5);
        \draw [double] (1,.5) -- (1.5,.5);
        \draw [thick] (1.5,.5) -- (2,1.5);
        \draw [thick] (-.2,2) to [out=0,in=150] (2,1.5);
        \draw [double] (2,1.5) -- (3,1.5);
        \draw [thick] (1.5,.5) to [out=-30,in=180] (4.2,0);
        \draw [thick] (3,1.5) to [out=-30,in=180] (4,1);
        \draw [thick] (3,1.5) to [out=30,in=180] (3.8,2);
  \end{tikzpicture}
  };
  \end{tikzpicture}
\label{R31}
\end{equation}

\begin{equation}
  \begin{tikzpicture}[anchorbase]
    \node at (-4,0) {
            \begin{tikzpicture}[anchorbase,scale=.5]
        \draw [line width=4, white] (1,0) to [out=90,in=-90] (2,1) to [out=90,in=-90] (1,2);
        \draw [thick, ->] (1,0) to [out=90,in=-90] (2,1) to [out=90,in=-90] (1,2); 
        \draw [line width=4, white] (2,0) -- (0,2);
        \draw [thick, ->] (2,0) -- (0,2);
        \draw [line width=4, white] (0,0) -- (2,2);
        \draw [double, ->] (0,0) -- (2,2);
        \node at (2,.5) {\tiny $c_1$};
        \node at (.5,1) {\tiny $c_2$};
        \node at (2,1.5) {\tiny $c_3$};
      \end{tikzpicture}
    };
    \node at (-3,0) {$\rightarrow$};
        \node (A) at (-2,0) {
          \begin{tikzpicture}[anchorbase,scale=.5]
            \draw [double] (0,0) -- (0,.7);
            \draw [thick] (1,0) to [out=90,in=-120] (2,.3);
            \draw [thick] (2,0) to [out=90,in=-90] (2,.3);
            \draw [double] (2,.3) -- (2,.7);
            \draw [thick] (2,.7) -- (1,1);
            \draw [thick] (0,.7) -- (1,1);
            \draw [double] (1,1) -- (1,1.4);
            \draw [thick] (0,.7) -- (0,2);
            \draw [thick] (1,1.4) -- (1,2);
            \draw [thick] (1,1.4) -- (2,1.6);
            \draw [thick] (2,.7) to [out=90,in=-90] (2,1.6);
            \draw [double] (2,1.6) -- (2,2);
      \end{tikzpicture}
    };
        \node (B) at (1,0) {
      \begin{tikzpicture}[anchorbase,scale=.5]
        \draw [double] (0,0) -- (0,.7);
        \draw [thick] (1,0) -- (1,1);
        \draw [thick] (2,0) -- (2,1.6);
        \draw [thick] (0,.7) -- (1,1);
        \draw [thick] (0,.7) -- (0,2);
        \draw [double] (1,1) -- (1,1.4);
        \draw [thick] (1,1.4) -- (2,1.6);
        \draw [thick] (1,1.4) -- (1,2);
        \draw [double] (2,1.6) -- (2,2);
      \end{tikzpicture}
    };
     \draw (A) -- (B) node[midway,above] {\tiny $\langle \cdot,c_1\rangle$};
    \node at (-4,-6) {
      \begin{tikzpicture}[anchorbase,scale=.5]
        \draw [line width=4, white] (1,0) to [out=90,in=-90] (0,1) to [out=90,in=-90] (1,2);
        \draw [thick, ->] (1,0) to [out=90,in=-90] (0,1) to [out=90,in=-90] (1,2); 
        \draw [line width=4, white] (2,0) -- (0,2);
        \draw [thick, ->] (2,0) -- (0,2);
        \draw [line width=4, white] (0,0) -- (2,2);
        \draw [double, ->] (0,0) -- (2,2);
        \node at (0,.5) {\tiny $c'_1$};
        \node at (1.5,1) {\tiny $c'_2$};
        \node at (0,1.5) {\tiny $c'_3$};
      \end{tikzpicture}
    };
    \node at (-3,-6) {$\rightarrow$};
        \node (A') at (-2,-6) {
      \begin{tikzpicture}[anchorbase,scale=.5]
        \draw [double] (0,0) -- (0,.4);
        \draw [thick] (0,.4) -- (1,.6);
        \draw [thick] (1,0) -- (1,.6);
        \draw [double] (1,.6) -- (1,1);
        \draw [thick] (1,1) -- (2,1.2);
        \draw [thick] (2,0) -- (2,1.2);
        \draw [double] (2,1.2) -- (2,2);
        \draw [thick] (1,1) -- (0,1.3);
        \draw [thick] (0,.4) -- (0,1.3);
        \draw [double] (0,1.3) -- (0,1.7);
        \draw [thick] (0,1.7) -- (0,2);
        \draw [thick] (0,1.7) to [out=60,in=-90] (1,2);
      \end{tikzpicture}
    };
        \node (B') at (1,-6) {
      \begin{tikzpicture}[anchorbase,scale=.5]
        \draw [double] (0,0) -- (0,.7);
        \draw [thick] (1,0) -- (1,1);
        \draw [thick] (2,0) -- (2,1.6);
        \draw [thick] (0,.7) -- (1,1);
        \draw [thick] (0,.7) -- (0,2);
        \draw [double] (1,1) -- (1,1.4);
        \draw [thick] (1,1.4) -- (2,1.6);
        \draw [thick] (1,1.4) -- (1,2);
        \draw [double] (2,1.6) -- (2,2);
      \end{tikzpicture}
    };
     \draw (A') -- (B') node[midway,above] {\tiny $\langle \cdot,c'_3\rangle$};
    \draw [green,->] (A) edge node [above,rotate=90] {\tiny $F\otimes \langle \cdot,c_1\wedge c_2\wedge c_3\rangle \wedge c'_1\wedge c'_2\wedge c'_3$} (A');
    \draw [green,->] (B) edge node [above,rotate=90] {\tiny $-\langle \cdot,c_2\wedge c_3\rangle c'_1\wedge c'_2$} (B');
        \node at (4,-2) {$F\leftrightarrow$  
      \begin{tikzpicture}[anchorbase,scale=.3]
        \draw [double] (-.2,2) -- (1,2);
        \draw [thick] (0,1) to [out=0,in=150] (.6,0);
        \draw [thick] (.2,0) -- (.6,0);
        \draw [double] (.6,0) -- (1,0);
        \draw [thick] (1,2) to [out=-30,in=150] (1.5,1);
        \draw [thick] (1,0) to [out=30,in=-150] (1.5,1);
        \draw [double] (1.5,1) -- (2.2,1);
        \draw [thick] (2.2,1) to [out=-30,in=150] (3,0);
        \draw [thick] (1,0) -- (3,0);
        \draw [double] (3,0) -- (4.2,0);
        \draw [thick] (1,2) -- (3.8,2);
        \draw [thick] (2.2,1) -- (4,1);
        \fill [blue,opacity=.6] (.8,6) -- (2.4,6) to [out=-90,in=0] (1.6,4.3) to [out=180,in=-90] (.8,6);
        \fill [yellow, opacity=.6] (-.2,6) -- (.8,6) to [out=-90,in=180] (1.6,4.3) to [out=0,in=-90] (2.4,6) -- (3,6) to [out=-90,in=90] (1,2) -- (-.2,2) -- (-.2,6);
        \fill [blue, opacity=.6] (3,6) -- (3.8,6) -- (3.8,2) -- (1,2) to [out=90,in=-90] (3,6);
        \draw (-.2,2) -- (-.2,6);
        \draw (3.8,2) -- (3.8,6);
        \fill [blue, opacity=.6] (2.4,6) to [out=-90,in=0] (1.6,4.3) to [out=180,in=-90] (.8,6) -- (1.2,5) to [out=-90,in=180] (1.6,4.6) to [out=0,in=-90] (2,5) -- (2.4,6);
        \fill [blue, opacity=.6] (4,5) to [out=180,in=-30] (3,6) to [out=-90,in=90] (1,2) to [out=-30,in=150] (1.5,1) to [out=90,in=180] (1.85,1.5) to [out=0,in=90] (2.2,1) -- (4,1) -- (4,5);
        \draw [red,thick]  (2.4,6) to [out=-90,in=0] (1.6,4.3) to [out=180,in=-90] (.8,6);
        \draw [red, thick] (3,6) to [out=-90,in=90] (1,2);
        \fill [yellow, opacity=.6] (1.2,5) to [out=-90,in=180] (1.6,4.6) to [out=0,in=-90] (2,5);
        \fill [yellow, opacity=.6] (1.5,1) -- (2.2,1) to [out=90,in=0] (1.85,1.5) to [out=180,in=90] (1.5,1);
        \draw (0,1) -- (0,5);
        \draw (4,1) -- (4,5);
        \fill [blue, opacity=.6] (0,5) -- (1.2,5) to [out=-90,in=180] (1.6,4.6) to [out=0, in=-90] (2,5) to [out=-30,in=150] (2.4,4) to [out=-90,in=90] (3,0) to [out=150,in=-30] (2.2,1) to [out=90,in=0] (1.85,1.5) to [out=180,in=90] (1.5,1) to [out=-150,in=30] (1,0) to [out=90,in=0] (.8,.3) to [out=180,in=90] (.6,0) to [out=150,in=0] (0,1) -- (0,5);
        \draw [red, thick] (1.2,5) to [out=-90,in=180] (1.6,4.6) to [out=0, in=-90] (2,5);
        \draw [red, thick] (2.2,1) to [out=90,in=0] (1.85,1.5) to [out=180,in=90] (1.5,1);
        \fill [blue, opacity=.6] (.2,0) -- (.6,0) to [out=90,in=180] (.8,.3) to [out=0,in=90] (1,0) -- (3,0) to [out=90,in=-90] (2.4,4) -- (.2,4) -- (.2,0);
        \fill [yellow, opacity=.6] (1,0) to [out=90,in=0] (.8,.3) to [out=180,in=90] (.6,0);
        \fill [yellow, opacity=.6] (3,0) to [out=90,in=-90] (2.4,4) -- (4.2,4) -- (4.2,0) -- (3,0); 
        \draw [red, thick] (1,0) to [out=90,in=0] (.8,.3) to [out=180,in=90] (.6,0);
        \draw [red, thick] (2.4,4) to [out=-90,in=90] (3,0);
        \draw (.2,0) -- (.2,4);
        \draw (4.2,0) -- (4.2,4);
        \draw [double] (-.2,6) -- (.8,6);
        \draw [thick] (0,5) -- (1.2,5);
        \draw [thick] (.8,6) -- (1.2,5);
        \draw [double] (1.2,5) -- (2,5);
        \draw [thick] (2,5) -- (2.4,6);
        \draw [double] (2.4,6) -- (3,6);
        \draw [thick] (.8,6) -- (2.4,6);
        \draw [thick] (3,6) to [out=-30,in=180] (4,5);
        \draw [thick] (2,5) to [out=-30,in=150] (2.4,4);
        \draw [thick] (.2,4) -- (2.4,4);
        \draw [double] (2.4,4) -- (4.2,4);
        \draw [thick] (3,6) -- (3.8,6);
      \end{tikzpicture}
  };
  \end{tikzpicture}
\label{R32}
\end{equation}

\begin{equation}
  \begin{tikzpicture}[anchorbase]
    \node at (-4,0) {
            \begin{tikzpicture}[anchorbase,scale=.5]
        \draw [ultra thick, white] (2,0) -- (0,2);
        \draw [thick, ->] (2,0) -- (0,2);
        \draw [ultra thick, white, double=white] (1,0) to [out=90,in=-90] (2,1) to [out=90,in=-90] (1,2);
        \draw [double, ->] (1,0) to [out=90,in=-90] (2,1) to [out=90,in=-90] (1,2); 
        \draw [ultra thick, white, double=white] (0,0) -- (2,2);
        \draw [double, ->] (0,0) -- (2,2);
        \node at (2.1,.5) {\tiny $c_1$};
        \node at (.5,1) {\tiny $c_2$};
      \end{tikzpicture}
    };
    \node at (-3,0) {$\rightarrow$};
        \node (A) at (0,0) {
      \begin{tikzpicture}[anchorbase,scale=.5]
        \draw [double] (0,0) -- (0,.7);
        \draw [double] (1,0) -- (1,.7);
        \draw [thick] (2,0) -- (2,1.3);
        \draw [thick] (0,.7) -- (1,1.3);
        \draw [thick] (1,.7) -- (1,1.3);
        \draw [thick] (1,.7) -- (2,1.3);
        \draw [thick] (0,.7) -- (0,2);
        \draw [double] (1,1.3) -- (1,2);
        \draw [double] (2,1.3) -- (2,2);
      \end{tikzpicture}
    };
    \node at (-4,-5) {
      \begin{tikzpicture}[anchorbase,scale=.5]
        \draw [ultra thick, white] (2,0) -- (0,2);
        \draw [thick, ->] (2,0) -- (0,2);
        \draw [ultra thick, white, double=white] (1,0) to [out=90,in=-90] (0,1) to [out=90,in=-90] (1,2);
        \draw [double, ->] (1,0) to [out=90,in=-90] (0,1) to [out=90,in=-90] (1,2); 
        \draw [ultra thick, white, double=white] (0,0) -- (2,2);
        \draw [double, ->] (0,0) -- (2,2);
        \node at (1.5,1) {\tiny $c'_1$};
        \node at (-.2,1.5) {\tiny $c'_2$};
      \end{tikzpicture}
    };
    \node at (-3,-5) {$\rightarrow$};
    \node (B) at (0,-5) {
      \begin{tikzpicture}[anchorbase,scale=.5]
        \draw [double] (0,0) -- (0,.7);
        \draw [double] (1,0) -- (1,.7);
        \draw [thick] (2,0) -- (2,1.3);
        \draw [thick] (0,.7) -- (1,1.3);
        \draw [thick] (1,.7) -- (1,1.3);
        \draw [thick] (1,.7) -- (2,1.3);
        \draw [thick] (0,.7) -- (0,2);
        \draw [double] (1,1.3) -- (1,2);
        \draw [double] (2,1.3) -- (2,2);
      \end{tikzpicture}
};
    \draw [green,->] (A) edge[bend right=30] node [above,rotate=90] {\tiny $\langle \cdot, c_1\wedge c_2\rangle \wedge c'_1\wedge c'_2$} (B);
    \draw [red,->] (B) edge[bend right=30] node [below,rotate=90] {\tiny $\langle \cdot, c'_1\wedge c'_2\rangle \wedge c_1\wedge c_2$} (A);
  \end{tikzpicture}
\label{R33}
\end{equation}

\begin{equation}
  \begin{tikzpicture}[anchorbase]
    \node at (-4,0) {
            \begin{tikzpicture}[anchorbase,scale=.5]
        \draw [ultra thick, white] (2,0) -- (0,2);
        \draw [thick, ->] (2,0) -- (0,2);
        \draw [ultra thick, white, double=white] (1,0) to [out=90,in=-90] (2,1) to [out=90,in=-90] (1,2);
        \draw [double, ->] (1,0) to [out=90,in=-90] (2,1) to [out=90,in=-90] (1,2); 
        \draw [ultra thick, white] (0,0) -- (2,2);
        \draw [thick, ->] (0,0) -- (2,2);
        \node at (2,.5) {\tiny $c_1$};
        \node at (.5,1) {\tiny $c_2$};
        \node at (2,1.5) {\tiny $c_3$};
      \end{tikzpicture}
    };
    \node at (-3,0) {$\rightarrow$};
        \node (A) at (0,0) {
      \begin{tikzpicture}[anchorbase,scale=.5]
        \draw [thick] (0,0) -- (0,2);
        \draw [double] (1,0) -- (1,.4);
        \draw [thick] (2,0) -- (2,.6);
        \draw [thick] (1,.4) -- (2,.6);
        \draw [double] (2,.6) -- (2,1.4);
        \draw [thick] (1,.4) -- (1,1.6);
        \draw [thick] (2,1.4) -- (2,2);
        \draw [thick] (2,1.4) -- (1,1.6);
        \draw [double] (1,1.6) -- (1,2);
      \end{tikzpicture}
    };
        \node (B) at (3,0) {
      \begin{tikzpicture}[anchorbase,scale=.5]
        \draw [thick] (0,0) -- (0,.6);
        \draw [double] (1,0) -- (1,.4);
        \draw [thick] (2,0) -- (2,.6);
        \draw [thick] (1,.4) -- (2,.6);
        \draw [double] (2,.6) -- (2,1.4);
        \draw [thick] (1,.4) -- (0,.6);
        \draw [double] (0,.6) -- (0,1.4);
        \draw [thick] (0,1.4) -- (1,1.6);
        \draw [thick] (2,1.4) -- (2,2);
        \draw [thick] (0,1.4) -- (0,2);
        \draw [thick] (2,1.4) -- (1,1.6);
        \draw [double] (1,1.6) -- (1,2);
      \end{tikzpicture}
    };
     \draw (A) -- (B) node[midway,above] {\tiny $\cdot\wedge c_2$};
    \node at (-4,-5) {
      \begin{tikzpicture}[anchorbase,scale=.5]
        \draw [ultra thick, white] (2,0) -- (0,2);
        \draw [thick, ->] (2,0) -- (0,2);
        \draw [ultra thick, white, double=white] (1,0) to [out=90,in=-90] (0,1) to [out=90,in=-90] (1,2);
        \draw [double, ->] (1,0) to [out=90,in=-90] (0,1) to [out=90,in=-90] (1,2); 
        \draw [ultra thick, white] (0,0) -- (2,2);
        \draw [thick, ->] (0,0) -- (2,2);
        \node at (-.2,.5) {\tiny $c'_1$};
        \node at (1.5,1) {\tiny $c'_2$};
        \node at (-.2,1.5) {\tiny $c'_3$};
      \end{tikzpicture}
    };
    \node at (-3,-5) {$\rightarrow$};
        \node (A') at (0,-5) {
      \begin{tikzpicture}[anchorbase,scale=.5]
        \draw [thick] (2,0) -- (2,2);
        \draw [double] (1,0) -- (1,.4);
        \draw [thick] (0,0) -- (0,.6);
        \draw [thick] (1,.4) -- (0,.6);
        \draw [double] (0,.6) -- (0,1.4);
        \draw [thick] (1,.4) -- (1,1.6);
        \draw [thick] (0,1.4) -- (0,2);
        \draw [thick] (0,1.4) -- (1,1.6);
        \draw [double] (1,1.6) -- (1,2);
      \end{tikzpicture}
    };
        \node (B') at (3,-5) {
      \begin{tikzpicture}[anchorbase,scale=.5]
        \draw [thick] (0,0) -- (0,.6);
        \draw [double] (1,0) -- (1,.4);
        \draw [thick] (2,0) -- (2,.6);
        \draw [thick] (1,.4) -- (2,.6);
        \draw [double] (2,.6) -- (2,1.4);
        \draw [thick] (1,.4) -- (0,.6);
        \draw [double] (0,.6) -- (0,1.4);
        \draw [thick] (0,1.4) -- (1,1.6);
        \draw [thick] (2,1.4) -- (2,2);
        \draw [thick] (0,1.4) -- (0,2);
        \draw [thick] (2,1.4) -- (1,1.6);
        \draw [double] (1,1.6) -- (1,2);
      \end{tikzpicture}
    };
     \draw (A') -- (B') node[midway,above] {\tiny $\cdot\wedge c'_2$};
    \draw [green,->] (A) edge [bend right=30] node [above,rotate=90] {\tiny $F \otimes \langle \cdot,c_1\wedge c_3\rangle \wedge c'_1\wedge c'_3$} (A');
    \draw [green,->] (B) edge [bend right=30] node [above,rotate=90] {\tiny $-\langle \cdot,c_1\wedge c_2\wedge c_3\rangle c'_1\wedge c'_2\wedge c'_3$} (B');
        \node at (6,-2) {$F\leftrightarrow$  
      \begin{tikzpicture}[anchorbase,scale=.3]
        \draw [double] (0,1) -- (1,1);
        \draw [thick] (.2,0) -- (1.5,0);
        \draw [thick] (-.2,2) -- (3.8,2);
        \draw [thick] (1,1) to [out=-30,in=150] (1.5,0);
        \draw [double] (1.5,0) -- (2.5,0);
        \draw [thick] (2.5,0) to [out=30,in=-150] (3,1);
        \draw [thick] (1,1) -- (3,1);
        \draw [double] (3,1) -- (4,1);
        \draw [thick] (2.5,0) -- (4.2,0);
        \fill [yellow,opacity=.6] (1.5,6) -- (2.5,6) to [out=-90,in=0] (2,5.2) to [out=180,in=-90] (1.5,6);
        \fill [blue, opacity=.6] (-.2,6) -- (1.5,6) to [out=-90,in=180] (2,5.2) to [out=0,in=-90] (2.5,6) -- (3.8,6) to [out=-90,in=90] (3.8,2) -- (-.2,2) -- (-.2,6);
        \draw (-.2,2) -- (-.2,6);
        \draw (3.8,2) -- (3.8,6);
        \fill [blue, opacity=.6] (2.5,6) to [out=-90,in=0] (2,5.2) to [out=180,in=-90] (1.5,6) to [out=-150,in=30] (1,5) to [out=-90,in=180] (2,3.8) to [out=0,in=-90] (3,5) to [out=150,in=-30] (2.5,6);
        \draw [red,thick]  (2.5,6) to [out=-90,in=0] (2,5.2) to [out=180,in=-90] (1.5,6);
        \fill [blue, opacity=.6] (1,5) to [out=-90,in=180] (2,3.8) to [out=0,in=-90] (3,5);
        \fill [blue, opacity=.6] (1,1) -- (3,1) to [out=90,in=0] (2,2.2) to [out=180,in=90] (1,1);
        \fill [yellow, opacity=.6] (0,5) -- (1,5) to [out=-90,in=180] (2,3.8) to [out=0,in=-90] (3,5) -- (4,5) -- (4,1) -- (3,1) to [out=90,in=0] (2,2.2) to [out=180,in=90] (1,1) -- (0,1) -- (0,5);
        \draw [thick, red] (1,5) to [out=-90,in=180] (2,3.8) to [out=0,in=-90] (3,5);
        \draw (0,1) -- (0,5);
        \draw (4,1) -- (4,5);
        \draw [red, thick]        (3,1) to [out=90,in=0] (2,2.2) to [out=180,in=90] (1,1);
        \fill [blue, opacity=.6] (1,1) to [out=90,in=180] (2,2.2) to [out=0,in=90] (3,1) to [out=-150,in=30] (2.5,0) to [out=90,in=0] (2,.8) to [out=180,in=90] (1.5,0) to [out=150,in=-30] (1,1);
        \fill [yellow, opacity=.6] (2.5,0) to [out=90,in=0] (2,.8) to [out=180,in=90] (1.5,0) -- (2.5,0);
        \fill [blue, opacity=.6] (.2,0) -- (1.5,0) to [out=90,in=180] (2,.8) to [out=0,in=90] (2.5,0) -- (4.2,0) -- (4.2,4) -- (.2,4) -- (.2,0);
        \draw [red, thick] (2.5,0) to [out=90,in=0] (2,.8) to [out=180,in=90] (1.5,0);
        \draw [double] (0,5) -- (1,5);
        \draw [thick] (-.2,6) -- (1.5,6);
        \draw [thick] (.2,4) -- (4.2,4);
        \draw [thick] (1,5) to [out=30,in=-150] (1.5,6);
        \draw [double] (1.5,6) -- (2.5,6);
        \draw [thick] (2.5,6) to [out=-30,in=150] (3,5);
        \draw [thick] (1,5) -- (3,5);
        \draw [double] (3,5) -- (4,5);
        \draw [thick] (2.5,6) -- (3.8,6);
      \end{tikzpicture}
  };
  \end{tikzpicture}
\label{R34}
\end{equation}

\begin{equation}
  \begin{tikzpicture}[anchorbase]
    \node at (-4,0) {
            \begin{tikzpicture}[anchorbase,scale=.5]
        \draw [ultra thick, white] (0,0) -- (2,2);
        \draw [thick, ->] (0,0) -- (2,2);
        \draw [ultra thick, white, double=white] (1,0) to [out=90,in=-90] (2,1) to [out=90,in=-90] (1,2);
        \draw [double, ->] (1,0) to [out=90,in=-90] (2,1) to [out=90,in=-90] (1,2); 
        \draw [ultra thick, white] (2,0) -- (0,2);
        \draw [thick, ->] (2,0) -- (0,2);
        \node at (2,.5) {\tiny $c_1$};
        \node at (.5,1) {\tiny $c_2$};
        \node at (2,1.5) {\tiny $c_3$};
      \end{tikzpicture}
    };
    \node at (-3,0) {$\rightarrow$};
        \node (A) at (0,0) {
      \begin{tikzpicture}[anchorbase,scale=.5]
        \draw [thick] (0,0) -- (0,.6);
        \draw [double] (1,0) -- (1,.4);
        \draw [thick] (2,0) -- (2,.6);
        \draw [thick] (1,.4) -- (2,.6);
        \draw [double] (2,.6) -- (2,1.4);
        \draw [thick] (1,.4) -- (0,.6);
        \draw [double] (0,.6) -- (0,1.4);
        \draw [thick] (0,1.4) -- (1,1.6);
        \draw [thick] (2,1.4) -- (2,2);
        \draw [thick] (0,1.4) -- (0,2);
        \draw [thick] (2,1.4) -- (1,1.6);
        \draw [double] (1,1.6) -- (1,2);
      \end{tikzpicture}
    };
        \node (B) at (3,0) {
      \begin{tikzpicture}[anchorbase,scale=.5]
        \draw [thick] (0,0) -- (0,2);
        \draw [double] (1,0) -- (1,.4);
        \draw [thick] (2,0) -- (2,.6);
        \draw [thick] (1,.4) -- (2,.6);
        \draw [double] (2,.6) -- (2,1.4);
        \draw [thick] (1,.4) -- (1,1.6);
        \draw [thick] (2,1.4) -- (2,2);
        \draw [thick] (2,1.4) -- (1,1.6);
        \draw [double] (1,1.6) -- (1,2);
      \end{tikzpicture}
    };
     \draw (A) -- (B) node[midway,above] {\tiny $\langle\cdot,c_2\rangle$};
    \node at (-4,-5) {
      \begin{tikzpicture}[anchorbase,scale=.5]
        \draw [ultra thick, white] (0,0) -- (2,2);
        \draw [thick, ->] (0,0) -- (2,2);
        \draw [ultra thick, white, double=white] (1,0) to [out=90,in=-90] (0,1) to [out=90,in=-90] (1,2);
        \draw [double, ->] (1,0) to [out=90,in=-90] (0,1) to [out=90,in=-90] (1,2); 
        \draw [ultra thick, white] (2,0) -- (0,2);
        \draw [thick, ->] (2,0) -- (0,2);
        \node at (-.2,.5) {\tiny $c'_1$};
        \node at (1.5,1) {\tiny $c'_2$};
        \node at (-.2,1.5) {\tiny $c'_3$};
      \end{tikzpicture}
    };
    \node at (-3,-5) {$\rightarrow$};
        \node (A') at (0,-5) {
      \begin{tikzpicture}[anchorbase,scale=.5]
        \draw [thick] (0,0) -- (0,.6);
        \draw [double] (1,0) -- (1,.4);
        \draw [thick] (2,0) -- (2,.6);
        \draw [thick] (1,.4) -- (2,.6);
        \draw [double] (2,.6) -- (2,1.4);
        \draw [thick] (1,.4) -- (0,.6);
        \draw [double] (0,.6) -- (0,1.4);
        \draw [thick] (0,1.4) -- (1,1.6);
        \draw [thick] (2,1.4) -- (2,2);
        \draw [thick] (0,1.4) -- (0,2);
        \draw [thick] (2,1.4) -- (1,1.6);
        \draw [double] (1,1.6) -- (1,2);
      \end{tikzpicture}
    };
        \node (B') at (3,-5) {
      \begin{tikzpicture}[anchorbase,scale=.5]
        \draw [thick] (2,0) -- (2,2);
        \draw [double] (1,0) -- (1,.4);
        \draw [thick] (0,0) -- (0,.6);
        \draw [thick] (1,.4) -- (0,.6);
        \draw [double] (0,.6) -- (0,1.4);
        \draw [thick] (1,.4) -- (1,1.6);
        \draw [thick] (0,1.4) -- (0,2);
        \draw [thick] (0,1.4) -- (1,1.6);
        \draw [double] (1,1.6) -- (1,2);
      \end{tikzpicture}
    };
     \draw (A') -- (B') node[midway,above] {\tiny $\langle \cdot,c'_2\rangle$};
    \draw [green,->] (A) edge [bend right=30] node [above,rotate=90] {\tiny $-\langle \cdot,c_1\wedge c_2\wedge c_3\rangle c'_1\wedge c'_2\wedge c'_3$} (A');
    \draw [green,->] (B) edge [bend right=30] node [above,rotate=90] {\tiny $F \otimes \langle \cdot,c_1\wedge c_3\rangle \wedge c'_1\wedge c'_3$} (B');
        \node at (6,-2) {$F\leftrightarrow$  

    };
     \draw (A') -- (B') node[midway,above] {\tiny $\cdot\wedge c'_1$};
    \draw [green,->] (A) edge node [above,rotate=90] {\tiny $-\langle \cdot,c_1\wedge c_2\rangle \wedge c'_2\wedge c'_3$} (A');
    \draw [green,->] (B) edge node [above,rotate=90] {\tiny $\text{cap/cup}\otimes \langle \cdot,c_1\wedge c_2\wedge c_3\rangle \wedge c'_1\wedge c'_2\wedge c'_3$} (B');
  \end{tikzpicture}
\label{R312}
\end{equation}

\begin{equation}
  \begin{tikzpicture}[anchorbase]
    \node at (-4,0) {
            \begin{tikzpicture}[anchorbase,scale=.5]
        \draw [line width=4, white] (2,0) -- (0,2);
        \draw [double, ->] (2,0) -- (0,2);
        \draw [line width=4, white] (0,0) -- (2,2);
        \draw [thick, ->] (0,0) -- (2,2);
        \draw [line width=4, white] (1,0) to [out=90,in=-90] (2,1) to [out=90,in=-90] (1,2);
        \draw [thick, ->] (1,0) to [out=90,in=-90] (2,1) to [out=90,in=-90] (1,2); 
        \node at (2,.5) {\tiny $c_1$};
        \node at (.5,1) {\tiny $c_2$};
        \node at (2,1.5) {\tiny $c_3$};
      \end{tikzpicture}
    };
    \node at (-3,0) {$\rightarrow$};
        \node (A) at (-2,0) {
      \begin{tikzpicture}[anchorbase,scale=.5,xscale=-1]
        \draw [double] (0,0) -- (0,.4);
        \draw [thick] (0,.4) -- (1,.6);
        \draw [thick] (1,0) -- (1,.6);
        \draw [double] (1,.6) -- (1,1);
        \draw [thick] (1,1) -- (2,1.2);
        \draw [thick] (2,0) -- (2,1.2);
        \draw [double] (2,1.2) -- (2,2);
        \draw [thick] (1,1) -- (0,1.3);
        \draw [thick] (0,.4) -- (0,1.3);
        \draw [double] (0,1.3) -- (0,1.7);
        \draw [thick] (0,1.7) -- (0,2);
        \draw [thick] (0,1.7) to [out=60,in=-90] (1,2);
      \end{tikzpicture}
    };
        \node (B) at (1,0) {
      \begin{tikzpicture}[anchorbase,scale=.5,xscale=-1]
        \draw [double] (0,0) -- (0,.7);
        \draw [thick] (1,0) -- (1,1);
        \draw [thick] (2,0) -- (2,1.6);
        \draw [thick] (0,.7) -- (1,1);
        \draw [thick] (0,.7) -- (0,2);
        \draw [double] (1,1) -- (1,1.4);
        \draw [thick] (1,1.4) -- (2,1.6);
        \draw [thick] (1,1.4) -- (1,2);
        \draw [double] (2,1.6) -- (2,2);
      \end{tikzpicture}
    };
     \draw (A) -- (B) node[midway,above] {\tiny $\langle \cdot, c_3\rangle$};
    \node at (-4,-6) {
      \begin{tikzpicture}[anchorbase,scale=.5]
        \draw [line width=4, white] (2,0) -- (0,2);
        \draw [double, ->] (2,0) -- (0,2);
        \draw [line width=4, white] (0,0) -- (2,2);
        \draw [thick, ->] (0,0) -- (2,2);
        \draw [line width=4, white] (1,0) to [out=90,in=-90] (0,1) to [out=90,in=-90] (1,2);
        \draw [thick, ->] (1,0) to [out=90,in=-90] (0,1) to [out=90,in=-90] (1,2); 
        \node at (0,.5) {\tiny $c'_1$};
        \node at (1.5,1) {\tiny $c'_2$};
        \node at (0,1.5) {\tiny $c'_3$};
      \end{tikzpicture}
    };
    \node at (-3,-6) {$\rightarrow$};
        \node (A') at (-2,-6) {
          \begin{tikzpicture}[anchorbase,scale=.5,xscale=-1]
            \draw [double] (0,0) -- (0,.7);
            \draw [thick] (1,0) to [out=90,in=-120] (2,.3);
            \draw [thick] (2,0) to [out=90,in=-90] (2,.3);
            \draw [double] (2,.3) -- (2,.7);
            \draw [thick] (2,.7) -- (1,1);
            \draw [thick] (0,.7) -- (1,1);
            \draw [double] (1,1) -- (1,1.4);
            \draw [thick] (0,.7) -- (0,2);
            \draw [thick] (1,1.4) -- (1,2);
            \draw [thick] (1,1.4) -- (2,1.6);
            \draw [thick] (2,.7) to [out=90,in=-90] (2,1.6);
            \draw [double] (2,1.6) -- (2,2);
      \end{tikzpicture}
    };
        \node (B') at (1,-6) {
      \begin{tikzpicture}[anchorbase,scale=.5,xscale=-1]
        \draw [double] (0,0) -- (0,.7);
        \draw [thick] (1,0) -- (1,1);
        \draw [thick] (2,0) -- (2,1.6);
        \draw [thick] (0,.7) -- (1,1);
        \draw [thick] (0,.7) -- (0,2);
        \draw [double] (1,1) -- (1,1.4);
        \draw [thick] (1,1.4) -- (2,1.6);
        \draw [thick] (1,1.4) -- (1,2);
        \draw [double] (2,1.6) -- (2,2);
      \end{tikzpicture}
    };
     \draw (A') -- (B') node[midway,above] {\tiny $\langle \cdot, c'_1\rangle$};
    \draw [green,->] (A) edge node [above,rotate=90] {\tiny $\text{cap/cup}\otimes \langle \cdot,c_1\wedge c_2\wedge c_3\rangle \wedge c'_1\wedge c'_2\wedge c'_3$} (A');
    \draw [green,->] (B) edge node [above,rotate=90] {\tiny $-\langle \cdot,c_1\wedge c_2\rangle \wedge c'_2\wedge c'_3$} (B');
  \end{tikzpicture}
\label{R313}
\end{equation}

\begin{equation}
  \begin{tikzpicture}[anchorbase]
    \node at (-4,0) {
            \begin{tikzpicture}[anchorbase,scale=.5]
        \draw [line width=4, white] (2,0) -- (0,2);
        \draw [double, ->] (2,0) -- (0,2);
        \draw [line width=4,white] (0,0) -- (2,2);
        \draw [double, ->] (0,0) -- (2,2);
        \draw [line width=4, white] (1,0) to [out=90,in=-90] (2,1) to [out=90,in=-90] (1,2);
        \draw [double, ->] (1,0) to [out=90,in=-90] (2,1) to [out=90,in=-90] (1,2); 
        \end{tikzpicture}
    };
    \node at (-3,0) {$\rightarrow$};
        \node (A) at (0,0) {
      \begin{tikzpicture}[anchorbase,scale=.5]
        \draw [double] (0,0) -- (0,2);
        \draw [double] (1,0) -- (1,2);
        \draw [double] (2,0) -- (2,2);
      \end{tikzpicture}
    };
    \node at (-4,-5) {
      \begin{tikzpicture}[anchorbase,scale=.5]
        \draw [line width=4, white] (2,0) -- (0,2);
        \draw [double, ->] (2,0) -- (0,2);
        \draw [line width=4,white] (0,0) -- (2,2);
        \draw [double, ->] (0,0) -- (2,2);
        \draw [line width=4, white] (1,0) to [out=90,in=-90] (0,1) to [out=90,in=-90] (1,2);
        \draw [double, ->] (1,0) to [out=90,in=-90] (0,1) to [out=90,in=-90] (1,2); 
        \end{tikzpicture}
    };
    \node at (-3,-5) {$\rightarrow$};
    \node (B) at (0,-5) {
  \begin{tikzpicture}[anchorbase,scale=.5]
        \draw [double] (0,0) -- (0,2);
        \draw [double] (1,0) -- (1,2);
        \draw [double] (2,0) -- (2,2);
      \end{tikzpicture}
  };
    \draw [green,->] (A) edge[bend right=30] node [above,rotate=90] {\tiny $\id$} (B);
    \draw [red,->] (B) edge[bend right=30] node [below,rotate=90] {\tiny $\id$} (A);
  \end{tikzpicture}
\label{R314}
\end{equation}

\begin{equation}
  \begin{tikzpicture}[anchorbase]
    \node at (-4,0) {
      \begin{tikzpicture}[anchorbase,scale=.5]
        \draw [line width=4, white] (2,0) -- (0,2);
        \draw [double, ->] (2,0) -- (0,2);
        \draw [line width=4, white] (1,0) to [out=90,in=-90] (0,1) to [out=90,in=-90] (1,2);
        \draw [double, ->] (1,0) to [out=90,in=-90] (0,1) to [out=90,in=-90] (1,2); 
        \draw [line width=4,white] (0,0) -- (2,2);
        \draw [thick, ->] (0,0) -- (2,2);
        \node at (-.2,.5) {\tiny $c_1$};
        \node at (1.7,1) {\tiny $c_2$};
      \end{tikzpicture}
    };
    \node at (-3,0) {$\rightarrow$};
        \node (A) at (0,0) {
      \begin{tikzpicture}[anchorbase,scale=.5,xscale=-1]
        \draw [double] (0,0) -- (0,.7);
        \draw [double] (1,0) -- (1,.7);
        \draw [thick] (2,0) -- (2,1.3);
        \draw [thick] (0,.7) -- (1,1.3);
        \draw [thick] (1,.7) -- (1,1.3);
        \draw [thick] (1,.7) -- (2,1.3);
        \draw [thick] (0,.7) -- (0,2);
        \draw [double] (1,1.3) -- (1,2);
        \draw [double] (2,1.3) -- (2,2);
      \end{tikzpicture}
    };
    \node at (-4,-5) {
            \begin{tikzpicture}[anchorbase,scale=.5]
        \draw [line width=4, white] (2,0) -- (0,2);
        \draw [double, ->] (2,0) -- (0,2);
        \draw [line width=4, white] (1,0) to [out=90,in=-90] (2,1) to [out=90,in=-90] (1,2);
        \draw [double, ->] (1,0) to [out=90,in=-90] (2,1) to [out=90,in=-90] (1,2); 
        \draw [line width=4,white] (0,0) -- (2,2);
        \draw [thick, ->] (0,0) -- (2,2);
        \node at (.4,1) {\tiny $c'_1$};
        \node at (2.1,1.5) {\tiny $c'_2$};
      \end{tikzpicture}
    };
    \node at (-3,-5) {$\rightarrow$};
    \node (B) at (0,-5) {
      \begin{tikzpicture}[anchorbase,scale=.5,xscale=-1]
        \draw [double] (0,0) -- (0,.7);
        \draw [double] (1,0) -- (1,.7);
        \draw [thick] (2,0) -- (2,1.3);
        \draw [thick] (0,.7) -- (1,1.3);
        \draw [thick] (1,.7) -- (1,1.3);
        \draw [thick] (1,.7) -- (2,1.3);
        \draw [thick] (0,.7) -- (0,2);
        \draw [double] (1,1.3) -- (1,2);
        \draw [double] (2,1.3) -- (2,2);
      \end{tikzpicture}
};
    \draw [green,->] (A) edge[bend right=30] node [above,rotate=90] {\tiny $\langle\cdot, c_1\wedge c_2\rangle \wedge c'_1\wedge c'_2$} (B);
    \draw [red,->] (B) edge[bend right=30] node [below,rotate=90] {\tiny $\langle \cdot, c'_1\wedge c'_2\rangle \wedge c_1\wedge c_2$} (A);
  \end{tikzpicture}
\label{R315}
\end{equation}

\begin{equation}
  \begin{tikzpicture}[anchorbase]
    \node at (-4,0) {
            \begin{tikzpicture}[anchorbase,scale=.5]
        \draw [line width=4,white] (0,0) -- (2,2);
        \draw [double, ->] (0,0) -- (2,2);
        \draw [line width=4, white] (2,0) -- (0,2);
        \draw [double, ->] (2,0) -- (0,2);
        \draw [line width=4, white] (1,0) to [out=90,in=-90] (2,1) to [out=90,in=-90] (1,2);
        \draw [thick, ->] (1,0) to [out=90,in=-90] (2,1) to [out=90,in=-90] (1,2); 
        \node at (2.1,.5) {\tiny $c_1$};
        \node at (2.1,1.5) {\tiny $c_2$};
      \end{tikzpicture}
    };
    \node at (-3,0) {$\rightarrow$};
        \node (A) at (0,0) {
      \begin{tikzpicture}[anchorbase,scale=.5]
        \draw [double] (0,0) -- (0,2);
        \draw [thick] (1,0) -- (1,.6);
        \draw [double] (2,0) -- (2,.4);
        \draw [thick] (1,.6) -- (2,.4);
        \draw [double] (1,.6) -- (1,1.4);
        \draw [thick] (2,.4) -- (2,1.6);
        \draw [thick] (1,1.4) -- (2,1.6);
        \draw [thick] (1,1.4) -- (1,2);
        \draw [double] (2,1.6) -- (2,2);
      \end{tikzpicture}
    };
    \node at (-4,-5) {
      \begin{tikzpicture}[anchorbase,scale=.5]
        \draw [line width=4,white] (0,0) -- (2,2);
        \draw [double, ->] (0,0) -- (2,2);
        \draw [line width=4, white] (2,0) -- (0,2);
        \draw [double, ->] (2,0) -- (0,2);
        \draw [line width=4, white] (1,0) to [out=90,in=-90] (0,1) to [out=90,in=-90] (1,2);
        \draw [thick, ->] (1,0) to [out=90,in=-90] (0,1) to [out=90,in=-90] (1,2); 
        \node at (-.2,.5) {\tiny $c'_1$};
        \node at (-.2,1.5) {\tiny $c'_2$};
      \end{tikzpicture}
    };
    \node at (-3,-5) {$\rightarrow$};
    \node (B) at (0,-5) {
      \begin{tikzpicture}[anchorbase,scale=.5]
        \draw [double] (2,0) -- (2,2);
        \draw [thick] (1,0) -- (1,.6);
        \draw [double] (0,0) -- (0,.4);
        \draw [thick] (1,.6) -- (0,.4);
        \draw [double] (1,.6) -- (1,1.4);
        \draw [thick] (0,.4) -- (0,1.6);
        \draw [thick] (1,1.4) -- (0,1.6);
        \draw [thick] (1,1.4) -- (1,2);
        \draw [double] (0,1.6) -- (0,2);
      \end{tikzpicture}
};
    \draw [green,->] (A) edge[bend right=30] node [above,rotate=90] {\tiny $\text{cap/cup}\otimes\langle \cdot, c_1\wedge c_2\rangle \wedge c'_1\wedge c'_2$} (B);
    \draw [red,->] (B) edge[bend right=30] node [below,rotate=90] {\tiny $\text{cap/cup}\otimes \langle \cdot, c'_1\wedge c'_2\rangle \wedge c_1\wedge c_2$} (A);
  \end{tikzpicture}
\label{R316}
\end{equation}

\begin{equation}
  \begin{tikzpicture}[anchorbase]
    \node at (-4,0) {
            \begin{tikzpicture}[anchorbase,scale=.5]
        \draw [line width=4, white] (2,0) -- (0,2);
        \draw [thick, ->] (2,0) -- (0,2);
        \draw [line width=4, white] (1,0) to [out=90,in=-90] (2,1) to [out=90,in=-90] (1,2);
        \draw [thick, ->] (1,0) to [out=90,in=-90] (2,1) to [out=90,in=-90] (1,2); 
        \draw [line width=4, white] (0,0) -- (2,2);
        \draw [double, ->] (0,0) -- (2,2);
        \node at (2,.5) {\tiny $c_1$};
        \node at (.5,1) {\tiny $c_2$};
        \node at (2,1.5) {\tiny $c_3$};
      \end{tikzpicture}
    };
    \node at (-3,0) {$\rightarrow$};
        \node (A) at (-2,0) {
      \begin{tikzpicture}[anchorbase,scale=.5]
        \draw [double] (0,0) -- (0,.7);
        \draw [thick] (1,0) -- (1,1);
        \draw [thick] (2,0) -- (2,1.6);
        \draw [thick] (0,.7) -- (1,1);
        \draw [thick] (0,.7) -- (0,2);
        \draw [double] (1,1) -- (1,1.4);
        \draw [thick] (1,1.4) -- (2,1.6);
        \draw [thick] (1,1.4) -- (1,2);
        \draw [double] (2,1.6) -- (2,2);
      \end{tikzpicture}
    };
        \node (B) at (1,0) {
          \begin{tikzpicture}[anchorbase,scale=.5]
            \draw [double] (0,0) -- (0,.7);
            \draw [thick] (1,0) to [out=90,in=-120] (2,.3);
            \draw [thick] (2,0) to [out=90,in=-90] (2,.3);
            \draw [double] (2,.3) -- (2,.7);
            \draw [thick] (2,.7) -- (1,1);
            \draw [thick] (0,.7) -- (1,1);
            \draw [double] (1,1) -- (1,1.4);
            \draw [thick] (0,.7) -- (0,2);
            \draw [thick] (1,1.4) -- (1,2);
            \draw [thick] (1,1.4) -- (2,1.6);
            \draw [thick] (2,.7) to [out=90,in=-90] (2,1.6);
            \draw [double] (2,1.6) -- (2,2);
      \end{tikzpicture}
    };
     \draw (A) -- (B) node[midway,above] {\tiny $\cdot \wedge c_1$};
    \node at (-4,-6) {
      \begin{tikzpicture}[anchorbase,scale=.5]
        \draw [line width=4, white] (2,0) -- (0,2);
        \draw [thick, ->] (2,0) -- (0,2);
        \draw [line width=4, white] (1,0) to [out=90,in=-90] (0,1) to [out=90,in=-90] (1,2);
        \draw [thick, ->] (1,0) to [out=90,in=-90] (0,1) to [out=90,in=-90] (1,2); 
        \draw [line width=4, white] (0,0) -- (2,2);
        \draw [double, ->] (0,0) -- (2,2);
        \node at (0,.5) {\tiny $c'_1$};
        \node at (1.5,1) {\tiny $c'_2$};
        \node at (0,1.5) {\tiny $c'_3$};
      \end{tikzpicture}
    };
    \node at (-3,-6) {$\rightarrow$};
        \node (A') at (-2,-6) {
      \begin{tikzpicture}[anchorbase,scale=.5]
        \draw [double] (0,0) -- (0,.7);
        \draw [thick] (1,0) -- (1,1);
        \draw [thick] (2,0) -- (2,1.6);
        \draw [thick] (0,.7) -- (1,1);
        \draw [thick] (0,.7) -- (0,2);
        \draw [double] (1,1) -- (1,1.4);
        \draw [thick] (1,1.4) -- (2,1.6);
        \draw [thick] (1,1.4) -- (1,2);
        \draw [double] (2,1.6) -- (2,2);
      \end{tikzpicture}
    };
        \node (B') at (1,-6) {
      \begin{tikzpicture}[anchorbase,scale=.5]
        \draw [double] (0,0) -- (0,.4);
        \draw [thick] (0,.4) -- (1,.6);
        \draw [thick] (1,0) -- (1,.6);
        \draw [double] (1,.6) -- (1,1);
        \draw [thick] (1,1) -- (2,1.2);
        \draw [thick] (2,0) -- (2,1.2);
        \draw [double] (2,1.2) -- (2,2);
        \draw [thick] (1,1) -- (0,1.3);
        \draw [thick] (0,.4) -- (0,1.3);
        \draw [double] (0,1.3) -- (0,1.7);
        \draw [thick] (0,1.7) -- (0,2);
        \draw [thick] (0,1.7) to [out=60,in=-90] (1,2);
      \end{tikzpicture}
    };
     \draw (A') -- (B') node[midway,above] {\tiny $\cdot\wedge c'_3$};
    \draw [green,->] (A) edge node [above,rotate=90] {\tiny $-\langle \cdot,c_2\wedge c_3\rangle c'_1\wedge c'_2$} (A');
    \draw [green,->] (B) edge node [above,rotate=90] {\tiny $\text{cap/cup}\otimes \langle \cdot,c_1\wedge c_2\wedge c_3\rangle \wedge c'_1\wedge c'_2\wedge c'_3$} (B');
  \end{tikzpicture}
\label{R317}
\end{equation}

\begin{equation}
  \begin{tikzpicture}[anchorbase]
    \node at (-4,0) {
      \begin{tikzpicture}[anchorbase,scale=.5]
        \draw [line width=4, white] (1,0) to [out=90,in=-90] (0,1) to [out=90,in=-90] (1,2);
        \draw [double, ->] (1,0) to [out=90,in=-90] (0,1) to [out=90,in=-90] (1,2); 
        \draw [line width=4, white] (0,0) -- (2,2);
        \draw [thick, ->] (0,0) -- (2,2);
        \draw [line width=4, white] (2,0) -- (0,2);
        \draw [thick, ->] (2,0) -- (0,2);
        \node at (-.2,.5) {\tiny $c_1$};
        \node at (1.5,1) {\tiny $c_2$};
        \node at (-.2,1.5) {\tiny $c_3$};
      \end{tikzpicture}
    };
    \node at (-3,0) {$\rightarrow$};
        \node (A) at (0,0) {
      \begin{tikzpicture}[anchorbase,scale=.5]
        \draw [thick] (0,0) -- (0,.6);
        \draw [double] (1,0) -- (1,.4);
        \draw [thick] (2,0) -- (2,.6);
        \draw [thick] (1,.4) -- (2,.6);
        \draw [double] (2,.6) -- (2,1.4);
        \draw [thick] (1,.4) -- (0,.6);
        \draw [double] (0,.6) -- (0,1.4);
        \draw [thick] (0,1.4) -- (1,1.6);
        \draw [thick] (2,1.4) -- (2,2);
        \draw [thick] (0,1.4) -- (0,2);
        \draw [thick] (2,1.4) -- (1,1.6);
        \draw [double] (1,1.6) -- (1,2);
      \end{tikzpicture}
    };
        \node (B) at (3,0) {
      \begin{tikzpicture}[anchorbase,scale=.5]
        \draw [thick] (0,0) -- (0,2);
        \draw [double] (1,0) -- (1,.4);
        \draw [thick] (2,0) -- (2,.6);
        \draw [thick] (1,.4) -- (2,.6);
        \draw [double] (2,.6) -- (2,1.4);
        \draw [thick] (1,.4) -- (1,1.6);
        \draw [thick] (2,1.4) -- (2,2);
        \draw [thick] (2,1.4) -- (1,1.6);
        \draw [double] (1,1.6) -- (1,2);
      \end{tikzpicture}
    };
     \draw (A) -- (B) node[midway,above] {\tiny $\langle\cdot,c_2\rangle$};
    \node at (-4,-5) {
            \begin{tikzpicture}[anchorbase,scale=.5]
        \draw [line width=4, white] (1,0) to [out=90,in=-90] (2,1) to [out=90,in=-90] (1,2);
        \draw [double, ->] (1,0) to [out=90,in=-90] (2,1) to [out=90,in=-90] (1,2); 
        \draw [line width=4, white] (0,0) -- (2,2);
        \draw [thick, ->] (0,0) -- (2,2);
        \draw [line width=4, white] (2,0) -- (0,2);
        \draw [thick, ->] (2,0) -- (0,2);
        \node at (2,.5) {\tiny $c'_1$};
        \node at (.5,1) {\tiny $c'_2$};
        \node at (2,1.5) {\tiny $c'_3$};
      \end{tikzpicture}
    };
    \node at (-3,-5) {$\rightarrow$};
        \node (A') at (0,-5) {
      \begin{tikzpicture}[anchorbase,scale=.5]
        \draw [thick] (0,0) -- (0,.6);
        \draw [double] (1,0) -- (1,.4);
        \draw [thick] (2,0) -- (2,.6);
        \draw [thick] (1,.4) -- (2,.6);
        \draw [double] (2,.6) -- (2,1.4);
        \draw [thick] (1,.4) -- (0,.6);
        \draw [double] (0,.6) -- (0,1.4);
        \draw [thick] (0,1.4) -- (1,1.6);
        \draw [thick] (2,1.4) -- (2,2);
        \draw [thick] (0,1.4) -- (0,2);
        \draw [thick] (2,1.4) -- (1,1.6);
        \draw [double] (1,1.6) -- (1,2);
      \end{tikzpicture}
    };
        \node (B') at (3,-5) {
      \begin{tikzpicture}[anchorbase,scale=.5]
        \draw [thick] (2,0) -- (2,2);
        \draw [double] (1,0) -- (1,.4);
        \draw [thick] (0,0) -- (0,.6);
        \draw [thick] (1,.4) -- (0,.6);
        \draw [double] (0,.6) -- (0,1.4);
        \draw [thick] (1,.4) -- (1,1.6);
        \draw [thick] (0,1.4) -- (0,2);
        \draw [thick] (0,1.4) -- (1,1.6);
        \draw [double] (1,1.6) -- (1,2);
      \end{tikzpicture}
    };
     \draw (A') -- (B') node[midway,above] {\tiny $\langle \cdot,c'_2\rangle$};
    \draw [green,->] (A) edge [bend right=30] node [above,rotate=90] {\tiny $\langle \cdot,c_1\wedge c_2\wedge c_3\rangle c'_1\wedge c'_2\wedge c'_3$} (A');
    \draw [green,->] (B) edge [bend right=30] node [above,rotate=90] {\tiny $-\text{cap/cup} \otimes \langle \cdot,c_1\wedge c_3\rangle \wedge c'_1\wedge c'_3$} (B');
  \end{tikzpicture}
\label{R318}
\end{equation}

\begin{equation}
  \begin{tikzpicture}[anchorbase]
    \node at (-4,0) {
      \begin{tikzpicture}[anchorbase,scale=.5]
        \draw [line width=4, white] (0,0) -- (2,2);
        \draw [thick, ->] (0,0) -- (2,2);
        \draw [line width=4, white] (2,0) -- (0,2);
        \draw [thick, ->] (2,0) -- (0,2);
        \draw [line width=4, white] (1,0) to [out=90,in=-90] (0,1) to [out=90,in=-90] (1,2);
        \draw [double, ->] (1,0) to [out=90,in=-90] (0,1) to [out=90,in=-90] (1,2); 
        \node at (-.2,.5) {\tiny $c_1$};
        \node at (1.5,1) {\tiny $c_2$};
        \node at (-.2,1.5) {\tiny $c_3$};
      \end{tikzpicture}
    };
    \node at (-3,0) {$\rightarrow$};
        \node (A) at (0,0) {
      \begin{tikzpicture}[anchorbase,scale=.5]
        \draw [thick] (0,0) -- (0,.6);
        \draw [double] (1,0) -- (1,.4);
        \draw [thick] (2,0) -- (2,.6);
        \draw [thick] (1,.4) -- (2,.6);
        \draw [double] (2,.6) -- (2,1.4);
        \draw [thick] (1,.4) -- (0,.6);
        \draw [double] (0,.6) -- (0,1.4);
        \draw [thick] (0,1.4) -- (1,1.6);
        \draw [thick] (2,1.4) -- (2,2);
        \draw [thick] (0,1.4) -- (0,2);
        \draw [thick] (2,1.4) -- (1,1.6);
        \draw [double] (1,1.6) -- (1,2);
      \end{tikzpicture}
    };
        \node (B) at (3,0) {
      \begin{tikzpicture}[anchorbase,scale=.5]
        \draw [thick] (0,0) -- (0,2);
        \draw [double] (1,0) -- (1,.4);
        \draw [thick] (2,0) -- (2,.6);
        \draw [thick] (1,.4) -- (2,.6);
        \draw [double] (2,.6) -- (2,1.4);
        \draw [thick] (1,.4) -- (1,1.6);
        \draw [thick] (2,1.4) -- (2,2);
        \draw [thick] (2,1.4) -- (1,1.6);
        \draw [double] (1,1.6) -- (1,2);
      \end{tikzpicture}
    };
     \draw (A) -- (B) node[midway,above] {\tiny $\langle\cdot,c_2\rangle$};
    \node at (-4,-5) {
            \begin{tikzpicture}[anchorbase,scale=.5]
        \draw [line width=4, white] (0,0) -- (2,2);
        \draw [thick, ->] (0,0) -- (2,2);
        \draw [line width=4, white] (2,0) -- (0,2);
        \draw [thick, ->] (2,0) -- (0,2);
        \draw [line width=4, white] (1,0) to [out=90,in=-90] (2,1) to [out=90,in=-90] (1,2);
        \draw [double, ->] (1,0) to [out=90,in=-90] (2,1) to [out=90,in=-90] (1,2); 
        \node at (2,.5) {\tiny $c'_1$};
        \node at (.5,1) {\tiny $c'_2$};
        \node at (2,1.5) {\tiny $c'_3$};
      \end{tikzpicture}
    };
    \node at (-3,-5) {$\rightarrow$};
        \node (A') at (0,-5) {
      \begin{tikzpicture}[anchorbase,scale=.5]
        \draw [thick] (0,0) -- (0,.6);
        \draw [double] (1,0) -- (1,.4);
        \draw [thick] (2,0) -- (2,.6);
        \draw [thick] (1,.4) -- (2,.6);
        \draw [double] (2,.6) -- (2,1.4);
        \draw [thick] (1,.4) -- (0,.6);
        \draw [double] (0,.6) -- (0,1.4);
        \draw [thick] (0,1.4) -- (1,1.6);
        \draw [thick] (2,1.4) -- (2,2);
        \draw [thick] (0,1.4) -- (0,2);
        \draw [thick] (2,1.4) -- (1,1.6);
        \draw [double] (1,1.6) -- (1,2);
      \end{tikzpicture}
    };
        \node (B') at (3,-5) {
      \begin{tikzpicture}[anchorbase,scale=.5]
        \draw [thick] (2,0) -- (2,2);
        \draw [double] (1,0) -- (1,.4);
        \draw [thick] (0,0) -- (0,.6);
        \draw [thick] (1,.4) -- (0,.6);
        \draw [double] (0,.6) -- (0,1.4);
        \draw [thick] (1,.4) -- (1,1.6);
        \draw [thick] (0,1.4) -- (0,2);
        \draw [thick] (0,1.4) -- (1,1.6);
        \draw [double] (1,1.6) -- (1,2);
      \end{tikzpicture}
    };
     \draw (A') -- (B') node[midway,above] {\tiny $\langle \cdot,c'_2\rangle$};
    \draw [green,->] (A) edge [bend right=30] node [above,rotate=90] {\tiny $\langle \cdot,c_1\wedge c_2\wedge c_3\rangle c'_1\wedge c'_2\wedge c'_3$} (A');
    \draw [green,->] (B) edge [bend right=30] node [above,rotate=90] {\tiny $-\text{cap/cup} \otimes \langle \cdot,c_1\wedge c_3\rangle \wedge c'_1\wedge c'_3$} (B');
  \end{tikzpicture}
\label{R319}
\end{equation}

\begin{equation}
  \begin{tikzpicture}[anchorbase]
    \node at (-4,0) {
            \begin{tikzpicture}[anchorbase,scale=.5]
        \draw [ultra thick, white, double=white] (0,0) -- (2,2);
        \draw [double, ->] (0,0) -- (2,2);
        \draw [line width=4, white] (1,0) to [out=90,in=-90] (2,1) to [out=90,in=-90] (1,2);
        \draw [double, ->] (1,0) to [out=90,in=-90] (2,1) to [out=90,in=-90] (1,2); 
        \draw [line width=4, white] (2,0) -- (0,2);
        \draw [thick, ->] (2,0) -- (0,2);
        \node at (2.1,.5) {\tiny $c_1$};
        \node at (.5,1) {\tiny $c_2$};
      \end{tikzpicture}
    };
    \node at (-3,0) {$\rightarrow$};
        \node (A) at (0,0) {
      \begin{tikzpicture}[anchorbase,scale=.5]
        \draw [double] (0,0) -- (0,.7);
        \draw [double] (1,0) -- (1,.7);
        \draw [thick] (2,0) -- (2,1.3);
        \draw [thick] (0,.7) -- (1,1.3);
        \draw [thick] (1,.7) -- (1,1.3);
        \draw [thick] (1,.7) -- (2,1.3);
        \draw [thick] (0,.7) -- (0,2);
        \draw [double] (1,1.3) -- (1,2);
        \draw [double] (2,1.3) -- (2,2);
      \end{tikzpicture}
    };
    \node at (-4,-5) {
      \begin{tikzpicture}[anchorbase,scale=.5]
        \draw [ultra thick, white, double=white] (0,0) -- (2,2);
        \draw [double, ->] (0,0) -- (2,2);
        \draw [line width=4, white] (1,0) to [out=90,in=-90] (0,1) to [out=90,in=-90] (1,2);
        \draw [double, ->] (1,0) to [out=90,in=-90] (0,1) to [out=90,in=-90] (1,2); 
        \draw [line width=4, white] (2,0) -- (0,2);
        \draw [thick, ->] (2,0) -- (0,2);
        \node at (1.5,1) {\tiny $c'_1$};
        \node at (-.2,1.5) {\tiny $c'_2$};
      \end{tikzpicture}
    };
    \node at (-3,-5) {$\rightarrow$};
    \node (B) at (0,-5) {
      \begin{tikzpicture}[anchorbase,scale=.5]
        \draw [double] (0,0) -- (0,.7);
        \draw [double] (1,0) -- (1,.7);
        \draw [thick] (2,0) -- (2,1.3);
        \draw [thick] (0,.7) -- (1,1.3);
        \draw [thick] (1,.7) -- (1,1.3);
        \draw [thick] (1,.7) -- (2,1.3);
        \draw [thick] (0,.7) -- (0,2);
        \draw [double] (1,1.3) -- (1,2);
        \draw [double] (2,1.3) -- (2,2);
      \end{tikzpicture}
};
    \draw [green,->] (A) edge[bend right=30] node [above,rotate=90] {\tiny $\langle \cdot, c_1\wedge c_2\rangle \wedge c'_1\wedge c'_2$} (B);
    \draw [red,->] (B) edge[bend right=30] node [below,rotate=90] {\tiny $\langle \cdot, c'_1\wedge c'_2\rangle \wedge c_1\wedge c_2$} (A);
  \end{tikzpicture}
\label{R320}
\end{equation}

\begin{equation}
  \begin{tikzpicture}[anchorbase]
    \node at (-4,0) {
            \begin{tikzpicture}[anchorbase,scale=.5]
        \draw [line width=4, white] (1,0) to [out=90,in=-90] (2,1) to [out=90,in=-90] (1,2);
        \draw [double, ->] (1,0) to [out=90,in=-90] (2,1) to [out=90,in=-90] (1,2); 
        \draw [ultra thick, white, double=white] (0,0) -- (2,2);
        \draw [double, ->] (0,0) -- (2,2);
        \draw [line width=4, white] (2,0) -- (0,2);
        \draw [thick, ->] (2,0) -- (0,2);
        \node at (2.1,.5) {\tiny $c_1$};
        \node at (.5,1) {\tiny $c_2$};
      \end{tikzpicture}
    };
    \node at (-3,0) {$\rightarrow$};
        \node (A) at (0,0) {
      \begin{tikzpicture}[anchorbase,scale=.5]
        \draw [double] (0,0) -- (0,.7);
        \draw [double] (1,0) -- (1,.7);
        \draw [thick] (2,0) -- (2,1.3);
        \draw [thick] (0,.7) -- (1,1.3);
        \draw [thick] (1,.7) -- (1,1.3);
        \draw [thick] (1,.7) -- (2,1.3);
        \draw [thick] (0,.7) -- (0,2);
        \draw [double] (1,1.3) -- (1,2);
        \draw [double] (2,1.3) -- (2,2);
      \end{tikzpicture}
    };
    \node at (-4,-5) {
      \begin{tikzpicture}[anchorbase,scale=.5]
        \draw [line width=4, white] (1,0) to [out=90,in=-90] (0,1) to [out=90,in=-90] (1,2);
        \draw [double, ->] (1,0) to [out=90,in=-90] (0,1) to [out=90,in=-90] (1,2); 
        \draw [ultra thick, white, double=white] (0,0) -- (2,2);
        \draw [double, ->] (0,0) -- (2,2);
        \draw [line width=4, white] (2,0) -- (0,2);
        \draw [thick, ->] (2,0) -- (0,2);
        \node at (1.5,1) {\tiny $c'_1$};
        \node at (-.2,1.5) {\tiny $c'_2$};
      \end{tikzpicture}
    };
    \node at (-3,-5) {$\rightarrow$};
    \node (B) at (0,-5) {
      \begin{tikzpicture}[anchorbase,scale=.5]
        \draw [double] (0,0) -- (0,.7);
        \draw [double] (1,0) -- (1,.7);
        \draw [thick] (2,0) -- (2,1.3);
        \draw [thick] (0,.7) -- (1,1.3);
        \draw [thick] (1,.7) -- (1,1.3);
        \draw [thick] (1,.7) -- (2,1.3);
        \draw [thick] (0,.7) -- (0,2);
        \draw [double] (1,1.3) -- (1,2);
        \draw [double] (2,1.3) -- (2,2);
      \end{tikzpicture}
};
    \draw [green,->] (A) edge[bend right=30] node [above,rotate=90] {\tiny $\langle \cdot, c_1\wedge c_2\rangle \wedge c'_1\wedge c'_2$} (B);
    \draw [red,->] (B) edge[bend right=30] node [below,rotate=90] {\tiny $\langle \cdot, c'_1\wedge c'_2\rangle \wedge c_1\wedge c_2$} (A);
  \end{tikzpicture}
\label{R321}
\end{equation}

\begin{equation}
  \begin{tikzpicture}[anchorbase]
    \node at (-4,0) {
      \begin{tikzpicture}[anchorbase,scale=.5]
        \draw [line width=4, white] (1,0) to [out=90,in=-90] (0,1) to [out=90,in=-90] (1,2);
        \draw [thick, ->] (1,0) to [out=90,in=-90] (0,1) to [out=90,in=-90] (1,2); 
        \draw [line width=4, white] (0,0) -- (2,2);
        \draw [double, ->] (0,0) -- (2,2);
        \draw [line width=4, white] (2,0) -- (0,2);
        \draw [thick, ->] (2,0) -- (0,2);
        \node at (0,.5) {\tiny $c_1$};
        \node at (1.5,1) {\tiny $c_2$};
        \node at (0,1.5) {\tiny $c_3$};
      \end{tikzpicture}
    };
    \node at (-3,0) {$\rightarrow$};
        \node (A) at (-2,0) {
      \begin{tikzpicture}[anchorbase,scale=.5]
        \draw [double] (0,0) -- (0,.4);
        \draw [thick] (0,.4) -- (1,.6);
        \draw [thick] (1,0) -- (1,.6);
        \draw [double] (1,.6) -- (1,1);
        \draw [thick] (1,1) -- (2,1.2);
        \draw [thick] (2,0) -- (2,1.2);
        \draw [double] (2,1.2) -- (2,2);
        \draw [thick] (1,1) -- (0,1.3);
        \draw [thick] (0,.4) -- (0,1.3);
        \draw [double] (0,1.3) -- (0,1.7);
        \draw [thick] (0,1.7) -- (0,2);
        \draw [thick] (0,1.7) to [out=60,in=-90] (1,2);
      \end{tikzpicture}
    };
        \node (B) at (1,0) {
      \begin{tikzpicture}[anchorbase,scale=.5]
        \draw [double] (0,0) -- (0,.7);
        \draw [thick] (1,0) -- (1,1);
        \draw [thick] (2,0) -- (2,1.6);
        \draw [thick] (0,.7) -- (1,1);
        \draw [thick] (0,.7) -- (0,2);
        \draw [double] (1,1) -- (1,1.4);
        \draw [thick] (1,1.4) -- (2,1.6);
        \draw [thick] (1,1.4) -- (1,2);
        \draw [double] (2,1.6) -- (2,2);
      \end{tikzpicture}
    };
     \draw (A) -- (B) node[midway,above] {\tiny $\langle \cdot,c_3\rangle$};
    \node at (-4,-6) {
            \begin{tikzpicture}[anchorbase,scale=.5]
        \draw [line width=4, white] (1,0) to [out=90,in=-90] (2,1) to [out=90,in=-90] (1,2);
        \draw [thick, ->] (1,0) to [out=90,in=-90] (2,1) to [out=90,in=-90] (1,2); 
        \draw [line width=4, white] (0,0) -- (2,2);
        \draw [double, ->] (0,0) -- (2,2);
        \draw [line width=4, white] (2,0) -- (0,2);
        \draw [thick, ->] (2,0) -- (0,2);
        \node at (2,.5) {\tiny $c'_1$};
        \node at (.5,1) {\tiny $c'_2$};
        \node at (2,1.5) {\tiny $c'_3$};
      \end{tikzpicture}
    };
    \node at (-3,-6) {$\rightarrow$};
        \node (A') at (-2,-6) {
          \begin{tikzpicture}[anchorbase,scale=.5]
            \draw [double] (0,0) -- (0,.7);
            \draw [thick] (1,0) to [out=90,in=-120] (2,.3);
            \draw [thick] (2,0) to [out=90,in=-90] (2,.3);
            \draw [double] (2,.3) -- (2,.7);
            \draw [thick] (2,.7) -- (1,1);
            \draw [thick] (0,.7) -- (1,1);
            \draw [double] (1,1) -- (1,1.4);
            \draw [thick] (0,.7) -- (0,2);
            \draw [thick] (1,1.4) -- (1,2);
            \draw [thick] (1,1.4) -- (2,1.6);
            \draw [thick] (2,.7) to [out=90,in=-90] (2,1.6);
            \draw [double] (2,1.6) -- (2,2);
      \end{tikzpicture}
    };
        \node (B') at (1,-6) {
      \begin{tikzpicture}[anchorbase,scale=.5]
        \draw [double] (0,0) -- (0,.7);
        \draw [thick] (1,0) -- (1,1);
        \draw [thick] (2,0) -- (2,1.6);
        \draw [thick] (0,.7) -- (1,1);
        \draw [thick] (0,.7) -- (0,2);
        \draw [double] (1,1) -- (1,1.4);
        \draw [thick] (1,1.4) -- (2,1.6);
        \draw [thick] (1,1.4) -- (1,2);
        \draw [double] (2,1.6) -- (2,2);
      \end{tikzpicture}
    };
     \draw (A') -- (B') node[midway,above] {\tiny $\langle \cdot,c'_1\rangle$};
    \draw [green,->] (A) edge node [above,rotate=90] {\tiny $-\text{cap/cup}\otimes \langle \cdot,c_1\wedge c_2\wedge c_3\rangle \wedge c'_1\wedge c'_2\wedge c'_3$} (A');
    \draw [green,->] (B) edge node [above,rotate=90] {\tiny $\langle \cdot,c_1\wedge c_2\rangle c'_2\wedge c'_3$} (B');
  \end{tikzpicture}
\label{R322}
\end{equation}

\begin{equation}
  \begin{tikzpicture}[anchorbase]
    \node at (-4,0) {
            \begin{tikzpicture}[anchorbase,scale=.5]
        \draw [line width=4, white] (2,0) -- (0,2);
        \draw [thick, ->] (2,0) -- (0,2);
        \draw [line width=4, white] (0,0) -- (2,2);
        \draw [double, ->] (0,0) -- (2,2);
        \draw [line width=4, white] (1,0) to [out=90,in=-90] (2,1) to [out=90,in=-90] (1,2);
        \draw [double, ->] (1,0) to [out=90,in=-90] (2,1) to [out=90,in=-90] (1,2); 
        \node at (2.1,.5) {\tiny $c_1$};
        \node at (.5,1) {\tiny $c_2$};
      \end{tikzpicture}
    };
    \node at (-3,0) {$\rightarrow$};
        \node (A) at (0,0) {
      \begin{tikzpicture}[anchorbase,scale=.5]
        \draw [double] (0,0) -- (0,.7);
        \draw [double] (1,0) -- (1,.7);
        \draw [thick] (2,0) -- (2,1.3);
        \draw [thick] (0,.7) -- (1,1.3);
        \draw [thick] (1,.7) -- (1,1.3);
        \draw [thick] (1,.7) -- (2,1.3);
        \draw [thick] (0,.7) -- (0,2);
        \draw [double] (1,1.3) -- (1,2);
        \draw [double] (2,1.3) -- (2,2);
      \end{tikzpicture}
    };
    \node at (-4,-5) {
      \begin{tikzpicture}[anchorbase,scale=.5]
        \draw [line width=4, white] (2,0) -- (0,2);
        \draw [thick, ->] (2,0) -- (0,2);
        \draw [line width=4, white] (0,0) -- (2,2);
        \draw [double, ->] (0,0) -- (2,2);
        \draw [line width=4, white] (1,0) to [out=90,in=-90] (0,1) to [out=90,in=-90] (1,2);
        \draw [double, ->] (1,0) to [out=90,in=-90] (0,1) to [out=90,in=-90] (1,2); 
        \node at (1.5,1) {\tiny $c'_1$};
        \node at (-.2,1.5) {\tiny $c'_2$};
      \end{tikzpicture}
    };
    \node at (-3,-5) {$\rightarrow$};
    \node (B) at (0,-5) {
      \begin{tikzpicture}[anchorbase,scale=.5]
        \draw [double] (0,0) -- (0,.7);
        \draw [double] (1,0) -- (1,.7);
        \draw [thick] (2,0) -- (2,1.3);
        \draw [thick] (0,.7) -- (1,1.3);
        \draw [thick] (1,.7) -- (1,1.3);
        \draw [thick] (1,.7) -- (2,1.3);
        \draw [thick] (0,.7) -- (0,2);
        \draw [double] (1,1.3) -- (1,2);
        \draw [double] (2,1.3) -- (2,2);
      \end{tikzpicture}
};
    \draw [green,->] (A) edge[bend right=30] node [above,rotate=90] {\tiny $\langle \cdot, c_1\wedge c_2\rangle \wedge c'_1\wedge c'_2$} (B);
    \draw [red,->] (B) edge[bend right=30] node [below,rotate=90] {\tiny $\langle \cdot, c'_1\wedge c'_2\rangle \wedge c_1\wedge c_2$} (A);
  \end{tikzpicture}
\label{R323}
\end{equation}

\begin{equation}
  \begin{tikzpicture}[anchorbase]
    \node at (-4,0) {
      \begin{tikzpicture}[anchorbase,scale=.5]
        \draw [line width=4,white] (0,0) -- (2,2);
        \draw [thick, ->] (0,0) -- (2,2);
        \draw [line width=4, white] (1,0) to [out=90,in=-90] (0,1) to [out=90,in=-90] (1,2);
        \draw [double, ->] (1,0) to [out=90,in=-90] (0,1) to [out=90,in=-90] (1,2); 
        \draw [line width=4, white] (2,0) -- (0,2);
        \draw [double, ->] (2,0) -- (0,2);
        \node at (-.2,.5) {\tiny $c_1$};
        \node at (1.7,1) {\tiny $c_2$};
      \end{tikzpicture}
    };
    \node at (-3,0) {$\rightarrow$};
        \node (A) at (0,0) {
      \begin{tikzpicture}[anchorbase,scale=.5,xscale=-1]
        \draw [double] (0,0) -- (0,.7);
        \draw [double] (1,0) -- (1,.7);
        \draw [thick] (2,0) -- (2,1.3);
        \draw [thick] (0,.7) -- (1,1.3);
        \draw [thick] (1,.7) -- (1,1.3);
        \draw [thick] (1,.7) -- (2,1.3);
        \draw [thick] (0,.7) -- (0,2);
        \draw [double] (1,1.3) -- (1,2);
        \draw [double] (2,1.3) -- (2,2);
      \end{tikzpicture}
    };
    \node at (-4,-5) {
            \begin{tikzpicture}[anchorbase,scale=.5]
        \draw [line width=4,white] (0,0) -- (2,2);
        \draw [thick, ->] (0,0) -- (2,2);
        \draw [line width=4, white] (1,0) to [out=90,in=-90] (2,1) to [out=90,in=-90] (1,2);
        \draw [double, ->] (1,0) to [out=90,in=-90] (2,1) to [out=90,in=-90] (1,2); 
        \draw [line width=4, white] (2,0) -- (0,2);
        \draw [double, ->] (2,0) -- (0,2);
        \node at (.4,1) {\tiny $c'_1$};
        \node at (2.1,1.5) {\tiny $c'_2$};
      \end{tikzpicture}
    };
    \node at (-3,-5) {$\rightarrow$};
    \node (B) at (0,-5) {
      \begin{tikzpicture}[anchorbase,scale=.5,xscale=-1]
        \draw [double] (0,0) -- (0,.7);
        \draw [double] (1,0) -- (1,.7);
        \draw [thick] (2,0) -- (2,1.3);
        \draw [thick] (0,.7) -- (1,1.3);
        \draw [thick] (1,.7) -- (1,1.3);
        \draw [thick] (1,.7) -- (2,1.3);
        \draw [thick] (0,.7) -- (0,2);
        \draw [double] (1,1.3) -- (1,2);
        \draw [double] (2,1.3) -- (2,2);
      \end{tikzpicture}
};
    \draw [green,->] (A) edge[bend right=30] node [above,rotate=90] {\tiny $\langle\cdot, c_1\wedge c_2\rangle \wedge c'_1\wedge c'_2$} (B);
    \draw [red,->] (B) edge[bend right=30] node [below,rotate=90] {\tiny $\langle \cdot, c'_1\wedge c'_2\rangle \wedge c_1\wedge c_2$} (A);
  \end{tikzpicture}
\label{R324}
\end{equation}

\begin{equation}
  \begin{tikzpicture}[anchorbase]
    \node at (-4,0) {
            \begin{tikzpicture}[anchorbase,scale=.5]
        \draw [line width=4, white] (2,0) -- (0,2);
        \draw [double, ->] (2,0) -- (0,2);
        \draw [line width=4,white] (0,0) -- (2,2);
        \draw [double, ->] (0,0) -- (2,2);
        \draw [line width=4, white] (1,0) to [out=90,in=-90] (2,1) to [out=90,in=-90] (1,2);
        \draw [thick, ->] (1,0) to [out=90,in=-90] (2,1) to [out=90,in=-90] (1,2); 
        \node at (2.1,.5) {\tiny $c_1$};
        \node at (2.1,1.5) {\tiny $c_2$};
      \end{tikzpicture}
    };
    \node at (-3,0) {$\rightarrow$};
        \node (A) at (0,0) {
      \begin{tikzpicture}[anchorbase,scale=.5]
        \draw [double] (0,0) -- (0,2);
        \draw [thick] (1,0) -- (1,.6);
        \draw [double] (2,0) -- (2,.4);
        \draw [thick] (1,.6) -- (2,.4);
        \draw [double] (1,.6) -- (1,1.4);
        \draw [thick] (2,.4) -- (2,1.6);
        \draw [thick] (1,1.4) -- (2,1.6);
        \draw [thick] (1,1.4) -- (1,2);
        \draw [double] (2,1.6) -- (2,2);
      \end{tikzpicture}
    };
    \node at (-4,-5) {
      \begin{tikzpicture}[anchorbase,scale=.5]
        \draw [line width=4, white] (2,0) -- (0,2);
        \draw [double, ->] (2,0) -- (0,2);
        \draw [line width=4,white] (0,0) -- (2,2);
        \draw [double, ->] (0,0) -- (2,2);
        \draw [line width=4, white] (1,0) to [out=90,in=-90] (0,1) to [out=90,in=-90] (1,2);
        \draw [thick, ->] (1,0) to [out=90,in=-90] (0,1) to [out=90,in=-90] (1,2); 
        \node at (-.2,.5) {\tiny $c'_1$};
        \node at (-.2,1.5) {\tiny $c'_2$};
      \end{tikzpicture}
    };
    \node at (-3,-5) {$\rightarrow$};
    \node (B) at (0,-5) {
      \begin{tikzpicture}[anchorbase,scale=.5]
        \draw [double] (2,0) -- (2,2);
        \draw [thick] (1,0) -- (1,.6);
        \draw [double] (0,0) -- (0,.4);
        \draw [thick] (1,.6) -- (0,.4);
        \draw [double] (1,.6) -- (1,1.4);
        \draw [thick] (0,.4) -- (0,1.6);
        \draw [thick] (1,1.4) -- (0,1.6);
        \draw [thick] (1,1.4) -- (1,2);
        \draw [double] (0,1.6) -- (0,2);
      \end{tikzpicture}
};
    \draw [green,->] (A) edge[bend right=30] node [above,rotate=90] {\tiny $\text{cap/cup}\otimes\langle \cdot, c_1\wedge c_2\rangle \wedge c'_1\wedge c'_2$} (B);
    \draw [red,->] (B) edge[bend right=30] node [below,rotate=90] {\tiny $\text{cap/cup}\otimes \langle \cdot, c'_1\wedge c'_2\rangle \wedge c_1\wedge c_2$} (A);
  \end{tikzpicture}
\label{R325}
\end{equation}

For the following versions of the ${\rm R}_3$ move, we only collect its effect on some states in fully deformed homology, as this is sufficient for our purpose.

\begin{equation}
            \begin{tikzpicture}[anchorbase,scale=.5]
        \draw [line width=4, white] (0,0) -- (2,2);
        \draw [thick, ->] (0,0) -- (2,2);
        \draw [line width=4, white] (2,0) -- (0,2);
        \draw [thick, ->] (2,0) -- (0,2);
        \draw [line width=4, white] (1,0) to [out=90,in=-90] (2,1) to [out=90,in=-90] (1,2);
        \draw [thick, ->] (1,0) to [out=90,in=-90] (2,1) to [out=90,in=-90] (1,2); 
        \node at (2,.5) {\tiny $c_1$};
        \node at (.5,1) {\tiny $c_2$};
        \node at (2,1.5) {\tiny $c_3$};
        \node at (0,-.3) {\tiny $\epsilon_1$};
        \node at (1,-.3) {\tiny $\epsilon_2$};
        \node at (2,-.3) {\tiny $\epsilon_3$};
      \end{tikzpicture}
      \quad \rightarrow \quad
          \begin{tikzpicture}[anchorbase,scale=.5]
        \draw [line width=4, white] (0,0) -- (2,2);
        \draw [thick, ->] (0,0) -- (2,2);
        \draw [line width=4, white] (2,0) -- (0,2);
        \draw [thick, ->] (2,0) -- (0,2);
        \draw [line width=4, white] (1,0) to [out=90,in=-90] (0,1) to [out=90,in=-90] (1,2);
        \draw [thick, ->] (1,0) to [out=90,in=-90] (0,1) to [out=90,in=-90] (1,2); 
        \node at (-.2,.5) {\tiny $c'_1$};
        \node at (1.5,1) {\tiny $c'_2$};
        \node at (-.2,1.5) {\tiny $c'_3$};
        \node at (0,-.3) {\tiny $\epsilon_1$};
        \node at (1,-.3) {\tiny $\epsilon_2$};
        \node at (2,-.3) {\tiny $\epsilon_3$};
      \end{tikzpicture}
  \label{R326}
\end{equation}

\begin{itemize}
\item if $(\epsilon_1,\epsilon_2,\epsilon_3)=(+,+,+)$ or $(-,-,-)$, then the map is $\id$;
\item if $(\epsilon_1,\epsilon_2,\epsilon_3)=(-,+,+)$, then the map is $-\id\otimes \langle \cdot, c_2\wedge c_3\rangle \wedge c'_1\wedge c'_2$.
\end{itemize}

\begin{equation}
          \begin{tikzpicture}[anchorbase,scale=.5]
        \draw [line width=4, white] (1,0) to [out=90,in=-90] (0,1) to [out=90,in=-90] (1,2);
        \draw [thick, ->] (1,0) to [out=90,in=-90] (0,1) to [out=90,in=-90] (1,2); 
        \draw [line width=4, white] (0,0) -- (2,2);
        \draw [thick, ->] (0,0) -- (2,2);
        \draw [line width=4, white] (2,0) -- (0,2);
        \draw [thick, ->] (2,0) -- (0,2);
        \node at (-.2,.5) {\tiny $c_1$};
        \node at (1.5,1) {\tiny $c_2$};
        \node at (-.2,1.5) {\tiny $c_3$};
        \node at (0,-.3) {\tiny $\epsilon_1$};
        \node at (1,-.3) {\tiny $\epsilon_2$};
        \node at (2,-.3) {\tiny $\epsilon_3$};
      \end{tikzpicture}
      \quad \rightarrow \quad
            \begin{tikzpicture}[anchorbase,scale=.5]
        \draw [line width=4, white] (1,0) to [out=90,in=-90] (2,1) to [out=90,in=-90] (1,2);
        \draw [thick, ->] (1,0) to [out=90,in=-90] (2,1) to [out=90,in=-90] (1,2); 
        \draw [line width=4, white] (0,0) -- (2,2);
        \draw [thick, ->] (0,0) -- (2,2);
        \draw [line width=4, white] (2,0) -- (0,2);
        \draw [thick, ->] (2,0) -- (0,2);
        \node at (2,.5) {\tiny $c'_1$};
        \node at (.5,1) {\tiny $c'_2$};
        \node at (2,1.5) {\tiny $c'_3$};
        \node at (0,-.3) {\tiny $\epsilon_1$};
        \node at (1,-.3) {\tiny $\epsilon_2$};
        \node at (2,-.3) {\tiny $\epsilon_3$};
      \end{tikzpicture}
  \label{R327}
\end{equation}

\begin{itemize}
\item if $(\epsilon_1,\epsilon_2,\epsilon_3)=(+,+,+)$, then the map is $\id$;
\item if $(\epsilon_1,\epsilon_2,\epsilon_3)=(+,-,+)$, then the map is $-\text{cap/cup}\otimes \langle \cdot, c_1\wedge c_3\rangle \wedge c'_1\wedge c'_3$;
\item if $(\epsilon_1,\epsilon_2,\epsilon_3)=(-,+,+)$ then the map is $-\id \otimes \langle-,c_1\wedge c_2\rangle \otimes \wedge c'_2\wedge c'_3$. 
\end{itemize}

\begin{equation}
          \begin{tikzpicture}[anchorbase,scale=.5]
        \draw [line width=4, white] (2,0) -- (0,2);
        \draw [thick, ->] (2,0) -- (0,2);
        \draw [line width=4, white] (1,0) to [out=90,in=-90] (0,1) to [out=90,in=-90] (1,2);
        \draw [thick, ->] (1,0) to [out=90,in=-90] (0,1) to [out=90,in=-90] (1,2); 
        \draw [line width=4, white] (0,0) -- (2,2);
        \draw [thick, ->] (0,0) -- (2,2);
        \node at (-.2,.5) {\tiny $c_1$};
        \node at (1.5,1) {\tiny $c_2$};
        \node at (-.2,1.5) {\tiny $c_3$};
        \node at (0,-.3) {\tiny $\epsilon_1$};
        \node at (1,-.3) {\tiny $\epsilon_2$};
        \node at (2,-.3) {\tiny $\epsilon_3$};
      \end{tikzpicture}
      \quad \rightarrow \quad
            \begin{tikzpicture}[anchorbase,scale=.5]
        \draw [line width=4, white] (2,0) -- (0,2);
        \draw [thick, ->] (2,0) -- (0,2);
        \draw [line width=4, white] (1,0) to [out=90,in=-90] (2,1) to [out=90,in=-90] (1,2);
        \draw [thick, ->] (1,0) to [out=90,in=-90] (2,1) to [out=90,in=-90] (1,2); 
        \draw [line width=4, white] (0,0) -- (2,2);
        \draw [thick, ->] (0,0) -- (2,2);
        \node at (2,.5) {\tiny $c'_1$};
        \node at (.5,1) {\tiny $c'_2$};
        \node at (2,1.5) {\tiny $c'_3$};
        \node at (0,-.3) {\tiny $\epsilon_1$};
        \node at (1,-.3) {\tiny $\epsilon_2$};
        \node at (2,-.3) {\tiny $\epsilon_3$};
      \end{tikzpicture}
  \label{R328}
\end{equation}

\begin{itemize}
\item if $(\epsilon_1,\epsilon_2,\epsilon_3)=(+,-,+)$, then the map is $-\text{cap/cup}\otimes \langle -,c_1\wedge c_3\rangle \wedge c'_1\wedge c'_3$.
\end{itemize}

\subsection{Forkslides}

\begin{equation}
  \begin{tikzpicture}[anchorbase]
    \node at (-4,0) {
      \begin{tikzpicture}[anchorbase,scale=.5]
        \draw [very thick, white] (0,0) to [out=90,in=-120] (.5,1);
        \draw [very thick, white] (1,0) to [out=90,in=-60] (.5,1);
        \draw [thick] (0,0) to [out=90,in=-120] (.5,1);
        \draw [thick] (1,0) to [out=90,in=-60] (.5,1);
        \draw [ultra thick, white, double, double=white] (.5,1) to [out=90,in=-90] (1.5,2);
        \draw [double,->] (.5,1) to [out=90,in=-90] (1.5,2);
        \draw [very thick, white, double=black, double distance=.8pt] (2,0) to [out=90,in=-90]  (.5,2);
        \draw [thick, ->] (.5,1.95) -- (.5,2);
        \node at (1.3,1.3) {\tiny $c_3$};
      \end{tikzpicture}
    };
    \node at (-3,0) {$\rightarrow$};
        \node (A) at (2.5,0) {
      \begin{tikzpicture}[anchorbase,scale=.5]
        \draw [thick] (0,0) to [out=90,in=-120] (.5,1);
        \draw [thick] (1,0) to [out=90,in=-60] (.5,1);
        \draw [double] (.5,1) to [out=90,in=-90] (1,1.3);
        \draw [thick] (1,1.3) to [out=120,in=-90] (.5,2);
        \draw [thick] (1,1.3) to [out=60,in=-120] (1.5,1.5);
        \draw [thick] (2,0) to [out=90,in=-60] (1.5,1.5);
        \draw [double] (1.5,1.5) -- (1.5,2);
      \end{tikzpicture}
    };
    \node at (-4,-5) {
      \begin{tikzpicture}[anchorbase,scale=.5]
        \draw [very thick, white] (0,0) to [out=90,in=-120] (1.5,1.5);
        \draw [very thick, white] (1,0) to [out=90,in=-60] (1.5,1.5);
        \draw [thick] (0,0) to [out=90,in=-120] (1.5,1.5);
        \draw [thick] (1,0) to [out=90,in=-60] (1.5,1.5);
        \draw [ultra thick, white, double, double=white] (1.5,1.5) to [out=90,in=-90] (1.5,2);
        \draw [double,->] (1.5,1.5) to [out=90,in=-90] (1.5,2);
        \draw [very thick, white, double=black, double distance=.8pt] (2,0) to [out=90,in=-90]  (.5,2);
        \draw [thick, ->] (.5,1.95) -- (.5,2);
        \node at (1.8,.9) {\tiny $c_1$};
        \node at (.5,1.2) {\tiny $c_2$};
      \end{tikzpicture}
    };
    \node at (-3,-5) {$\rightarrow$};
    \node (B) at (0,-5) {
      \begin{tikzpicture}[anchorbase,scale=.5]
        \node (C1) at (1.5,.3) {};
        \node (C2) at (1.5,.7) {};
        \node (C3) at (.8,1.1) {};
        \node (C4) at (.8,1.5) {};
        \node (C5) at (1.5,1.7) {};
        \draw [thick] (0,0) to [out=90,in=-120] (C3.center);
        \draw [thick] (1,0) to [out=90,in=-120] (C1.center);
        \draw [thick] (2,0) to [out=90,in=-60] (C1.center);
        \draw [double] (C1.center) -- (C2.center);
        \draw [thick] (C2.center) to [out=120,in=-60] (C3.center);
        \draw [thick] (C2.center) to [out=60,in=-60] (C5.center);
        \draw [double] (C3.center) -- (C4.center);
        \draw [thick] (C4.center) to [out=120,in=-90] (.5,2);
        \draw [thick] (C4.center) to [out=60,in=-120] (C5.center);
        \draw [double] (C5.center) to [out=90,in=-90] (1.5,2);
      \end{tikzpicture}};
    \node (C) at (2.5,-4) {
      \begin{tikzpicture}[anchorbase,scale=.5]
        \node (C1) at (1.5,.3) {};
        \node (C2) at (1.5,.7) {};
        \node (C3) at (.8,1.1) {};
        \node (C4) at (.8,1.5) {};
        \node (C5) at (1.5,1.7) {};
        \draw [thick] (0,0) to [out=90,in=-120] (C3.center);
        \draw [thick] (1,0) to [out=90,in=-60] (C3.center);
        \draw [thick] (2,0) to [out=90,in=-60] (C5.center);
        \draw [double] (C3.center) -- (C4.center);
        \draw [thick] (C4.center) to [out=120,in=-90] (.5,2);
        \draw [thick] (C4.center) to [out=60,in=-120] (C5.center);
        \draw [double] (C5.center) to [out=90,in=-90] (1.5,2);
      \end{tikzpicture}};
    \node (D) at (2.5,-6) {
      \begin{tikzpicture}[anchorbase,scale=.5]
        \node (C1) at (1.5,.3) {};
        \node (C2) at (1.5,.7) {};
        \node (C3) at (.8,1.1) {};
        \node (C4) at (.8,1.5) {};
        \node (C5) at (1.5,1.7) {};
        \draw [thick] (0,0) to [out=90,in=-90] (.5,2);
        \draw [thick] (1,0) to [out=90,in=-120] (C1.center);
        \draw [thick] (2,0) to [out=90,in=-60] (C1.center);
        \draw [double] (C1.center) -- (C2.center);
        \draw [thick] (C2.center) to [out=120,in=-90] (1.1,1.2) to [out=90,in=-120] (C5.center);
        \draw [thick] (C2.center) to [out=60,in=-90] (1.9,1.2) to [out=90,in=-60] (C5.center);
        \draw [double] (C5.center) to [out=90,in=-90] (1.5,2);
      \end{tikzpicture}};
    \node (E) at (5,-5) {
      \begin{tikzpicture}[anchorbase,scale=.5]
        \node (C1) at (1.5,.3) {};
        \node (C2) at (1.5,.7) {};
        \node (C3) at (.8,1.1) {};
        \node (C4) at (.8,1.5) {};
        \node (C5) at (1.5,1.7) {};
        \draw [thick] (0,0) to [out=90,in=-90] (.5,2);
        \draw [thick] (1,0) to [out=90,in=-120] (C5.center);
        \draw [thick] (2,0) to [out=90,in=-60] (C5.center);
        \draw [double] (C5.center) to [out=90,in=-90] (1.5,2);
      \end{tikzpicture}};
    \draw (B) -- (C) node[midway,above] {\tiny $\langle \cdot,c_1\rangle$};
    \draw (B) -- (D) node[midway,above] {\tiny $\langle \cdot,c_2\rangle$};
    \draw (C) -- (E) node[midway,above] {\tiny $\langle \cdot,c_2\rangle$};
    \draw (D) -- (E) node[midway,above] {\tiny $\langle \cdot,c_1\rangle$};
    \draw [green,->] (A) edge[bend right=30] node [below,rotate=90] {\tiny $a_{11}\langle \cdot,c_3\rangle \wedge c_2$} (C);
    \draw [green,->] (A) edge[bend right=40] node [above,rotate=90] {\tiny $-a_{11}F\otimes \langle \cdot,c_3\rangle \wedge c_1$} (D);
    \draw [red,->] (C) edge[bend right=30] node [above,rotate=90] {\tiny $a_{11}\langle \cdot,c_2\rangle \wedge c_3$} (A);
    \draw [red,->] (D) edge[bend right=40] node [below,rotate=90] {\tiny $-a_{11}F\otimes \langle \cdot,c_1\rangle \wedge c_3$} (A);
        \node at (8,-2) {$F\leftrightarrow$  
      \begin{tikzpicture}[anchorbase,scale=.3]
        \draw [double] (-.2,1.5)  -- (1,1.5);
        \draw [thick] (0,.5) -- (1.5,.5);
        \draw [thick] (1,1.5) -- (1.5,.5);
        \draw [double] (1.5,.5) -- (3,.5);
        \draw [thick] (1,1.5) to [out=0,in=180] (3.8,2);
        \draw [thick] (3,.5) to [out=30,in=180] (4,1);
        \draw [thick] (3,.5) to [out=-30,in=180] (4.2,0);
        \fill [blue, opacity=.6]  (.5,5.5) to [out=60,in=180] (1.25,6) to [out=0,in=120] (2,5.5) to [out=-90,in=0] (1.25,3.8) to [out=180,in=-90] (.5,5.5);
        \fill [blue, opacity=.6]  (.5,5.5) to [out=-60,in=180] (1.25,5) to [out=0,in=-120] (2,5.5) to [out=-90,in=0] (1.25,3.8) to [out=180,in=-90] (.5,5.5);
        \draw [thick,red] (2,5.5) to [out=-90,in=0] (1.25,3.8) to [out=180,in=-90] (.5,5.5);
        \fill [yellow, opacity=.6] (-.2,1.5) --(1,1.5) to [out=90,in=-90] (3,5.5) -- (2,5.5)  to [out=-90,in=0] (1.25,3.8) to [out=180,in=-90] (.5,5.5) -- (-.2,5.5) -- (-.2,1.5);
        \draw (-.2,1.5) -- (-.2,5.5); 
        \fill [blue, opacity=.6] (1,1.5) to [out=0,in=180] (3.8,2) -- (3.8,6) to [out=180,in=30] (3,5.5) to [out=-90,in=90] (1,1.5);
        \draw (3.8,2) -- (3.8,6);
        \fill [blue, opacity=.6] (4,1) -- (4,5) to [out=180,in=-30] (3,5.5) to [out=-90,in=90] (1,1.5) -- (1.5,.5) to [out=90,in=180] (2.25,1.3) to [out=0,in=90] (3,.5) to [out=30,in=180] (4,1);
        \draw [thick, red] (3,5.5) to [out=-90,in=90] (1,1.5);
        \draw (4,1) -- (4,5);
        \fill [blue,opacity=.6] (0,.5) -- (1.5,.5) to [out=90,in=180] (2.25,1.3) to [out=0,in=90] (3,.5) to [out=-30,in=180] (4.2,0) -- (4.2,4) to [out=180,in=0] (0,4.5) -- (0,.5);
        \draw (0,4.5) -- (0,.5);
        \draw (4.2,0) -- (4.2,4);
        \fill [yellow, opacity=.6] (1.5,.5) -- (3,.5) to [out=90,in=0] (2.25,1.3) to [out=180,in=90] (1.5,.5);
        \draw [thick, red] (3,.5) to [out=90,in=0] (2.25,1.3) to [out=180,in=90] (1.5,.5);
        \draw [double] (-.2,5.5)  -- (.5,5.5);
        \draw [thick] (.5,5.5) to [out=60,in=180] (1.25,6) to [out=0,in=120] (2,5.5);
        \draw [thick] (.5,5.5) to [out=-60,in=180] (1.25,5) to [out=0,in=-120] (2,5.5);
        \draw [thick] (0,4.5) to [out=0,in=180] (4.2,4);
        \draw [double] (2,5.5) -- (3,5.5);
        \draw [thick] (3,5.5) to [out=30,in=180] (3.8,6);
        \draw [thick] (3,5.5) to [out=-30,in=180] (4,5);
  \end{tikzpicture}
  };
\end{tikzpicture}
\label{FS1}
\end{equation}

\begin{equation}
  \begin{tikzpicture}[anchorbase]
    \node at (-4,0) {
      \begin{tikzpicture}[anchorbase,scale=.5]
        \draw [very thick, white] (2,0) to [out=90,in=-60] (1.5,1);
        \draw [very thick, white] (1,0) to [out=90,in=-120] (1.5,1);
        \draw [thick] (2,0) to [out=90,in=-60] (1.5,1);
        \draw [thick] (1,0) to [out=90,in=-120] (1.5,1);
        \draw [ultra thick, white, double, double=white] (1.5,1) to [out=90,in=-90] (.5,2);
        \draw [double,->] (1.5,1) to [out=90,in=-90] (.5,2);
        \draw [very thick, white, double=black, double distance=.8pt] (0,0) to [out=90,in=-90]  (1.5,2);
        \draw [thick, ->] (1.5,1.95) -- (1.5,2);
        \node at (.7,1.3) {\tiny $c_3$};
      \end{tikzpicture}
    };
    \node at (-3,0) {$\rightarrow$};
        \node (A) at (2.5,0) {
      \begin{tikzpicture}[anchorbase,scale=.5,xscale=-1]
        \draw [thick] (0,0) to [out=90,in=-120] (.5,1);
        \draw [thick] (1,0) to [out=90,in=-60] (.5,1);
        \draw [double] (.5,1) to [out=90,in=-90] (1,1.3);
        \draw [thick] (1,1.3) to [out=120,in=-90] (.5,2);
        \draw [thick] (1,1.3) to [out=60,in=-120] (1.5,1.5);
        \draw [thick] (2,0) to [out=90,in=-60] (1.5,1.5);
        \draw [double] (1.5,1.5) -- (1.5,2);
      \end{tikzpicture}
    };
    \node at (-4,-5) {
      \begin{tikzpicture}[anchorbase,scale=.5]
        \draw [very thick, white] (2,0) to [out=90,in=-60] (.5,1.5);
        \draw [very thick, white] (1,0) to [out=90,in=-120] (.5,1.5);
        \draw [thick] (2,0) to [out=90,in=-60] (.5,1.5);
        \draw [thick] (1,0) to [out=90,in=-120] (.5,1.5);
        \draw [ultra thick, white, double, double=white] (.5,1.5) to [out=90,in=-90] (.5,2);
        \draw [double,->] (.5,1.5) to [out=90,in=-90] (.5,2);
        \draw [very thick, white, double=black, double distance=.8pt] (0,0) to [out=90,in=-90]  (1.5,2);
        \draw [thick, ->] (1.5,1.95) -- (1.5,2);
        \node at (.2,.9) {\tiny $c_1$};
        \node at (1.5,1.2) {\tiny $c_2$};
      \end{tikzpicture}
    };
    \node at (-3,-5) {$\rightarrow$};
    \node (B) at (0,-5) {
      \begin{tikzpicture}[anchorbase,scale=.5,xscale=-1]
        \node (C1) at (1.5,.3) {};
        \node (C2) at (1.5,.7) {};
        \node (C3) at (.8,1.1) {};
        \node (C4) at (.8,1.5) {};
        \node (C5) at (1.5,1.7) {};
        \draw [thick] (0,0) to [out=90,in=-90] (.5,2);
        \draw [thick] (1,0) to [out=90,in=-120] (C5.center);
        \draw [thick] (2,0) to [out=90,in=-60] (C5.center);
        \draw [double] (C5.center) to [out=90,in=-90] (1.5,2);
      \end{tikzpicture}};
    \node (C) at (2.5,-4) {
      \begin{tikzpicture}[anchorbase,scale=.5,xscale=-1]
        \node (C1) at (1.5,.3) {};
        \node (C2) at (1.5,.7) {};
        \node (C3) at (.8,1.1) {};
        \node (C4) at (.8,1.5) {};
        \node (C5) at (1.5,1.7) {};
        \draw [thick] (0,0) to [out=90,in=-120] (C3.center);
        \draw [thick] (1,0) to [out=90,in=-60] (C3.center);
        \draw [thick] (2,0) to [out=90,in=-60] (C5.center);
        \draw [double] (C3.center) -- (C4.center);
        \draw [thick] (C4.center) to [out=120,in=-90] (.5,2);
        \draw [thick] (C4.center) to [out=60,in=-120] (C5.center);
        \draw [double] (C5.center) to [out=90,in=-90] (1.5,2);
      \end{tikzpicture}};
    \node (D) at (2.5,-6) {
      \begin{tikzpicture}[anchorbase,scale=.5,xscale=-1]
        \node (C1) at (1.5,.3) {};
        \node (C2) at (1.5,.7) {};
        \node (C3) at (.8,1.1) {};
        \node (C4) at (.8,1.5) {};
        \node (C5) at (1.5,1.7) {};
        \draw [thick] (0,0) to [out=90,in=-90] (.5,2);
        \draw [thick] (1,0) to [out=90,in=-120] (C1.center);
        \draw [thick] (2,0) to [out=90,in=-60] (C1.center);
        \draw [double] (C1.center) -- (C2.center);
        \draw [thick] (C2.center) to [out=120,in=-90] (1.1,1.2) to [out=90,in=-120] (C5.center);
        \draw [thick] (C2.center) to [out=60,in=-90] (1.9,1.2) to [out=90,in=-60] (C5.center);
        \draw [double] (C5.center) to [out=90,in=-90] (1.5,2);
      \end{tikzpicture}};
    \node (E) at (5,-5) {
      \begin{tikzpicture}[anchorbase,scale=.5,xscale=-1]
        \node (C1) at (1.5,.3) {};
        \node (C2) at (1.5,.7) {};
        \node (C3) at (.8,1.1) {};
        \node (C4) at (.8,1.5) {};
        \node (C5) at (1.5,1.7) {};
        \draw [thick] (0,0) to [out=90,in=-120] (C3.center);
        \draw [thick] (1,0) to [out=90,in=-120] (C1.center);
        \draw [thick] (2,0) to [out=90,in=-60] (C1.center);
        \draw [double] (C1.center) -- (C2.center);
        \draw [thick] (C2.center) to [out=120,in=-60] (C3.center);
        \draw [thick] (C2.center) to [out=60,in=-60] (C5.center);
        \draw [double] (C3.center) -- (C4.center);
        \draw [thick] (C4.center) to [out=120,in=-90] (.5,2);
        \draw [thick] (C4.center) to [out=60,in=-120] (C5.center);
        \draw [double] (C5.center) to [out=90,in=-90] (1.5,2);
      \end{tikzpicture}};
    \draw (B) -- (C) node[midway,above] {\tiny $\cdot \wedge c_2$};
    \draw (B) -- (D) node[midway,above] {\tiny $\cdot\wedge c_1$};
    \draw (C) -- (E) node[midway,above] {\tiny $\cdot\wedge c_1$};
    \draw (D) -- (E) node[midway,above] {\tiny $\cdot\wedge c_2$};
    \draw [green,->] (A) edge[bend right=30] node [below,rotate=90] {\tiny $a_3a_{11}\langle \cdot,c_3\rangle \wedge c_2$} (C);
    \draw [green,->] (A) edge[bend right=40] node [above,rotate=90] {\tiny $-a_3a_{11}F\otimes \langle \cdot,c_3\rangle \wedge c_1$} (D);
    \draw [red,->] (C) edge[bend right=30] node [above,rotate=90] {\tiny $a_3a_{11}\langle \cdot,c_2\rangle \wedge c_3$} (A);
    \draw [red,->] (D) edge[bend right=40] node [below,rotate=90] {\tiny $-a_3a_{11}F\otimes \langle \cdot,c_1\rangle \wedge c_3$} (A);
        \node at (8,-2) {$F\leftrightarrow$  
      \begin{tikzpicture}[anchorbase,scale=.3]
        \draw [double] (0,.5)  -- (1,.5);
        \draw [thick] (-.2,1.5) -- (1.5,1.5);
        \draw [thick] (1,.5) -- (1.5,1.5);
        \draw [double] (1.5,1.5) -- (3,1.5);
        \draw [thick] (1,.5) to [out=0,in=180] (4.2,0);
        \draw [thick] (3,1.5) to [out=-30,in=180] (4,1);
        \draw [thick] (3,1.5) to [out=30,in=180] (3.8,2);
        \fill [blue,opacity=.6] (-.2,1.5) -- (1.5,1.5) to [out=90,in=180] (2.25,2.3) to [out=0,in=90] (3,1.5) to [out=-30,in=180] (3.8,2) -- (3.8,6) to [out=180,in=0] (-.2,5.5) -- (-.2,1.5);
        \draw (-.2,5.5) -- (-.2,1.5);
        \draw (3.8,2) -- (3.8,6);
        \fill [yellow, opacity=.6] (1.5,1.5) -- (3,1.5) to [out=90,in=0] (2.25,2.3) to [out=180,in=90] (1.5,1.5);
        \draw [thick, red] (3,1.5) to [out=90,in=0] (2.25,2.3) to [out=180,in=90] (1.5,1.5);
        \fill [blue, opacity=.6] (4,1) -- (4,5) to [out=180,in=30] (3,4.5) to [out=-90,in=90] (1,.5) -- (1.5,1.5) to [out=90,in=180] (2.25,2.3) to [out=0,in=90] (3,1.5) to [out=-30,in=180] (4,1);
        \draw (4,1) -- (4,5);
        \fill [yellow, opacity=.6] (0,.5) --(1,.5) to [out=90,in=-90] (3,4.5) -- (2,4.5)  to [out=-90,in=0] (1.25,2.8) to [out=180,in=-90] (.5,4.5) -- (0,4.5) -- (0,.5);
        \draw (0,.5) -- (0,4.5); 
        \fill [blue, opacity=.6] (1,.5) to [out=0,in=180] (4.2,0) -- (4.2,4) to [out=180,in=-30] (3,4.5) to [out=-90,in=90] (1,.5);
        \draw (4.2,0) -- (4.2,4);
        \draw [thick, red] (3,4.5) to [out=-90,in=90] (1,.5);
        \fill [blue, opacity=.6]  (.5,4.5) to [out=60,in=180] (1.25,5) to [out=0,in=120] (2,4.5) to [out=-90,in=0] (1.25,2.8) to [out=180,in=-90] (.5,4.5);
        \fill [blue, opacity=.6]  (.5,4.5) to [out=-60,in=180] (1.25,4) to [out=0,in=-120] (2,4.5) to [out=-90,in=0] (1.25,2.8) to [out=180,in=-90] (.5,4.5);
        \draw [thick,red] (2,4.5) to [out=-90,in=0] (1.25,2.8) to [out=180,in=-90] (.5,4.5);
        \draw [double] (0,4.5)  -- (.5,4.5);
        \draw [thick] (.5,4.5) to [out=60,in=180] (1.25,5) to [out=0,in=120] (2,4.5);
        \draw [thick] (.5,4.5) to [out=-60,in=180] (1.25,4) to [out=0,in=-120] (2,4.5);
        \draw [thick] (-.2,5.5) to [out=0,in=180] (3.8,6);
        \draw [double] (2,4.5) -- (3,4.5);
        \draw [thick] (3,4.5) to [out=-30,in=180] (4.2,4);
        \draw [thick] (3,4.5) to [out=30,in=180] (4,5);
  \end{tikzpicture}
  };
  \end{tikzpicture}
  \label{FS2}
\end{equation}

\begin{equation}
  \begin{tikzpicture}[anchorbase]
    \node at (-4,0) {
      \begin{tikzpicture}[anchorbase,scale=.5]
        \draw [very thick, white, double=black, double distance=.8pt] (2,0) to [out=90,in=-90]  (.5,2);
        \draw [thick, ->] (.5,1.95) -- (.5,2);
        \draw [very thick, white] (0,0) to [out=90,in=-120] (.5,1);
        \draw [very thick, white] (1,0) to [out=90,in=-60] (.5,1);
        \draw [thick] (0,0) to [out=90,in=-120] (.5,1);
        \draw [thick] (1,0) to [out=90,in=-60] (.5,1);
        \draw [ultra thick, white, double, double=white] (.5,1) to [out=90,in=-90] (1.5,2);
        \draw [double,->] (.5,1) to [out=90,in=-90] (1.5,2);
        \node at (1.3,1.3) {\tiny $c_3$};
      \end{tikzpicture}
    };
    \node at (-3,0) {$\rightarrow$};
        \node (A) at (2.5,0) {
      \begin{tikzpicture}[anchorbase,scale=.5]
        \draw [thick] (0,0) to [out=90,in=-120] (.5,1);
        \draw [thick] (1,0) to [out=90,in=-60] (.5,1);
        \draw [double] (.5,1) to [out=90,in=-90] (1,1.3);
        \draw [thick] (1,1.3) to [out=120,in=-90] (.5,2);
        \draw [thick] (1,1.3) to [out=60,in=-120] (1.5,1.5);
        \draw [thick] (2,0) to [out=90,in=-60] (1.5,1.5);
        \draw [double] (1.5,1.5) -- (1.5,2);
      \end{tikzpicture}
    };
    \node at (-4,-5) {
      \begin{tikzpicture}[anchorbase,scale=.5]
        \draw [very thick, white, double=black, double distance=.8pt] (2,0) to [out=90,in=-90]  (.5,2);
        \draw [thick, ->] (.5,1.95) -- (.5,2);
        \draw [ultra thick, white] (0,0) to [out=90,in=-120] (1.5,1.5);
        \draw [ultra thick, white] (1,0) to [out=90,in=-60] (1.5,1.5);
        \draw [thick] (0,0) to [out=90,in=-120] (1.5,1.5);
        \draw [thick] (1,0) to [out=90,in=-60] (1.5,1.5);
        \draw [ultra thick, white, double, double=white] (1.5,1.5) to [out=90,in=-90] (1.5,2);
        \draw [double,->] (1.5,1.5) to [out=90,in=-90] (1.5,2);
        \node at (1.8,.9) {\tiny $c_1$};
        \node at (.5,1.2) {\tiny $c_2$};
      \end{tikzpicture}
    };
    \node at (-3,-5) {$\rightarrow$};
    \node (B) at (0,-5) {
      \begin{tikzpicture}[anchorbase,scale=.5]
        \node (C1) at (1.5,.3) {};
        \node (C2) at (1.5,.7) {};
        \node (C3) at (.8,1.1) {};
        \node (C4) at (.8,1.5) {};
        \node (C5) at (1.5,1.7) {};
        \draw [thick] (0,0) to [out=90,in=-90] (.5,2);
        \draw [thick] (1,0) to [out=90,in=-120] (C5.center);
        \draw [thick] (2,0) to [out=90,in=-60] (C5.center);
        \draw [double] (C5.center) to [out=90,in=-90] (1.5,2);
      \end{tikzpicture}};
    \node (C) at (2.5,-4) {
      \begin{tikzpicture}[anchorbase,scale=.5]
        \node (C1) at (1.5,.3) {};
        \node (C2) at (1.5,.7) {};
        \node (C3) at (.8,1.1) {};
        \node (C4) at (.8,1.5) {};
        \node (C5) at (1.5,1.7) {};
        \draw [thick] (0,0) to [out=90,in=-120] (C3.center);
        \draw [thick] (1,0) to [out=90,in=-60] (C3.center);
        \draw [thick] (2,0) to [out=90,in=-60] (C5.center);
        \draw [double] (C3.center) -- (C4.center);
        \draw [thick] (C4.center) to [out=120,in=-90] (.5,2);
        \draw [thick] (C4.center) to [out=60,in=-120] (C5.center);
        \draw [double] (C5.center) to [out=90,in=-90] (1.5,2);
      \end{tikzpicture}};
    \node (D) at (2.5,-6) {
      \begin{tikzpicture}[anchorbase,scale=.5]
        \node (C1) at (1.5,.3) {};
        \node (C2) at (1.5,.7) {};
        \node (C3) at (.8,1.1) {};
        \node (C4) at (.8,1.5) {};
        \node (C5) at (1.5,1.7) {};
        \draw [thick] (0,0) to [out=90,in=-90] (.5,2);
        \draw [thick] (1,0) to [out=90,in=-120] (C1.center);
        \draw [thick] (2,0) to [out=90,in=-60] (C1.center);
        \draw [double] (C1.center) -- (C2.center);
        \draw [thick] (C2.center) to [out=120,in=-90] (1.1,1.2) to [out=90,in=-120] (C5.center);
        \draw [thick] (C2.center) to [out=60,in=-90] (1.9,1.2) to [out=90,in=-60] (C5.center);
        \draw [double] (C5.center) to [out=90,in=-90] (1.5,2);
      \end{tikzpicture}};
    \node (E) at (5,-5) {
      \begin{tikzpicture}[anchorbase,scale=.5]
        \node (C1) at (1.5,.3) {};
        \node (C2) at (1.5,.7) {};
        \node (C3) at (.8,1.1) {};
        \node (C4) at (.8,1.5) {};
        \node (C5) at (1.5,1.7) {};
        \draw [thick] (0,0) to [out=90,in=-120] (C3.center);
        \draw [thick] (1,0) to [out=90,in=-120] (C1.center);
        \draw [thick] (2,0) to [out=90,in=-60] (C1.center);
        \draw [double] (C1.center) -- (C2.center);
        \draw [thick] (C2.center) to [out=120,in=-60] (C3.center);
        \draw [thick] (C2.center) to [out=60,in=-60] (C5.center);
        \draw [double] (C3.center) -- (C4.center);
        \draw [thick] (C4.center) to [out=120,in=-90] (.5,2);
        \draw [thick] (C4.center) to [out=60,in=-120] (C5.center);
        \draw [double] (C5.center) to [out=90,in=-90] (1.5,2);
      \end{tikzpicture}};
    \draw (B) -- (C) node[midway,above] {\tiny $\cdot\wedge c_2$};
    \draw (B) -- (D) node[midway,above] {\tiny $\cdot\wedge c_1$};
    \draw (C) -- (E) node[midway,above] {\tiny $\cdot\wedge c_1$};
    \draw (D) -- (E) node[midway,above] {\tiny $\cdot \wedge c_2$};
    \draw [green,->] (A) edge[bend right=30] node [below,rotate=90] {\tiny $a_{11}\langle \cdot,c_3\rangle \wedge c_2$} (C);
    \draw [green,->] (A) edge[bend right=40] node [above,rotate=90] {\tiny $-a_{11}F\otimes \langle \cdot,c_3\rangle \wedge c_1$} (D);
    \draw [red,->] (C) edge[bend right=30] node [above,rotate=90] {\tiny $a_{11}\langle \cdot,c_2\rangle \wedge c_3$} (A);
    \draw [red,->] (D) edge[bend right=40] node [below,rotate=90] {\tiny $-a_{11}F\otimes \langle \cdot,c_1\rangle \wedge c_3$} (A);
        \node at (8,-2) {$F\leftrightarrow$  
      \begin{tikzpicture}[anchorbase,scale=.3]
        \draw [double] (-.2,1.5)  -- (1,1.5);
        \draw [thick] (0,.5) -- (1.5,.5);
        \draw [thick] (1,1.5) -- (1.5,.5);
        \draw [double] (1.5,.5) -- (3,.5);
        \draw [thick] (1,1.5) to [out=0,in=180] (3.8,2);
        \draw [thick] (3,.5) to [out=30,in=180] (4,1);
        \draw [thick] (3,.5) to [out=-30,in=180] (4.2,0);
        \fill [blue, opacity=.6]  (.5,5.5) to [out=60,in=180] (1.25,6) to [out=0,in=120] (2,5.5) to [out=-90,in=0] (1.25,3.8) to [out=180,in=-90] (.5,5.5);
        \fill [blue, opacity=.6]  (.5,5.5) to [out=-60,in=180] (1.25,5) to [out=0,in=-120] (2,5.5) to [out=-90,in=0] (1.25,3.8) to [out=180,in=-90] (.5,5.5);
        \draw [thick,red] (2,5.5) to [out=-90,in=0] (1.25,3.8) to [out=180,in=-90] (.5,5.5);
        \fill [yellow, opacity=.6] (-.2,1.5) --(1,1.5) to [out=90,in=-90] (3,5.5) -- (2,5.5)  to [out=-90,in=0] (1.25,3.8) to [out=180,in=-90] (.5,5.5) -- (-.2,5.5) -- (-.2,1.5);
        \draw (-.2,1.5) -- (-.2,5.5); 
        \fill [blue, opacity=.6] (1,1.5) to [out=0,in=180] (3.8,2) -- (3.8,6) to [out=180,in=30] (3,5.5) to [out=-90,in=90] (1,1.5);
        \draw (3.8,2) -- (3.8,6);
        \fill [blue, opacity=.6] (4,1) -- (4,5) to [out=180,in=-30] (3,5.5) to [out=-90,in=90] (1,1.5) -- (1.5,.5) to [out=90,in=180] (2.25,1.3) to [out=0,in=90] (3,.5) to [out=30,in=180] (4,1);
        \draw [thick, red] (3,5.5) to [out=-90,in=90] (1,1.5);
        \draw (4,1) -- (4,5);
        \fill [blue,opacity=.6] (0,.5) -- (1.5,.5) to [out=90,in=180] (2.25,1.3) to [out=0,in=90] (3,.5) to [out=-30,in=180] (4.2,0) -- (4.2,4) to [out=180,in=0] (0,4.5) -- (0,.5);
        \draw (0,4.5) -- (0,.5);
        \draw (4.2,0) -- (4.2,4);
        \fill [yellow, opacity=.6] (1.5,.5) -- (3,.5) to [out=90,in=0] (2.25,1.3) to [out=180,in=90] (1.5,.5);
        \draw [thick, red] (3,.5) to [out=90,in=0] (2.25,1.3) to [out=180,in=90] (1.5,.5);
        \draw [double] (-.2,5.5)  -- (.5,5.5);
        \draw [thick] (.5,5.5) to [out=60,in=180] (1.25,6) to [out=0,in=120] (2,5.5);
        \draw [thick] (.5,5.5) to [out=-60,in=180] (1.25,5) to [out=0,in=-120] (2,5.5);
        \draw [thick] (0,4.5) to [out=0,in=180] (4.2,4);
        \draw [double] (2,5.5) -- (3,5.5);
        \draw [thick] (3,5.5) to [out=30,in=180] (3.8,6);
        \draw [thick] (3,5.5) to [out=-30,in=180] (4,5);
  \end{tikzpicture}
  };
  \end{tikzpicture}
  \label{FS3}
\end{equation}

\begin{equation}
  \begin{tikzpicture}[anchorbase]
    \node at (-4,0) {
      \begin{tikzpicture}[anchorbase,scale=.5]
        \draw [very thick, white, double=black, double distance=.8pt] (0,0) to [out=90,in=-90]  (1.5,2);
        \draw [thick, ->] (1.5,1.95) -- (1.5,2);
        \draw [ultra thick, white] (2,0) to [out=90,in=-60] (1.5,1);
        \draw [ultra thick, white] (1,0) to [out=90,in=-120] (1.5,1);
        \draw [thick] (2,0) to [out=90,in=-60] (1.5,1);
        \draw [thick] (1,0) to [out=90,in=-120] (1.5,1);
        \draw [ultra thick, white, double, double=white] (1.5,1) to [out=90,in=-90] (.5,2);
        \draw [double,->] (1.5,1) to [out=90,in=-90] (.5,2);
        \node at (.7,1.3) {\tiny $c_3$};
      \end{tikzpicture}
    };
    \node at (-3,0) {$\rightarrow$};
        \node (A) at (2.5,0) {
      \begin{tikzpicture}[anchorbase,scale=.5,xscale=-1]
        \draw [thick] (0,0) to [out=90,in=-120] (.5,1);
        \draw [thick] (1,0) to [out=90,in=-60] (.5,1);
        \draw [double] (.5,1) to [out=90,in=-90] (1,1.3);
        \draw [thick] (1,1.3) to [out=120,in=-90] (.5,2);
        \draw [thick] (1,1.3) to [out=60,in=-120] (1.5,1.5);
        \draw [thick] (2,0) to [out=90,in=-60] (1.5,1.5);
        \draw [double] (1.5,1.5) -- (1.5,2);
      \end{tikzpicture}
    };
    \node at (-4,-5) {
      \begin{tikzpicture}[anchorbase,scale=.5]
        \draw [very thick, white, double=black, double distance=.8pt] (0,0) to [out=90,in=-90]  (1.5,2);
        \draw [thick, ->] (1.5,1.95) -- (1.5,2);
        \draw [ultra thick, white] (2,0) to [out=90,in=-60] (.5,1.5);
        \draw [ultra thick, white] (1,0) to [out=90,in=-120] (.5,1.5);
        \draw [thick] (2,0) to [out=90,in=-60] (.5,1.5);
        \draw [thick] (1,0) to [out=90,in=-120] (.5,1.5);
        \draw [ultra thick, white, double, double=white] (.5,1.5) to [out=90,in=-90] (.5,2);
        \draw [double,->] (.5,1.5) to [out=90,in=-90] (.5,2);
        \node at (.2,.9) {\tiny $c_1$};
        \node at (1.5,1.2) {\tiny $c_2$};
      \end{tikzpicture}
    };
    \node at (-3,-5) {$\rightarrow$};
    \node (B) at (0,-5) {
      \begin{tikzpicture}[anchorbase,scale=.5,xscale=-1]
        \node (C1) at (1.5,.3) {};
        \node (C2) at (1.5,.7) {};
        \node (C3) at (.8,1.1) {};
        \node (C4) at (.8,1.5) {};
        \node (C5) at (1.5,1.7) {};
        \draw [thick] (0,0) to [out=90,in=-120] (C3.center);
        \draw [thick] (1,0) to [out=90,in=-120] (C1.center);
        \draw [thick] (2,0) to [out=90,in=-60] (C1.center);
        \draw [double] (C1.center) -- (C2.center);
        \draw [thick] (C2.center) to [out=120,in=-60] (C3.center);
        \draw [thick] (C2.center) to [out=60,in=-60] (C5.center);
        \draw [double] (C3.center) -- (C4.center);
        \draw [thick] (C4.center) to [out=120,in=-90] (.5,2);
        \draw [thick] (C4.center) to [out=60,in=-120] (C5.center);
        \draw [double] (C5.center) to [out=90,in=-90] (1.5,2);
      \end{tikzpicture}};
    \node (C) at (2.5,-4) {
      \begin{tikzpicture}[anchorbase,scale=.5,xscale=-1]
        \node (C1) at (1.5,.3) {};
        \node (C2) at (1.5,.7) {};
        \node (C3) at (.8,1.1) {};
        \node (C4) at (.8,1.5) {};
        \node (C5) at (1.5,1.7) {};
        \draw [thick] (0,0) to [out=90,in=-120] (C3.center);
        \draw [thick] (1,0) to [out=90,in=-60] (C3.center);
        \draw [thick] (2,0) to [out=90,in=-60] (C5.center);
        \draw [double] (C3.center) -- (C4.center);
        \draw [thick] (C4.center) to [out=120,in=-90] (.5,2);
        \draw [thick] (C4.center) to [out=60,in=-120] (C5.center);
        \draw [double] (C5.center) to [out=90,in=-90] (1.5,2);
      \end{tikzpicture}};
    \node (D) at (2.5,-6) {
      \begin{tikzpicture}[anchorbase,scale=.5,xscale=-1]
        \node (C1) at (1.5,.3) {};
        \node (C2) at (1.5,.7) {};
        \node (C3) at (.8,1.1) {};
        \node (C4) at (.8,1.5) {};
        \node (C5) at (1.5,1.7) {};
        \draw [thick] (0,0) to [out=90,in=-90] (.5,2);
        \draw [thick] (1,0) to [out=90,in=-120] (C1.center);
        \draw [thick] (2,0) to [out=90,in=-60] (C1.center);
        \draw [double] (C1.center) -- (C2.center);
        \draw [thick] (C2.center) to [out=120,in=-90] (1.1,1.2) to [out=90,in=-120] (C5.center);
        \draw [thick] (C2.center) to [out=60,in=-90] (1.9,1.2) to [out=90,in=-60] (C5.center);
        \draw [double] (C5.center) to [out=90,in=-90] (1.5,2);
      \end{tikzpicture}};
    \node (E) at (5,-5) {
      \begin{tikzpicture}[anchorbase,scale=.5,xscale=-1]
        \node (C1) at (1.5,.3) {};
        \node (C2) at (1.5,.7) {};
        \node (C3) at (.8,1.1) {};
        \node (C4) at (.8,1.5) {};
        \node (C5) at (1.5,1.7) {};
        \draw [thick] (0,0) to [out=90,in=-90] (.5,2);
        \draw [thick] (1,0) to [out=90,in=-120] (C5.center);
        \draw [thick] (2,0) to [out=90,in=-60] (C5.center);
        \draw [double] (C5.center) to [out=90,in=-90] (1.5,2);
      \end{tikzpicture}};
    \draw (B) -- (C) node[midway,above] {\tiny $\langle \cdot, c_1\rangle$};
    \draw (B) -- (D) node[midway,above] {\tiny $\langle \cdot,c_2\rangle$};
    \draw (C) -- (E) node[midway,above] {\tiny $\langle \cdot, c_2\rangle$};
    \draw (D) -- (E) node[midway,above] {\tiny $\langle\cdot,c_1\rangle$};
    \draw [green,->] (A) edge[bend right=30] node [below,rotate=90] {\tiny $-a_4a_{11}\langle \cdot,c_3\rangle \wedge c_2$} (C);
    \draw [green,->] (A) edge[bend right=40] node [above,rotate=90] {\tiny $a_4a_{11}F\otimes \langle \cdot,c_3\rangle \wedge c_1$} (D);
    \draw [red,->] (C) edge[bend right=30] node [above,rotate=90] {\tiny $-a_4a_{11}\langle \cdot,c_2\rangle \wedge c_3$} (A);
    \draw [red,->] (D) edge[bend right=40] node [below,rotate=90] {\tiny $a_4a_{11}F\otimes \langle \cdot,c_1\rangle \wedge c_3$} (A);
        \node at (8,-2) {$F\leftrightarrow$  
};
    \draw (B) -- (C) node[midway,above] {\tiny $\langle \cdot,c_2\rangle$};
    \draw (B) -- (D) node[midway,above] {\tiny $\langle \cdot,c_1\rangle$};
    \draw (C) -- (E) node[midway,above] {\tiny $\langle \cdot,c_1\rangle$};
    \draw (D) -- (E) node[midway,above] {\tiny $\langle \cdot,c_2\rangle$};
    \draw [green,->] (A) edge[bend right=30] node [below,rotate=90] {\tiny $-a_{11}\langle \cdot,c_3\rangle \wedge c_1$} (C);
    \draw [green,->] (A) edge[bend right=40] node [above,rotate=90] {\tiny $a_{11}F\otimes \langle \cdot,c_3\rangle \wedge c_2$} (D);
    \draw [red,->] (C) edge[bend right=30] node [above,rotate=90] {\tiny $-a_{11}\langle \cdot,c_1\rangle \wedge c_3$} (A);
    \draw [red,->] (D) edge[bend right=40] node [below,rotate=90] {\tiny $a_{11}F\otimes \langle \cdot,c_2\rangle \wedge c_3$} (A);
        \node at (8,-2) {$F\leftrightarrow$  
      \begin{tikzpicture}[anchorbase,scale=-.3]
        \draw [double] (-.2,5.5)  -- (.5,5.5);
        \draw [thick] (.5,5.5) to [out=60,in=180] (1.25,6) to [out=0,in=120] (2,5.5);
        \draw [thick] (.5,5.5) to [out=-60,in=180] (1.25,5) to [out=0,in=-120] (2,5.5);
        \draw [thick] (0,4.5) to [out=0,in=180] (4.2,4);
        \draw [double] (2,5.5) -- (3,5.5);
        \draw [thick] (3,5.5) to [out=30,in=180] (3.8,6);
        \draw [thick] (3,5.5) to [out=-30,in=180] (4,5);
        \fill [blue,opacity=.6] (0,.5) -- (1.5,.5) to [out=90,in=180] (2.25,1.3) to [out=0,in=90] (3,.5) to [out=-30,in=180] (4.2,0) -- (4.2,4) to [out=180,in=0] (0,4.5) -- (0,.5);
        \draw (0,4.5) -- (0,.5);
        \draw (4.2,0) -- (4.2,4);
        \fill [yellow, opacity=.6] (1.5,.5) -- (3,.5) to [out=90,in=0] (2.25,1.3) to [out=180,in=90] (1.5,.5);
        \draw [thick, red] (3,.5) to [out=90,in=0] (2.25,1.3) to [out=180,in=90] (1.5,.5);
        \fill [blue, opacity=.6] (4,1) -- (4,5) to [out=180,in=-30] (3,5.5) to [out=-90,in=90] (1,1.5) -- (1.5,.5) to [out=90,in=180] (2.25,1.3) to [out=0,in=90] (3,.5) to [out=30,in=180] (4,1);
        \draw [thick, red] (3,5.5) to [out=-90,in=90] (1,1.5);
        \draw (4,1) -- (4,5);
        \fill [yellow, opacity=.6] (-.2,1.5) --(1,1.5) to [out=90,in=-90] (3,5.5) -- (2,5.5)  to [out=-90,in=0] (1.25,3.8) to [out=180,in=-90] (.5,5.5) -- (-.2,5.5) -- (-.2,1.5);
        \draw (-.2,1.5) -- (-.2,5.5); 
        \fill [blue, opacity=.6] (1,1.5) to [out=0,in=180] (3.8,2) -- (3.8,6) to [out=180,in=30] (3,5.5) to [out=-90,in=90] (1,1.5);
        \draw (3.8,2) -- (3.8,6);
        \fill [blue, opacity=.6]  (.5,5.5) to [out=60,in=180] (1.25,6) to [out=0,in=120] (2,5.5) to [out=-90,in=0] (1.25,3.8) to [out=180,in=-90] (.5,5.5);
        \fill [blue, opacity=.6]  (.5,5.5) to [out=-60,in=180] (1.25,5) to [out=0,in=-120] (2,5.5) to [out=-90,in=0] (1.25,3.8) to [out=180,in=-90] (.5,5.5);
        \draw [thick,red] (2,5.5) to [out=-90,in=0] (1.25,3.8) to [out=180,in=-90] (.5,5.5);
        \draw [double] (-.2,1.5)  -- (1,1.5);
        \draw [thick] (0,.5) -- (1.5,.5);
        \draw [thick] (1,1.5) -- (1.5,.5);
        \draw [double] (1.5,.5) -- (3,.5);
        \draw [thick] (1,1.5) to [out=0,in=180] (3.8,2);
        \draw [thick] (3,.5) to [out=30,in=180] (4,1);
        \draw [thick] (3,.5) to [out=-30,in=180] (4.2,0);
  \end{tikzpicture}
  };
  \end{tikzpicture}
  \label{FS9}
\end{equation}

\begin{equation}
  \begin{tikzpicture}[anchorbase]
    \node at (-4,0) {
      \begin{tikzpicture}[anchorbase,scale=.5, rotate=180]
        \draw [very thick, white] (2,0) to [out=90,in=-60] (1.5,1);
        \draw [very thick, white] (1,0) to [out=90,in=-120] (1.5,1);
        \draw [<-,thick] (2,0) to [out=90,in=-60] (1.5,1);
        \draw [<-,thick] (1,0) to [out=90,in=-120] (1.5,1);
        \draw [ultra thick, white, double, double=white] (1.5,1) to [out=90,in=-90] (.5,2);
        \draw [double] (1.5,1) to [out=90,in=-90] (.5,2);
        \draw [very thick, white, double=black, double distance=.8pt] (0,0) to [out=90,in=-90]  (1.5,2);
        \draw [thick, ->] (0,0.05) -- (0,0);
        \node at (.7,1.3) {\tiny $c_3$};
      \end{tikzpicture}
    };
    \node at (-3,0) {$\rightarrow$};
        \node (A) at (2.5,0) {
      \begin{tikzpicture}[anchorbase,scale=.5,xscale=-1, rotate=180]
        \draw [thick] (0,0) to [out=90,in=-120] (.5,1);
        \draw [thick] (1,0) to [out=90,in=-60] (.5,1);
        \draw [double] (.5,1) to [out=90,in=-90] (1,1.3);
        \draw [thick] (1,1.3) to [out=120,in=-90] (.5,2);
        \draw [thick] (1,1.3) to [out=60,in=-120] (1.5,1.5);
        \draw [thick] (2,0) to [out=90,in=-60] (1.5,1.5);
        \draw [double] (1.5,1.5) -- (1.5,2);
      \end{tikzpicture}
    };
    \node at (-4,-5) {
      \begin{tikzpicture}[anchorbase,scale=.5, rotate=180]
        \draw [very thick, white] (2,0) to [out=90,in=-60] (.5,1.5);
        \draw [very thick, white] (1,0) to [out=90,in=-120] (.5,1.5);
        \draw [<-,thick] (2,0) to [out=90,in=-60] (.5,1.5);
        \draw [<-,thick] (1,0) to [out=90,in=-120] (.5,1.5);
        \draw [ultra thick, white, double, double=white] (.5,1.5) to [out=90,in=-90] (.5,2);
        \draw [double] (.5,1.5) to [out=90,in=-90] (.5,2);
        \draw [very thick, white, double=black, double distance=.8pt] (0,0) to [out=90,in=-90]  (1.5,2);
        \draw [thick, ->] (0,0.05) -- (0,0);
        \node at (.2,.9) {\tiny $c_2$};
        \node at (1.5,1.2) {\tiny $c_1$};
      \end{tikzpicture}
    };
    \node at (-3,-5) {$\rightarrow$};
    \node (B) at (0,-5) {
      \begin{tikzpicture}[anchorbase,scale=.5,xscale=-1, rotate=180]
        \node (C1) at (1.5,.3) {};
        \node (C2) at (1.5,.7) {};
        \node (C3) at (.8,1.1) {};
        \node (C4) at (.8,1.5) {};
        \node (C5) at (1.5,1.7) {};
        \draw [thick] (0,0) to [out=90,in=-90] (.5,2);
        \draw [thick] (1,0) to [out=90,in=-120] (C5.center);
        \draw [thick] (2,0) to [out=90,in=-60] (C5.center);
        \draw [double] (C5.center) to [out=90,in=-90] (1.5,2);
      \end{tikzpicture}};
    \node (C) at (2.5,-4) {
      \begin{tikzpicture}[anchorbase,scale=.5,xscale=-1, rotate=180]
        \node (C1) at (1.5,.3) {};
        \node (C2) at (1.5,.7) {};
        \node (C3) at (.8,1.1) {};
        \node (C4) at (.8,1.5) {};
        \node (C5) at (1.5,1.7) {};
        \draw [thick] (0,0) to [out=90,in=-120] (C3.center);
        \draw [thick] (1,0) to [out=90,in=-60] (C3.center);
        \draw [thick] (2,0) to [out=90,in=-60] (C5.center);
        \draw [double] (C3.center) -- (C4.center);
        \draw [thick] (C4.center) to [out=120,in=-90] (.5,2);
        \draw [thick] (C4.center) to [out=60,in=-120] (C5.center);
        \draw [double] (C5.center) to [out=90,in=-90] (1.5,2);
      \end{tikzpicture}};
    \node (D) at (2.5,-6) {
      \begin{tikzpicture}[anchorbase,scale=.5,xscale=-1, rotate=180]
        \node (C1) at (1.5,.3) {};
        \node (C2) at (1.5,.7) {};
        \node (C3) at (.8,1.1) {};
        \node (C4) at (.8,1.5) {};
        \node (C5) at (1.5,1.7) {};
        \draw [thick] (0,0) to [out=90,in=-90] (.5,2);
        \draw [thick] (1,0) to [out=90,in=-120] (C1.center);
        \draw [thick] (2,0) to [out=90,in=-60] (C1.center);
        \draw [double] (C1.center) -- (C2.center);
        \draw [thick] (C2.center) to [out=120,in=-90] (1.1,1.2) to [out=90,in=-120] (C5.center);
        \draw [thick] (C2.center) to [out=60,in=-90] (1.9,1.2) to [out=90,in=-60] (C5.center);
        \draw [double] (C5.center) to [out=90,in=-90] (1.5,2);
      \end{tikzpicture}};
    \node (E) at (5,-5) {
      \begin{tikzpicture}[anchorbase,scale=.5,xscale=-1, rotate=180]
        \node (C1) at (1.5,.3) {};
        \node (C2) at (1.5,.7) {};
        \node (C3) at (.8,1.1) {};
        \node (C4) at (.8,1.5) {};
        \node (C5) at (1.5,1.7) {};
        \draw [thick] (0,0) to [out=90,in=-120] (C3.center);
        \draw [thick] (1,0) to [out=90,in=-120] (C1.center);
        \draw [thick] (2,0) to [out=90,in=-60] (C1.center);
        \draw [double] (C1.center) -- (C2.center);
        \draw [thick] (C2.center) to [out=120,in=-60] (C3.center);
        \draw [thick] (C2.center) to [out=60,in=-60] (C5.center);
        \draw [double] (C3.center) -- (C4.center);
        \draw [thick] (C4.center) to [out=120,in=-90] (.5,2);
        \draw [thick] (C4.center) to [out=60,in=-120] (C5.center);
        \draw [double] (C5.center) to [out=90,in=-90] (1.5,2);
      \end{tikzpicture}};
    \draw (B) -- (C) node[midway,above] {\tiny $\cdot \wedge c_1$};
    \draw (B) -- (D) node[midway,above] {\tiny $\cdot\wedge c_2$};
    \draw (C) -- (E) node[midway,above] {\tiny $\cdot\wedge c_2$};
    \draw (D) -- (E) node[midway,above] {\tiny $\cdot\wedge c_1$};
    \draw [green,->] (A) edge[bend right=30] node [below,rotate=90] {\tiny $-a_3a_{11}\langle \cdot,c_3\rangle \wedge c_1$} (C);
    \draw [green,->] (A) edge[bend right=40] node [above,rotate=90] {\tiny $a_3a_{11}F\otimes \langle \cdot,c_3\rangle \wedge c_2$} (D);
    \draw [red,->] (C) edge[bend right=30] node [above,rotate=90] {\tiny $-a_3a_{11}\langle \cdot,c_1\rangle \wedge c_3$} (A);
    \draw [red,->] (D) edge[bend right=40] node [below,rotate=90] {\tiny $a_3a_{11}F\otimes \langle \cdot,c_2\rangle \wedge c_3$} (A);
        \node at (8,-2) {$F\leftrightarrow$  
      \begin{tikzpicture}[anchorbase,scale=-.3 ]
        \draw [double] (0,4.5)  -- (.5,4.5);
        \draw [thick] (.5,4.5) to [out=60,in=180] (1.25,5) to [out=0,in=120] (2,4.5);
        \draw [thick] (.5,4.5) to [out=-60,in=180] (1.25,4) to [out=0,in=-120] (2,4.5);
        \draw [thick] (-.2,5.5) to [out=0,in=180] (3.8,6);
        \draw [double] (2,4.5) -- (3,4.5);
        \draw [thick] (3,4.5) to [out=-30,in=180] (4.2,4);
        \draw [thick] (3,4.5) to [out=30,in=180] (4,5);
        \fill [blue, opacity=.6]  (.5,4.5) to [out=60,in=180] (1.25,5) to [out=0,in=120] (2,4.5) to [out=-90,in=0] (1.25,2.8) to [out=180,in=-90] (.5,4.5);
        \fill [blue, opacity=.6]  (.5,4.5) to [out=-60,in=180] (1.25,4) to [out=0,in=-120] (2,4.5) to [out=-90,in=0] (1.25,2.8) to [out=180,in=-90] (.5,4.5);
        \draw [thick,red] (2,4.5) to [out=-90,in=0] (1.25,2.8) to [out=180,in=-90] (.5,4.5);
        \fill [yellow, opacity=.6] (0,.5) --(1,.5) to [out=90,in=-90] (3,4.5) -- (2,4.5)  to [out=-90,in=0] (1.25,2.8) to [out=180,in=-90] (.5,4.5) -- (0,4.5) -- (0,.5);
        \draw (0,.5) -- (0,4.5); 
        \fill [blue, opacity=.6] (1,.5) to [out=0,in=180] (4.2,0) -- (4.2,4) to [out=180,in=-30] (3,4.5) to [out=-90,in=90] (1,.5);
        \draw (4.2,0) -- (4.2,4);
        \draw [thick, red] (3,4.5) to [out=-90,in=90] (1,.5);
        \fill [blue, opacity=.6] (4,1) -- (4,5) to [out=180,in=30] (3,4.5) to [out=-90,in=90] (1,.5) -- (1.5,1.5) to [out=90,in=180] (2.25,2.3) to [out=0,in=90] (3,1.5) to [out=-30,in=180] (4,1);
        \draw (4,1) -- (4,5);
        \fill [blue,opacity=.6] (-.2,1.5) -- (1.5,1.5) to [out=90,in=180] (2.25,2.3) to [out=0,in=90] (3,1.5) to [out=-30,in=180] (3.8,2) -- (3.8,6) to [out=180,in=0] (-.2,5.5) -- (-.2,1.5);
        \draw (-.2,5.5) -- (-.2,1.5);
        \draw (3.8,2) -- (3.8,6);
        \fill [yellow, opacity=.6] (1.5,1.5) -- (3,1.5) to [out=90,in=0] (2.25,2.3) to [out=180,in=90] (1.5,1.5);
        \draw [thick, red] (3,1.5) to [out=90,in=0] (2.25,2.3) to [out=180,in=90] (1.5,1.5);
        \draw [double] (0,.5)  -- (1,.5);
        \draw [thick] (-.2,1.5) -- (1.5,1.5);
        \draw [thick] (1,.5) -- (1.5,1.5);
        \draw [double] (1.5,1.5) -- (3,1.5);
        \draw [thick] (1,.5) to [out=0,in=180] (4.2,0);
        \draw [thick] (3,1.5) to [out=-30,in=180] (4,1);
        \draw [thick] (3,1.5) to [out=30,in=180] (3.8,2);
  \end{tikzpicture}
  };
  \end{tikzpicture}
  \label{FS10}
\end{equation}

\begin{equation}
  \begin{tikzpicture}[anchorbase]
    \node at (-4,0) {
      \begin{tikzpicture}[anchorbase,scale=.5, rotate=180]
        \draw [very thick, white, double=black, double distance=.8pt] (2,0) to [out=90,in=-90]  (.5,2);
        \draw [thick, ->] (2,0.05) -- (2,0);
        \draw [very thick, white] (0,0) to [out=90,in=-120] (.5,1);
        \draw [very thick, white] (1,0) to [out=90,in=-60] (.5,1);
        \draw [<-,thick] (0,0) to [out=90,in=-120] (.5,1);
        \draw [<-,thick] (1,0) to [out=90,in=-60] (.5,1);
        \draw [ultra thick, white, double, double=white] (.5,1) to [out=90,in=-90] (1.5,2);
        \draw [double] (.5,1) to [out=90,in=-90] (1.5,2);
        \node at (1.3,1.3) {\tiny $c_3$};
      \end{tikzpicture}
    };
    \node at (-3,0) {$\rightarrow$};
        \node (A) at (2.5,0) {
      \begin{tikzpicture}[anchorbase,scale=.5, rotate=180]
        \draw [thick] (0,0) to [out=90,in=-120] (.5,1);
        \draw [thick] (1,0) to [out=90,in=-60] (.5,1);
        \draw [double] (.5,1) to [out=90,in=-90] (1,1.3);
        \draw [thick] (1,1.3) to [out=120,in=-90] (.5,2);
        \draw [thick] (1,1.3) to [out=60,in=-120] (1.5,1.5);
        \draw [thick] (2,0) to [out=90,in=-60] (1.5,1.5);
        \draw [double] (1.5,1.5) -- (1.5,2);
      \end{tikzpicture}
    };
    \node at (-4,-5) {
      \begin{tikzpicture}[anchorbase,scale=.5, rotate=180]
        \draw [very thick, white, double=black, double distance=.8pt] (2,0) to [out=90,in=-90]  (.5,2);
        \draw [thick, ->] (2,0.05) -- (2,0);
        \draw [ultra thick, white] (0,0) to [out=90,in=-120] (1.5,1.5);
        \draw [ultra thick, white] (1,0) to [out=90,in=-60] (1.5,1.5);
        \draw [<-,thick] (0,0) to [out=90,in=-120] (1.5,1.5);
        \draw [<-,thick] (1,0) to [out=90,in=-60] (1.5,1.5);
        \draw [ultra thick, white, double, double=white] (1.5,1.5) to [out=90,in=-90] (1.5,2);
        \draw [double] (1.5,1.5) to [out=90,in=-90] (1.5,2);
        \node at (1.8,.9) {\tiny $c_2$};
        \node at (.5,1.2) {\tiny $c_1$};
      \end{tikzpicture}
    };
    \node at (-3,-5) {$\rightarrow$};
    \node (B) at (0,-5) {
      \begin{tikzpicture}[anchorbase,scale=.5, rotate=180]
        \node (C1) at (1.5,.3) {};
        \node (C2) at (1.5,.7) {};
        \node (C3) at (.8,1.1) {};
        \node (C4) at (.8,1.5) {};
        \node (C5) at (1.5,1.7) {};
        \draw [thick] (0,0) to [out=90,in=-90] (.5,2);
        \draw [thick] (1,0) to [out=90,in=-120] (C5.center);
        \draw [thick] (2,0) to [out=90,in=-60] (C5.center);
        \draw [double] (C5.center) to [out=90,in=-90] (1.5,2);
      \end{tikzpicture}};
    \node (C) at (2.5,-4) {
      \begin{tikzpicture}[anchorbase,scale=.5, rotate=180]
        \node (C1) at (1.5,.3) {};
        \node (C2) at (1.5,.7) {};
        \node (C3) at (.8,1.1) {};
        \node (C4) at (.8,1.5) {};
        \node (C5) at (1.5,1.7) {};
        \draw [thick] (0,0) to [out=90,in=-120] (C3.center);
        \draw [thick] (1,0) to [out=90,in=-60] (C3.center);
        \draw [thick] (2,0) to [out=90,in=-60] (C5.center);
        \draw [double] (C3.center) -- (C4.center);
        \draw [thick] (C4.center) to [out=120,in=-90] (.5,2);
        \draw [thick] (C4.center) to [out=60,in=-120] (C5.center);
        \draw [double] (C5.center) to [out=90,in=-90] (1.5,2);
      \end{tikzpicture}};
    \node (D) at (2.5,-6) {
      \begin{tikzpicture}[anchorbase,scale=.5, rotate=180]
        \node (C1) at (1.5,.3) {};
        \node (C2) at (1.5,.7) {};
        \node (C3) at (.8,1.1) {};
        \node (C4) at (.8,1.5) {};
        \node (C5) at (1.5,1.7) {};
        \draw [thick] (0,0) to [out=90,in=-90] (.5,2);
        \draw [thick] (1,0) to [out=90,in=-120] (C1.center);
        \draw [thick] (2,0) to [out=90,in=-60] (C1.center);
        \draw [double] (C1.center) -- (C2.center);
        \draw [thick] (C2.center) to [out=120,in=-90] (1.1,1.2) to [out=90,in=-120] (C5.center);
        \draw [thick] (C2.center) to [out=60,in=-90] (1.9,1.2) to [out=90,in=-60] (C5.center);
        \draw [double] (C5.center) to [out=90,in=-90] (1.5,2);
      \end{tikzpicture}};
    \node (E) at (5,-5) {
      \begin{tikzpicture}[anchorbase,scale=.5, rotate=180]
        \node (C1) at (1.5,.3) {};
        \node (C2) at (1.5,.7) {};
        \node (C3) at (.8,1.1) {};
        \node (C4) at (.8,1.5) {};
        \node (C5) at (1.5,1.7) {};
        \draw [thick] (0,0) to [out=90,in=-120] (C3.center);
        \draw [thick] (1,0) to [out=90,in=-120] (C1.center);
        \draw [thick] (2,0) to [out=90,in=-60] (C1.center);
        \draw [double] (C1.center) -- (C2.center);
        \draw [thick] (C2.center) to [out=120,in=-60] (C3.center);
        \draw [thick] (C2.center) to [out=60,in=-60] (C5.center);
        \draw [double] (C3.center) -- (C4.center);
        \draw [thick] (C4.center) to [out=120,in=-90] (.5,2);
        \draw [thick] (C4.center) to [out=60,in=-120] (C5.center);
        \draw [double] (C5.center) to [out=90,in=-90] (1.5,2);
      \end{tikzpicture}};
    \draw (B) -- (C) node[midway,above] {\tiny $\cdot\wedge c_1$};
    \draw (B) -- (D) node[midway,above] {\tiny $\cdot\wedge c_2$};
    \draw (C) -- (E) node[midway,above] {\tiny $\cdot\wedge c_2$};
    \draw (D) -- (E) node[midway,above] {\tiny $\cdot \wedge c_1$};
    \draw [green,->] (A) edge[bend right=30] node [below,rotate=90] {\tiny $-a_{11}\langle \cdot,c_3\rangle \wedge c_1$} (C);
    \draw [green,->] (A) edge[bend right=40] node [above,rotate=90] {\tiny $a_{11}F\otimes \langle \cdot,c_3\rangle \wedge c_2$} (D);
    \draw [red,->] (C) edge[bend right=30] node [above,rotate=90] {\tiny $-a_{11}\langle \cdot,c_1\rangle \wedge c_3$} (A);
    \draw [red,->] (D) edge[bend right=40] node [below,rotate=90] {\tiny $a_{11}F\otimes \langle \cdot,c_2\rangle \wedge c_3$} (A);
        \node at (8,-2) {$F\leftrightarrow$  
      \begin{tikzpicture}[anchorbase,scale=-.3]
        \draw [double] (-.2,5.5)  -- (.5,5.5);
        \draw [thick] (.5,5.5) to [out=60,in=180] (1.25,6) to [out=0,in=120] (2,5.5);
        \draw [thick] (.5,5.5) to [out=-60,in=180] (1.25,5) to [out=0,in=-120] (2,5.5);
        \draw [thick] (0,4.5) to [out=0,in=180] (4.2,4);
        \draw [double] (2,5.5) -- (3,5.5);
        \draw [thick] (3,5.5) to [out=30,in=180] (3.8,6);
        \draw [thick] (3,5.5) to [out=-30,in=180] (4,5);
        \fill [blue,opacity=.6] (0,.5) -- (1.5,.5) to [out=90,in=180] (2.25,1.3) to [out=0,in=90] (3,.5) to [out=-30,in=180] (4.2,0) -- (4.2,4) to [out=180,in=0] (0,4.5) -- (0,.5);
        \draw (0,4.5) -- (0,.5);
        \draw (4.2,0) -- (4.2,4);
        \fill [yellow, opacity=.6] (1.5,.5) -- (3,.5) to [out=90,in=0] (2.25,1.3) to [out=180,in=90] (1.5,.5);
        \draw [thick, red] (3,.5) to [out=90,in=0] (2.25,1.3) to [out=180,in=90] (1.5,.5);
        \fill [blue, opacity=.6] (4,1) -- (4,5) to [out=180,in=-30] (3,5.5) to [out=-90,in=90] (1,1.5) -- (1.5,.5) to [out=90,in=180] (2.25,1.3) to [out=0,in=90] (3,.5) to [out=30,in=180] (4,1);
        \draw [thick, red] (3,5.5) to [out=-90,in=90] (1,1.5);
        \draw (4,1) -- (4,5);
        \fill [yellow, opacity=.6] (-.2,1.5) --(1,1.5) to [out=90,in=-90] (3,5.5) -- (2,5.5)  to [out=-90,in=0] (1.25,3.8) to [out=180,in=-90] (.5,5.5) -- (-.2,5.5) -- (-.2,1.5);
        \draw (-.2,1.5) -- (-.2,5.5); 
        \fill [blue, opacity=.6] (1,1.5) to [out=0,in=180] (3.8,2) -- (3.8,6) to [out=180,in=30] (3,5.5) to [out=-90,in=90] (1,1.5);
        \draw (3.8,2) -- (3.8,6);
        \fill [blue, opacity=.6]  (.5,5.5) to [out=60,in=180] (1.25,6) to [out=0,in=120] (2,5.5) to [out=-90,in=0] (1.25,3.8) to [out=180,in=-90] (.5,5.5);
        \fill [blue, opacity=.6]  (.5,5.5) to [out=-60,in=180] (1.25,5) to [out=0,in=-120] (2,5.5) to [out=-90,in=0] (1.25,3.8) to [out=180,in=-90] (.5,5.5);
        \draw [thick,red] (2,5.5) to [out=-90,in=0] (1.25,3.8) to [out=180,in=-90] (.5,5.5);
        \draw [double] (-.2,1.5)  -- (1,1.5);
        \draw [thick] (0,.5) -- (1.5,.5);
        \draw [thick] (1,1.5) -- (1.5,.5);
        \draw [double] (1.5,.5) -- (3,.5);
        \draw [thick] (1,1.5) to [out=0,in=180] (3.8,2);
        \draw [thick] (3,.5) to [out=30,in=180] (4,1);
        \draw [thick] (3,.5) to [out=-30,in=180] (4.2,0);
  \end{tikzpicture}
  };
  \end{tikzpicture}
  \label{FS11}
\end{equation}

\begin{equation}
  \begin{tikzpicture}[anchorbase]
    \node at (-4,0) {
      \begin{tikzpicture}[anchorbase,scale=.5, rotate=180]
        \draw [very thick, white, double=black, double distance=.8pt] (0,0) to [out=90,in=-90]  (1.5,2);
        \draw [thick, ->] (0,0.05) -- (0,0);
        \draw [ultra thick, white] (2,0) to [out=90,in=-60] (1.5,1);
        \draw [ultra thick, white] (1,0) to [out=90,in=-120] (1.5,1);
        \draw [<-,thick] (2,0) to [out=90,in=-60] (1.5,1);
        \draw [<-,thick] (1,0) to [out=90,in=-120] (1.5,1);
        \draw [ultra thick, white, double, double=white] (1.5,1) to [out=90,in=-90] (.5,2);
        \draw [double] (1.5,1) to [out=90,in=-90] (.5,2);
        \node at (.7,1.3) {\tiny $c_3$};
      \end{tikzpicture}
    };
    \node at (-3,0) {$\rightarrow$};
        \node (A) at (2.5,0) {
      \begin{tikzpicture}[anchorbase,scale=.5,xscale=-1, rotate=180]
        \draw [thick] (0,0) to [out=90,in=-120] (.5,1);
        \draw [thick] (1,0) to [out=90,in=-60] (.5,1);
        \draw [double] (.5,1) to [out=90,in=-90] (1,1.3);
        \draw [thick] (1,1.3) to [out=120,in=-90] (.5,2);
        \draw [thick] (1,1.3) to [out=60,in=-120] (1.5,1.5);
        \draw [thick] (2,0) to [out=90,in=-60] (1.5,1.5);
        \draw [double] (1.5,1.5) -- (1.5,2);
      \end{tikzpicture}
    };
    \node at (-4,-5) {
      \begin{tikzpicture}[anchorbase,scale=.5, rotate=180]
        \draw [very thick, white, double=black, double distance=.8pt] (0,0) to [out=90,in=-90]  (1.5,2);
        \draw [thick, ->] (0,0.05) -- (0,0);
        \draw [ultra thick, white] (2,0) to [out=90,in=-60] (.5,1.5);
        \draw [ultra thick, white] (1,0) to [out=90,in=-120] (.5,1.5);
        \draw [<-,thick] (2,0) to [out=90,in=-60] (.5,1.5);
        \draw [<-,thick] (1,0) to [out=90,in=-120] (.5,1.5);
        \draw [ultra thick, white, double, double=white] (.5,1.5) to [out=90,in=-90] (.5,2);
        \draw [double] (.5,1.5) to [out=90,in=-90] (.5,2);
        \node at (.2,.9) {\tiny $c_2$};
        \node at (1.5,1.2) {\tiny $c_1$};
      \end{tikzpicture}
    };
    \node at (-3,-5) {$\rightarrow$};
    \node (B) at (0,-5) {
      \begin{tikzpicture}[anchorbase,scale=.5,xscale=-1, rotate=180]
        \node (C1) at (1.5,.3) {};
        \node (C2) at (1.5,.7) {};
        \node (C3) at (.8,1.1) {};
        \node (C4) at (.8,1.5) {};
        \node (C5) at (1.5,1.7) {};
        \draw [thick] (0,0) to [out=90,in=-120] (C3.center);
        \draw [thick] (1,0) to [out=90,in=-120] (C1.center);
        \draw [thick] (2,0) to [out=90,in=-60] (C1.center);
        \draw [double] (C1.center) -- (C2.center);
        \draw [thick] (C2.center) to [out=120,in=-60] (C3.center);
        \draw [thick] (C2.center) to [out=60,in=-60] (C5.center);
        \draw [double] (C3.center) -- (C4.center);
        \draw [thick] (C4.center) to [out=120,in=-90] (.5,2);
        \draw [thick] (C4.center) to [out=60,in=-120] (C5.center);
        \draw [double] (C5.center) to [out=90,in=-90] (1.5,2);
      \end{tikzpicture}};
    \node (C) at (2.5,-4) {
      \begin{tikzpicture}[anchorbase,scale=.5,xscale=-1, rotate=180]
        \node (C1) at (1.5,.3) {};
        \node (C2) at (1.5,.7) {};
        \node (C3) at (.8,1.1) {};
        \node (C4) at (.8,1.5) {};
        \node (C5) at (1.5,1.7) {};
        \draw [thick] (0,0) to [out=90,in=-120] (C3.center);
        \draw [thick] (1,0) to [out=90,in=-60] (C3.center);
        \draw [thick] (2,0) to [out=90,in=-60] (C5.center);
        \draw [double] (C3.center) -- (C4.center);
        \draw [thick] (C4.center) to [out=120,in=-90] (.5,2);
        \draw [thick] (C4.center) to [out=60,in=-120] (C5.center);
        \draw [double] (C5.center) to [out=90,in=-90] (1.5,2);
      \end{tikzpicture}};
    \node (D) at (2.5,-6) {
      \begin{tikzpicture}[anchorbase,scale=.5,xscale=-1, rotate=180]
        \node (C1) at (1.5,.3) {};
        \node (C2) at (1.5,.7) {};
        \node (C3) at (.8,1.1) {};
        \node (C4) at (.8,1.5) {};
        \node (C5) at (1.5,1.7) {};
        \draw [thick] (0,0) to [out=90,in=-90] (.5,2);
        \draw [thick] (1,0) to [out=90,in=-120] (C1.center);
        \draw [thick] (2,0) to [out=90,in=-60] (C1.center);
        \draw [double] (C1.center) -- (C2.center);
        \draw [thick] (C2.center) to [out=120,in=-90] (1.1,1.2) to [out=90,in=-120] (C5.center);
        \draw [thick] (C2.center) to [out=60,in=-90] (1.9,1.2) to [out=90,in=-60] (C5.center);
        \draw [double] (C5.center) to [out=90,in=-90] (1.5,2);
      \end{tikzpicture}};
    \node (E) at (5,-5) {
      \begin{tikzpicture}[anchorbase,scale=.5,xscale=-1, rotate=180]
        \node (C1) at (1.5,.3) {};
        \node (C2) at (1.5,.7) {};
        \node (C3) at (.8,1.1) {};
        \node (C4) at (.8,1.5) {};
        \node (C5) at (1.5,1.7) {};
        \draw [thick] (0,0) to [out=90,in=-90] (.5,2);
        \draw [thick] (1,0) to [out=90,in=-120] (C5.center);
        \draw [thick] (2,0) to [out=90,in=-60] (C5.center);
        \draw [double] (C5.center) to [out=90,in=-90] (1.5,2);
      \end{tikzpicture}};
    \draw (B) -- (C) node[midway,above] {\tiny $\langle \cdot, c_2\rangle$};
    \draw (B) -- (D) node[midway,above] {\tiny $\langle \cdot,c_1\rangle$};
    \draw (C) -- (E) node[midway,above] {\tiny $\langle \cdot, c_1\rangle$};
    \draw (D) -- (E) node[midway,above] {\tiny $\langle\cdot,c_2\rangle$};
    \draw [green,->] (A) edge[bend right=30] node [below,rotate=90] {\tiny $a_{4}a_{11}\langle \cdot,c_3\rangle \wedge c_1$} (C);
    \draw [green,->] (A) edge[bend right=40] node [above,rotate=90] {\tiny $-a_{4}a_{11}F\otimes \langle \cdot,c_3\rangle \wedge c_2$} (D);
    \draw [red,->] (C) edge[bend right=30] node [above,rotate=90] {\tiny $a_4a_{11}\langle \cdot,c_1\rangle \wedge c_3$} (A);
    \draw [red,->] (D) edge[bend right=40] node [below,rotate=90] {\tiny $-a_4a_{11}F\otimes \langle \cdot,c_2\rangle \wedge c_3$} (A);
        \node at (8,-2) {$F\leftrightarrow$  

  \label{FS16}
\end{equation}

\subsection{Trivalent twists} \label{sec:trivtwists}

\begin{equation} \label{TT1}
  %
\end{equation}

In fully deformed homology, on both surviving generators, the maps equal:
\begin{itemize}
  \item $4a_3 \zeta \rm{cap}^2 \otimes\langle \cdot,c_1\wedge c_2 \rangle$ going downward; 
  \item $-a_3\zeta\rm{cup}^2\otimes \wedge c_1\wedge c_2$ going upward.
\end{itemize}

\begin{equation} \label{TT2}
  %
\end{equation}

In fully deformed homology, on both surviving generators, the maps equal:
\begin{itemize}
  \item $4a_3\zeta \rm{cap}^2 \otimes\langle \cdot,c_1\wedge c_2 \rangle$ going downward;
  \item $-a_3\zeta\rm{cup}^2\otimes \wedge c_1\wedge c_2$ going upward.
\end{itemize}

\begin{equation} \label{TT3}
  %
\end{equation}

In fully deformed homology, on both surviving generators, the maps equal:
\begin{itemize}
  \item $-a_3\zeta\rm{cap}^2\otimes \langle \cdot, c_1\wedge c_2\rangle$ going downward.
  \item $4a_3\zeta \rm{cup}^2 \otimes \wedge c_1\wedge c_2 $ going upward;
\end{itemize}

\begin{equation} \label{TT4}
  %
};
  \end{tikzpicture}
\end{equation}

In fully deformed homology, on both surviving generators, the maps equal:
\begin{itemize}
  \item $-a_3\zeta\rm{cap}^2\otimes \langle \cdot, c_1\wedge c_2\rangle$ going downward.
  \item $4a_3\zeta \rm{cup}^2 \otimes \wedge c_1\wedge c_2 $ going upward;
\end{itemize}


\end{document}